\documentclass{amsart}

\usepackage{amsmath}
\usepackage{amsthm}
\usepackage{amssymb}
\usepackage{amsfonts}
\usepackage{amscd}

\makeindex
\makeatletter

\begin{document}

\bibliographystyle{amsalpha}
\makeatletter

\def\comm#1{{\bf #1}}
\newcommand{\e}{\epsilon}
\newcommand{\0}{{\bold 0}}
\newcommand{\w}{{\bold w}}
\newcommand{\y}{{\bold y}}
\newcommand{\p}{{\bold p}}
\newcommand{\q}{{\bold q}}
\newcommand{\ba}{{\bold a}}
\newcommand{\bc}{{\bold c}}
\newcommand{\be}{{\bold e}}

\newcommand{\balpha}{{\boldsymbol \alpha}}
\newcommand{\bbeta}{{\boldsymbol \beta}}
\newcommand{\blambda}{{\boldsymbol \lambda}}
\newcommand{\bt}{{\bold t}}
\newcommand{\bl}{{\bf l}}
\newcommand{\z}{{\bold z}}
\newcommand{\x}{{\bold x}}
\newcommand{\bH}{{\mathbf H}}
\newcommand{\N}{{\bold N}}
\newcommand{\Z}{{\bold Z}}
\newcommand{\F}{{\bold F}}
\newcommand{\R}{{\bold R}}
\newcommand{\Q}{{\bold Q}}
\newcommand{\C}{{\bold C}}
\newcommand{\BP}{{\mathbf P}}
\newcommand{\BF}{{\mathbf F}}
\newcommand{\cA}{{\mathcal A}}
\newcommand{\cG}{{\mathcal G}}
\newcommand{\cO}{{\mathcal O}}
\newcommand{\cX}{{\mathcal X}}
\newcommand{\cH}{{\mathcal H}}
\newcommand{\cM}{{\mathcal M}}
\newcommand{\cD}{{\mathcal D}}
\newcommand{\cB}{{\mathcal B}}
\newcommand{\cC}{{\mathcal C}}
\newcommand{\cT}{{\mathcal T}}
\newcommand{\cI}{{\mathcal I}}
\newcommand{\cJ}{{\mathcal J}}
\newcommand{\cS}{{\mathcal S}}
\newcommand{\cF}{{\mathcal F}}
\newcommand{\cE}{{\mathcal E}}
\newcommand{\cP}{{\mathcal P}}
\newcommand{\cR}{{\mathcal R}}
\newcommand{\E}{{\mathbf E}}
\newcommand{\V}{{\mathbf V}}

\newcommand{\cL}{{\mathcal L}}
\newcommand{\cY}{{\mathcal Y}}
\newcommand{\sB}{{\sf B}}
\newcommand{\sE}{{\sf E}}

\newcommand{\sA}{{\sf A}}
\newcommand{\ga}{{\sf a}}
\newcommand{\es}{{\sf s}}
\newcommand{\m}{{\bold m}}
\newcommand{\bS}{{\bold S}}

\newcommand{\ovf}{{\overline{f}}}

\newcommand{\ihra}{\stackrel{i}{\hookrightarrow}}
\newcommand\rank{\mathop{\rm rank}\nolimits}
\newcommand\im{\mathop{\rm Im}\nolimits}
\newcommand\coker{\mathop{\rm coker}\nolimits}
\newcommand{\Aut}{\mathop{\rm Aut}\nolimits}
\newcommand{\Tr}{\mathop{\rm Tr}\nolimits}
\newcommand\Li{\mathop{\rm Li}\nolimits}
\newcommand\NS{\mathop{\rm NS}\nolimits}
\newcommand\Hom{\mathop{\rm Hom}\nolimits}
\newcommand\Ext{\mathop{\rm Ext}\nolimits}
\newcommand\End{\mathop{\rm End}\nolimits}
\newcommand\Pic{\mathop{\rm Pic}\nolimits}
\newcommand\Spec{\mathop{\rm Spec}\nolimits}
\newcommand\Hilb{\mathop{\rm Hilb}\nolimits}
\newcommand\pardeg{\mathop{\rm pardeg}\nolimits}
\newcommand\Ker{\mathop{\rm Ker}\nolimits}
\newcommand\RH{\mathop{\bf RH}\nolimits}
\newcommand{\length}{\mathop{\rm length}\nolimits}
\newcommand{\res}{\mathop{\sf res}\nolimits}
\newcommand\Quot{\mathop{\rm Quot}\nolimits}
\newcommand\Grass{\mathop{\rm Grass}\nolimits}
\newcommand\Proj{\mathop{\rm Proj}\nolimits}
\newcommand{\codim}{\mathop{\rm codim}\nolimits}

\newcommand\lra{\longrightarrow}
\newcommand\ra{\rightarrow}
\newcommand\la{\leftarrow}
\newcommand\JG{J_{\Gamma}}
\newcommand{\wvskp}{\vspace{1cm}}
\newcommand{\vskp}{\vspace{5mm}}
\newcommand{\nvskp}{\vspace{1mm}}
\newcommand{\nid}{\noindent}
\newcommand{\new}{\nvskp \nid}
\newtheorem{Assumption}{Assumption}[section]
\newtheorem{Theorem}{Theorem}[section]
\newtheorem{Lemma}{Lemma}[section]
\newtheorem{Remark}{Remark}[section]
\newtheorem{Corollary}{Corollary}[section]
\newtheorem{Conjecture}{Conjecture}[section]
\newtheorem{Proposition}{Proposition}[section]
\newtheorem{Example}{Example}[section]
\newtheorem{Definition}{Definition}[section]
\newtheorem{Question}{Question}[section]
\newtheorem{Fact}{Fact}[section]

\renewcommand{\thefootnote}{\fnsymbol{footnote}}

\baselineskip=18 pt

\title[Moduli of Stable Parabolic Connections]
{Moduli of Stable Parabolic Connections, 
Riemann-Hilbert correspondence and   
Geometry of Painlev\'{e} equation of type $VI$, Part I}
\thanks{Partly supported by Grant-in Aid
for Scientific Research  (B-12440008), (B-12440043),   (Houga-138740023),   (Wakate-B-15740018)
the Ministry of Education, Science and Culture and JSPS-NWO exchange program. }
\keywords{Stable parabolic connections , Representation of fundamental 
groups, Riemann-Hilbert correspondences, Symplectic structure, Painlev\'e equations, Garnier equations}
\author{Michi-aki Inaba}
\author{Katsunori Iwasaki}
\author{Masa-Hiko Saito}
\dedicatory{Dedicated to Professor Kyoichi Takano on his 60th birthday}
\address{Faculty of Mathematics, Kyushu University, 6-10-1, Hakozaki, 
Higashi-ku, Fukuoka 812-8581, Japan}
\email{inaba@math.kyushu-u.ac.jp}
\email{iwasaki@math.kyushu-u.ac.jp}
\address{Department of Mathematics, Faculty of Science, 
Kobe University, Kobe, Rokko, 657-8501, Japan}
\email{mhsaito@math.kobe-u.ac.jp}
\subjclass[2000]{34M55, 14D20, 32G34, 34G20, 58F05}
\footnote[0]{Communicated by T. Kawai. Received March 15, 2005}

\begin{abstract} In this paper, we will give 
a complete geometric background for the geometry of 
Painlev\'e $VI$ and Garnier equations.  By geometric invariant theory, 
we will construct a smooth fine moduli space $M_n^{\balpha}(\bt, \blambda, L) $ 
 of stable parabolic connections  
on $\BP^1$ with logarithmic poles 
at $D(\bt) = t_1 + \cdots + t_n$  as well as its natural compactification.  
Moreover the moduli space $\cR(\cP_{n, \bt})_{\ba}$ 
 of Jordan equivalence classes of $SL_2(\C)$-representations of the fundamental 
group $\pi_1(\BP^1 \setminus D(\bt),\ast)$ are defined as the categorical 
quotient.  We define the Riemann-Hilbert 
correspondence $\RH: M_n^{\balpha}(\bt, \blambda, L) \lra \cR(\cP_{n, \bt})_{\ba}$ 
and prove that $\RH$ is a bimeromorphic proper surjective analytic map. 
Painlev\'e and Garnier equations can be derived from the isomonodromic flows 
and Painlev\'e property of these equations are easily derived from the properties of 
$\RH$.  We also prove that the smooth parts of 
both moduli spaces have natural symplectic structures and 
$\RH$ is a symplectic resolution of singularities of 
$\cR(\cP_{n, \bt})_{\ba}$, from which one can give geometric backgrounds for 
other interesting  phenomena, like Hamiltonian structures, 
 B\"acklund transformations, special solutions of 
these equations.  
\end{abstract}

\maketitle
%\tableofcontents

\section{Introduction}

\subsection{The purpose}

\quad 

The purpose of the series of  papers  is to give  a complete 
geometric background for Painlev\'e equations of type $VI$ or
more generally for the so-called Garnier equations.  

As is well-known, these nonlinear differential 
equations have the Painlev\'e property which means that 
 {\em generic}  solutions of  these equations have 
 no movable singularity except 
for poles so that solutions have the
analytic continuations  on whole of the universal covering 
of the space of time variables.    

Besides the Painlev\'e property, there are several interesting 
phenomena related to these equations which  have been investigated by many 
authors.  
\begin{itemize}
\item Each of these equations can be written in a Hamiltonian system 
by a natural  symplectic coordinate system (\cite{Mal}, \cite{O3}, \cite{Iwa91},  \cite{Iwa92}, 
 \cite{K}, \cite{ST}).

\item These equations have  natural parameters 
 $\blambda = (\lambda_1, \ldots, \lambda_n) \in \C^n  $. 
 Moreover there exist birational symmetries of these equations, called 
 {\em B\"{a}cklund transformations of these equations}, 
 which act on both of variables and the parameters  and preserve the 
 equations.  (\cite{O4}).

 \item In Painlev\'e $VI$ case, the group  of  all B\"{a}cklund 
   transformations is isomorphic  to 
   the  affine Weyl group $W(D_{4}^{(1)})$ of the type $D_{4}^{(1)}$.  
   (\cite{O4}, \cite{Sakai}, \cite{AL2}, \cite{NY}, \cite{IIS-0}).

 \item In Painlev\'e $VI$ case,  if  $\blambda \in \C^4$ lies on a reflection 
 hyperplane of a reflection in  $W(D_{4}^{(1)})$, then the 
 corresponding equation has one parameter family of  Riccati solutions. 
 (\cite{LY}, \cite{FA}, \cite{W}, \cite{STe02}, \cite{SU01}).

 \item A natural compactification of each  space of initial conditions for $P_{VI}$, 
 introduced by Okamoto \cite{O1}, can be obtained by a series of 
  explicit blowing-ups of 
 $\BP^1_{\C} \times \BP^1_{\C}$ or $\F_2$.  
   The compactification is given by a smooth projective rational surface $S$ and it 
has a  unique anti-canonical divisors $-K_S=Y$ 
such that $S \setminus Y_{red}$ is the space of initial conditions for $P_{VI}$.   
The pair $(S, Y)$ becomes 
an Okamoto-Painlev\'e pair of type $D_4^{(1)}$ 
in the sense of \cite{STT02}.  (See also \cite{Sakai}).   
 \end{itemize}

Though these phenomena are discussed and  investigated 
by many authors, the intrinsic mathematical 
 background for these facts remains to be understood.   
  Therefore, for example,  it is worthwhile to ask 
the following fundamental questions:
$$
\fbox{\begin{tabular}{ll}
\\
{$\bullet$} & What is the geometric meaning of Painlev\'e property for these equations? \\ 
\\

{$\bullet$}  & What is the geometric meaning of the symplectic structure ?  \\
\\
{$\bullet$} & What is the geometric origin of 
 B\"acklund transformations? \\
\\
{$\bullet$}  &  Why Riccati solutions or some classical solutions  appear 
for the parameters   \\  
& on the  reflection hyperplanes of the B\"acklund transformations ?\\
\\
\end{tabular}}
$$

In the series of the papers, the authors will give  answers  
to these questions in a natural intrinsic framework.

\subsection{Natural Framework}
\label{subsec:naturalfram}
\quad

It is already known (cf. \cite{Fuchs}, \cite{Garnier}, \cite{Sch}, \cite{JMU}, \cite{O3},  \cite{Iwa91} and 
\cite{Iwa92} )
that these equations can be derived from the isomonodromic 
deformation of  the systems of linear equations of rank 
2 with regular singularities  over $\BP^1$ 
 or equivalently linear connections on vector bundles  of rank 2 with 
logarithmic poles  over $\BP^1$.  
Although we will  follow this line in  this paper,  for several essential  reasons,  
we have to  introduce  a slight generalization of  linear connections 
which will be explained as follows.

Let $n \geq 3$ and  let us set  $T_n = \{ \bt =(t_1, \ldots, t_n) \in (\BP^1_{\C})^n \ | \ t_i \not= t_j, ( i \not= j) \}  $ , $\Lambda_n = \{ \blambda = (\lambda_1, \ldots, \lambda_n) \in \C^n \} $.   Fix a data 
$ (\bt, \blambda) \in T_n \times \Lambda_n $ and set $D(\bt) 
= t_1 + \cdots + t_n $.  We also fix  a line bundle $L$ on 
$\BP^1_{\C}$ with a logarithmic connection $ \nabla_L : L \lra L \otimes 
\Omega_{\BP^1_{\C}}^1(D(\bt))$.  

A quadruple $(E, \nabla, \varphi, l = \{ l_i \}_{i=1}^n )$ 
consisting  of:
\begin{enumerate}
\item 
a rank 2 vector bundle 
$E$ on $\BP^1$, 

\item 
a logarithmic connection $\nabla:E \lra E 
\otimes \Omega_{\BP^1}^1(D(\bt)) $  

\item 
a bundle isomorphism $ \varphi: \wedge^2 E \stackrel{\simeq}{\lra} L  $ and 

\item  one dimensional subspace $l_i$ of the fiber $ E_{t_i}$  of $E$ 
at $t_i$,   $l_i \subset E_{t_i}$,   $ i = 1, \ldots, n$,   
\end{enumerate}
is called a {\em 
$(\bt, \blambda)$-parabolic connection } with the 
determinant $(L, \nabla_L)$  if they satisfy the 
following conditions:
 \begin{enumerate}
\item  for any local sections $s_1, s_2$ of $E$, 
$$
(\varphi \otimes {\rm id}) 
(\nabla s_1 \wedge s_2 + s_1 \wedge \nabla s_2) = \nabla_L(\varphi(s_1 \wedge s_2)),  
$$

\item  
$l_i \subset  \Ker (\res_{t_i}(\nabla) - \lambda_i)$, that is, $\lambda_i$ is an  
eigenvalue of the residue $\res_{t_i}(\nabla)$ of $\nabla$ at $t_i$ and 
$l_i$ is a corresponding one-dimensional eigensubspace of $\res_{t_i}(\nabla)$.   
\end{enumerate}

We introduce a series of  rational  numbers 
$\balpha=(\alpha_1, \ldots, \alpha_{2n})$ such that $0 \leq \alpha_1 < \ldots < \alpha_{2n} < 1$, 
which is called {\em a weight}.  By using a weight $\balpha$, 
one can define  parabolic degrees for  
$(\bt, \blambda)$-parabolic connections  $(E, \nabla, \varphi, l)$
and introduce the notion of the parabolic  stability.  Let 
$
M_n^{\balpha}(\bt, \blambda, L)
$
be  the coarse moduli space of stable 
$(\bt, \blambda)$-parabolic connections on  $\BP^1$ 
with the determinant $(L, \nabla_L)$.  
Considering the relative setting over the parameter space $ T_n \times \Lambda_n= \{ (\bt, \blambda) \} $, 
we can construct a family of moduli spaces 
\begin{equation}\label{eq:family-con}
\pi_n:M_n^{\balpha} (L) \lra  T_n \times \Lambda_n 
\end{equation}
such that $\pi_n^{-1}(\bt, \blambda) \simeq M_n^{\balpha}(\bt, \blambda, L)$.  
Later, we have  to extend the family by  a finite \'etale covering 
$T'_n \lra T_n$, and for simplicity,  we denote it  also by $\pi_n: M_n^{\balpha}(L) \lra T'_n \times \Lambda_n$.

Next, let us fix   $\bt \in T_n$ and   consider a representation  
$
\rho:\pi_1(\BP^1_{\C} \setminus D(\bt), *) \lra SL_2(\C)
$
of the fundamental group  $\pi_1(\BP^1_{\C} \setminus D(\bt), \ast)$ 
with a fixed base point $* \in \BP^1_{\C}$.   
Two representations $\rho_1 $ and $\rho_2$ 
 are said to be equivalent 
if there exists an element $P \in SL_2(\C)$ such 
that $\rho_2 = P^{-1} \rho_1 P$. 
To each representation $\rho$, one can associate
a local system 
${\bf E}_{\rho}$ of rank $2$ on $\BP^1_{\C} \setminus D(\bt)$ 
with an isomorphism 
$\wedge^2 {\bf E}_{\rho} \simeq \C_{\BP^1_{\C} \setminus D(\bt)}$. 
Moreover, two representations $\rho_1$ and $\rho_2$ are equivalent to each 
other if and only if  ${\bf E}_{\rho_1}$ and ${\bf E}_{\rho_2}$ are isomorphic as 
local systems.  Hence the moduli space of the isomorphism classes of local systems on $\BP^1 \setminus D(\bt)$ with  trivial determinants 
is isomorphic to the moduli space of equivalence classes of the representations. 

Since $\pi_1(\BP^1_{\C} \setminus D(\bt), *)$ is a free group 
generated by $\gamma_i $ for $1 \leq i \leq n-1$ where 
$\gamma_i$ is a loop around the point $t_i$, such a 
representation can be determined by $M_i = \rho(\gamma_i) \in SL_2(\C) $ 
for $ 1 \leq i \leq n-1$.  Therefore the moduli space should be a quotient space of  
$SL_2(\C)^{n-1}$ by a diagonal 
adjoint action of $SL_2(\C)$. 
 
However there is no canonical way to give a scheme structure on the set of 
equivalence classes of the representations.  In this sense,   
we have to  introduce a stronger equivalence relation.  
Two $SL_2(\C)$-representations $\rho_1$ and $\rho_2$ 
of $\pi_1(\BP^1_{\C}\setminus D(\bt), *)$ 
are said to be {\em Jordan equivalent} if their 
semisimplifications are equivalent. This means that if a local system 
${\bf E}_{\rho}$ is an extension of rank one local systems $L_1$ and $L_2$
one can not distinguish the extension classes.   
As is shown by 
Simpson \cite{Simp-II}, the set of {\em the Jordan equivalence classes} of the 
local systems or representations  is equal to the set of closed points of the categorical quotient 
$$
\cR(\cP_{n, \bt}) = SL_2(\C)^{n-1}/Ad(SL_2(\C)), 
$$
of $SL_2(\C)^{n-1}$ by the diagonal adjoint action of $SL_2(\C)$. The categorical quotients is defined as the affine scheme of the ring of 
invariant functions on   $SL_2(\C)^{n-1}$ by the action of  $SL_2(\C)$.   (Cf.\  \S \ref{sec:monodromy}) .   

Fixing the canonical generators $\gamma_i $ ( $1 \leq i \leq n$) of 
$\pi_1(\BP^1_{\C} \setminus D(\bt), *)$, to each representation 
$\rho: \pi_1(\BP^1_{\C} \setminus D(\bt), *) \lra SL_2(\C)$, we can associate   $n$-algebraic functions on $SL_2(\C)^{n-1}$ 
$$
\Tr (\rho(\gamma_i) ) = a_i,  \quad \Tr(\rho( (\gamma_1 \cdots \gamma_{n-1})^{-1}))  
= \Tr (\rho(\gamma_n))  =a_n
$$
which are clearly invariant under the adjoint action.  Setting 
$\cA_n = \Spec \C[a_1, \ldots, a_n] \simeq \C^n $, we obtain a natural morphism 
$$
p_n: \cR(\cP_{n, \bt})  \lra \cA_n. 
$$
For a fixed closed point $\ba=(a_1, \ldots, a_n) \in \cA_n$, let us  denote by 
$\cR(\cP_{n, \bt})_{\ba} = p_n^{-1}(\ba) $ the closed fiber at $\ba$, that is, 
we set 
$$
\cR(\cP_{n, \bt})_{\ba} = \{   \  [\rho] \in \cR(\cP_{n,\bt}) \ | \ 
\Tr (\rho(\gamma_i))  = a_i, \ 1 \leq i \leq n \}. 
$$
Moreover, taking   a finite \'etale covering  $T'_n \lra T_n $  
we can obtain a family of moduli spaces 
\begin{equation}\label{eq:family-rep-1}
\phi_n : \cR_n \lra T'_n \times \cA_n 
\end{equation}
such that   $\phi_n^{-1}(\bt, \ba) = \cR(\cP_{n, \bt})_{\ba}$ (cf.  \S \ref{sec:monodromy}).

Now,  we have obtained two kinds of moduli spaces 
$M_n^{\balpha} (\bt, \blambda,  L)$ and $\cR(\cP_{n,\bt})_{\ba} $  for fixed 
$(\bt, \blambda) \in T'_n \times \Lambda_n$ and $(\bt, \ba) \in T'_n \times \cA_n$.  
Moreover we have two families of moduli spaces 
as in (\ref{eq:family-con}) and (\ref{eq:family-rep-1}).  
(Note that we have already pulled back the family in (1) by the finite covering $T'_n \lra T_n $.)

Next, 
 let us assume that  eigenvalues of $\res_{t_i}(\nabla_L)$ are  
 integers for all $1 \leq i \leq n$. 
Then  we can 
define the {\em Riemann-Hilbert correspondence} ${\bf RH}_n :
M_n^{\balpha}(L) \lra \cR_n$ such that  the following diagram commutes:

\begin{equation}\label{eq:rh-intro}
\begin{CD} 
M_n^{\balpha}(L) @>{\bf RH}_n>>   \cR_n  \\
 @V \pi_n VV   @VV \phi_n V  \\
 T'_n \times \Lambda_n   @> (1 \times \mu_n) >>     T'_n \times \ \cA_n.  \\  
\end{CD}
\end{equation}

Here, the map $ 1 \times \mu_n$ in the bottom row  in (\ref{eq:rh-intro}) 
is given by the map $
(1 \times \mu_n)( \bt, \blambda) = (\bt, \ba)$ where 
\begin{equation}\label{eq:exponent-rel-intro}
\fbox{ \ $a_i = 2 \cos 2 \pi  \lambda_i $ \ }  \quad \mbox{ for $ 1 \leq i \leq n$}.
\end{equation}
 Under these relations,  ${\bf RH}_n$ induces the 
analytic morphism of the fibers for each $(\bt, \blambda) \in T_n' \times \Lambda_n$:
\begin{equation}\label{eq:fiber-cor}
{\bf RH}_{\bt, \blambda}: M_{n}^{\balpha}(\bt, \blambda, L) \lra \cR(\cP_{n, \bt})_{\ba}. 
\end{equation} 
To define the correspondence, 
take a stable $(\bt, \blambda)$-parabolic 
connection  $(E, \nabla, \varphi, \{ l_i \}) $. 
Then restricting the connection $ \nabla $ to 
$ \BP^1_{\C} \setminus D(\bt)$, 
define the local system on $\BP^1_{\C} \setminus D(\bt)$ by
\begin{equation}\label{eq:localsystem}
{\bf E}(\nabla)  := \ker \left(\nabla_{|\BP^1_{\C} \setminus D(\bt) }\right)^{an}.
\end{equation}
(Here $\left(\nabla_{|\BP^1_{\C} \setminus D(\bt) }\right)^{an}$ denotes 
the analytic connection associated to $\nabla_{|\BP^1_{\C} \setminus D(\bt) }$). 
Then it is easy to see that 
the map $(E, \nabla, \varphi, \{l_i\}) \mapsto {\bf E}(\nabla)$ 
induces the correspondence in (\ref{eq:rh-intro}) or (\ref{eq:fiber-cor}).  
Basically, our framework for understanding the Painlev\'e or Garnier equations 
is the Riemann-Hilbert correspondences in  (\ref{eq:rh-intro}) and (\ref{eq:fiber-cor}).

There exists one more thing  which we should mention here.  
Let $\beta_1, \beta_2$ be positive integers, 
$\balpha' = (\alpha'_1, \ldots, \alpha'_{2n})$ a  series of 
rational numbers with $0 \leq \alpha'_1 < \ldots < \alpha'_{2n} < 1 $ and set
$\bbeta = (\beta_1, \beta_2)$. 
Setting $\balpha = \balpha' \frac{\beta_1}{\beta_1 + \beta_2}$, we obtain a 
weight $\balpha$ for $(\bt, \blambda)$-parabolic connections.  (Note that 
since $\alpha_{2n} = \alpha'_{2n} \frac{\beta_1}{\beta_1+ \beta_2} < \frac{\beta_1}{\beta_1+ \beta_2}$, 
this gives a restriction for the weight $\balpha$).  For the weight $\balpha$, we consider the family of 
moduli spaces  $M_{n}^{\balpha}(L) \lra T'_n \times \Lambda_n$.  
On the other hand,  we will introduce the notion of 
$(\balpha', \bbeta)$-stable $(\bt, \blambda)$-parabolic $\phi$-connection which  is a  
generalization of  $\balpha$-stable 
$(\bt, \blambda)$-parabolic connections.  The moduli space 
$\overline{M_n^{\balpha' \bbeta}}(\bt, \blambda, L) $ contains the moduli space 
$M_n^{\balpha}(\bt, \blambda, L) $   as a Zariski  open set.  
Moreover we can construct the family of the moduli spaces such that 
the following diagram commutes:
\begin{equation}\label{eq:compact-intro}
\begin{CD} 
M_n^{\balpha}(L) & \stackrel{\iota}{\hookrightarrow} & \overline{M_n^{\balpha'\bbeta}}(L)   \\
 @V \pi_n VV   @VV \overline{\pi_n} V  \\
 T'_n \times \Lambda_n  @=    T'_n \times \Lambda_n.  \\  
\end{CD}
\end{equation}

\subsection{Main Results}
\quad

In the framework as above, 
we can state  our main results in this paper as follows.

%%%%%%%%%%%%%%%%%% MAIN THEOREM %%%%%%%%%%%%%%%%%%%%%%%%%%%%%%%%%%%%%%%%%%%%%%

\subsubsection{Projectivity of the moduli space   
$\overline{M_n^{\balpha'\bbeta}}(\bt, \blambda, L)$,  Smoothness,  Irreducibility and the Symplectic Structure 
of $M_n^{\balpha}(\bt, \blambda, L)$}

We  first  prove that  the moduli space  
$\overline{M_n^{\balpha'\bbeta}}(\bt, \blambda, L)$ is a projective scheme.  
Moreover one can show that  the moduli space  
$M_n^{\balpha}(\bt, \blambda, L)$ for each $(\bt, \blambda) \in T_n \times \Lambda_n$ 
is smooth and endowed with a natural intrinsic symplectic structure induced by Serre duality of  tangent complexes. 
The irreduciblity of the moduli space $M_n^{\balpha}(\bt, \blambda, L)$ for each 
 $(\bt, \blambda) \in T_n \times \Lambda_n$ follows from the irreduciblity of $\cR(\cP_{n, \bt})_{\ba}$ 
 via the Riemann-Hilbert correspondece (\ref{eq:fiber-cor}).

\begin{Theorem}\label{thm:fund-intro} $($Cf. Theorem \ref{thm:fund},  Theorem \ref{thm:exist-proj}, 
Proposition \ref{prop:smoothness} and 
Proposition \ref{prop:stable-mod-irr}$)$. 
\begin{enumerate}
\item For a generic  weight $(\balpha', \bbeta)$, 
$\overline{\pi}_n:\overline{M_n^{\balpha'\bbeta}}(L)  \lra T'_n \times \Lambda_n$
is a {\em projective} morphism. In particular, the moduli space 
$\overline{M_n^{\balpha'\bbeta}}(\bt, \blambda, L)$ is a {\em projective} 
algebraic scheme for all 
$(\bt, \blambda) \in T'_n \times \Lambda_n$.
\item For a generic  weight $\balpha$, 
$\pi_n:M_n^{\balpha}(L)  \lra T'_n \times \Lambda_n$
is a {\em smooth morphism} of relative dimension $2n-6$ with irreducible closed fibers. 
Therefore, the moduli space 
$M_n^{\balpha}(\bt, \blambda, L)$ is a {\em smooth, irreducible} algebraic 
variety of dimension $2n -6$  for all 
$(\bt, \blambda) \in T'_n \times \Lambda_n$.  
\end{enumerate}
\end{Theorem}

\begin{Theorem}\label{thm:symplectic-intro}
 $($Cf. Proposition \ref{prop:symplectic}$)$. 
 There exists a  global relative 2-form
\begin{equation}\label{eq:symplectic-intro}
\Omega \in H^0(M^{\balpha}_n(L), \Omega^2_{M^{\balpha}_n(L)/T_n\times\Lambda_n} ). 
\end{equation}
which induces a symplectic structure on each fiber of $\pi_n $. Consequently, 
for each $(\bt, \blambda)$, the moduli space $M_n^{\balpha}(\bt, \blambda, L)$ 
becomes a smooth  symplectic algebraic variety.  
\end{Theorem}

%%%%%%%%%%%%%%%%%% END OF  MAIN THEOREM %%%%%%%%%%%%%%%%%%%%%%%%%%%%%%%%%%%%%%%%%%%%%%

\subsubsection{Irreducibility,  symplectic structure and singularities  
of $\cR(\cP_{n, \bt})_{\ba}$}

\quad

Let us  call a data $\blambda \in \Lambda_n$  a set of 
local exponents of  connections.  

\begin{Definition}\label{def:exponents-intro}
{\rm 
\begin{enumerate}
\item A set of local exponents 
 $ \blambda =(\lambda_1, \ldots, \lambda_n) \in \Lambda_n $ is said to be  {\em special} if  
\begin{enumerate}
\item  $\blambda$ is {\em resonant}, that is, for some $ 1 \leq i \leq n$, 
\begin{equation} \label{eq:special1-intro}
  2 \lambda_i \in \Z,
\end{equation}  
\item or $\blambda$ is {\em reducible}, that is,  
for some $ (\epsilon_1, \ldots, \epsilon_n) \in \{ \pm 1 \}^n $ 
\begin{equation} \label{eq:special2-intro}
\sum_{i=1}^{n} \epsilon_i  \lambda_i  \in \Z 
\end{equation}
\end{enumerate}
\item If $\blambda \in \Lambda_n$ is not special, $\blambda$ is said to be 
{\em generic}. 
\item The data $\ba= (a_1, \ldots, a_n)  \in \cA_n$ is said to be 
{\em special} if   $\mu_n(\blambda) = \ba$ for some special $\blambda 
\in \Lambda_n$.   
\end{enumerate}}
\end{Definition}

For a monodromy representation $\rho: \pi_1(\BP^1 \setminus D(\bt), *) \lra SL_2(\C)$, 
set $M_i = \rho(\gamma_i) \in SL_2(\C) $ for $ 1 \leq i \leq n$.  
We consider the following conditions which are  invariant under the adjoint action of 
$SL_2(\C)$.   
\begin{equation}\label{eq:smooth-cond1-intro}
\mbox{The representation $\rho$ is irreducible.} 
\end{equation}
\begin{equation}\label{eq:smooth-cond2-intro}
\mbox{ For all $i, 1 \leq i \leq n$, the local monodromy matrix $M_i$ around 
$t_i$ is {\em not} equal  to $\pm I_2$.}  
\end{equation}

\begin{Theorem} \label{thm:moduli-rep-intro} 
$($Cf.\ Proposition \ref{prop:irr-p2}, Proposition \ref{prop:symplectic-rep-1} and 
Theorem \ref{thm:RH} $)$. 
Assume that $ n\geq 4$.  
\begin{enumerate}
\item 
For any $\ba \in \cA_n$, the moduli space  
$\cR(\cP_{n, \bt})_{\ba}$ is an irreducible affine scheme.

\item Let  $\cR(\cP_{n, \bt})^{\sharp}_{\ba}$ be the Zariski dense open
subset of $\cR(\cP_{n, \bt})_{\ba}$ whose closed points satisfy the conditions 
(\ref{eq:smooth-cond1-intro}) and (\ref{eq:smooth-cond2-intro}). Then $\cR(\cP_{n, \bt})^{\sharp}_{\ba}$
is smooth and there exists a natural symplectic form $\Omega_1 $ on $\cR(\cP_{n, \bt})^{\sharp}_{\ba}$.  

\item    The codimension of the locus 
 $\cR(\cP_{n, \bt})_{\ba}^{sing} := \cR(\cP_{n, \bt})_{\ba} \setminus   \cR(\cP_{n, \bt})^{\sharp}_{\ba}$ 
 is at least $2$.
\end{enumerate}
\end{Theorem}

\subsubsection{Surjectivity and Properness of the Riemann-Hilbert correspondence}
\quad

Next, the most important result for the Riemann-Hilbert correspondence is the 
surjectivity and the properness.  One can show that the correspondence
${\bf RH}_{\bt,\blambda}$ 
in (\ref{eq:fiber-cor}) 
gives an analytic isomorphism between two moduli spaces
if $\blambda \in \Lambda_n$ is generic (i.e. non-special). 
However, for a special $\blambda \in \Lambda_n$, one can see that 
the map  (\ref{eq:fiber-cor}) contracts some subvarieties of  
$M_n^{\balpha}(\bt, \blambda, L)$ to singular locus  of $\cR(\cP_{n, \bt})_{\ba}$.  Note that 
since the correspondence is not an algebraic morphism, 
one  can not directly apply the valuative criterion for the proof of the properness.

\begin{Theorem} \label{thm:RH-intro} $($Cf. Theorem \ref{thm:RH}$)$.  
Under the notation above and assume that $n \geq 4$ and 
$\balpha$ is general.  For all $(\bt, \blambda) \in T'_n \times \Lambda_n$, 
the Riemann-Hilbert correspondence 
\begin{equation} \label{eq:rh-surj}
\RH_{\bt,\blambda}: M_n^{\balpha}(\bt, \blambda, L) \lra \cR(\cP_{n, \bt})_{\ba}
\end{equation} 
is a bimeromorphic proper  surjective morphism.   
\end{Theorem}

\subsubsection{The Riemann-Hilbert correspondence as a symplectic resolution of singularities of 
$\cR(\cP_{n, \bt})_{\ba} $}

\quad

Moreover, we can introduce the natural intrinsic symplectic structure on  the 
smooth part $\cR(\cP_{n, \bt})^{\sharp}_{\ba}$ of the moduli spaces $\cR(\cP_{n, \bt})_{\ba}$.  
Together with the 
natural symplectic structure of the moduli space $M_n^{\balpha}(\bt, \blambda, L)$, 
the map ${\bf RH}_{\bt, \blambda}$ gives a symplectic map, which means that 
the pullback of the symplectic structure on the smooth part of 
$\cR(\cP_{n, \bt})_{\ba}$ coincides with the symplectic structure on $M_n^{\balpha}(\bt, \blambda, L)$. 
This identification will be given by a kind of 
infinitesimal Riemann-Hilbert correspondence (cf. Lemma \ref{lem:inf-rh}).  
  Together with the surjectivity,  the properness of 
${\bf RH}_{\bt, \blambda}$ and the fact that 
$M_n^{\balpha}(\bt, \blambda, L)$ is smooth, 
we can say that ${\bf RH}_{\bt, \blambda}$ gives 
an {\em analytic} symplectic resolution of the singularities of 
$\cR(\cP_{n, \bt})_{\ba}$.  Moreover, we can say that the map 
${\bf RH}_n $ in \eqref{eq:rh-intro}  gives a {\em simultaneous 
resolution of the family}  $\phi_n:\cR_n \lra T'_n \times \cA_n$ with 
the base extension $1 \times \mu_n :T'_n \times \Lambda_n 
\lra T'_n \times \cA_n$.  (For definition, see [Definition 4.26, \cite{KM}]).

\begin{Theorem} \label{thm:RH2-intro}  $($Theorem \ref{thm:RH}, Lemma 
\ref{lem:inf-rh}$)$. 
Under the assumption of Theorem (\ref{thm:RH-intro}), we have the following.

\begin{enumerate}
\item 
For any $(\bt, \blambda)$, let $\cR(\cP_{n, \bt})^{\sharp}_{\ba}$ be as in Theorem \ref{thm:moduli-rep-intro}, 
and set $M_n^{\balpha}(\bt, \blambda, L)^{\sharp} = 
\RH_{\bt, \blambda}^{-1}( \cR(\cP_{n, \bt})^{\sharp}_{\ba})$. 
Then the Riemann-Hilbert correspondence gives an analytic 
isomorphism 
\begin{equation}\label{eq:restricted-rh} 
{\bf RH}_{\bt,\blambda, | M_n^{\balpha}(\bt, \blambda, L)^{\sharp}} : 
M_n^{\balpha}(\bt, \blambda, L)^{\sharp}  \stackrel{\simeq}{\lra} \cR(\cP_{n, \bt})^{\sharp}_{\ba} .
\end{equation}  
(Note that if $\blambda$ is not special (cf. Definition \ref{def:exponents}, (\ref{eq:special1}),  
(\ref{eq:special2})),    
$\cR(\cP_{n, \bt})^{\sharp}_{\ba} = \cR(\cP_{n, \bt})_{\ba}$, hence $\RH_{\bt, \blambda}$ 
gives an analytic isomorphism between 
$M_n^{\balpha}(\bt, \blambda, L)$ and $\cR(\cP_{n, \bt})_{\ba}$.) 

\item The symplectic structures  $\Omega$ restricted to 
$M_n^{\balpha}(\bt, \blambda, L)^{\sharp}$  and 
$ \Omega_1$ on $\cR(\cP_{n, \bt})^{\sharp}_{\ba}$ can be identified 
with each other 
via $\RH_{\bt, \blambda}$, that is, 
\begin{equation} \label{eq:pull-back-intro}
\Omega_{|M_n^{\balpha}(\bt, \blambda, L)^{\sharp}} = 
\RH_{\bt, \blambda | M_n^{\balpha}(\bt, \blambda, L)^{\sharp}}^*(\Omega_1) \quad 
\mbox{on}  \ \  M_n^{\balpha}(\bt, \blambda, L)^{\sharp}. 
\end{equation} 
\item Putting together all results, the correspondence ${\bf RH}_n$ in (\ref{eq:rh-intro}) 
gives an analytic simultaneous symplectic resolution of 
singularities after the base extension 
$1 \times \mu_n : T'_n \times \Lambda_n \ra T'_n \times \cA_n$.  
\end{enumerate}
\end{Theorem}

\subsection{Painlev\'e and Garnier equations and their Painlev\'e property}
\quad

In the framework of this paper, 
we can derive the Painlev\'e and Garnier equations as follows.  
Take the universal covering map  $ \tilde{T_n} \lra   T'_n \lra T_n $ and pull back the diagram 
(\ref{eq:rh-intro})  to obtain the following commutative diagram:  
\begin{equation}\label{eq:rh-un-intro}
\begin{CD} 
\tilde{M}_n^{\balpha}(L) 
@>{\bf RH}_n>>   \tilde{\cR}_n  \\
 @V \tilde{\pi}_n VV   @VV \tilde{\phi}_n V  \\
 \tilde{T_n} \times \Lambda_n   @> (1 \times \mu_n) >>  \tilde{T_n} \times \ \cA_n.  \\  
\end{CD}
\end{equation}

\subsubsection{The case of  generic $\blambda$ }
\quad

Now let us fix  $\blambda \in \Lambda_n$  and  set $\ba = \mu_n(\blambda) $.  
First, assume that $\blambda$ is {\em generic}.   We denote by 
$ (\pi_n)_{\blambda}: \tilde{M}_n^{\balpha}(\blambda, L) \lra \tilde{T}_n $ 
and 
$(\phi_n)_{\ba}: (\tilde{\cR}_n)_{\ba}  \lra \tilde{T}_n$ the families 
obtained by  restricting the families in (\ref{eq:rh-un-intro}) 
to $\tilde{T}_n\times \{\blambda\}$  and $\tilde{ T}_n \times \{ \ba \}$.  
Moreover we denote by 
$
{\bf RH}_{\blambda}: \tilde{M}_n^{\balpha}(\blambda, L) 
\lra (\tilde{\cR}_n)_{\ba} 
$ 
the restriction of ${\bf RH}_n$ to the restricted families.  
Since $\blambda$ is generic, 
${\bf RH}_{\blambda}$ 
induces  an analytic isomorphism between 
$\tilde{M}_n^{\balpha}(\blambda, L)$ and $ (\tilde{\cR}_n)_{\ba} $.  
Fix a point $\bt_0 \in T'_n$.  Since the original fibration 
$(\phi_n)_{\ba}: (\cR_n)_{\ba}  \lra T'_n \times \{ \ba  \}$ is locally trivial, 
we can obtain an isomorphism 
$ (\tilde{\cR}_n)_{\ba} \simeq \cR(\cP_{n, \bt_0})_{\ba}  \times \tilde{T}_n$ and 
the following commutative diagram 
for fixed $\blambda$ and $\ba$.  
\begin{equation}\label{eq:rh2-un-intro}
\begin{CD} 
\tilde{M}_n^{\balpha}(\blambda, L) 
@>{\bf RH}_{\blambda}   >\simeq>   (\tilde{\cR})_{\ba}   & 
\simeq \cR(\cP_{n, \bt_0})_{\ba} \times \tilde{T}_n \\
 @V (\tilde{\pi}_n)_{\blambda} VV   @VV (\tilde{\phi}_n)_{\ba} V  \\
 \tilde{T}_n \times \{ \blambda \}   @> = >>   
 \tilde{T}_n \times \ \{ \ba \} . \\
\end{CD}
\end{equation}
By using this global trivialization, {for} each closed point 
$\x \in  \cR(\cP_{n, \bt_0})_{\ba} $, we can define the unique constant section 
$s_{\x}: \tilde{T}_n  \lra \cR(\cP_{n, \bt_0})_{\ba}  \times \tilde{T}_n$ for $(\phi_n)_{\ba}$
by the formula 
 $s_{\x}(\bt) = (\x, \bt)$.  Pulling back this constant section $s_{\x}$ via 
 ${\bf RH}_{\blambda}$  we obtain the global analytic section $ \tilde{s}_{\x}$ 
  for  the morphism  $(\pi_n)_{\blambda}$.  
  Varying the initial points $\x$, we obtain 
  the family of  constant sections 
  $\{ s_{\x}\}_{\x \in \cR(\cP_{n, \bt_0})_{\ba}}$  
  of  
  $\cR(\cP_{n, \bt_0})_{\ba}  \times \tilde{T}_n \lra \tilde{T}_n $ 
  and also the family  of pullback sections 
  $\{ \tilde{s}_{\x} \}_{\x \in \cR(\cP_{n, \bt_0})_{\ba}}$ 
  for $\tilde{M}_n^{\balpha}(\blambda, L)$.  
  
\begin{figure}
%WinTpicVersion2.15
\unitlength 0.1in
\begin{picture}(62.10,41.10)(0.90,-43.40)
% VECTOR 1 0 3 0
% 2 3594 3022 5386 3030
% 
\special{pn 13}%
\special{pa 3594 2622}%
\special{pa 5386 2630}%
\special{fp}%
\special{sh 1}%
\special{pa 5386 2630}%
\special{pa 5319 2610}%
\special{pa 5333 2630}%
\special{pa 5319 2650}%
\special{pa 5386 2630}%
\special{fp}%
% POLYGON 1 0 3 0
% 6 3562 1406 3570 2454 4346 2846 4346 2838 4346 1734 3562 1406
% 
\special{pn 13}%
\special{pa 3562 1006}%
\special{pa 3570 2054}%
\special{pa 4346 2446}%
\special{pa 4346 2438}%
\special{pa 4346 1334}%
\special{pa 3562 1006}%
\special{fp}%
% POLYGON 1 0 3 0
% 6 4602 1398 4610 2446 5386 2838 5386 2830 5386 1726 4602 1398
% 
\special{pn 13}%
\special{pa 4602 998}%
\special{pa 4610 2046}%
\special{pa 5386 2438}%
\special{pa 5386 2430}%
\special{pa 5386 1326}%
\special{pa 4602 998}%
\special{fp}%
% POLYGON 1 0 3 0
% 6 730 1406 738 2454 1514 2846 1514 2838 1514 1734 730 1406
% 
\special{pn 13}%
\special{pa 730 1006}%
\special{pa 738 2054}%
\special{pa 1514 2446}%
\special{pa 1514 2438}%
\special{pa 1514 1334}%
\special{pa 730 1006}%
\special{fp}%
% POLYGON 1 0 3 0
% 6 1770 1398 1778 2446 2554 2838 2554 2830 2554 1726 1770 1398
% 
\special{pn 13}%
\special{pa 1770 998}%
\special{pa 1778 2046}%
\special{pa 2554 2438}%
\special{pa 2554 2430}%
\special{pa 2554 1326}%
\special{pa 1770 998}%
\special{fp}%
% DOT 0 0 3 0
% 2 3986 3038 3986 3030
% 
\special{pn 20}%
\special{sh 1}%
\special{ar 3986 2638 10 10 0  6.28318530717959E+0000}%
\special{sh 1}%
\special{ar 3986 2630 10 10 0  6.28318530717959E+0000}%
% STR 2 0 3 0
% 3 3834 3150 3834 3230 2 0
% $\bt_0$
\put(38.3400,-28.3000){\makebox(0,0)[lb]{$\bt_0$}}%
% STR 2 0 3 0
% 3 5090 3140 5090 3220 2 0
% $\bt$
\put(50.9000,-28.2000){\makebox(0,0)[lb]{$\bt$}}%
% DOT 0 0 3 0
% 2 5058 3038 5058 3030
% 
\special{pn 20}%
\special{sh 1}%
\special{ar 5058 2638 10 10 0  6.28318530717959E+0000}%
\special{sh 1}%
\special{ar 5058 2630 10 10 0  6.28318530717959E+0000}%
% STR 2 0 3 0
% 3 1082 3166 1082 3246 2 0
% $\bt_0$
\put(10.8200,-28.4600){\makebox(0,0)[lb]{$\bt_0$}}%
% STR 2 0 3 0
% 3 2300 3160 2300 3240 2 0
% $\bt$
\put(23.0000,-28.4000){\makebox(0,0)[lb]{$\bt$}}%
% VECTOR 1 0 3 0
% 2 850 3006 2642 3014
% 
\special{pn 13}%
\special{pa 850 2606}%
\special{pa 2642 2614}%
\special{fp}%
\special{sh 1}%
\special{pa 2642 2614}%
\special{pa 2575 2594}%
\special{pa 2589 2614}%
\special{pa 2575 2634}%
\special{pa 2642 2614}%
\special{fp}%
% DOT 0 0 3 0
% 2 1242 3022 1242 3014
% 
\special{pn 20}%
\special{sh 1}%
\special{ar 1242 2622 10 10 0  6.28318530717959E+0000}%
\special{sh 1}%
\special{ar 1242 2614 10 10 0  6.28318530717959E+0000}%
% DOT 0 0 3 0
% 2 2314 3022 2314 3014
% 
\special{pn 20}%
\special{sh 1}%
\special{ar 2314 2622 10 10 0  6.28318530717959E+0000}%
\special{sh 1}%
\special{ar 2314 2614 10 10 0  6.28318530717959E+0000}%
% VECTOR 1 0 3 0
% 2 2850 1926 3298 1926
% 
\special{pn 13}%
\special{pa 2850 1526}%
\special{pa 3298 1526}%
\special{fp}%
\special{sh 1}%
\special{pa 3298 1526}%
\special{pa 3231 1506}%
\special{pa 3245 1526}%
\special{pa 3231 1546}%
\special{pa 3298 1526}%
\special{fp}%
% VECTOR 1 0 3 0
% 2 3714 1878 5082 1878
% 
\special{pn 13}%
\special{pa 3714 1478}%
\special{pa 5082 1478}%
\special{fp}%
\special{sh 1}%
\special{pa 5082 1478}%
\special{pa 5015 1458}%
\special{pa 5029 1478}%
\special{pa 5015 1498}%
\special{pa 5082 1478}%
\special{fp}%
% VECTOR 1 0 3 0
% 2 3730 2150 5098 2150
% 
\special{pn 13}%
\special{pa 3730 1750}%
\special{pa 5098 1750}%
\special{fp}%
\special{sh 1}%
\special{pa 5098 1750}%
\special{pa 5031 1730}%
\special{pa 5045 1750}%
\special{pa 5031 1770}%
\special{pa 5098 1750}%
\special{fp}%
% VECTOR 1 0 3 0
% 2 3770 2414 5138 2414
% 
\special{pn 13}%
\special{pa 3770 2014}%
\special{pa 5138 2014}%
\special{fp}%
\special{sh 1}%
\special{pa 5138 2014}%
\special{pa 5071 1994}%
\special{pa 5085 2014}%
\special{pa 5071 2034}%
\special{pa 5138 2014}%
\special{fp}%
% SPLINE 1 0 3 0
% 6 970 1894 1274 1806 1394 1822 1658 1982 2122 1814 2122 1814
% 
\special{pn 13}%
\special{pa 970 1494}%
\special{pa 1000 1481}%
\special{pa 1030 1469}%
\special{pa 1060 1457}%
\special{pa 1091 1445}%
\special{pa 1121 1435}%
\special{pa 1152 1426}%
\special{pa 1183 1418}%
\special{pa 1214 1412}%
\special{pa 1246 1408}%
\special{pa 1278 1406}%
\special{pa 1310 1406}%
\special{pa 1342 1409}%
\special{pa 1373 1416}%
\special{pa 1404 1426}%
\special{pa 1433 1439}%
\special{pa 1462 1455}%
\special{pa 1490 1474}%
\special{pa 1517 1493}%
\special{pa 1543 1513}%
\special{pa 1570 1532}%
\special{pa 1596 1550}%
\special{pa 1623 1566}%
\special{pa 1649 1579}%
\special{pa 1677 1588}%
\special{pa 1704 1593}%
\special{pa 1733 1595}%
\special{pa 1762 1593}%
\special{pa 1791 1588}%
\special{pa 1820 1581}%
\special{pa 1850 1571}%
\special{pa 1881 1558}%
\special{pa 1911 1544}%
\special{pa 1942 1528}%
\special{pa 1973 1510}%
\special{pa 2004 1491}%
\special{pa 2036 1471}%
\special{pa 2067 1451}%
\special{pa 2099 1430}%
\special{pa 2122 1414}%
\special{sp}%
% VECTOR 1 0 3 0
% 2 2130 1806 2202 1758
% 
\special{pn 13}%
\special{pa 2130 1406}%
\special{pa 2202 1358}%
\special{fp}%
\special{sh 1}%
\special{pa 2202 1358}%
\special{pa 2135 1378}%
\special{pa 2158 1388}%
\special{pa 2158 1412}%
\special{pa 2202 1358}%
\special{fp}%
% STR 2 0 3 0
% 3 2620 2680 2620 2760 2 0
% $M_n^{\balpha}(\bt, \blambda, L)$
\put(26.2000,-23.6000){\makebox(0,0)[lb]{$M_n^{\balpha}(\bt, \blambda, L)$}}%
% STR 2 0 3 0
% 3 5470 2290 5470 2370 2 0
% $\cR(\cP_{n, \bt_0})_{\ba} $
\put(54.7000,-19.7000){\makebox(0,0)[lb]{$\cR(\cP_{n, \bt_0})_{\ba} $}}%
% STR 2 0 3 0
% 3 3090 1660 3090 1740 5 0
% ${\bf RH}_{\lambda}$
\put(30.9000,-13.4000){\makebox(0,0){${\bf RH}_{\lambda}$}}%
% STR 2 0 3 0
% 3 4720 1060 4720 1140 5 0
% constant flows = monodromy is consant
\put(47.2000,-7.4000){\makebox(0,0){constant flows = monodromy is constant}}%
% STR 2 0 3 0
% 3 200 1130 200 1210 2 0
% Isomondromic flows = Painlev\'{e} or Garnier flows
\put(2.0000,-8.1000){\makebox(0,0)[lb]{Isomonodromic flows = Painlev\'{e} or Garnier flows}}%
% STR 2 0 3 0
% 3 3020 2990 3020 3090 2 0
% $=$
\put(30.2000,-26.9000){\makebox(0,0)[lb]{$=$}}%
% VECTOR 2 0 3 0
% 2 1640 1280 1630 1920
% 
\special{pn 8}%
\special{pa 1640 880}%
\special{pa 1630 1520}%
\special{fp}%
\special{sh 1}%
\special{pa 1630 1520}%
\special{pa 1651 1454}%
\special{pa 1631 1467}%
\special{pa 1611 1453}%
\special{pa 1630 1520}%
\special{fp}%
% VECTOR 2 0 3 0
% 2 4440 1310 4450 1810
% 
\special{pn 8}%
\special{pa 4440 910}%
\special{pa 4450 1410}%
\special{fp}%
\special{sh 1}%
\special{pa 4450 1410}%
\special{pa 4469 1343}%
\special{pa 4449 1357}%
\special{pa 4429 1344}%
\special{pa 4450 1410}%
\special{fp}%
% STR 2 0 3 0
% 3 5430 1880 5430 1960 2 0
% $\cR(\cP_{n, \bt})_{\ba} $
\put(54.3000,-15.6000){\makebox(0,0)[lb]{$\cR(\cP_{n, \bt})_{\ba} $}}%
% STR 2 0 3 0
% 3 5710 2070 5710 2170 2 0
% $||$
\put(57.1000,-17.7000){\makebox(0,0)[lb]{$||$}}%
% STR 2 0 3 0
% 3 180 2660 180 2740 2 0
% $M_n^{\balpha}(\bt_0, \blambda, L)$
\put(1.8000,-23.4000){\makebox(0,0)[lb]{$M_n^{\balpha}(\bt_0, \blambda, L)$}}%
% STR 2 0 3 0
% 3 1600 3150 1600 3250 2 0
% $\tilde{T}_n \times \{\blambda \} $
\put(16.0000,-28.5000){\makebox(0,0)[lb]{$\tilde{T}_n \times \{\blambda \} $}}%
% STR 2 0 3 0
% 3 4360 3160 4360 3260 2 0
% $\tilde{T}_n \times \{ \ba \}$
\put(43.6000,-28.6000){\makebox(0,0)[lb]{$\tilde{T}_n \times \{ \ba \}$}}%
% SPLINE 1 0 3 0
% 6 1030 2190 1334 2102 1454 2118 1718 2278 2182 2110 2182 2110
% 
\special{pn 13}%
\special{pa 1030 1790}%
\special{pa 1060 1777}%
\special{pa 1090 1765}%
\special{pa 1120 1753}%
\special{pa 1151 1741}%
\special{pa 1181 1731}%
\special{pa 1212 1722}%
\special{pa 1243 1714}%
\special{pa 1274 1708}%
\special{pa 1306 1704}%
\special{pa 1338 1702}%
\special{pa 1370 1702}%
\special{pa 1402 1705}%
\special{pa 1433 1712}%
\special{pa 1464 1722}%
\special{pa 1493 1735}%
\special{pa 1522 1751}%
\special{pa 1550 1770}%
\special{pa 1577 1789}%
\special{pa 1603 1809}%
\special{pa 1630 1828}%
\special{pa 1656 1846}%
\special{pa 1683 1862}%
\special{pa 1709 1875}%
\special{pa 1737 1884}%
\special{pa 1764 1889}%
\special{pa 1793 1891}%
\special{pa 1822 1889}%
\special{pa 1851 1884}%
\special{pa 1880 1877}%
\special{pa 1910 1867}%
\special{pa 1941 1854}%
\special{pa 1971 1840}%
\special{pa 2002 1824}%
\special{pa 2033 1806}%
\special{pa 2064 1787}%
\special{pa 2096 1767}%
\special{pa 2127 1747}%
\special{pa 2159 1726}%
\special{pa 2182 1710}%
\special{sp}%
% VECTOR 1 0 3 0
% 2 2180 2100 2252 2052
% 
\special{pn 13}%
\special{pa 2180 1700}%
\special{pa 2252 1652}%
\special{fp}%
\special{sh 1}%
\special{pa 2252 1652}%
\special{pa 2185 1672}%
\special{pa 2208 1682}%
\special{pa 2208 1706}%
\special{pa 2252 1652}%
\special{fp}%
% SPLINE 1 0 3 0
% 6 1050 2412 1354 2324 1474 2340 1738 2500 2202 2332 2202 2332
% 
\special{pn 13}%
\special{pa 1050 2012}%
\special{pa 1080 1999}%
\special{pa 1110 1987}%
\special{pa 1140 1975}%
\special{pa 1171 1963}%
\special{pa 1201 1953}%
\special{pa 1232 1944}%
\special{pa 1263 1936}%
\special{pa 1294 1930}%
\special{pa 1326 1926}%
\special{pa 1358 1924}%
\special{pa 1390 1924}%
\special{pa 1422 1927}%
\special{pa 1453 1934}%
\special{pa 1484 1944}%
\special{pa 1513 1957}%
\special{pa 1542 1973}%
\special{pa 1570 1992}%
\special{pa 1597 2011}%
\special{pa 1623 2031}%
\special{pa 1650 2050}%
\special{pa 1676 2068}%
\special{pa 1703 2084}%
\special{pa 1729 2097}%
\special{pa 1757 2106}%
\special{pa 1784 2111}%
\special{pa 1813 2113}%
\special{pa 1842 2111}%
\special{pa 1871 2106}%
\special{pa 1900 2099}%
\special{pa 1930 2089}%
\special{pa 1961 2076}%
\special{pa 1991 2062}%
\special{pa 2022 2046}%
\special{pa 2053 2028}%
\special{pa 2084 2009}%
\special{pa 2116 1989}%
\special{pa 2147 1969}%
\special{pa 2179 1948}%
\special{pa 2202 1932}%
\special{sp}%
% VECTOR 1 0 3 0
% 2 2200 2322 2272 2274
% 
\special{pn 13}%
\special{pa 2200 1922}%
\special{pa 2272 1874}%
\special{fp}%
\special{sh 1}%
\special{pa 2272 1874}%
\special{pa 2205 1894}%
\special{pa 2228 1904}%
\special{pa 2228 1928}%
\special{pa 2272 1874}%
\special{fp}%
% STR 2 0 3 0
% 3 3020 2010 3020 2110 2 0
% $\simeq$
\put(30.2000,-17.1000){\makebox(0,0)[lb]{$\simeq$}}%
% LINE 2 5 3 0
% 2 3090 4450 4830 4450
% 
\special{pn 8}%
\special{pa 3090 4050}%
\special{pa 4830 4050}%
\special{ip}%
% BOX 1 0 3 0
% 2 6300 4740 90 630
% 
\special{pn 13}%
\special{pa 6300 4340}%
\special{pa 90 4340}%
\special{pa 90 230}%
\special{pa 6300 230}%
\special{pa 6300 4340}%
\special{fp}%
% STR 2 0 3 0
% 3 1520 3820 1520 3920 2 0
% Isomonodromic Flows and Painlev\'e or Garnier Flows
\put(15.2000,-35.2000){\makebox(0,0)[lb]{Isomonodromic Flows and Painlev\'e or Garnier Flows}}%
\end{picture}%
\caption{Riemann-Hilbert correspondence and isomonodromic flows for generic $\blambda$}
\label{fig:isom}
\end{figure}
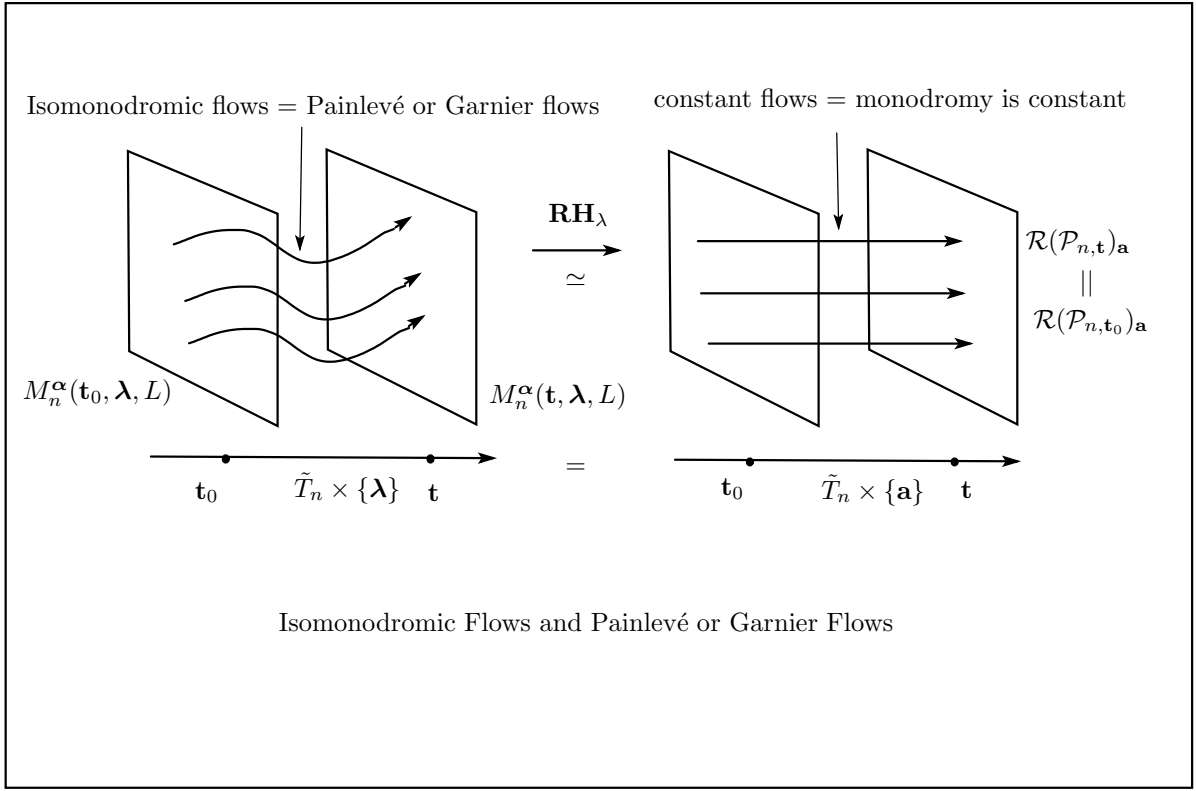

The family of sections  
$\{ \tilde{s}_{\x} \}_{\x \in \cR(\cP_{n, \bt_0})_{\ba}}$ 
gives the splitting homomorphism 
\begin{equation}\label{eq:splitting}
\tilde{v}_{\blambda} : (\pi_n)_{\blambda}^{*} ( \Theta_{\tilde{T}_n \times \{\blambda \} }) 
\lra \Theta_{\tilde{M}_n^{\balpha}(\blambda, L)} 
\end{equation}
for the natural 
surjective homomorphism 
$
\Theta_{\tilde{M}_n^{\balpha}(\blambda, L)} 
\lra (\pi_n)_{\blambda}^{*}( \Theta_{\tilde{T}_n \times \{\blambda \} })
$.
Consider the following commutative diagram:

\begin{equation}\label{eq:covering-com}
\begin{CD} 
\tilde{M}_n^{\balpha}(\blambda, L) 
@> \tilde{u}   >  >   M_n^{\balpha} (\blambda, L)   \\
 @V (\tilde{\pi}_n)_{\blambda} VV   @VV (\pi_n)_{\blambda} V  \\
 \tilde{T}_n \times \{ \blambda \}   @> u > \mbox{\tiny universal covering } >  T_n \times \ \{ \blambda \} .  \\\end{CD}
\end{equation}
We can see that the splitting homomorphism (\ref{eq:splitting})  descends to a splitting homomorphism
\begin{equation}\label{eq:descend-split}
v_{\blambda} : (\pi_n)_{\blambda}^{*} ( \Theta_{T_n \times \{\blambda \} }) 
\lra \Theta_{M_n^{\balpha}(\blambda, L)}. 
\end{equation}
(One can show that this splitting is an algebraic homomorphism). Therefore, 
each  algebraic vector field  $\theta $ on $T_n \times \{\blambda \}$ determines 
an algebraic  vector field $v_{\blambda}(\theta)$ on $M_n^{\balpha}(\blambda, L)$.  
The natural generators of the tangent sheaf of $T_n \times \{\blambda \}$ can be given by 
$$
 \langle  \frac{\partial}{\partial t_1}, \ldots,   \frac{\partial}{\partial t_n}   
 \rangle. 
$$ 
Defining 
\begin{equation}\label{eq:def-vf}
v_i(\blambda) = v_{\blambda} (\frac{\partial}{\partial t_i})  
\in H^0( M_n^{\balpha}(\blambda, L), \Theta_{M_n^{\balpha}(\blambda, L)})
\end{equation}
we obtain 
the differential system  
\begin{equation}\label{eq:diff-sys}
\langle v_1 (\blambda), \ldots, v_n(\blambda)  \rangle
\end{equation}
on $M_n^{\balpha}(\blambda, L)$.    
 From the construction, it is obvious that these 
vector fields $\{v_i (\blambda) \}_{1 \leq i \leq n }$ 
commute to each other, that is,  the differential systems are integrable. 
Since $\tilde{u}: \tilde{M}_n^{\balpha}(\blambda, L) \lra   
M_n^{\balpha} (\blambda, L)$ in (\ref{eq:covering-com}) 
is also a covering map, each section $ \tilde{s}_{\x} : \tilde{T}_n \times \{\blambda\} 
\lra \tilde{M}_n^{\balpha}(\blambda, L) $ 
defines a multi-section for $ M_n^{\balpha}(\blambda, L) \lra T_n \times \{ \blambda \}$,  which 
gives an integral submanifold of $ M_n^{\balpha}(\blambda, L)$  
for the differential system (\ref{eq:diff-sys}) at least locally. 
Hence   the  submanifold   $ \tilde{s}_{\x}( \tilde{T}_n \times \{\blambda \} )$  of 
$\tilde{M}_n^{\balpha}(\blambda, L)$  
given by the image of the section $\tilde{s}_{\x}$  can be considered as the integral submanifold (or 
a solution submanifold) for (\ref{eq:diff-sys}) over the universal covering space $\tilde{T_n}$. 
 (It is natural to call the  submanifold   $ \tilde{s}_{\x}( \tilde{T}_n \times \{\blambda \} )$    
an {\em isomonodromic flow}).  Since the integral submanifold $ \tilde{s}_{\x}( \tilde{T}_n \times \{\blambda \} )$ 
is isomorphic 
to the parameter space $ \tilde{T}_n \times \{\blambda \}$ and the morphism 
  $\tilde{\pi}_n : \tilde{M}_n^{\balpha}(\blambda, L) \lra \tilde{T}_n \times \{ \blambda \}$ 
  is algebraic, we can conculde that
\begin{equation}
\fbox{  the diffrential system  
$\{ v_i(\blambda)\}_{1 \leq i \leq n} $ on $M_n^{\balpha}(\blambda, L) $ 
has Painlev\'e property.  (See Figure \ref{fig:isom}). } 
\end{equation}
Actually, the dynamical system on $M_n^{\balpha}(\blambda, L) $ detemined by 
$\{ v_i(\blambda)\}_{1 \leq i \leq n} $ has geometric Painlev\'e property 
in the sense of \cite{IISA} (cf. [Definition 2.2, \cite{IISA}]). 
The differential system  $\{v_i(\blambda)\}_{1 \leq i \leq n} $ in (\ref{eq:diff-sys}) is 
called {\em Painlev\'e $VI$ system} for $n = 4$ and 
{\em Garnier system }  for $n \geq 5$.  (Moreover we call each vector field 
 $v_i(\blambda)$  Painlev\'e or Garnier vector field).

By using a suitable algebraic 
local coordinate system for $M_n^{\alpha}(\blambda, L)$, one can 
write down the differential equations associated to $v_i (\blambda)$ and  
see that these differential systems are equivalent to known 
Painlev\'e $VI$  systems and Garnier systems.   
(It is possible to reduce the number of the time variables $t_i$ applying the 
automorphism of $\BP^1_{\C}$  from
 $n$ to $n-3$).  
Moreover, one can apply a standard argument to show that the vector fields $v_i (\blambda)$ are  
algebraic vector fields on $M_n^{\alpha}(\blambda, L)$.  

\subsubsection{The case of special $\blambda$} 
\quad 

Next, let us consider the case when $\blambda$ is special. 
We have the same commutative diagram as (\ref{eq:rh2-un-intro}), 
however we encounter the following new phenomena. 
\begin{enumerate}
\item Although the moduli space $M_n^{\balpha}(\bt, \blambda, L)$ is nonsingular, the moduli space $\cR(\cP_{n, \bt_0})_{\ba}$ has  singularities.
\item
The Riemann-Hilbert correspondence $
{\bf RH}_{\blambda}: \tilde{M}_n^{\balpha}(\blambda, L) 
\lra (\tilde{\cR}_n)_{\ba} $ (or $
{\bf RH}_{\bt, \blambda}: \tilde{M}_n^{\balpha}(\bt, \blambda, L) 
\lra  \cR(\cP_{n, \bt_0})_{\ba} $ for a fixed $\bt$ )  is still a 
bimeromorphic proper surjective  map,  but it contracts  
some  families of compact  subvarieties to singular locus of $\cR(\cP_{n, \bt_0})_{\ba}$.  
\end{enumerate}
For example, in case when $n = 4$ (Painlev\'e $VI$  case) and $\blambda$ is special,  
$M_n^{\balpha}(\bt, \blambda, L)$ contains at least one $(-2)$-rational curve. 
For simplicity, assume that there is a unique $(-2)$-rational curve on  
$\tilde{M}_4^{\balpha}(\bt_0, \blambda, L)$. 
Since $\cR(\cP_{n, \bt_0})_{\ba}$ is an irreducible affine scheme, it 
cannot contain complete subvarieties of  positive dimension, and  
hence ${\bf RH}_{\bt_0, \blambda}$  has to contract the $(-2)$-rational 
curve onto a singular point of type $A_1$.  (See Figure \ref{fig:isom2}). 
Let us define  the subset  $\tilde{M}_n^{\balpha}(\blambda, L)^{\sharp}$ 
the complement of the subvarieties contracted by 
${\bf RH}_{\blambda}$ in 
$ \tilde{M}_n^{\balpha}(\blambda, L) $ and set 
$(\tilde{\cR}_n)_{\ba}^{\sharp}:
={\bf RH}_{\blambda}(\tilde{M}_n^{\balpha}(\blambda, L)^{\sharp} ) $
 so that 
${\bf RH}_{\blambda|\tilde{M}_n^{\balpha}(\blambda, L)^{\sharp} }:  
\tilde{M}_n^{\balpha}(\blambda, L)^{\sharp} 
\lra (\tilde{\cR}_n)_{\ba}^{\sharp}$ is an analytic isomorphism. 
For any $ n \geq 4$, 
we can pull back the constant 
sections $s_{\x}$  by 
${\bf RH}_{\blambda}$ for  $ \x \in (\tilde{\cR}_n)_{\ba}^{\sharp}  $ 
and obtain analytic sections $\tilde{s}_{\x}$, 
for $(\pi_n)_{\blambda}$.  
Now consider the family 
$ (\pi_n)_{\blambda}: M_n^{\balpha}(\blambda, L) \lra T_n \times \{\blambda\} $ 
over $T_n \times \{ \blambda \} $ and 
define $ M_n^{\balpha}(\blambda, L)^{\sharp} \subset  M_n^{\balpha}(\blambda, L)$  as above.  Then we can also obtain mutually commutative 
Pailev\'e $VI $ or Garnier  vector fields 
$v_i (\blambda)$ for $ 1 \leq i \leq n $ 
 on  $M_n^{\balpha}(\blambda, L)^{\sharp}$, and    
 $\{ v_i (\blambda) \}_{1 \leq i \leq n}$ 
defines an integrable differential system 
on $M_n^{\balpha}(\blambda, L)^{\sharp}$.   
Varying $\blambda$, we obtain the  set of algebraic  vector fields 
$\{ v_i \}_{1 \leq i \leq n} $ 
on $M_n^{\balpha}(L)^{\sharp} $ over $T_n \times \Lambda_n$. 
Since the codimension of $M_n^{\balpha}(L) \setminus M_n^{\balpha}(L)^{\sharp} $
 in  $M_n^{\balpha}(L) $ 
is greater than $2$, one  can extend the algebraic vector field $v_i $ to 
$M_n^{\balpha}(L) $.  
Hence $v_i (\blambda)$ can also be extended to the total  
space of the family of the moduli spaces $(\pi_n)_{\blambda}: M_n^{\balpha}(\blambda, L) \lra T_n \times \{ \blambda \} $.  
From the properness of the  Riemann-Hilbert correspondence 
${\bf RH}_{\blambda}: \tilde{M}_n^{\balpha}(\blambda, L) 
\lra (\tilde{\cR}_n)_{\ba} $, we can  conclude that the differential system 
$\{ v_i (\blambda) \}_{1 \leq i \leq n} $ also has the geometric Painlv\'e 
property (cf. \cite{IISA}).

The extended vector fields should be tangent to the family of contracted subvarieties (see Figure \ref{fig:isom2}). The restriction of Painlev\'e $VI$ or Garnier vector fields $\{ v_i (\blambda) \}_{1 \leq i \leq n} $ to the family of the contracted subvarieties yields integrable differential systems on the subvarieties whose solutions are given by  a family of classical solutions like Riccati solutions for Painlev\'e $VI$ system.   For example, in the Painlev\'e $VI$ case, we can observe the following correspondence (cf. \cite{STe02}, \cite{IISA}).  
 (See \cite{Iwa02-2} or \cite{IISA} for the meaning of the nonlinear monodromy group 
 for Painlev\'e $VI$).  
 \begin{equation}
 \begin{array}{ccc}
    M_4^{\balpha}(\bt, \blambda, L)   &   &    \cR(\cP_{4, \bt})_{\ba} \\
        &   &  \\
\fbox{  $(-2)$ rational curves in  $M_4^{\balpha}(\bt, \blambda, L) $} 
&  \stackrel{{\bf RH}_{\blambda}}{\Longleftrightarrow}   &  \fbox{ Rational double points on $\cR(\cP_{4, \bt})_{\ba}$} \\
  &  &  \\
\Updownarrow &    & \Updownarrow \\
  &   &  \\
\fbox{Riccati solutions for $P_{VI}$ }  & \Longleftrightarrow  & \fbox{Fixed points of the nonlinear 
monodromies}              \\
 \end{array}
 \end{equation}

In Garnier case $(n \geq 5)$, when $\blambda$ is reducible (\ref{eq:special2-intro}), one can obtain a special classical solution of the 
equation integrated by hypergeometric functions $F_D$ of Lauricella (cf. [Proposition 1.7 \cite{K}]).   
One can see that these  
classical solutions of Garnier systems ${\mathcal G}_n$  
correspond to the subvariety isomorphic to $\BP^{n-3} $ which parametrizes 
reducible stable   parabolic connections.  Moreover when $\blambda$ is resonant (\ref{eq:special1-intro}), the Garnier system  
${\mathcal G}_n$ 
degenerates  into a  Riccati system over a Garnier system ${\mathcal G}_{n-1}$.   A subvariety which can be contracted by ${\bf RH}_{\bt, \blambda}$ is  
isomorphic to $\BP^1$-bundle over $M_{n-1}^{\balpha}(\bt', 
\blambda', L')$ at a generic point of the contracted subvatiety.

\subsubsection{Painlev\'e $VI$ or Garnier equations parametrized by $\blambda \in \Lambda_n$}
\quad

In the above formulation, for each fixed local exponent $\blambda \in \Lambda_n$, 
we obtain the Painlev\'e or Garnier vector 
fields $v_i (\blambda) $ for  $i, 1 \leq i \leq n$  as in (\ref{eq:def-vf})  
such that $\{ v_i (\blambda)\}_{1 \leq i \leq n} $ forms an  integrable 
differential system.   Moreover the solution manifold for the differential system 
can be given by the isomonodromic flows.  
 Varying the data $\blambda$, we obtain vector fields  
\begin{equation}\label{eq:painleve-vf1-intro}
\fbox{ $ v_i    
\in H^0(M_n^{\balpha}(L),  \Theta_{ M_n^{\balpha}(L)/ \Lambda_n})  $,  $ 1 \leq i \leq n $    }  
\end{equation}
for  $M_n^{\balpha}(L) \lra T_n \times \Lambda_n $ such that 
$v_{i| M_n^{\balpha}(\blambda, L)} = v_i (\blambda)$.

%%%%%%%%%%% FIGURE 1  =ISOMONODROMIC DEFORMATION %%%%%%%%%%%%%

%%%%%%%%%%%END OF FIGURE 1 %%%%%%%%%%%%%%%%%%%%%%%%%%%%%%

%%%%%%%%%%% FIGURE 2  =ISOMONODROMIC DEFORMATION %%%%%%%%%%%%%

\begin{center}
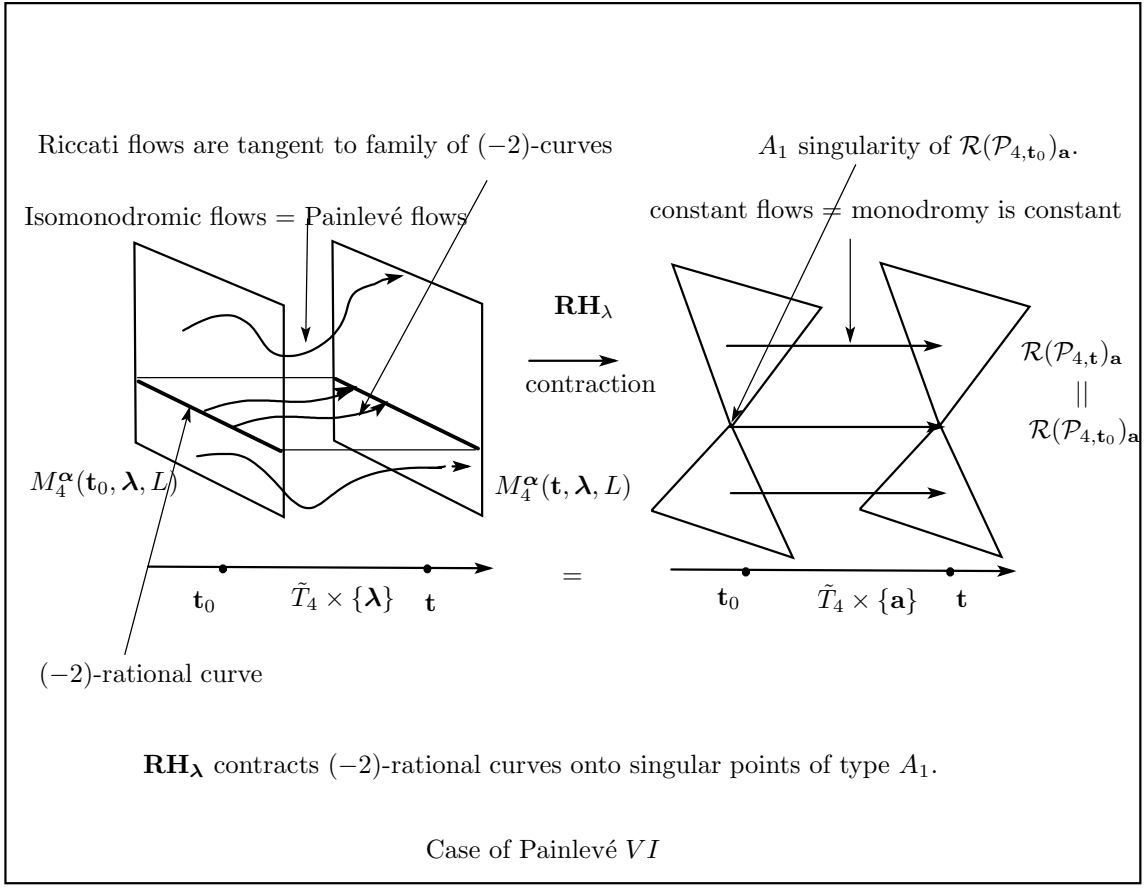
\begin{figure}

%WinTpicVersion2.15
\unitlength 0.1in
\begin{picture}(59.40,46.10)(0.10,-42.40)
% VECTOR 1 0 3 0
% 2 3496 2998 5283 3006
% 
\special{pn 13}%
\special{pa 3496 2598}%
\special{pa 5283 2606}%
\special{fp}%
\special{sh 1}%
\special{pa 5283 2606}%
\special{pa 5216 2586}%
\special{pa 5230 2606}%
\special{pa 5216 2626}%
\special{pa 5283 2606}%
\special{fp}%
% POLYGON 1 0 3 0
% 6 686 1286 693 2332 1468 2722 1468 2715 1468 1614 686 1286
% 
\special{pn 13}%
\special{pa 686 886}%
\special{pa 693 1932}%
\special{pa 1468 2322}%
\special{pa 1468 2315}%
\special{pa 1468 1214}%
\special{pa 686 886}%
\special{fp}%
% POLYGON 1 0 3 0
% 6 1723 1278 1731 2323 2505 2715 2505 2706 2505 1605 1723 1278
% 
\special{pn 13}%
\special{pa 1723 878}%
\special{pa 1731 1923}%
\special{pa 2505 2315}%
\special{pa 2505 2306}%
\special{pa 2505 1205}%
\special{pa 1723 878}%
\special{fp}%
% DOT 0 0 3 0
% 2 3886 3014 3886 3006
% 
\special{pn 20}%
\special{sh 1}%
\special{ar 3886 2614 10 10 0  6.28318530717959E+0000}%
\special{sh 1}%
\special{ar 3886 2606 10 10 0  6.28318530717959E+0000}%
% STR 2 0 3 0
% 3 3735 3126 3735 3206 2 0
% $\bt_0$
\put(37.3500,-28.0600){\makebox(0,0)[lb]{$\bt_0$}}%
% STR 2 0 3 0
% 3 4987 3116 4987 3196 2 0
% $\bt$
\put(49.8700,-27.9600){\makebox(0,0)[lb]{$\bt$}}%
% DOT 0 0 3 0
% 2 4956 3014 4956 3006
% 
\special{pn 20}%
\special{sh 1}%
\special{ar 4956 2614 10 10 0  6.28318530717959E+0000}%
\special{sh 1}%
\special{ar 4956 2606 10 10 0  6.28318530717959E+0000}%
% STR 2 0 3 0
% 3 990 3141 990 3221 2 0
% $\bt_0$
\put(9.9000,-28.2100){\makebox(0,0)[lb]{$\bt_0$}}%
% STR 2 0 3 0
% 3 2205 3136 2205 3216 2 0
% $\bt$
\put(22.0500,-28.1600){\makebox(0,0)[lb]{$\bt$}}%
% VECTOR 1 0 3 0
% 2 758 2982 2546 2990
% 
\special{pn 13}%
\special{pa 758 2582}%
\special{pa 2546 2590}%
\special{fp}%
\special{sh 1}%
\special{pa 2546 2590}%
\special{pa 2479 2570}%
\special{pa 2493 2590}%
\special{pa 2479 2610}%
\special{pa 2546 2590}%
\special{fp}%
% DOT 0 0 3 0
% 2 1149 2998 1149 2990
% 
\special{pn 20}%
\special{sh 1}%
\special{ar 1149 2598 10 10 0  6.28318530717959E+0000}%
\special{sh 1}%
\special{ar 1149 2590 10 10 0  6.28318530717959E+0000}%
% DOT 0 0 3 0
% 2 2219 2998 2219 2990
% 
\special{pn 20}%
\special{sh 1}%
\special{ar 2219 2598 10 10 0  6.28318530717959E+0000}%
\special{sh 1}%
\special{ar 2219 2590 10 10 0  6.28318530717959E+0000}%
% VECTOR 1 0 3 0
% 2 2753 1904 3200 1904
% 
\special{pn 13}%
\special{pa 2753 1504}%
\special{pa 3200 1504}%
\special{fp}%
\special{sh 1}%
\special{pa 3200 1504}%
\special{pa 3133 1484}%
\special{pa 3147 1504}%
\special{pa 3133 1524}%
\special{pa 3200 1504}%
\special{fp}%
% STR 2 0 3 0
% 3 2570 2557 2570 2637 2 0
% $M_4^{\balpha}(\bt, \blambda, L)$
\put(25.7000,-22.3700){\makebox(0,0)[lb]{$M_4^{\balpha}(\bt, \blambda, L)$}}%
% STR 2 0 3 0
% 3 5367 2267 5367 2347 2 0
% $\cR(\cP_{4, \bt_0})_{\ba} $
\put(53.6700,-19.4700){\makebox(0,0)[lb]{$\cR(\cP_{4, \bt_0})_{\ba} $}}%
% STR 2 0 3 0
% 3 3040 1540 3040 1620 5 0
% ${\bf RH}_{\lambda}$
\put(30.4000,-12.2000){\makebox(0,0){${\bf RH}_{\lambda}$}}%
% STR 2 0 3 0
% 3 4619 1041 4619 1121 5 0
% constant flows = monodromy is consant
\put(46.1900,-7.2100){\makebox(0,0){constant flows = monodromy is constant}}%
% STR 2 0 3 0
% 3 110 1110 110 1190 2 0
% Isomondromic flows = Painlev\'{e} flows
\put(1.1000,-7.9000){\makebox(0,0)[lb]{Isomonodromic flows = Painlev\'{e} flows}}%
% STR 2 0 3 0
% 3 2923 2967 2923 3065 2 0
% $=$
\put(29.2300,-26.6500){\makebox(0,0)[lb]{$=$}}%
% VECTOR 2 0 3 0
% 2 1593 1161 1582 1799
% 
\special{pn 8}%
\special{pa 1593 761}%
\special{pa 1582 1399}%
\special{fp}%
\special{sh 1}%
\special{pa 1582 1399}%
\special{pa 1603 1333}%
\special{pa 1583 1346}%
\special{pa 1563 1332}%
\special{pa 1582 1399}%
\special{fp}%
% STR 2 0 3 0
% 3 5327 1859 5327 1939 2 0
% $\cR(\cP_{4, \bt})_{\ba} $
\put(53.2700,-15.3900){\makebox(0,0)[lb]{$\cR(\cP_{4, \bt})_{\ba} $}}%
% STR 2 0 3 0
% 3 5607 2049 5607 2148 2 0
% $||$
\put(56.0700,-17.4800){\makebox(0,0)[lb]{$||$}}%
% STR 2 0 3 0
% 3 137 2537 137 2617 2 0
% $M_4^{\balpha}(\bt_0, \blambda, L)$
\put(1.3700,-22.1700){\makebox(0,0)[lb]{$M_4^{\balpha}(\bt_0, \blambda, L)$}}%
% STR 2 0 3 0
% 3 1506 3126 1506 3225 2 0
% $\tilde{T}_4 \times \{\blambda \} $
\put(15.0600,-28.2500){\makebox(0,0)[lb]{$\tilde{T}_4 \times \{\blambda \} $}}%
% STR 2 0 3 0
% 3 4260 3136 4260 3235 2 0
% $\tilde{T}_4 \times \{ \ba \}$
\put(42.6000,-28.3500){\makebox(0,0)[lb]{$\tilde{T}_4 \times \{ \ba \}$}}%
% STR 2 0 3 0
% 3 730 4001 730 4100 2 0
% ${\bf RH}_{\blambda}$ contracts $(-2)$-rational curves onto singular points of type $A_1$.  
\put(7.3000,-37.0000){\makebox(0,0)[lb]{${\bf RH}_{\blambda}$ contracts $(-2)$-rational curves onto singular points of type $A_1$.  }}%
% LINE 0 0 3 0
% 2 1739 1999 2484 2364
% 
\special{pn 20}%
\special{pa 1739 1599}%
\special{pa 2484 1964}%
\special{fp}%
% LINE 0 0 3 0
% 2 708 2012 1453 2378
% 
\special{pn 20}%
\special{pa 708 1612}%
\special{pa 1453 1978}%
\special{fp}%
% LINE 2 0 3 0
% 2 702 1992 1733 1992
% 
\special{pn 8}%
\special{pa 702 1592}%
\special{pa 1733 1592}%
\special{fp}%
% LINE 2 0 3 0
% 2 1440 2371 2471 2371
% 
\special{pn 8}%
\special{pa 1440 1971}%
\special{pa 2471 1971}%
\special{fp}%
% POLYGON 1 0 3 0
% 5 3502 1400 4280 1626 3814 2245 3807 2245 3502 1400
% 
\special{pn 13}%
\special{pa 3502 1000}%
\special{pa 4280 1226}%
\special{pa 3814 1845}%
\special{pa 3807 1845}%
\special{pa 3502 1000}%
\special{fp}%
% POLYGON 1 0 3 0
% 5 3807 2238 3395 2684 3395 2690 4133 2936 3807 2238
% 
\special{pn 13}%
\special{pa 3807 1838}%
\special{pa 3395 2284}%
\special{pa 3395 2290}%
\special{pa 4133 2536}%
\special{pa 3807 1838}%
\special{fp}%
% POLYGON 1 0 3 0
% 5 4592 1393 5370 1620 4905 2238 4898 2238 4592 1393
% 
\special{pn 13}%
\special{pa 4592 993}%
\special{pa 5370 1220}%
\special{pa 4905 1838}%
\special{pa 4898 1838}%
\special{pa 4592 993}%
\special{fp}%
% POLYGON 1 0 3 0
% 5 4898 2231 4486 2677 4486 2684 5224 2930 4898 2231
% 
\special{pn 13}%
\special{pa 4898 1831}%
\special{pa 4486 2277}%
\special{pa 4486 2284}%
\special{pa 5224 2530}%
\special{pa 4898 1831}%
\special{fp}%
% VECTOR 1 0 3 0
% 2 3801 2251 4911 2251
% 
\special{pn 13}%
\special{pa 3801 1851}%
\special{pa 4911 1851}%
\special{fp}%
\special{sh 1}%
\special{pa 4911 1851}%
\special{pa 4844 1831}%
\special{pa 4858 1851}%
\special{pa 4844 1871}%
\special{pa 4911 1851}%
\special{fp}%
% VECTOR 1 0 3 0
% 2 3787 1826 4898 1826
% 
\special{pn 13}%
\special{pa 3787 1426}%
\special{pa 4898 1426}%
\special{fp}%
\special{sh 1}%
\special{pa 4898 1426}%
\special{pa 4831 1406}%
\special{pa 4845 1426}%
\special{pa 4831 1446}%
\special{pa 4898 1426}%
\special{fp}%
% VECTOR 1 0 3 0
% 2 3807 2597 4918 2597
% 
\special{pn 13}%
\special{pa 3807 2197}%
\special{pa 4918 2197}%
\special{fp}%
\special{sh 1}%
\special{pa 4918 2197}%
\special{pa 4851 2177}%
\special{pa 4865 2197}%
\special{pa 4851 2217}%
\special{pa 4918 2197}%
\special{fp}%
% SPLINE 1 0 3 0
% 7 941 1746 1201 1659 1320 1733 1420 1859 1812 1679 1866 1540 1985 1480
% 
\special{pn 13}%
\special{pa 941 1346}%
\special{pa 972 1328}%
\special{pa 1002 1310}%
\special{pa 1033 1294}%
\special{pa 1063 1280}%
\special{pa 1094 1268}%
\special{pa 1124 1260}%
\special{pa 1154 1256}%
\special{pa 1184 1256}%
\special{pa 1214 1262}%
\special{pa 1243 1274}%
\special{pa 1271 1290}%
\special{pa 1298 1310}%
\special{pa 1321 1334}%
\special{pa 1342 1361}%
\special{pa 1361 1389}%
\special{pa 1379 1415}%
\special{pa 1399 1440}%
\special{pa 1421 1459}%
\special{pa 1446 1473}%
\special{pa 1475 1481}%
\special{pa 1506 1484}%
\special{pa 1539 1482}%
\special{pa 1573 1474}%
\special{pa 1607 1463}%
\special{pa 1642 1448}%
\special{pa 1675 1429}%
\special{pa 1707 1406}%
\special{pa 1736 1382}%
\special{pa 1763 1354}%
\special{pa 1786 1325}%
\special{pa 1805 1294}%
\special{pa 1819 1261}%
\special{pa 1829 1228}%
\special{pa 1838 1197}%
\special{pa 1849 1168}%
\special{pa 1864 1143}%
\special{pa 1885 1122}%
\special{pa 1912 1107}%
\special{pa 1943 1094}%
\special{pa 1976 1083}%
\special{pa 1985 1080}%
\special{sp}%
% VECTOR 1 0 3 0
% 2 1985 1480 2078 1460
% 
\special{pn 13}%
\special{pa 1985 1080}%
\special{pa 2078 1060}%
\special{fp}%
\special{sh 1}%
\special{pa 2078 1060}%
\special{pa 2009 1054}%
\special{pa 2026 1071}%
\special{pa 2017 1094}%
\special{pa 2078 1060}%
\special{fp}%
% SPLINE 1 0 3 0
% 6 1014 2404 1267 2404 1600 2684 1999 2484 2318 2457 2311 2457
% 
\special{pn 13}%
\special{pa 1014 2004}%
\special{pa 1047 1998}%
\special{pa 1080 1992}%
\special{pa 1113 1988}%
\special{pa 1145 1986}%
\special{pa 1177 1985}%
\special{pa 1209 1988}%
\special{pa 1240 1995}%
\special{pa 1270 2005}%
\special{pa 1299 2020}%
\special{pa 1327 2039}%
\special{pa 1355 2061}%
\special{pa 1382 2085}%
\special{pa 1408 2111}%
\special{pa 1433 2138}%
\special{pa 1458 2165}%
\special{pa 1481 2190}%
\special{pa 1504 2215}%
\special{pa 1527 2237}%
\special{pa 1549 2255}%
\special{pa 1570 2270}%
\special{pa 1590 2281}%
\special{pa 1610 2286}%
\special{pa 1629 2285}%
\special{pa 1649 2280}%
\special{pa 1668 2270}%
\special{pa 1689 2257}%
\special{pa 1711 2242}%
\special{pa 1734 2223}%
\special{pa 1760 2204}%
\special{pa 1787 2183}%
\special{pa 1818 2162}%
\special{pa 1852 2142}%
\special{pa 1889 2123}%
\special{pa 1931 2105}%
\special{pa 1977 2090}%
\special{pa 2027 2078}%
\special{pa 2082 2069}%
\special{pa 2137 2062}%
\special{pa 2190 2058}%
\special{pa 2240 2055}%
\special{pa 2282 2055}%
\special{pa 2315 2055}%
\special{pa 2335 2055}%
\special{pa 2341 2056}%
\special{pa 2330 2057}%
\special{pa 2311 2057}%
\special{sp}%
% VECTOR 1 0 3 0
% 2 2351 2457 2431 2457
% 
\special{pn 13}%
\special{pa 2351 2057}%
\special{pa 2431 2057}%
\special{fp}%
\special{sh 1}%
\special{pa 2431 2057}%
\special{pa 2364 2037}%
\special{pa 2378 2057}%
\special{pa 2364 2077}%
\special{pa 2431 2057}%
\special{fp}%
% VECTOR 2 0 3 0
% 2 635 3442 974 2152
% 
\special{pn 8}%
\special{pa 635 3042}%
\special{pa 974 1752}%
\special{fp}%
\special{sh 1}%
\special{pa 974 1752}%
\special{pa 938 1811}%
\special{pa 960 1804}%
\special{pa 976 1822}%
\special{pa 974 1752}%
\special{fp}%
% STR 2 0 3 0
% 3 186 3558 186 3624 2 0
% $(-2)$-rational curve
\put(1.8600,-32.2400){\makebox(0,0)[lb]{$(-2)$-rational curve}}%
% SPLINE 1 0 3 0
% 5 1201 2251 1420 2191 1726 2198 1972 2145 1972 2145
% 
\special{pn 13}%
\special{pa 1201 1851}%
\special{pa 1232 1840}%
\special{pa 1262 1829}%
\special{pa 1293 1819}%
\special{pa 1324 1810}%
\special{pa 1355 1802}%
\special{pa 1386 1796}%
\special{pa 1417 1791}%
\special{pa 1449 1789}%
\special{pa 1480 1788}%
\special{pa 1512 1789}%
\special{pa 1544 1791}%
\special{pa 1577 1793}%
\special{pa 1609 1795}%
\special{pa 1641 1797}%
\special{pa 1673 1799}%
\special{pa 1705 1799}%
\special{pa 1737 1797}%
\special{pa 1768 1794}%
\special{pa 1800 1789}%
\special{pa 1831 1783}%
\special{pa 1862 1776}%
\special{pa 1894 1768}%
\special{pa 1925 1759}%
\special{pa 1956 1750}%
\special{pa 1972 1745}%
\special{sp}%
% VECTOR 1 0 3 0
% 2 1965 2152 1999 2125
% 
\special{pn 13}%
\special{pa 1965 1752}%
\special{pa 1999 1725}%
\special{fp}%
\special{sh 1}%
\special{pa 1999 1725}%
\special{pa 1934 1751}%
\special{pa 1957 1758}%
\special{pa 1959 1782}%
\special{pa 1999 1725}%
\special{fp}%
% SPLINE 1 0 3 0
% 5 1054 2165 1274 2105 1580 2112 1826 2058 1826 2058
% 
\special{pn 13}%
\special{pa 1054 1765}%
\special{pa 1085 1754}%
\special{pa 1115 1743}%
\special{pa 1146 1733}%
\special{pa 1177 1724}%
\special{pa 1208 1716}%
\special{pa 1239 1710}%
\special{pa 1270 1705}%
\special{pa 1302 1703}%
\special{pa 1333 1702}%
\special{pa 1365 1703}%
\special{pa 1397 1705}%
\special{pa 1430 1707}%
\special{pa 1462 1709}%
\special{pa 1494 1711}%
\special{pa 1526 1713}%
\special{pa 1558 1713}%
\special{pa 1590 1711}%
\special{pa 1621 1708}%
\special{pa 1653 1703}%
\special{pa 1684 1697}%
\special{pa 1715 1690}%
\special{pa 1747 1681}%
\special{pa 1778 1672}%
\special{pa 1809 1663}%
\special{pa 1826 1658}%
\special{sp}%
% VECTOR 0 0 3 0
% 2 1754 2076 1820 2049
% 
\special{pn 20}%
\special{pa 1754 1676}%
\special{pa 1820 1649}%
\special{fp}%
\special{sh 1}%
\special{pa 1820 1649}%
\special{pa 1751 1656}%
\special{pa 1771 1669}%
\special{pa 1766 1693}%
\special{pa 1820 1649}%
\special{fp}%
% STR 2 0 3 0
% 3 180 784 180 850 2 0
% Riccati flows are tangent to family of $(-2)$-curves
\put(1.8000,-4.5000){\makebox(0,0)[lb]{Riccati flows are tangent to family of $(-2)$-curves}}%
% VECTOR 2 0 3 0
% 2 2531 961 1866 2185
% 
\special{pn 8}%
\special{pa 2531 561}%
\special{pa 1866 1785}%
\special{fp}%
\special{sh 1}%
\special{pa 1866 1785}%
\special{pa 1915 1736}%
\special{pa 1891 1738}%
\special{pa 1880 1717}%
\special{pa 1866 1785}%
\special{fp}%
% STR 2 0 3 0
% 3 3948 764 3948 860 2 0
% $A_1$ singularity of $\cR(\cP_{4, \bt_0})_{\ba}$. 
\put(39.4800,-4.6000){\makebox(0,0)[lb]{$A_1$ singularity of $\cR(\cP_{4, \bt_0})_{\ba}$. }}%
% VECTOR 2 0 3 0
% 2 4461 878 3815 2218
% 
\special{pn 8}%
\special{pa 4461 478}%
\special{pa 3815 1818}%
\special{fp}%
\special{sh 1}%
\special{pa 3815 1818}%
\special{pa 3862 1767}%
\special{pa 3838 1770}%
\special{pa 3826 1749}%
\special{pa 3815 1818}%
\special{fp}%
% VECTOR 2 0 3 0
% 2 4432 1268 4432 1800
% 
\special{pn 8}%
\special{pa 4432 868}%
\special{pa 4432 1400}%
\special{fp}%
\special{sh 1}%
\special{pa 4432 1400}%
\special{pa 4452 1333}%
\special{pa 4432 1347}%
\special{pa 4412 1333}%
\special{pa 4432 1400}%
\special{fp}%
% STR 2 0 3 0
% 3 2732 1971 2732 2066 2 0
% contraction
\put(27.3200,-16.6600){\makebox(0,0)[lb]{contraction}}%
% BOX 1 0 3 0
% 2 5950 4640 10 30
% 
\special{pn 13}%
\special{pa 5950 4240}%
\special{pa 10 4240}%
\special{pa 10 -370}%
\special{pa 5950 -370}%
\special{pa 5950 4240}%
\special{fp}%
% STR 2 0 3 0
% 3 2210 4410 2210 4510 2 0
% Case of Painlev\'e $VI$
\put(22.1000,-41.1000){\makebox(0,0)[lb]{Case of Painlev\'e $VI$}}%
\end{picture}%
\caption{Riemann-Hilbert correspondence and isomonodromic flows for special 
 $\blambda$}
\label{fig:isom2}
\end{figure}
\end{center}

%%%%%%%%%%%END OF FIGURE 1 %%%%%%%%%%%%%%%%%%%%%%%%%%%%%%

\subsection{The  Hamiltonian system} 
\quad 

It is well-known that the Painlev\'e and Garnier equations can be written in 
the Hamiltonian systems.   Now we can explain this as follows.  Since 
the constant flows on $
(\phi_n)_{\ba}: (\tilde{\cR}_n)_{\ba}  \lra \tilde{T}_n$ preserve the natural 
symplectic form $\Omega_1$  
on the fiber $\cR(\cP_{n, \bt})_{\ba}$ and the pullback of $\Omega_1$  by 
Riemann-Hilbert correspondence coincides with the symplectic structure 
 $\Omega$ on $M_n^{\balpha}(\bt, \blambda, L)$, Painlev\'e 
 or  Garnier vector fields  preserve the symplectic structure $\Omega$.  
 Therefore, we can write the differential equations in the Hamiltonian systems 
by using suitable cannonical coordinate systems.  Then an argument  shows that  
such vector fields are actually regular algebraic, hence the Hamiltonians are 
given by  regular algebraic functions.

\subsection{The relation of the space of initial conditions of Okamoto or 
Okamoto-Painlev\'e pairs for $P_{VI}$}  
\quad

In the case of $P_{VI}$, Okamoto  \cite{O1}  constructed  
the spaces of initial conditions by blowing up  the accessible singularities 
of 4 parameter family of  Painlev\'e $VI$ equations.  They are open algebraic surfaces which are 
complements of the anti-canonical divisors of projective rational surfaces 
obtained by the 8-point blowing-ups of $\BP^1_{\C} \times \BP^1_{\C}$ or 
$\F_2$ .  In \cite{Sakai}, \cite{STT02}, 
the notion of the pairs of projective rational surfaces and its effective anti-canonical divisors  
with suitable conditions was introduced and its relation to 
Painlev\'e equation was revealed.  In \cite{STT02}, such a pair is called    an {\em Okamoto--Painlev\'e pair}.  
Okamoto-Painlev\'e pairs of type $D_4^{(1)}$ 
correspond to Painlev\'e $VI$ equations.  A semiuniversal family of Okamoto--Painlev\'e pairs is a  family of projective surfaces $\pi: \overline{\cS} \lra T'_4 \times \Lambda_4 $ with the effective relative anticanonical divisor  $\cY$
such that the configuration of the anticanonical divisor $\cY_{\bt, \blambda}$ 
is of type $D_4^{(1)}$.  Then family of  spaces of the initial conditions 
of Okamoto can be obtained as an open subset $ \cS := \overline{\cS} \setminus 
\cY$.   

   In the second part of this paper\cite{IIS-2}, 
 we will show that the family of Okamoto-Pailev\'e pairs 
 $ \overline{\cS} \lra T_4 \times \Lambda_4 $ can be identified with 
 the family of the moduli spaces 
 $ \overline{M_4^{\balpha' \bbeta}}(\cO_{\BP^1}(-1))\lra T_4 \times \Lambda_4 $, 
 while  $ \cS \lra T_4 \times \Lambda_4$ can be identified with 
 $ M_4^{\balpha} ( \cO_{\BP^1}(-1))  \lra T_4 \times \Lambda_4 $. 
 (In this case, we will take $\beta_1 = \beta_2 = 1$, 
 hence $\balpha = \balpha'\frac{\beta_1}{\beta_1+ \beta_2}=\balpha'/2$). 
 So our constructions of the moduli spaces give an intrinsic meaning of 
Okamoto's explicit hard calculations in \cite{O1}.

\subsection{The B\"acklund transformations--Symmetries of the equations}
\label{subsection:symmetry}
\quad

In our framework,  B\"acklund transformations for 
the Painlev\'e equations or Garnier equations can be defined as follows.  
Consider the Painlev\'e $VI$ or Garnier system $\{ v_i \}_{1 \leq i \leq n } $ 
defined in (\ref{eq:painleve-vf1-intro}) and the family of 
moduli spaces $\pi_n :  M_n^{\balpha}(L) \lra T'_n \times \Lambda_n$.

\begin{Definition}  \label{def:backlund-intro} 
{\rm The pair $(\tilde{s}, s)$  of 
a birational map $\tilde{s}: M_n^{\balpha}(L) 
\cdots \ra M_n^{\balpha}(L) $ (or $\tilde{s}: \overline{M_n^{\balpha'\bbeta}}(L)\cdots \ra \overline{M_n^{\balpha'\bbeta}}(L) $  )
and  an affine transformation 
$
s: \Lambda_n \lra \Lambda_n$ is said to be a 
{\em B\"acklund transformation of 
the differential system}   $\{v_i \}_{1 \leq i \leq n } $  or $
\{ v_i( \blambda) \}_{ 1 \leq i \leq n, \blambda \in \Lambda_n }  $  
if they make the following diagram commutative:
\begin{equation}\label{eq:backlund1-intro}
\begin{CD}
M_n^{\balpha}(L)  & \stackrel{\tilde{s}}{\cdots \ra } &  
M_n^{\balpha}(L)  \\
\quad \downarrow \pi_n   & & \quad \downarrow \pi_n \\
T'_n \times \Lambda_n  &  \stackrel{1 \times s}{\lra} &  T'_n \times \Lambda_n, \end{CD}
\end{equation}
and  it satisfies the condition:
\begin{equation}\label{eq:backlund2-intro}
\fbox{\quad $
\tilde{s}_{*} (v_i )  = v_i $,   or equivalently   $\tilde{s}_{*} (v_i(\blambda) )  = 
v_i( s(\blambda)) $ \quad}. 
\end{equation}}
\end{Definition}

There exists a natural class of B\"acklund transformations of 
$M_n^{\balpha}(L)$ for any $n \geq 4$ 
which are   induced by {\em elementary transformations} of  
stable parabolic connections (cf.\ \S \ref{sec:elemnt}). 
Since  such transformations induce the 
identity on the moduli space of the monodromy representations via  
Riemann-Hilbert correspondence, we can conclude that the  
transformations preserve the vector field as in (\ref{eq:backlund2-intro}).  
(This notion  is 
equivalent to the rational  gauge transformation or Schlesinger transformation  of  connections ).    
In \S \ref{sec:elemnt}, we will list up these kinds of B\"acklund transformations.

\begin{figure}
\begin{center}
%WinTpicVersion2.15
\unitlength 0.1in
\begin{picture}(20.39,26.50)(20.85,-36.95)
% CIRCLE 2 0 3 0
% 4 3125 2681 3205 2841 3205 2841 3285 3021
% 
\special{pn 8}%
\special{ar 3125 2281 179 179  1.1309537 6.2831853}%
\special{ar 3125 2281 179 179  0.0000000 1.1071487}%
% CIRCLE 2 0 3 0
% 4 3945 1881 4025 2041 4025 2041 4105 2221
% 
\special{pn 8}%
\special{ar 3945 1481 179 179  1.1309537 6.2831853}%
\special{ar 3945 1481 179 179  0.0000000 1.1071487}%
% CIRCLE 2 0 3 0
% 4 3945 3481 4025 3641 4025 3641 4105 3821
% 
\special{pn 8}%
\special{ar 3945 3081 179 179  1.1309537 6.2831853}%
\special{ar 3945 3081 179 179  0.0000000 1.1071487}%
% CIRCLE 2 0 3 0
% 4 2325 3461 2405 3621 2405 3621 2485 3801
% 
\special{pn 8}%
\special{ar 2325 3061 179 179  1.1309537 6.2831853}%
\special{ar 2325 3061 179 179  0.0000000 1.1071487}%
% CIRCLE 2 0 3 0
% 4 2325 1881 2405 2041 2405 2041 2485 2221
% 
\special{pn 8}%
\special{ar 2325 1481 179 179  1.1309537 6.2831853}%
\special{ar 2325 1481 179 179  0.0000000 1.1071487}%
% LINE 2 0 3 0
% 2 3805 1981 3261 2549
% 
\special{pn 8}%
\special{pa 3805 1581}%
\special{pa 3261 2149}%
\special{fp}%
% LINE 2 0 3 0
% 2 3009 2809 2489 3369
% 
\special{pn 8}%
\special{pa 3009 2409}%
\special{pa 2489 2969}%
\special{fp}%
% LINE 2 0 3 0
% 2 3265 2821 3825 3381
% 
\special{pn 8}%
\special{pa 3265 2421}%
\special{pa 3825 2981}%
\special{fp}%
% STR 2 0 3 0
% 3 3185 4061 3185 4115 5 0
% Root system $D_4^{(1)}$
\put(31.8500,-37.1500){\makebox(0,0){Root system $D_4^{(1)}$}}%
% STR 2 0 3 0
% 3 3355 2515 3355 2715 2 0
% $\alpha_0$
\put(33.5500,-23.1500){\makebox(0,0)[lb]{$\alpha_0$}}%
% STR 2 0 3 0
% 3 2105 1421 2105 1621 2 0
% $\alpha_1$
\put(21.0500,-12.2100){\makebox(0,0)[lb]{$\alpha_1$}}%
% STR 2 0 3 0
% 3 3905 1415 3905 1615 2 0
% $\alpha_2$
\put(39.0500,-12.1500){\makebox(0,0)[lb]{$\alpha_2$}}%
% STR 2 0 3 0
% 3 2535 3635 2535 3835 5 0
% $\alpha_3$
\put(25.3500,-34.3500){\makebox(0,0){$\alpha_3$}}%
% STR 2 0 3 0
% 3 3770 3710 3770 3910 2 0
% $\alpha_4$
\put(37.7000,-35.1000){\makebox(0,0)[lb]{$\alpha_4$}}%
% LINE 2 0 3 0
% 2 2461 2009 3001 2541
% 
\special{pn 8}%
\special{pa 2461 1609}%
\special{pa 3001 2141}%
\special{fp}%
% STR 2 0 3 0
% 3 2565 3995 2565 4095 5 0
% 
\put(25.6500,-36.9500){\makebox(0,0){}}%
\end{picture}%
\end{center}
\caption{Dynkin diagram $D_4^{(1)}$}
\label{figure:dynkin}
\end{figure}
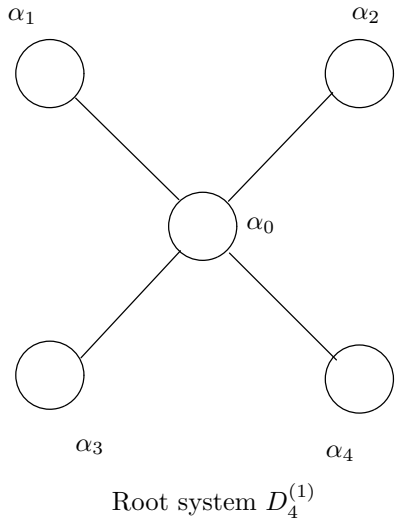

As for Painlev\'e $VI$ equations, the group of the B\"{a}cklund transformation 
in the  above sense is isomorphic to the affine Weyl group $W(D_4^{(1)})$ of 
type $D_4^{(1)}$,   (cf. e.g. \cite{O4}, \cite{IIS-0}).  
The affine Weyl group $W(D_4^{(1)})$  is generated by 
5  reflections $s_i$, $ i = 0, 1, \ldots, 4$
corresponding to the simple roots in the Dynkin diagram of type $D_4^{(1)}$
(see Figure \ref{figure:dynkin}). 
A  natural faithful affine action of $W(D_4^{(1)})$  to  $\Lambda_4 = \C^4 \ni \blambda = (\lambda_j)$ 
can be given by  
\begin{equation}\label{eq:action}
\begin{array}{rcl}
s_i(\lambda_j) & = & (-1)^{\delta_{ij}} \lambda_j, \quad i = 1, \ldots, 4 \\
s_0(\lambda_i) & = & \lambda_i - \frac{1}{2} \left( 
\sum_{k=1}^4 \lambda_k \right)+ \frac{1}{2}. 
\end{array}
\end{equation}

Recalling the identification of the family 
$ \overline{M_4^{\balpha'\bbeta}}(\cO_{\BP^1}(-1)) \lra T'_4 \times \Lambda_4 $ 
with the family of Okamoto-Painlev\'e 
pairs $\overline{\cS} \lra T'_4 \times \Lambda_4 $,  
one can see that the actions of $ W(D_4^{(1)}) $ in (\ref{eq:action}) 
can be lifted to  birational actions of
the total space of the 
family  $ \overline{\pi}:\overline{\mathcal S} \lra  
T'_4 \times \Lambda$, that is, for each $ s \in W(D_4^{(1)})$, there
 exists a commutative diagram   
\begin{equation}
\begin{CD}
\overline{\mathcal S} & \stackrel{\tilde{s}}{\cdots \ra } &  
\overline{\mathcal S} \\
\quad \downarrow \overline{\pi}   & & \quad \downarrow \overline{\pi} \\
T'_4 \times \Lambda_4  &  \stackrel{1 \times s}{\lra} &  T'_4 \times \Lambda_4. \end{CD}
\end{equation}
Moreover it is known \cite{O4} that 
the actions preserve the Painlev\'e vector field $v_i$ 
 in (\ref{eq:painleve-vf1-intro}).  
That is,  for each  $s \in W(D_4^{(1)})$,   we have  
\begin{equation}\label{eq:preserve-eq-intro}
\tilde{s}_{*}(v_i (\blambda)) = v_i(s(\blambda)) \quad \mbox{for}  \quad 1 \leq i \leq 4 .  
\end{equation}

In our framework, it is easy to give an 
intrinsic reason why  $\tilde{s_i}$ for $ 1 \leq i \leq 4$ preserve
the vector field. It is simply  
because these come from  elementary transformations. 
However, the origin of the transformation $\tilde{s_0}$ is still mysterious, and we cannot see any simple reason why $\tilde{s}_0$ preserves the vector field.    

As some experts suggested to us, 
it may be  plausible 
to believe that  $\tilde{s_0}$ is induced by   Laplace transformations of 
the stable  connection. (The authors were informed by H. Sakai that 
M. Mazzocco gives  some explanations for this fact  on this line).  
For simplicity, let us call 
the Laplace transform of the original connection the 
{\em dual} of the connection. 
In general, the dual  of logarithmic connections of rank $2$
becomes a connection of   higher rank which may not be 
logarithmic, so it is not so easy to identify the dual of the connection to 
the original one.   Only in the case of $n = 4$ (Painlev\'e case), 
we may miraculously identify the original connection with its dual or a further transformed object, 
so we have the extra B\"acklund transformation like $\tilde{s}_0$.  
It may be reasonable to consider 
the original connection and its  dual at once.  
Then we may include the Laplace transformation as a part of the B\"acklund transformations.

After we have finished 
the first version of this paper, Philip Boalch informed us 
that he can obtain $\tilde{s}_0$ using the method of 
\cite{Boalch} as follows. 
One can embed a rank $2$ connection with 
$4$-regular singular points over $\BP^1$ 
into a rank $3$ reducible connection.  Then there is a simple operation for shifting the eigenvalues of 
the rank $3$ connection.  For a special value of shifting, one can obtain a rank $2$ subconnection 
or a rank $1$ subconnection  in the shifted rank $3$ connection, then take the 
rank $2$ connection or the quotient of rank $1$ subconnection.  
This gives a transformation from a 
rank $2$ connection to another rank $2$ connection 
whose transformation on $\blambda$ gives $s_0$. 
Note that this transformation only works for the case of $n =4$.  
By using this result and  another result in \cite{Boalch}, 
he also gave a different proof of a result in 
\cite{IIS-0}.

Besides these stories, we should mention about the relation of the birational 
geometry and the B\"acklund transformations.   
As Saito and Umemura pointed out in \cite{SU01},   
 B\"acklund transformations of Painlev\'e equations which are reflections with respect to  roots of an affine root system are nothing but {\em flops} 
  corresponding to $(-2)$-rational curves in 
  Okamoto spaces of initial conditions. 

{F}rom the definition of elementary transformations, we can  easily see that 
the locus of indeterminacy of  birational transformations correspond to the subvarieties which are contracted by the Riemann-Hilbert correspondence.  Since  the Riemann-Hilbert correspondence gives a simultaneous symplectic resolution of the singularities of the family $\phi_n:\cR_n \lra T'_n \times 
\cA_n $, it is now obvious that those B\"acklund transformations are flops.   
(For definition and fundamental facts on flops, see [\S 6,1 \cite{KM}]).

\subsubsection{B\"acklund transformations and the Riemann-Hilbert correspondence}
\quad 

In \cite{IIS-0}, we have proved that all of 
the B\"acklund transformations  in $W(D_4^{(1)}) $  on 
$\overline{M_4^{\balpha}}(L) \lra T'_4 \times \Lambda_4$ induce
essentially identity on the moduli space $\cR_4$ after we 
take a finite quotient of $\cA_n$.  
(Note that this is nontrivial only for 
$\tilde{s_0}$).  
Therefore  in this sense, the group of the 
B\"acklund transformations $W(D_4^{(1)}) $ can be considered as 
the Galois group of the monodromy representations.  

\subsection{Related works}
\quad 

It is worthwhile to discuss about some works  related  to this paper and 
to clarify  what   are really new in this paper.  

The notion of $(\bt, \blambda)$-parabolic connection on $\BP^1$ 
is essentially introduced by Arinkin-Lysenko in \cite{AL1} 
as a 
quasiparabolic $SL_2$-bundle.  In \cite{AL1}, 
they also discussed about the 
moduli problem for quasiparabolic $SL_2$-bundles and 
consider the moduli space as an algebraic stack.  
In the case of $n=4$, under the assumption that 
$\blambda$ is generic  (cf. Definition \ref{def:exponents-intro}),  
they give an  explicit description of the coarse moduli space. 
 Moreover, in \cite{AL2}, by using the explicit  descriptions 
 of the moduli spaces, they describe the group of automorphisms of 
 the family of moduli spaces by using an explicit 
geometry of surfaces.   
Later, in \cite{A}, Arinkin introduced   
a notion of $\epsilon$-bundle, generalizing Deligne-Simpson's 
$\tau$-connections in \cite{Simp-I}.   
Again under the assumption that $\blambda$ is generic 
Arinkin gives a compactification of the moduli  space of quasiparabolic $SL_2$-bundles. 
Although the basic notions are introduced in their works, from the viewpoint 
of geometric background for Painlev\'e or Garnier equations, it is really necessary to 
construct the  moduli spaces even for special $\blambda$. 
For example, as we pointed out in \ref{subsection:symmetry} (cf. \cite{SU01}), 
some B\"acklund transformations of these equations are induced by flops in the terminology 
of  the modern birational 
geometry and the center of flops are lying over the special 
parameter $\blambda$.

In this sense, the advantage of introducing the notion of the stability for 
$(\bt, \blambda)$-parabolic connection  is obvious.  In the GIT setting 
(\cite{Mum:GIT}), 
despite   considerable careful computations of stability, 
we can construct the fine moduli space 
of stable objects as smooth irreducible schemes 
even for special  $ \blambda$.  Moreover, we introduce the notion of 
parabolic $\phi$-connections 
which is a generalization of the notion of $\epsilon$-connections 
due to Arinkin-Deligne-Simpson  and define the stability for them. 
One can understand the powers of these notions in  
Theorem \ref{thm:fund} and Theorem \ref{thm:n=4}.

The construction of the family of moduli spaces in \eqref{eq:family-rep-1}  
of $SL_2(\C)$-monodromy representations are 
essentially due to Simpson \cite{Simp-II}. 
However a systematic treatment of nonlinear monodromies of the
braid group is given by Dubrovin and Mazzocco [DM] for a 
 special case of $n = 4$, and by Iwasaki \cite{Iwa02-1}, 
 \cite{Iwa02-2} for the general case of $n=4$, and our 
 construction of the family in this paper is taking care of 
 the action of nonlinear monodromies of the braid groups to the
 moduli spaces.

Next, we would like to emphasize that only after 
we establish the natural setting in 
\ref{subsec:naturalfram} it becomes  possible 
to give a precise formulation of Hilbert's 21th problem for these
cases.   
In our setting,  the affirmative answer to the problem is 
equivalent to  the surjectivity of the Riemann-Hilbert 
correspondences $\RH_n$ and   $\RH_{\bt, \blambda}$ 
in \eqref{eq:rh-intro} and  \eqref{eq:fiber-cor}.  
As one can imagine easily, only reasonable result 
which one  can apply to proof of 
the surjectivity is Deligne's theorem in \cite{Deligne:70}.   

Moreover the properness of  $\RH_{\bt, \blambda}$ is also a new result. 
In the process of the proof, we need some analysis of the 
contraction induced by $\RH_{\bt, \blambda}$ and a 
technical lemma due to Professor A. Fujiki.  
The symplectic nature of the moduli spaces is 
discussed by many authors.  
(See for example \cite{Goldman:84}, \cite{Iwa91} \cite{Iwa92}).  
Iwasaki gave intrinsic symplectic structures on the moduli spaces 
of irreducible  logarithmic connections on a nonsingular 
complete curve and show that they are  obtained as 
the pullback of the symplectic structures 
on the moduli space of the irreducible representations.
   
Again, in this paper, we extended the symplectic structure  
to the whole moduli space of the stable 
parabolic connections and the 
smooth part of the moduli of the representations. 
Then one can show that these 
symplectic structures are identified via  $\RH_{\bt, \blambda}$.  
Our proof here is based on some complexes of sheaves 
whose hypercohomologies  
describe the tangent spaces to the moduli spaces.  
Together with the surjectivity and 
the properness of $\RH_{\bt, \blambda}$ 
these results can be understood as 
$\RH_{\bt, \blambda}$ gives an analytic  symplectic resolution 
of singularities in the sense of \cite{Bea}. 
These kinds of viewpoints seem to be new,  
and this gives a clear explanation that a simple reflection 
in the group of B\"acklund transformations is nothing 
but a flop with respect to this resolution.

The derivation of Painlev\'e equations from the isomonodromic deformation 
of the linear connections is 
well-known. (See for example \cite{JMU}, \cite{JM} and \cite{O3}).  
However in most cases, one 
first takes a normalized linear connection written in 
certain coordinate systems and 
then writes up the Painlev\'e equations as the  compatibility 
conditions for the extended linear connections. 
For a normalization, one has to assume that 
the vector bundle $E$ of rank $2$ and  degree $0$  
is isomorphic to $\cO_{\BP^1} \oplus \cO_{\BP^1}$, which is not true in general.  
In fact, the natural subscheme 
\begin{equation}\label{eq:taudivisor}
Z_{\bt} = \{ (E, \nabla, \varphi, l) \in \cM_{n}^{\balpha}(\bt, \blambda, L); 
E \not\simeq \cO_{\BP^1} \oplus \cO_{\BP^1} \}.
\end{equation} 
of $\cM_{n}^{\balpha}(\bt, \blambda, L)$ is a non-empty divisor.  
We note that the isomonodromic flow starting from some point 
$ p \in \cM_{n}^{\balpha}(\bt_0, \blambda, L) \setminus Z_{\bt_0}$ 
does not stay inside the open subset 
$\bigcup_{\bt} \left(\cM_{n}^{\balpha}(\bt, \blambda, L) \setminus Z_{\bt} \right)$, 
that is, the flow intersects with  $Z_{\bt}$ for some $\bt$. 
Therefore, in order to
prove the Painlev\'e property of Painlev\'e VI or Garnier equations completely, we have to 
consider the whole space $\cM_{n}^{\balpha}(\bt, \blambda, L)$ 
and the properness of 
the Riemann-Hilbert correspondence is essential 
for our proof of Painlev\'e property.  
(For a discussion of various definitions of Painlev\'e property, 
see \cite{IISA}. 
Moreover,  for some proofs of analytic Painlev\'e property of isomonodromic deformations, see \cite{Mal} and \cite{Miw}). 
Moreover most former approaches avoid dealing  with the case when 
$\blambda$ is special, because one has to introduce the notion of 
the stability of the parabolic connections to obtain a 
good moduli space which is smooth and Hausdorff.

In our framework, we can also discuss the Painlev\'e
or Garnier equations for special $\blambda$ 
in a natural framework.  
Interestingly enough, the classical solutions for 
these equations can be derived from the family of 
subvarieties contracted by $\RH_{\bt, \blambda}$. 
Now, the geometric meaning of these facts becomes
very clear.  
(For more detailed treatment in Painlev\'e $VI$  case, 
see \cite{W}, \cite{STe02} and \cite{SU01}).

We should mention that Nakajima \cite{N} obtained a 
smooth moduli space of stable parabolic connections 
as the moduli space of 
filtered regular $D$-modules by the 
technique of  the hyper-K\"{a}hler quotients of moment maps. 
Then he showed that the moduli space is diffeomorphic 
to the moduli space of parabolic Higgs bundles.  
Nitsure \cite{Ni} also constructed the moduli space 
of the stable logarithmic connection 
without parabolic structures in GIT setting.

\tableofcontents

\subsection{Acknowledgements}
First of all, the authors  would like to thank 
Akira Fujiki who kindly informed them 
of the proof of Lemma \ref{lem:fujiki}.
During the work, he has been  kind enough 
to send several careful replies to our questions via emails.  
The authors are grateful to Hidetaka Sakai, Kazuo Okamoto and 
 Kota Yoshioka for their  valuable suggestions.  
The last author gives his sincere thanks to Eiichi Sato for his  
supports for travels from Kobe to Fukuoka, 
which made these fruitful joint works possible.   
 The last author would like to thank JSPS-NWO  exchange program 
for their financial support for 
visiting Utrecht University in September--October 2003. 
He also would like to thank  
Jan Stienstra for his interest to our work and 
his hospitality in Utrecht.   
Last but not least, it is pleasure of three authors 
 to dedicate this paper to  Kyoichi Takano on 
his 60th birthday with the great respect to his 
leadership on the studies  of Painlev\'e 
equations as well as his 
friendships for younger generations.

%%%%%%%%%%%%%%%%%%% sec:moduli-spc %%%%%%%%%%%%%%

\section{Moduli spaces of stable parabolic 
connections on $\BP^1$ and their compactifications. }
\label{sec:moduli-spc}
%%%%%%%%%%%%%%%%%%%%%%

\subsection{Parabolic connections on $\BP^1$.}

Let $n \geq 3$ and set
\begin{equation}\label{eq:config-space}
T_n = \{ ( t_1, \ldots, t_n) \in (\BP^1)^n \quad | \quad t_i \not=t_j, (i \not= j) 
\}, 
\end{equation}
\begin{equation}\label{eq:exponents-space}
\Lambda_{n}  = \{ \blambda = (\lambda_1, \ldots, \lambda_n) \in  \C^{n} \}. 
\end{equation}  
Fixing a data  $(\bt, {\blambda}) = (t_1, \ldots, t_n, \lambda_1, \ldots, \lambda_n)  
\in T_n \times \Lambda_n$, we define a reduced divisor on $\BP^1$ as 
\begin{equation} \label{eq:divisor}
D(\bt) = t_1 + \cdots + t_n.
\end{equation}  
Moreover we fix a line bundle $L$ on $\BP^1$ with a logarithmic 
connection $\nabla_L: L 
\lra L \otimes \Omega^1_{\BP^1}(D(\bt))$.  

\begin{Definition}\label{def:parabolic} \rm 
A $($rank $2$ $)$ $(\bt, \blambda)$-parabolic connection on $\BP^1$ with the determinant 
$(L, \nabla_L)$ is a quadruplet
$(E, \nabla, \varphi, \{l_i\}_{1 \leq i \leq n}  )$ 
which consists of  
\begin{enumerate}
\item 
a rank 2 vector bundle 
$E$ on $\BP^1$,

\item 
a logarithmic connection $\nabla:E \lra E 
\otimes \Omega^1_{\BP^1}(D(\bt)) $

\item 
a bundle  isomorphism $ \varphi: \wedge^2 E \stackrel{\simeq}{\lra} L  $ 

\item  one dimensional subspace $l_i$ of the fiber $ E_{t_i}$  of $E$ 
at $t_i$,   $l_i \subset E_{t_i}$,   $ i = 1, \ldots, n$,   
such that 

\begin{enumerate}
\item  for any local sections $s_1, s_2$ of $E$, 
$$
\varphi \otimes id ( \nabla s_1 \wedge s_2 + s_1 \wedge \nabla s_2)  = 
\nabla_L(\varphi(s_1 \wedge s_2)), 
$$

\item  
$l_i \subset  \Ker (\res_{t_i}(\nabla) - \lambda_i)$, that is, $\lambda_i$ is aneigenvalue of the residue $\res_{t_i}(\nabla)$ of $\nabla$ at $t_i$ and 
$l_i$ is a one-dimensional eigensubspace of $\res_{t_i}(\nabla)$.   
\end{enumerate}
\end{enumerate}
\end{Definition}

\begin{Definition}\label{def:isom-spb} \rm 
  Two $(\bt, \blambda)$-parabolic connections 
$(E_1, \nabla_1, \varphi, \{l_i\}_{1 \leq i \leq n}  )$
$(E_2, \nabla_2, \varphi', \{l'_i\}_{1 \leq i \leq n}  )$
 on $\BP^1$ with the determinant 
$(L, \nabla_L)$ are isomorphic to each other if 
there is an isomorphism $\sigma: E_1 \stackrel{\sim}{\lra} E_2 $ and $c \in \C^{\times}$
 such that the diagrams
 \begin{equation}
\begin{CD}
  E_1 @>\nabla_1>> E_1\otimes\Omega^1_{\BP^1}(D(\bt)) \\
 @V\sigma V\cong V  @V\cong V \sigma \otimes\mathrm{id}V \\
  E_2 @>\nabla_2>> E_2\otimes\Omega^1_{\BP^1}(D(\bt)) 
 \end{CD}
 \hspace{50pt}
 \begin{CD}
  \bigwedge^2 E_1 @>\varphi>\cong > L \\
  @V\wedge^2\sigma V\cong V @VcV\cong V \\
  \bigwedge^2 E_2 @>\varphi'>\cong > L 
 \end{CD}
\end{equation}
commute and $(\sigma)_{t_i}(l_i)=l'_i$ for $i=1,\ldots,n$. 
\end{Definition}

\subsection{The set of local exponents  $\blambda \in \Lambda_n$}

Note that a data 
$ \blambda = (\lambda_1, \ldots, \lambda_n) 
\in \Lambda_n \simeq \C^n $
specifies the set of eigenvalues of the residue 
matrix of a connection $\nabla$ at $\bt=(t_1, \ldots, t_n)$, 
which will be called a set of   {\em local exponents} of $\nabla$. 

\begin{Definition}\label{def:exponents}
{\rm 
A set of local exponents 
 $ \blambda =(\lambda_1, \ldots, \lambda_n) \in \Lambda_n $ is called {\em special} if  
\begin{enumerate}
\item $\blambda$ is {\em resonant}, that is, for some $ 1 \leq i \leq n$, 
\begin{equation} \label{eq:special1}
  2 \lambda_i \in \Z, 
\end{equation}  
\item or $\blambda$ is {\em reducible}, that is, for some 
$ (\epsilon_1, \ldots, \epsilon_n) \in \{ \pm 1 \}^n $ 
\begin{equation} 
\label{eq:special2}\sum_{i=1}^{n} \epsilon_i  \lambda_i  \in \Z. 
\end{equation}
\end{enumerate}

If $\blambda \in \Lambda_n$ is not special,  
$\blambda$ is said to be  
{\em generic}. 
}
\end{Definition}

\begin{Lemma} \label{lem:special2}
Let  
 $(E, \nabla, \varphi, l = \{ l_i \} ) $ be  
  a $(\bt, \blambda)$-parabolic connection on $\BP^1$ with the 
  determinant $(L, \nabla_L)$. Assume that  eigenvalues of 
  $\res_{t_i}(\nabla_L) $ are  integers for $1 \leq i \leq n$.   
  Suppose that 
  there exists a subline bundle $F \subset E$ such that 
  $ \nabla F \subset F \otimes \Omega^1_{\BP^1}(D(\bt))$. Then 
  $\blambda$ is reducible, that is, 
$\blambda$ satisfies the condition (\ref{eq:special2}).  
\end{Lemma}

\begin{proof}
(Cf.\ [Proposition 1.1, \cite{AL1}]). 
Since  we have  a horizontal bundle 
isomorphism $\varphi:\wedge^2 E \simeq L $ 
with respect to the connections,  the eigenvalues of 
the residue matrix $ \res_{t_i} \nabla$ at $t_i$ are given by 
 $\lambda_i$ and 
 $ \res_{t_i} (\nabla_L) - \lambda_i$.  
Since $\nabla F \subset F  \otimes \Omega^1_{\BP^1}(D(\bt))$, 
the subspace 
$F_{t_i} \subset E_{t_i}$ is an eigenspace of $\res_{t_i} (\nabla)$.  
Therefore the eigenvalue of  $\res_{t_i} (\nabla_{|F})$ is congruent to 
 $\epsilon_i  \lambda_i$ 
modulo $\Z$  for $\epsilon_i = 1$ or $-1$.     
The residue 
theorem says that 
$$
\sum_{i=1}^n \res_{t_i} (\nabla_{|F})  \equiv -
 \deg F \equiv 0 \quad \mbox{mod} \quad  \Z 
$$
hence we have the lemma.
\end{proof}

\begin{Remark} For $n=4$, the data $\blambda \in \Lambda_4$ 
is special if and only if $\blambda \in \Lambda_4$ lies on  
a reflection hyperplane of a reflection $s \in W(D_4^{(1)})$.  
\end{Remark}

\subsection{Parabolic degrees}

Let us fix a series of positive rational  numbers $\balpha = (\alpha_1, \alpha_2, 
\ldots, \alpha_{2n} )$, which is called  
{\em a weight},   such that 
\begin{equation}
0 \leq \alpha_1 < \alpha_2 < \cdots< \alpha_{i} < \cdots < \alpha_{2n} < \alpha_{2n+1}=1.  
\end{equation} 
For a $(\bt, \blambda)$-parabolic connection on $\BP^1$ with the determinant 
$(L, \nabla_L)$, we can define the parabolic degree of $(E, \nabla, \varphi, l)$ 
 with respect to the weight $\balpha$ by  
\begin{eqnarray}\label{eq:pardeg}
\pardeg_{\balpha} E = 
\pardeg_{\balpha} (E, \nabla, \varphi, l) &  = &  \deg E + 
\sum_{i=1}^n \left(\alpha_{2i -1} \dim E_{t_i}/l_{i} + 
\alpha_{2i} \dim l_i \right) \\ 
& =  &\deg L + \sum_{i=1}^n 
(\alpha_{2i-1} + \alpha_{2i} ).  \nonumber
\end{eqnarray}
Let $F \subset E$ be a rank 1 subbundle of $E$ such that $ 
\nabla F \subset F \otimes \Omega^1_{\BP^1}(D(\bt))$.  
We define the parabolic degree of $(F, \nabla_{|F})$ by
\begin{equation}\label{eq:subpardeg} 
\pardeg_{\balpha} F  =   
\deg F +  \sum_{i=1}^n \left(\alpha_{2i-1} \dim F_{t_i}/l_{i} \cap F_{t_i} + 
\alpha_{2i} \dim l_i \cap F_{t_i}
\right)
\end{equation}

\begin{Definition} {\rm Fix a weight $\balpha$.  A $(\bt, \blambda)$-parabolic 
connection $(E, \nabla, \varphi, l)$ on $\BP^1$ with the determinant $(L, \nabla_L)$  
  is said to be  {\em $\balpha$-stable}
(resp. {\em $\balpha$-semistable} ) if 
for every rank-1 subbundle $F$ with $\nabla(F) \subset F \otimes \Omega^1_{\BP^1}
(D(\bt) )$ 
\begin{equation}\label{eq:stability}
\pardeg_{\balpha} F < \frac{\pardeg_{\balpha} E}{2}, \quad (\mbox{resp.} 
\pardeg_{\balpha} F \leq  \frac{\pardeg_{\balpha} E}{2} ).
\end{equation}}
(For simplicity, ``$\balpha$-stable" will be abbreviated to ``stable").
\end{Definition}

We define the coarse moduli space  by  

\begin{equation} \label{eq:modulispace}
M_n^{\balpha} (\bt, \blambda, L)= 
\left\{ (E, \nabla, \varphi, l ); 
\begin{array}{l}
\mbox{an $\balpha$-stable 
$(\bt, \blambda)$-parabolic connection } \\
\mbox{with the determinant $(L, \nabla_L)$ }
\end{array}
\right\}  / \mbox{isom.}  
\end{equation}

\subsection{Stable parabolic $\phi$-connections}

If $n \geq 4$, the moduli 
space $M_n^{\balpha} (\bt, \blambda, L)$ never becomes projective 
nor complete. 
In order to obtain a compactification of the moduli 
space $M_n^{\balpha} (\bt, \blambda, L)$, we will introduce  
the notion of a stable parabolic $\phi$-connection, or equivalently, 
a stable parabolic $\Lambda$-triple.  
Again, let us fix $(\bt,\blambda)\in T_n\times\Lambda_n$
and a line bundle $L$ on ${\bf P}^1$ with a connection
$\nabla_L:L\ra L\otimes\Omega^1_{\BP^1}(D(\bt))$.

%%%%%%%%%%%%%%%%%%%%%%%%%%%%%% Definition %%%%%%%%%%%%%%%%%%%%%%%%%%%%%%%%%%%
\begin{Definition}\rm The data 
 $(E_1,E_2,\phi,\nabla,\varphi,\{l_i\}_{i=1}^n)$
 is said to be a $(\bt,\blambda)$-parabolic $\phi$-connection
 of rank $2$ with the determinant $(L,\nabla_L)$ if
 $E_1,E_2$ are rank $2$ vector bundles on $\BP^1$
 with $\deg E_1=\deg L$,
 $\phi:E_1\ra E_2$,
 $\nabla:E_1\ra E_2\otimes\Omega^1_{\BP^1}(D(\bt))$
 are morphisms of sheaves,
 $\varphi:\bigwedge^2 E_2\stackrel{\sim}\lra L$
 is an isomorphism
 and $l_i\subset (E_1)_{t_i}$ are one dimensional subspaces
 for $i=1,\ldots,n$ such that
\begin{enumerate}
 \item $\phi(fa)=f\phi(a)$ and
  $\nabla(fa)=\phi(a)\otimes df+f\nabla(a)$
  for $f\in\cO_{\BP^1}$, $a\in E_1$,
 \item $(\varphi\otimes\mathrm{id})
 (\nabla(s_1)\wedge \phi(s_2)+\phi(s_1)\wedge\nabla(s_2))
 =\nabla_L(\varphi(\phi(s_1)\wedge\phi(s_2)))$
 for $s_1,s_2\in E_1$ and
 \item $(\res_{t_i}(\nabla)-\lambda_i\phi_{t_i})|_{l_i}=0$
  for $i=1,\ldots,n$.
\end{enumerate}
\end{Definition}

%%%%%%%%%%%%%%%%%%%%%%%%%%% Remark %%%%%%%%%%%%%%%%%%%%%%%%%%%%%%%%%%%%%%%%%%%
\begin{Definition}\rm
\item Two $(\bt,\blambda)$ parabolic $\phi$-connections
 $(E_1,E_2,\phi,\nabla,\varphi,\{l_i\})$,
 $(E'_1,E'_2,\phi',\nabla',\varphi',\{l'_i\})$
 are said  to be isomorphic to each other  if
 there are isomorphisms $\sigma_1:E_1\stackrel{\sim}\lra E'_1$,
 $\sigma_2:E_2\stackrel{\sim}\lra E'_2$ and $c\in\C\setminus\{0\}$
 such that the diagrams
\[
 \begin{CD}
  E_1 @>\phi>> E_2 \\
  @V\sigma_1V\cong V  @V\cong V\sigma_2V \\
  E'_1 @>\phi'>> E'_2
 \end{CD}
 \hspace{50pt}
 \begin{CD}
  E_1 @>\nabla>> E_2\otimes\Omega^1_{\BP^1}(D(\bt)) \\
  @V\sigma_1V\cong V  @V\cong V\sigma_2\otimes\mathrm{id}V \\
  E'_1 @>\nabla'>> E'_2\otimes\Omega^1_{\BP^1}(D(\bt)) 
 \end{CD}
 \hspace{50pt}
 \begin{CD}
  \bigwedge^2 E_2 @>\varphi>\cong > L \\
  @V\wedge^2\sigma_2V\cong V @VcV\cong V \\
  \bigwedge^2 E'_2 @>\varphi'>\cong > L 
 \end{CD}
\]
commute and $(\sigma_1)_{t_i}(l_i)=l'_i$ for $i=1,\ldots,n$.
\end{Definition}

\begin{Remark} \rm 
Assume that two vector bundles $E_1,E_2$ and 
 morphisms $\phi:E_1\ra E_2$,
 $\nabla:E_1\ra E_2\otimes\Omega^1_{\BP^1}(D(\bt))$ satisfying
 $\phi(fa)=f\phi(a)$, $\nabla(fa)=\phi(a)\otimes df+f\nabla(a)$
 for $f\in\cO_{\BP^1}$, $a\in E_1$ are given.
 If $\phi$ is an isomorphism, then
 $(\phi\otimes\mathrm{id})^{-1}\circ\nabla:
 E_1\ra E_1\otimes\Omega^1_{\BP^1}(D(\bt))$
 becomes a connection on $E_1$.
\end{Remark}

Fix rational numbers $ \alpha'_1, \alpha'_2, \ldots, \alpha'_{2n}, \alpha'_{2n+1} $
satisfying
$$
 0 \leq \alpha'_1 <  \alpha'_2  <  \cdots  <   \alpha'_{2n} <  \alpha'_{2n+1}=1
$$
and positive  integers $\beta_1, \beta_2$. 
 Setting $\balpha'=(\alpha'_1, \ldots, \alpha'_{2n}), \bbeta =(\beta_1, \beta_2)$, 
we obtain  a {\em weight}  $(\balpha', \bbeta) $ for  parabolic $\phi$-connections.   

%%%%%%%%%%%%%%%%%%%%%%%%%% Definition %%%%%%%%%%%%%%%%%%%%%%%%%%%%%%%%%%%%%

\begin{Definition}\rm Fix a sufficiently large integer $\gamma $.  
 A parabolic $\phi$-connection $(E_1,E_2,\phi,\nabla,\varphi,\{l_i\}_{i=1}^n)$
 is said to be $(\balpha', \bbeta)$-stable (resp. $(\balpha', \bbeta)$-semistable) if for any subbundles
 $F_1\subset E_1$, $F_2\subset E_2$ satisfying $\phi(F_1)\subset F_2$,
 $\nabla(F_1)\subset F_2\otimes\Omega^1_{\BP^1}(D(\bt))$ and
 $(F_1,F_2)\neq (E_1,E_2),(0,0)$,
 the inequality
 \begin{gather*}\label{mu-ineq}
  \frac{\beta_1(\deg F_1(-D(\bt))) + \beta_2 ( \deg F_2-\gamma\rank(F_2))+
  \sum_{i=1}^n \beta_1 (\alpha'_{2i-1}d_{2i-1}(F_1)+\alpha'_{2i}d_{2i}(F_1))}
  {\beta_1\rank(F_1)+\beta_2 \rank(F_2)} \\
  \genfrac{}{}{0pt}{}{<}{\text{(resp.\ $\leq$)}}
  \frac{\beta_1(\deg E_1(-D(\bt))) +\beta_2(\deg E_2-\gamma\rank(E_2))+
  \sum_{i=1}^n \beta_1(\alpha'_{2i-1}d_{2i-1}(E_1)+\alpha'_{2i}d_{2i}(E_1))}
  {\beta_1\rank(E_1)+\beta_2 \rank(E_2)}
 \end{gather*}
 holds, where $d_{2i-1}(F)=\dim ((F_1)_{t_i}/l_i\cap (F_1)_{t_i})$,
 $d_{2i}(F_1)=\dim ((F_1)_{t_i}\cap l_i)$,
 $d_{2i-1}(E_1)=\dim ((E_1)_{t_i}/l_i)(=1)$ and
 $d_{2i}(E_1)=\dim l_i(=1)$.
\end{Definition}

Define the coarse moduli space by  
\begin{equation}\label{eq:coarse-moduli-compact}
 \overline{M_n^{\balpha' \bbeta}}(\bt,\blambda,L):=
 \left\{ (E_1,E_2,\phi,\nabla,\varphi,\{l_i\})
 ; \begin{array}{l}
 \text{a $(\balpha',\bbeta)$-stable $(\bt,\blambda)$-parabolic $\phi$-connection} \\
 \text{with the determinant $(L,\nabla_L)$}
 \end{array} \right\}/\mathrm{isom}.
\end{equation}
For a given weight $(\balpha', \bbeta)$ and $1 \leq i \leq 2n$, 
define a rational number  $\alpha_i $ by 
\begin{equation}\label{eq:rel-weight}
\alpha_{i}=  \frac{\beta_1}{\beta_1 + \beta_2}\alpha'_i.
\end{equation}
Then $\balpha=(\alpha_i) $ satisfies the condition
\begin{equation}\label{eq:inequality}
0 \leq \alpha_1 < \alpha_2 < \cdots < \alpha_{2n}<  \frac{\beta_1}{(\beta_1+ \beta_2)} < 1, 
\end{equation}
hence $\balpha$ defines a weight for parabolic connections.
It is easy to see that if we take $\gamma $ sufficiently large
 $(E,\nabla,\varphi,\{l_i\})$ is $\balpha$-stable if and only if the 
 associated parabolic $\phi$-connection $(E,E,\mathrm{id}_E,\nabla,\varphi,\{l_i\})$
 is stable with respect to $(\balpha', \bbeta)$.  
 Therefore we  see that the natural map 
\begin{equation}\label{eq:map}
(E,\nabla,\varphi,\{l_i\})  \mapsto 
 (E,E,\mathrm{id}_E,\nabla,\varphi,\{l_i\})
\end{equation}
induces  an injection
\begin{equation}\label{eq:injection}
M_n^{\balpha}(\bt,\blambda,L)  \hookrightarrow 
\overline{M_n^{\balpha' \bbeta}}(\bt,\blambda,L).
\end{equation}
Conversely, assuming  that $\bbeta =(\beta_1, \beta_2)$ are given,  
for a weight $\balpha=(\alpha_i)$  satisfying the condition (\ref{eq:inequality}), we can define 
$\alpha'_i = \alpha_i \frac{\beta_1 + \beta_2}{\beta_1}$ for $1 \leq i \leq 2n$.  Since  
$0 \leq \alpha'_1 < \alpha'_2 < \cdots < \alpha'_{2n} = \alpha_{2n} \frac{\beta_1 + \beta_2}{\beta_1}< 1$, 
 $(\balpha', \bbeta)$ give a weight for parabolic $\phi$-connections.

Moreover, considering the relative setting over 
$T_n \times \Lambda_n$, we can  define  two families of 
the moduli spaces
\begin{equation}\label{eq:fam-moduli}
\overline{\pi}_n:
\overline{M_n^{\balpha'\bbeta}}(L) \lra T_n \times \Lambda_n, \quad 
\pi_n: M_n^{\balpha}(L) \lra T_n \times \Lambda_n
\end{equation}
such that the following diagram commutes;
\begin{equation}\label{eq:family-compact}
\begin{CD} 
M_n^{\balpha}(L) & \stackrel{\iota}{\hookrightarrow} & \overline{M_n^{\balpha' \bbeta }}(L)   \\
 @V \pi_n VV   @VV \overline{\pi}_n V  \\
 T_n \times \Lambda_n  @=    T_n \times \Lambda_n.  \\  
\end{CD}
\end{equation}
Here the fibers of $\pi_n$ and $\overline{\pi}_n $ 
over $(\bt, \blambda) \in T_n \times 
\Lambda_n $ are 
\begin{equation}
\pi_n^{-1}(\bt, \blambda) = M^{\balpha} (\bt, \blambda, L), \quad 
\overline{\pi}_n^{-1}(\bt, \blambda) = \overline{M^{\balpha' \bbeta}} (\bt, \blambda, L). 
\end{equation}

\subsection{The existence of moduli spaces and their properties}
\quad

The following theorem is one of our fundamental  results in this article which 
shows that the moduli spaces $\overline{M_n^{\balpha'\bbeta}}(\bt, \blambda, L)$ and 
$M_n^{\balpha}(\bt, \blambda, L)$ exist and they have good properties.

%%%%%%%%%%%%%%%%%% MAIN THEOREM %%%%%%%%%%%%%%%%%%%%%%%%%%%%%%%%%%%%%%%%%%%%%%
\begin{Theorem} \label{thm:fund}
\begin{enumerate}
\item Fix a weight $\bbeta =(\beta_1, \beta_2)$.  
For a generic weight $\balpha'$,
$\overline{\pi_n}:\overline{M_n^{\balpha' \bbeta }}(L)  \lra T_n \times \Lambda_n$
is a {\em projective} morphism. In particular, the moduli space 
$\overline{M^{\balpha'\bbeta}}(\bt, \blambda, L)$ is a {\em projective} 
algebraic scheme for all 
$(\bt, \blambda) \in T_n \times \Lambda_n$.

\item For a generic  weight $\balpha$, 
$\pi_n:M_n^{\balpha}(L)  \lra T_n \times \Lambda_n$
is a {\em smooth morphism} of relative dimension $2n-6$ with irreducible closed fibers. 
Therefore, the moduli space 
$M^{\balpha}_n(\bt, \blambda, L)$ is a {\em smooth, irreducible} algebraic 
variety of dimension $2n -6$  for all 
$(\bt, \blambda) \in T_n \times \Lambda_n$.  
\end{enumerate}
\end{Theorem}
%%%%%%%%%%%%%%%%%% END OF  MAIN THEOREM %%%%%%%%%%%%%%%%%%%%%%%%%%%%%%%%%%%%%%%%%%%%%%

The proof of Theorem \ref{thm:fund} can be separated into 3 parts.  
The construction  of the coarse moduli space 
of the parabolic $\phi$-connections over a projective smooth curve 
will be treated in Section \ref{sec:proof}. 
We deal with the relative settings  
and prove the 
projectivity of the morphism 
$ \overline{\pi_n}:\overline{M_n^{\balpha' \bbeta }}(L)  \lra T_n \times \Lambda_n$.  
(Cf. Theorem \ref{thm:exist-proj}).  Since we have a naturel embedding $M^{\balpha}_n(L) 
\hookrightarrow \overline{M_n^{\balpha' \bbeta }}(L) $, 
the existence of the moduli 
space $M^{\balpha}_n(L)$ easily follows from the first assertion. 
The smoothness of the morphism 
$\pi_n:M_n^{\balpha}(L)  \lra T_n \times \Lambda_n$ follows from 
Proposition \ref{prop:smoothness}. Finally, the irreducibility of the moduli 
space $M^{\balpha}_n(\bt, \blambda, L)$ is proved in Section \ref{sec:irr-stable}, 
(cf. Proposition \ref{prop:stable-mod-irr}),  
based on the irreducibility of  the moduli space 
$\cR(\cP_{n, \bt})_{\ba}$ proved in Proposition \ref{prop:irr-p2}.

\begin{Remark}\label{rem:fund-thm} 
\rm 
\begin{enumerate}
\item As we mentioned in Introduction, we sometimes extend the base by 
an \'etale covering $T'_n \lra T_n$ in Theorem \ref{thm:fund}, which causes
no change in the proof. 

\item The structure of 
moduli spaces $M^{\balpha}_n(L)$ and $\overline{M^{\balpha' \bbeta}_n}(L)$ may 
depend on the weight $\balpha$ and $\deg L$.

\item The moduli spaces $M^{\balpha}_n(L)$ and $\overline{M^{\balpha' \bbeta}_n}(L)$ 
are fine moduli spaces.  In fact, we have the universal families over these moduli 
spaces. See \S \ref{sec:proof}.

\item When we describe the explicit algebraic or geometric structure of the moduli spaces 
$M^{\balpha}_n(L)$ and $\overline{M^{\balpha' \bbeta}_n}(L)$, 
it is convenient to fix a 
determinant line bundle $(L, \nabla_L)$.  
As a typical example of the determinant bundle is 
\begin{equation}\label{eq:typical}
(L, \nabla_L) = (\cO_{\BP^1}(-t_n),  d) 
\end{equation}
where the connection is given by 
\begin{equation}
\nabla_L (z - t_n) = d(z- t_n) = (z - t_n) \otimes \frac{dz}{z - t_n}. 
\end{equation}
Here   $z$  is an inhomogeneous coordinate of $\BP^1 = \Spec \C[z] \cup \{ \infty\}$.  
For this $(L, \nabla_L) = (\cO_{\BP^1}(-t_n), d)$, 
we set 
$$
M_n^{\balpha} (\bt, \blambda, -1) = M_n^{\balpha} (\bt, \blambda, L), \quad  (\mbox{resp.} \  
\overline{M_n^{\balpha'\bbeta}} (\bt, \blambda,-1)= \overline{M_n^{\balpha'\bbeta}}
 (\bt, \blambda, L) \quad ). 
$$

\end{enumerate}
\end{Remark}

%%%%%%%%%%%%%%%%%subsection n=4 %%%%%%%%%%%%%%%%
\subsection{The case of  $n = 4$ (Painlev\'e $VI$ case).}
\label{subsec:n=4}
\quad

We will deal with the case of $n=4$ which corresponds to Painlev\'e $VI$ equation. 
Let us fix a sufficiently large intger $\gamma$ and take a weight $(\balpha', \bbeta)$ 
for parabolic $\phi$-connections where  $\balpha' = 
(\alpha'_1, \ldots, \alpha'_8)$, $\bbeta =(\beta_1, \beta_2), \gamma$  and fix
$(\bt, \blambda) = (t_1, \ldots, t_4, \lambda_1, \ldots,  \lambda_4)
\in T_4 \times \Lambda_4$.

Then the corresponding weight $\alpha = (\alpha_1, \ldots, \alpha_8)$ for parabolic 
connections can be given by 
 $$
\alpha_i = \alpha'_i \frac{\beta_1}{\beta_1 + \beta_2} \quad 1 \leq i  \leq 8. 
$$
Later, for simplicity,  we will assume that $\beta_1 = \beta_2$, hence 
$ \balpha = \balpha'/2$.  We also assume  $(L, \nabla_l) =   (\cO_{\BP^1}(-t_n), d)$
and in this case, we set  
$$
\overline{M^{\balpha'}_4}(\bt, \blambda, -1) = \overline{M^{\balpha' \bbeta}_4}(\bt, \blambda, L), 
\quad 
\overline{M^{\balpha'}_4}( -1) = \overline{M^{\balpha' \bbeta}_4}(L). 
$$

By Theorem \ref{thm:fund}, 
we can  obtain the commutative diagram:
%%%%%%%%%%%%%%%%%%%%%%%%%%%%%%%%%%%%%
\begin{equation}\label{eq:n=4family}
\begin{CD} 
M_4^{\balpha}(-1) & \stackrel{\iota}{\hookrightarrow} & 
\overline{M_4^{\balpha'}}(-1)   \\
 @V \pi_4 VV   @VV \overline{\pi}_4 V  \\
 T_4 \times \Lambda_4  @=    T_4 \times \Lambda_4,   \\  
\end{CD}
\end{equation}
%%%%%%%%%%%%%%%%%%%%%%%%%%%%%%%%%%%%%
such that  
$\pi_4^{-1}((\bt, \blambda)) \simeq M_4^{\balpha}(\bt, \blambda, -1 )$ 
and 
$ \overline{\pi}_4^{-1}(\bt, \blambda) \simeq \overline{M_4^{\balpha'}}(\bt, \blambda, -1 )$. 
(Note that $\balpha = \balpha'/2$). 
{F}rom Theorem \ref{thm:fund}, we see that for a generic weight 
$\balpha'$, $\overline{\pi}_4$ is a 
projective morphism and $\pi_4$ is a smooth morphism of relative dimension 
$2$.  In Part II, \cite{IIS-2}, we will give detailed descriptions of the moduli 
spaces $M_4^{\balpha}(\bt, \blambda, -1) $ and $\overline{M_4^{\balpha'}}(\bt, \blambda, -1 )$.
The following theorem shows that our family of the moduli space
$\overline{M_4^{\balpha'}}(-1)  \lra T_4 \times \Lambda_4$  
can be identified with the family of Okamoto-Painlev\'e pairs constructed by  
Okamoto \cite{O1}.  (See also \cite{Sakai},  \cite{STT02}).   
Note also that  
 Arinkin and Lysenko \cite{AL1} give   isomorphisms between
their  moduli spaces and  Okamoto spaces for generic $\blambda$.

%%%%%%%%%%%%%%%%%%Proposition%%%%%%%%%%%%%%%%%%%%%%%%
\begin{Theorem}\label{thm:n=4} (Cf. \cite{IIS-2}). 
\begin{enumerate}
\item
For a suitable choice of a weight $\balpha'$, 
the morphism 
$$
\overline{\pi}_4:\overline{M_4^{\balpha'}}(-1) \lra T_4 \times \Lambda_4
$$ 
is \underline{projective and  smooth} . 
Moreover for any $(\bt, \blambda)\in T_4 \times \Lambda_4$ the fiber 
$\overline{\pi}_4^{-1}(\bt, \blambda) := \overline{M_4^{\balpha'}}(\bt, \blambda, -1 )$ 
is irreducible, hence   a smooth projective surface.  
\item Let  $\cD = \overline{M_4^{\balpha'}}(-1) \setminus 
M_4^{\balpha}(-1) $ be the 
complement of $M_4^{\balpha}(-1)$ in $\overline{M_4^{\balpha'}}(-1)$. (Note that $\balpha = \balpha'/2$). 
Then $\cD$ is a flat reduced divisor over $T_4 \times \Lambda_4$.  

\item 
For each $(\bt, \blambda)$,  set
$$
\overline{S}_{\bt, \blambda} := \overline{\pi}_4^{-1}(\bt, \blambda) := \overline{M_4^{\balpha'}}(\bt, \blambda, -1 ).
$$
Then $\overline{S}_{\bt, \blambda}$ is  a smooth projective surface 
which can be  obtained by  blowing-ups at 
$8$ points of the Hirzebruch surface $\F_2 = \Proj(\cO_{\BP^1}(-2) \oplus \cO_{\BP^1})$
of degree $2$.  The surface  has a unique effective anti-canonical divisor 
$-K_{S_{\bt, \blambda}} = \cY_{\bt, \blambda}$ whose support is $\cD_{\bt, \blambda}$. 
Then the pair  
\begin{equation}
     (\overline{S}_{\bt, \blambda} , \cY_{\bt, \blambda})
\end{equation}
is an Okamoto-Painlev\'e pair of type $D_4^{(1)}$.  That is, the anti-canonical divisor 
$\cY_{\bt, \blambda}$ consists of $5$-nodal rational curves whose configuration 
is same as Kodaira--N\'eron degenerate elliptic curves of type $D_4^{(1)}$ (=Kodaira type 
$I_0^*$). 
Moreover we have  $(M_4^{\balpha}(-1))_{\bt, \blambda} =
(\overline{M_4^{\balpha'}}(-1))_{\bt, \blambda} \setminus \cY_{\bt, \blambda}$.
\end{enumerate}
\end{Theorem}

\section{Elementary transformation of 
parabolic connections}
\label{sec:elemnt}

In this section, we will give basic definitions and some  calculations of 
elementary transformations of stable parabolic connections.

\subsection{Definition}

Let us fix a line bundle $L$ with a connection $\nabla_L:
L \lra  L \otimes \Omega^1_{\BP^1} (D(\bt))$ and we set 
\begin{equation}
\mu_i = \res_{t_i} ( \nabla_{L}) \quad \mbox{ for }  1 \leq i \leq n.  
\end{equation}
The residue theorem implies that 
$ \sum_{i=1}^n \mu_i = - \deg L \in \Z$. 

For each  $i$, $ 1 \leq i \leq n$, we set  
$L(t_i) = L \otimes \cO_{\BP^1}(t_i) $, $L(-t_i) =   L \otimes \cO_{\BP^1}(-t_i) $ and so on.  
We will define two elementary transformations 
which induce morphisms of moduli spaces. 
\begin{eqnarray} \label{eq:upper-elm}
Elm_{t_i}^+: M_n^{\balpha}(\bt, \blambda, L) &  \lra & 
 M_n^{\balpha}(\bt, \blambda', L(t_i))  \\
\label{eq:lower-elm}
 Elm_{t_i}^-: M_n^{\balpha}(\bt, \blambda, L) &  \lra & 
 M_n^{\balpha}(\bt, \blambda", L(-t_i))  
\end{eqnarray}
Let $(E, \nabla_E, \varphi, \{l_j\}_{1 \leq j \leq n}  )$ be a 
$(\bt, \blambda)$-parabolic connection on $\BP^1$ with the determinant 
$(L, \nabla_L)$. 
Note that the eigenvalues of $\res_{t_i} (\nabla) $ are given by the following table.
\begin{equation}
E : 
\left(
\begin{array}{r|rrrrr|c}
   &t_1 & t_2 & \cdots & t_{n-1} & t_n & \wedge^2 E \\ \hline 
l_j=l_j^+& \lambda_1 & \lambda_2 & \cdots & \lambda_{n-1} & \lambda_n & L \\
E_{t_j}/l_j \simeq l_j^- & \mu_1- \lambda_1 & \mu_2- \lambda_2 & \cdots & 
\mu_{n-1}- \lambda_{n-1} & \mu_{n} - \lambda_n &  \\
\end{array}\right)
\end{equation}

\subsubsection{Definition of  $Elm_{t_i}^+$}
Take a subsheaf  $F_i$ as 
\begin{equation}\label{eq:subsheaf-1}
E(-t_i) \subset F_i \subset E  \mbox{ such that } \   
l_i = F_i /E(-t_i)  \subset E_{t_i} \quad  \mbox{and} \quad 
l_i(t_i) = F(t_i)/E   
\end{equation}
and define  
\begin{equation}
E^+_i = F_i (t_i)  =   \Ker \left[E(t_i) \lra E(t_i)/F(t_i) = E(t_i)_{t_i}/l_i(t_i) \right].\end{equation} 
Since $l_i$ is an eigenspace of $\res_{t_i} (\nabla_E)$, it is easy to see that  
$\nabla_E$ induces a connection  
\begin{equation}
\nabla_{E^+_i}:E^+_i \lra E^+_i  \otimes \Omega^1_{\BP^1}(D(\bt))
\end{equation} 
and $\varphi:\wedge^2 E \lra L $ 
induces a horizontal isomorphism $\varphi':\wedge^2 E^+_i \lra L(t_i)$. 
Moreover, one can see that the subspace 
$l_i' = E_{t_i}/l_i \subset( E^+_i)_{t_i} $
defines a new parabolic structure $\{ l'_j \}_{j=1}^n $  
with $l'_j = l_j $ for $ j \not= i$. 
Now we  define 
\begin{equation}
Elm_{t_i}^+( E) = (E^+_i, \nabla_{E^+_i}, \varphi', \{ l'_i \} ),  
\end{equation}
which is called  {\em an upper elementary transformation} of $E$ at $t_i$. 
Since $l_i' \simeq  E_{t_i}/l_i, (E^+_i)_{t_i}/l_i' \simeq l_i \otimes \cO(t_i) $, 
we see that 
$
 (\res_{t_i} (\nabla))_{|l_i'} = \mu_i - \lambda_i, \ 
 (\res_{t_i} (\nabla))_{|(E^+_i)_{t_i}/l_i'} = 
 -1 + \lambda_i.  
$
Therefore the eigenvalues of the residues of 
$\nabla_{E^+_i} $ on $Elm_{t_i}^+ (E) = E^+_i$  and the determinant $\wedge^2 E^+_i$ 
are given as follows.   
\begin{equation}
Elm_{t_i}^+ (E):
\left(
\begin{array}{c|rrrrr|c}
    t_j  &t_1 & \cdots &  t_i & \cdots  & t_n & \wedge^2 E^+_i \\ \hline 
l'_j = l_j^{'+} & \lambda_1 &  \cdots & \mu_i - \lambda_i & \cdots & \lambda_n & L(t_i)  \\
 l_j^{'-} &\mu_1 - \lambda_1 & \cdots &-1 + \lambda_{i}& 
\cdots  & \mu_n - \lambda_n &  \\
\end{array}\right). 
\end{equation}

\subsubsection{Definition of $Elm_{t_i}^{-}$}
By using (\ref{eq:subsheaf-1}) subsheaf $F_i \subset E$  we also define a filtration of sheaves
\begin{equation}
E^-_{i} =  F_i \supset E(-t_i) \supset  F_i(-t_i) 
\end{equation}
which defines a parabolic connection 
$(E^-_{i}, \nabla_{E^-_i}, \varphi', l' )$ such that
 $$
 l'_i = E(-t_i)/F_i(-t_i) = (E_{t_i}/l_i) \otimes \cO_{\BP^1}(-t_i).
 $$ 
This is called a {\em lower elementary transformation} of $E$ at $t_i$ 
and will be denoted by 
$$
Elm^-_{t_i}(E)  :=(E^-_{i}, \nabla_{E^-_{i}}, \varphi', l' ).
$$ 
Note that one has a horizontal isomorphism 
$ \varphi':\wedge^2 E^-_{i} 
\stackrel{\simeq}{\lra} L(-t_i) $ and 
the eigenvalues of the residues of  $\nabla_{E^-_i}$ on  
$Elm_{t_i}^- (E) =E^-_{i}$ are given as follows. 

\begin{equation}Elm_{t_i}^- (E): 
\left(\begin{array}{c|rrrrr|c}
     &t_1 & \cdots & t_i & \cdots & t_n & \wedge^2 E^-_i\\ \hline 
l_j'=l_j^{'+}  & \lambda_1 &   \cdots & 1+ \mu_i -\lambda_i &\cdots &  \lambda_n
& L(-t_i)  \\
l_j^{'-} & \mu_1- \lambda_1 & \cdots &  \lambda_i & \cdots & \mu_n- \lambda_n &  \\
\end{array}\right). 
\end{equation}

\subsubsection{Tensoring a line bundle $L_1$}

Let $L_1$ be a line bundle with a logarithmic connection $
\nabla_{L_1} $ and    set $
\nu_j = \res_{t_j} (\nabla_{L_1})$ for $ 1 \leq j \leq n$.  
We can define a transformation $ \otimes (L_1, \nabla_{L_1}) $ by 
\begin{equation}\label{eq:tensor}
(E, \nabla_E, \varphi, \{ l_j \} ) 
\mapsto  (E \otimes L_1, \nabla_{E \otimes L_1}, 
\varphi',  \{  l_j \otimes L_1 \} )  
\end{equation}
which induces a morphism of moduli spaces
\begin{equation}
  \otimes (L_1, \nabla_{L_1}) : M_n^{\balpha}(\bt, \blambda, L) 
  \lra M_n^{\balpha}(\bt, \blambda', L \otimes (L_1)^{\otimes 2}) . 
\end{equation}
The set of eigenvalues of new connection can be given as follows. 

\begin{equation} E \otimes L_1: 
\left(\begin{array}{c|rrrrr|c}
     &t_1 & \cdots & t_i & \cdots & t_n & \wedge^2 (E \otimes L_1)  \\ \hline 
l_j'=l_j^{'+}  & \nu_1+ \lambda_1 &   \cdots & \nu_i + \lambda_i &\cdots & \nu_n +  \lambda_n
& L \otimes( L_1)^{\otimes 2}  \\
 l_j^{'-} & \nu_1+ \mu_1- \lambda_1 & \cdots & \nu_i + \mu_i - \lambda_i & \cdots & \nu_n + \mu_n- \lambda_n &  \\
\end{array}\right). 
\end{equation}

\subsubsection{$R_i$: Interchanging the eigenspaces} 
\quad 

Under the assumption 
\begin{equation}\label{eq:non-resonance-1}
\lambda_i  \not= \mu_i - \lambda_i,  
\end{equation}
we see that there are unique eigenspaces $l^+_i = l_i $ and $l_i^{-}$    
of $\res_{t_i}(\nabla_E)$  with the eigenvalues $\lambda_i$ and $\mu_i 
- \lambda_i$ respectively.  Interchanging the 
eigenspaces $l_i^+$ and $l_i^{-}$ and 
keeping the other eigenspaces $l_j$ $j \not=i$ unchanged, 
we obtain a new parabolic connection
\begin{equation}\label{eq:eigen-change}
R_i(E) = (E, \nabla_E, \varphi, \{ l'_i \} ). 
\end{equation}
If $ \lambda_i = \mu_i -\lambda_i $, 
let us define $R_i(E) = (E, \nabla_E, \varphi, \{ l_i \} )$, that is, $R_i = Id$.   

The set of eigenvalues of new connection can be given as follows. 
\begin{equation} R_i(E) :  
\left(\begin{array}{r|rrrrr|c}
     &t_1 & \cdots & t_i & \cdots & t_n & \wedge^2  E \\ \hline 
l_j'=l_j^{'+}  &  \lambda_1 &   \cdots & \mu_i- \lambda_i &\cdots &   \lambda_n
& L  \\
E^-_{j, t_j}/l_j' \simeq l_j^{'-} &  \mu_1- \lambda_1 & \cdots & \lambda_i & \cdots &  \mu_n- \lambda_n &  \\
\end{array}\right). 
\end{equation}
  
Now assume that $\res_{t_i} (\nabla_L) \in \Z$ for all $ 1 \leq i \leq n$.   
 
\begin{Lemma}\label{lem:reflection} 
Assume that $\blambda$ is not reducible (cf. Definition \ref{def:exponents}).  
Then  $R_i$ induces an isomorphism 
\begin{equation}
R_i: M_n^{\balpha} (\bt, \blambda, L) \stackrel{\simeq}{\lra}  
M_n^{\alpha} (\bt, \blambda', L).
\end{equation}
\end{Lemma}
\begin{proof}  Since $\blambda$ is not reducible, any $(E, \nabla_E, \varphi, 
\{ l_i \}) \in M_n^{\balpha} (\bt, \blambda, L) $ are irreducible (Lemma \ref{lem:special2}), so is 
$R_i(E)$.  In particular $R_i(E)$ is $\balpha$-stable.  
Therefore it induces a morphism of moduli spaces.  Moreover it is 
obvious that $R_i^2 = Id$, so it must be an  isomorphism.  
\end{proof}

Later we will extend $R_i$ a birational map of the moduli spaces.

\subsection{Birational transformations arising from elementary 
transformations} 
\quad

\begin{Definition}
{\rm  Assume that $\balpha$ is generic.  
 {\em An affine  birational transformation} of 
the family of moduli spaces  $\pi_n: M_n^{\balpha}(L) \lra T'_n \times \Lambda_n$ is a pair of maps  $(\tilde{s}, s) $ consisting of  a birational map 
$\tilde{s}: M_n^{\balpha}(L) \cdots \lra M_n^{\balpha}(L)$ and an affine transformation $s : \Lambda_n \lra \Lambda_n$
such that the following diagram commutes:
\begin{equation}\label{eq:birational-trans}
\begin{CD}
M_n^{\balpha}(L)  & \stackrel{\tilde{s}}{\cdots \ra } &  
M_n^{\balpha}(L)  \\
\quad \downarrow \pi_n   & & \quad \downarrow \pi_n \\
T'_n \times \Lambda_n  &  \stackrel{1 \times s}{\lra} &  T'_n \times \Lambda_n.\end{CD}
\end{equation}
}
\end{Definition}

\subsubsection{The group $BL_n$}
\quad

Now we fix a determinant line bundle $(L, \nabla_L) = (\cO_{\BP^1}(-t_n), d)$ as in Remark \ref{rem:fund-thm} and consider the family of the 
 moduli spaces 
 $\pi_n: M_n^{\balpha}(\cO_{\BP^1}(-t_n)) \lra T'_n \times \Lambda_n $. 
 Let $\be_i \in \Lambda_n $ 
be the $i$-th standard base of $\Lambda_n \simeq \C^n$ and set
$\blambda = (\lambda_1, \ldots, \lambda_n) \in \Lambda_n$.  
We define a group $BL_n$  generated by the following affine automorphisms of 
$\Lambda_n$. 
\begin{equation} \left\{
 \begin{array}{rcl}
t^+_{i}(\blambda) & =  &   \blambda + \be_i 
= (\lambda_1, \ldots, \lambda_i +1,  \ldots, 
\lambda_n)  \\
t^+_{i,j}(\blambda) & = &   \blambda + \frac{1}{2}( \be_i + \be_j) 
= (\lambda_1, \ldots, \lambda_i + \frac{1}{2}, \ldots, \lambda_j + \frac{1}{2}, \ldots, 
\lambda_n) \\
t^-_{i,j}(\blambda) &  =  & 
(\lambda_1, \ldots, -\lambda_i + \frac{1}{2}, \ldots, -\lambda_j + \frac{1}{2}, \ldots, 
\lambda_n) \quad  ( 1 \leq i < j \leq  n-1) \\
t^-_{i,n}(\blambda) & = &  (\lambda_1, \ldots, -\lambda_i + \frac{1}{2}, \ldots, \ldots, \ldots, 
 - \lambda_n +  \frac{3}{2}) \\
r_i(\blambda)   &  = &    
( \lambda_1, \ldots, -\lambda_i, \ldots,  \lambda_n ) \quad
  ( 1 \leq i \leq n-1) \\
r_n(\blambda)   &  =  &   ( \lambda_1, \ldots, \lambda_i, \ldots,  1- \lambda_n ).
\end{array}\right.
\end{equation}
We can easily see the following relations. 
\begin{equation}
t^+_{i} = (t^{-}_{i,j} r_i)^2, \quad t^+_{i,j} = t^{-}_{i,j} r_i r_j.  
\end{equation}
Therefore we can define the group $BL_n$ as 
\begin{equation}\label{eq:backlund-group}
BL_n = \langle t^-_{i, j}, (1 \leq i < j \leq n),  \  r_k, ( 1 \leq k \leq n) \rangle. 
\end{equation}

In \cite{IIS-2}, we will show the following 

\begin{Proposition}\label{prop:backlund-elem}
Every element $s$  of the 
 group $BL_n$ of affine transformations of $\Lambda_n$ can 
 be lifted to a birational transformation 
 \begin{equation}
 \tilde{s}: M_n^{\balpha}(\cO_{\BP^1}(-t_n)) \cdots \lra 
 M_n^{\balpha}(\cO_{\BP^1}(-t_n))
 \end{equation}
 such that the pair $(\tilde{s}, s)$ becomes 
an affine  birational transformation of 
the family of moduli spaces.  
\end{Proposition}
\pagebreak

\section{Moduli of representations of fundamental groups}
\label{sec:monodromy}

\subsection{The family of punctured projective lines and their fundamental groups}

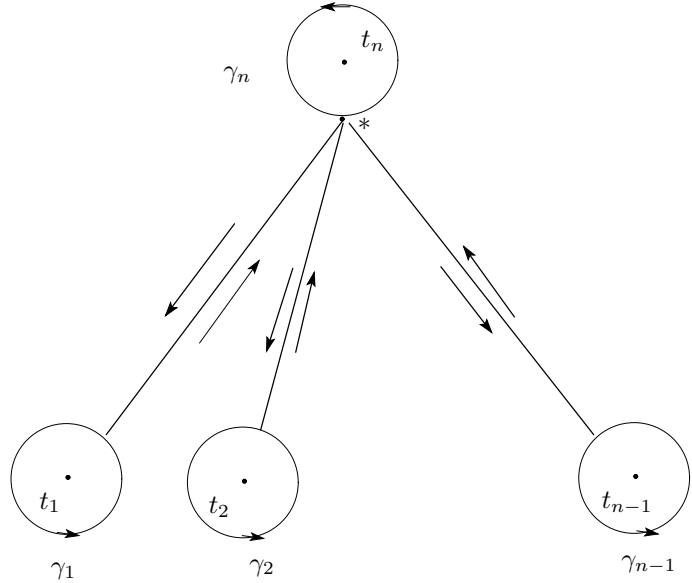
\begin{figure}
\begin{center}%WinTpicVersion2.15
\unitlength 0.1in
\begin{picture}(37.28,28.80)(12.50,-30.79)
% DOT 3 0 3 0
% 2 1724 3076 1724 3076
% 
\special{pn 4}%
\special{sh 1}%
\special{ar 1724 2676 10 10 0  6.28318530717959E+0000}%
\special{sh 1}%
\special{ar 1724 2676 10 10 0  6.28318530717959E+0000}%
% CIRCLE 3 0 3 0
% 4 1716 3084 1564 2836 1972 3084 1972 3084
% 
\special{pn 4}%
\special{ar 1716 2684 291 291  0.0000000 6.2831853}%
% DOT 3 0 3 0
% 2 2647 3095 2647 3095
% 
\special{pn 4}%
\special{sh 1}%
\special{ar 2647 2695 10 10 0  6.28318530717959E+0000}%
\special{sh 1}%
\special{ar 2647 2695 10 10 0  6.28318530717959E+0000}%
% CIRCLE 3 0 3 0
% 4 2639 3103 2487 2855 2895 3103 2895 3103
% 
\special{pn 4}%
\special{ar 2639 2703 291 291  0.0000000 6.2831853}%
% DOT 3 0 3 0
% 2 4695 3071 4695 3071
% 
\special{pn 4}%
\special{sh 1}%
\special{ar 4695 2671 10 10 0  6.28318530717959E+0000}%
\special{sh 1}%
\special{ar 4695 2671 10 10 0  6.28318530717959E+0000}%
% CIRCLE 3 0 3 0
% 4 4687 3079 4535 2831 4943 3079 4943 3079
% 
\special{pn 4}%
\special{ar 4687 2679 291 291  0.0000000 6.2831853}%
% STR 2 0 3 0
% 3 1572 3180 1572 3260 2 0
% $t_1$
\put(15.7200,-28.6000){\makebox(0,0)[lb]{$t_1$}}%
% STR 2 0 3 0
% 3 2460 3204 2460 3284 2 0
% $t_2$
\put(24.6000,-28.8400){\makebox(0,0)[lb]{$t_2$}}%
% STR 2 0 3 0
% 3 4516 3188 4516 3268 2 0
% $t_{n-1}$
\put(45.1600,-28.6800){\makebox(0,0)[lb]{$t_{n-1}$}}%
% DOT 3 0 3 0
% 2 3159 1199 3159 1199
% 
\special{pn 4}%
\special{sh 1}%
\special{ar 3159 799 10 10 0  6.28318530717959E+0000}%
\special{sh 1}%
\special{ar 3159 799 10 10 0  6.28318530717959E+0000}%
% STR 2 0 3 0
% 3 3320 710 3320 790 5 0
% $t_{n} $
\put(33.2000,-3.9000){\makebox(0,0){$t_{n} $}}%
% LINE 2 0 3 0
% 2 3164 1204 1924 2852
% 
\special{pn 8}%
\special{pa 3164 804}%
\special{pa 1924 2452}%
\special{fp}%
% LINE 2 0 3 0
% 2 3164 1220 2732 2828
% 
\special{pn 8}%
\special{pa 3164 820}%
\special{pa 2732 2428}%
\special{fp}%
% LINE 2 0 3 0
% 2 3196 1220 4476 2852
% 
\special{pn 8}%
\special{pa 3196 820}%
\special{pa 4476 2452}%
\special{fp}%
% VECTOR 2 0 3 0
% 2 2596 1748 2236 2228
% 
\special{pn 8}%
\special{pa 2596 1348}%
\special{pa 2236 1828}%
\special{fp}%
\special{sh 1}%
\special{pa 2236 1828}%
\special{pa 2292 1787}%
\special{pa 2268 1785}%
\special{pa 2260 1763}%
\special{pa 2236 1828}%
\special{fp}%
% VECTOR 3 0 3 0
% 2 1668 3364 1780 3380
% 
\special{pn 4}%
\special{pa 1668 2964}%
\special{pa 1780 2980}%
\special{fp}%
\special{sh 1}%
\special{pa 1780 2980}%
\special{pa 1717 2951}%
\special{pa 1727 2972}%
\special{pa 1711 2990}%
\special{pa 1780 2980}%
\special{fp}%
% VECTOR 3 0 3 0
% 2 2412 2372 2724 1940
% 
\special{pn 4}%
\special{pa 2412 1972}%
\special{pa 2724 1540}%
\special{fp}%
\special{sh 1}%
\special{pa 2724 1540}%
\special{pa 2669 1582}%
\special{pa 2693 1583}%
\special{pa 2701 1606}%
\special{pa 2724 1540}%
\special{fp}%
% VECTOR 3 0 3 0
% 2 2636 3380 2748 3396
% 
\special{pn 4}%
\special{pa 2636 2980}%
\special{pa 2748 2996}%
\special{fp}%
\special{sh 1}%
\special{pa 2748 2996}%
\special{pa 2685 2967}%
\special{pa 2695 2988}%
\special{pa 2679 3006}%
\special{pa 2748 2996}%
\special{fp}%
% VECTOR 3 0 3 0
% 2 4700 3356 4812 3372
% 
\special{pn 4}%
\special{pa 4700 2956}%
\special{pa 4812 2972}%
\special{fp}%
\special{sh 1}%
\special{pa 4812 2972}%
\special{pa 4749 2943}%
\special{pa 4759 2964}%
\special{pa 4743 2982}%
\special{pa 4812 2972}%
\special{fp}%
% STR 2 0 3 0
% 3 1700 3484 1700 3564 5 0
% $\gamma_1$
\put(17.0000,-31.6400){\makebox(0,0){$\gamma_1$}}%
% STR 2 0 3 0
% 3 2740 3468 2740 3548 5 0
% $\gamma_2$
\put(27.4000,-31.4800){\makebox(0,0){$\gamma_2$}}%
% STR 2 0 3 0
% 3 4764 3452 4764 3532 5 0
% $\gamma_{n-1}$
\put(47.6400,-31.3200){\makebox(0,0){$\gamma_{n-1}$}}%
% VECTOR 2 0 3 0
% 2 2900 1980 2764 2412
% 
\special{pn 8}%
\special{pa 2900 1580}%
\special{pa 2764 2012}%
\special{fp}%
\special{sh 1}%
\special{pa 2764 2012}%
\special{pa 2803 1954}%
\special{pa 2780 1961}%
\special{pa 2765 1942}%
\special{pa 2764 2012}%
\special{fp}%
% VECTOR 2 0 3 0
% 2 2916 2420 3012 2004
% 
\special{pn 8}%
\special{pa 2916 2020}%
\special{pa 3012 1604}%
\special{fp}%
\special{sh 1}%
\special{pa 3012 1604}%
\special{pa 2978 1664}%
\special{pa 3000 1656}%
\special{pa 3016 1673}%
\special{pa 3012 1604}%
\special{fp}%
% VECTOR 2 0 3 0
% 2 3676 1972 3940 2316
% 
\special{pn 8}%
\special{pa 3676 1572}%
\special{pa 3940 1916}%
\special{fp}%
\special{sh 1}%
\special{pa 3940 1916}%
\special{pa 3915 1851}%
\special{pa 3908 1874}%
\special{pa 3884 1875}%
\special{pa 3940 1916}%
\special{fp}%
% VECTOR 2 0 3 0
% 2 4060 2236 3796 1868
% 
\special{pn 8}%
\special{pa 4060 1836}%
\special{pa 3796 1468}%
\special{fp}%
\special{sh 1}%
\special{pa 3796 1468}%
\special{pa 3819 1534}%
\special{pa 3827 1511}%
\special{pa 3851 1511}%
\special{pa 3796 1468}%
\special{fp}%
% CIRCLE 3 0 3 0
% 4 3160 890 3008 642 3416 890 3416 890
% 
\special{pn 4}%
\special{ar 3160 490 291 291  0.0000000 6.2831853}%
% DOT 3 0 3 0
% 2 3170 900 3170 900
% 
\special{pn 4}%
\special{sh 1}%
\special{ar 3170 500 10 10 0  6.28318530717959E+0000}%
\special{sh 1}%
\special{ar 3170 500 10 10 0  6.28318530717959E+0000}%
% VECTOR 3 0 3 0
% 2 3200 610 3050 610
% 
\special{pn 4}%
\special{pa 3200 210}%
\special{pa 3050 210}%
\special{fp}%
\special{sh 1}%
\special{pa 3050 210}%
\special{pa 3117 230}%
\special{pa 3103 210}%
\special{pa 3117 190}%
\special{pa 3050 210}%
\special{fp}%
% STR 2 0 3 0
% 3 2610 890 2610 970 5 0
% $\gamma_{n}$
\put(26.1000,-5.7000){\makebox(0,0){$\gamma_{n}$}}%
% STR 2 0 3 0
% 3 3240 1110 3240 1190 1 0
% $\ast$
\put(32.4000,-7.9000){\makebox(0,0)[lt]{$\ast$}}%
\end{picture}%
\end{center}
\caption{Canonical generators of $\pi_1(\BP^1 \setminus D(\bt), \ast)$.}
\label{fig:fund}
\end{figure}

For $n \geq 3$, let us consider the space $
T_n = \{ (t_1, \ldots, t_n) \in (\BP^1)^n \ | \ t_i \not=t_j, (i \not= j) 
\}$
and its open subset  
\begin{equation}
W_{n} = \{ (t_1, \ldots, t_n) \in \C^n \quad | \quad t_i \not= t_j, ( i \not= j) \}. 
\end{equation}
Setting   $D(\bt) = t_1 + \cdots + t_n$ for each  
$ {\bt} = (t_1, \ldots, t_n) \in T_{n}$, we denote by  
\begin{equation}
\Gamma_{n, \bt}:= \pi_1(\BP^1 \setminus D(\bt), \ast),  
\end{equation}
the fundamental group of $\BP^1 \setminus D(\bt)$ 
with the base point $\ast$ which we take  very near to  $t_{n}$.  
It is easy to see that $ \Gamma_{n, \bt}$ is generated by $\gamma_1, \ldots, \gamma_{n-1}, \gamma_{n}$ 
in Figure \ref{fig:fund}  with one relation $\gamma_1 \gamma_2 \cdots \gamma_n = 1$.  
 This set of generators $\gamma_1, \ldots, \gamma_n$ is called  {\em canonical generators}  of $\Gamma_{n, \bt}$ 
 with respect to the ordered $n$-points $\bt$.

For each $i, 1 \leq i \leq n$, we define a divisor $\Sigma_{n,i}$ of $\BP^1 \times T_n$ as 
\begin{equation}\label{eq:divisor-section}
\Sigma_{n, i} = \{  (z, (t_1, \ldots, t_n) ) \in \BP^1\times T_n \  | \ z = t_i \ \}. 
\end{equation}
Setting $
{\mathcal P}_{n} := \left( \BP^1 \times T_{n} \right) \setminus 
\left( \cup_{i=1}^{n} \Sigma_{n, i} \right) \simeq T_{n+1}, 
$
we obtain a natural projection map which induces a smooth morphism 
\begin{equation}\label{eq:universal-fam}
\tau_{n} :{\mathcal P}_n \lra T_{n} 
\end{equation}
whose fiber ${\mathcal P}_{n, {\bt}}$ 
over  $ \bt = (t_1, \ldots, t_n)$ 
is  $ \BP^1 \setminus D(\bt) $.  
The family $\tau_{n} :{\mathcal P}_n \lra T_{n}$ 
in (\ref{eq:universal-fam}) is called  the {\em universal family of 
$n$-punctured lines}.  

By the universal covering map $\tilde{T_n} \lra T_{n}$, we can extend the family 
\begin{equation}\label{eq:pull-back-fam}
\begin{array}{ccc}
\tilde{\cP}_n & \lra & \cP_n  \\
\tilde{\tau}_n \downarrow \quad \ &    & \tau_n \downarrow  \quad \\
\tilde{T_n} & {\lra} & T_n,  \\
\end{array}
\end{equation} 
where we set $ \tilde{\mathcal P}_n = \cP_n \times_{T_n} \tilde{T_n}  $. 

Fix  a base point  $\bt_0 \in T_n$ and consider  the fundamental group 
$\pi_1( T_n, \bt_0)$.    The natural $n$-th projection 
$h_n: T_n  \lra  \BP^1$ ($ (t_1, \ldots, t_{n}) \mapsto t_n $) gives a structure of 
fiber bundle over $\BP^1$  whose  fiber at $t_n = \infty$ is isomorphic to $W_{n-1}$.   
By using the exact sequence of fundamental groups for fiber bundles, one can see 
that there exists an isomorphism
\begin{equation}
\pi_1(T_n, \bt_0) \simeq \pi_1(W_{n-1}, \bt_0).  
\end{equation}
On the other hand, it is well known that 
 the fundamental group $\pi_1(W_{n-1}, \bt_0) $ is isomorphic to the pure braid group 
 $PB_{n-1}$ of $n-1$ strings.   Therefore the pure braid group $PB_{n-1}$ acts on 
 the universal covering $\tilde{T}_n$ and also the typical fiber $\cP_{n, \bt_0}$  
 of $\tilde{\tau}_n$ in 
 (\ref{eq:pull-back-fam}).

Moreover the fiber bundle  $ \hat{\pi}_n: \tilde{\cP_n} \lra \tilde{T}_n $ becomes 
trivial, that is,  there exists a diffeomorphism 
$\tilde{\cP}_{n}\stackrel{\simeq}{\lra} \cP_{n, \tilde{\bt^0}} \times \tilde{T}_n$
such that the following diagram commutes:
\begin{equation}\label{eq:trivialization}
\begin{array}{ccccc} 
\tilde{\cP}_{n}& & \stackrel{\simeq}{\lra}& & \cP_{n, \bt_0} \times \tilde{T}_n \\
  & \tilde{\tau}_n \searrow &    & \swarrow &  \\ 
         &  &  \tilde{T}_n. & & \\
\end{array}
\end{equation}
By using the isomorphism, for every $\tilde{\bt} \in \tilde{T}_n$,  we can obtain  
the isomorphism of fundamental groups    
\begin{equation}\label{eq:canonical-isom}
\pi_1(\tilde{\cP}_{n, \tilde{\bt}} \ast)  \simeq \pi_1( \cP_{n, \bt_0}, \ast) = \Gamma_{n, 
\bt_0}
\end{equation}
as well as  the identification of  canonical generators  $\gamma_1, \ldots, \gamma_{n}$ in 
Figure \ref{fig:fund} . The action of 
the  pure braid group $PB_{n-1}$ on the fiber bundle  
$ \hat{\pi}_n: \tilde{\cP_n} \lra \tilde{T_n} $ induces an action  on canonical generators of 
$ \Gamma_{n, \bt_0}$, which  can be  written    in a very explicit way.  
(For example for the case of $n = 4$, see \cite{Iwa02-1}, \cite{Iwa02-2}).

\subsection{The moduli space of $SL_2(\C)$-representations}

\begin{Definition} {\rm \label{def:representation}
An  $SL_2(\C)$-representation of the fundamental group 
$\Gamma_{n, {\bt}} = \pi_1( \cP_{n, \bt}, \ast)$ of $\cP_{n, \bt} = \BP^1 \setminus  D(\bt)$ 
 is a group homomorphism 
\begin{equation}\label{eq:rep}
\rho: \Gamma_{n, \bt} = \pi_1(\cP_{n, \bt}, \ast)  \lra  SL_2(\C). 
\end{equation}
We denote by $\Hom(\Gamma_{n, \bt}, SL_2(\C))$ 
the set of all $SL_2(\C)$-representations of $\Gamma_{n, \bt}$. 
If we fix   a set of 
canonical  generators $\gamma_1, \ldots, \gamma_n $ of 
$ \Gamma_{n, \bt}$ as in Figure  \ref{fig:fund}, 
we have the identification   
$$
\Hom(\Gamma_{n, \bt}, SL_2(\C)) =  SL_2(\C)^{n-1}
$$
given by $\rho \mapsto (\rho(\gamma_i))$ for $ i = 1, \ldots, n-1$.  
}
\end{Definition}

\begin{Definition}{\rm \label{def:isom-equivalent}
\begin{enumerate}
\item 
Two $SL_2(\C)$-representations 
$\rho_1, \rho_2 $ are isomorphic to each other, if and 
only if 
there exists a matrix $P \in SL_2(\C)$ such that 
$$
\rho_2 (\gamma) = P^{-1} \cdot  \rho_2(\gamma) \cdot  P \quad \mbox{for all 
$ \gamma \in \pi_1(\tilde{\cP}_{n,\bt}, *) $ }.  
$$
\item A {\em semisimplification of a representation} $\rho$ 
is an associated graded of the composition series of $\rho$.

\item Two $SL_2(\C)$-representation is said to be {\em Jordan equivalent} if 
their semisimplifications are isomorphic. 
\end{enumerate}}
\end{Definition}

Fixing  $\bt_0 \in T_n $ and canonical generators 
$\gamma_1, \ldots, \gamma_n$ of $\Gamma_{n, \bt_0}$ and using the 
isomorphism in (\ref{eq:canonical-isom}),  for any $\bt \in \tilde{T}_n$,  we fix an identification  
\begin{equation}\label{eq:identify}
\Hom(\Gamma_{n, \bt}, SL_2(\C)) \stackrel{\simeq}{\lra} SL_2(\C)^{n-1}
\end{equation}  
by $\rho \mapsto (\rho(\gamma_1), \ldots, \rho(\gamma_{n-1}))$.  

Let $R_{n-1}$ denote the affine coordinate ring of $SL_2(\C)^{n-1}$ and consider the simultaneous
action of $SL_2(\C)$ on $SL_2(\C)^{n-1}$ as 
$$
(M_1, \cdots, M_{n-1}) \mapsto (P^{-1} M_1 P, \cdots, P^{-1} M_{n-1} P).  
$$
Hilbert shows that 
the ring of invariants, denoted by $\left(R_{n-1} \right)^{Ad(SL_2(\C))} $, is finitely generated.   
The following lemma is due to Simpson \cite{Simp-II},

\begin{Lemma}\label{lem:cat-quot}
 $($\cite{Mum:GIT}, $[${\rm Proposition 6.1}, \cite{Simp-II}$]$). 
For any $\bt \in \tilde{T}_n$,   under the identification (\ref{eq:identify}), 
there exists the universal categorical quotient map 
$$
\Phi_n: \Hom(\Gamma_{n, \bt}, SL_2(\C)) \simeq  SL_2(\C)^{n-1} \lra  \cR(\cP_{n, \bt})
= SL_2(\C)^{n-1}/Ad(SL_2(\C))
$$
where 
\begin{equation}\label{eq:categorical-quotient-1} 
\cR(\cP_{n, \bt}) = \Spec[ \left( R_{n-1} \right)^{Ad(SL_2(\C))}].
\end{equation}
The closed points of $\cR(\cP_{n, \bt})$ represent the Jordan equivalence classes of 
$SL_2(\C)$-representations of $\Gamma_{n, \bt}$. 
We say that $\cR\cP_n = \cR(\cP_{n, \bt})$ is the {\em moduli space} of 
$SL_2(\C)$-representation of 
$\pi_{1}(\BP^1 \setminus \Sigma(\bt))$. 
\end{Lemma}

\begin{Remark} {\rm \label{rem:structure}
Lemma \ref{lem:cat-quot} says that the set $\cR(\cP_{n, \bt})$ 
of Jordan equivalence classes  
of $SL_2(\C)$-representations admits a natural structure of an 
affine scheme. Moreover, it is easy to see that 
the moduli stack of isomorphism classes of 
$SL_2(\C)$-representations has no natural scheme structure.}
\end{Remark}

\begin{Remark} \rm It is obvious that the algebraic structure or complex structure of the  
moduli space $\cR(\cP_{n, \bt})$ does not depend on $ \bt \in \tilde{T}_n $.  However in order 
to define the isomorphism $\Hom(\Gamma_{n, \bt}, SL_2(\C)) \simeq  SL_2(\C)^{n-1}$
we have to fix canonical generators of 
 $\Gamma_{n, \bt} = \pi_1(\BP^1 \setminus D(\bt))$.  Since the pure braid group 
 $PB_{n-1}:=\pi_1(T_n, \ast)$ acts on 
 the sets of generators of $\Gamma_{n, \bt}$ and hence acts on $\cR(\cP_{n, \bt})$. 
 This action is called the \underline{\em  topological nonlinear monodromy action  } of the 
 pure braid group  $PB_{n-1}:=\pi_1(T_n, \ast)$.  (Cf. \cite{DM}, \cite{Iwa02-1}, \cite{Iwa02-2}). 
\end{Remark}

In our case, we can describe the categorical quotient 
$   \Spec[ \left( R_{n-1} \right)^{Ad(SL_2(\C))}]$ 
more explicitly. Denote the coordinate ring $R_{n-1}$ of $SL_2(\C)^{n-1}$ by 
\begin{equation}
 R_{n-1} = \C[a_{i}, b_{i}, c_{i}, d_{i} ]/(a_i d_i - b_i c_i -1) \ \ i =1, \ldots, n-1 
 \end{equation}
where 
$M_i = \left( \begin{array}{cc} a_i & b_i \\ c_i & d_i \end{array} \right) $. 

The following Proposition follows from  the  fundamental theorem  for matrix invariants. 
(See [Theorem 2, Theorem 7, \cite{For}] or [Theorem 1.3, \cite{P}]).  

\begin{Proposition}
\begin{equation}
\left( R_{n-1} \right) ^{{\rm Ad}(SL_2(\C))} = \C [ {\rm Tr}(M_{i_1} M_{i_2} \cdots M_{i_k}), 
 1 \leq i_1, \ldots, i_k \leq n-1 ].
\end{equation}
Moreover, the elements 
${\rm Tr}(M_{i_1} M_{i_2} \cdots M_{i_k})$ of degree $k \leq 3$ generate
the invariant ring, that is,  
\begin{equation}
\left(R_{n-1}\right)^{{\rm Ad}(SL_2(\C))} = \C [ {\rm Tr}(M_{i}), 
{\rm Tr}( M_{i}M_{j}), {\rm Tr}(M_{i}M_jM_k) \ | \ 
 1 \leq i, j , k  \leq n-1 ]. 
\end{equation}
\end{Proposition}

Let us  set 
\begin{equation}\label{eq:trace}
a_i = \Tr (M_i) \quad \mbox{for} \quad 1 \leq i \leq n,  
\end{equation}
which are elements of $\left(R_{n-1}\right)^{Ad(SL_2(\C))}$ 
and consider the 
subring $A_{n} = \C[a_1, \ldots, a_n]$ of $\left(R_{n-1}\right)^{Ad(SL_2(\C)}$. 
We have a natural morphism 
\begin{equation}\label{eq:isospect}
p_n: \cR(\cP_{n, \bt})  = 
\Spec \left[ \left(R_{n-1}\right)^{Ad(SL_2(\C))} \right] \lra \cA_n = 
\Spec \left[ A_n \right].
\end{equation}

\subsection{Construction of the family of moduli spaces $\phi_n : \cR_n \lra T_n' \times \cA_n $} 
\quad

Fix  ${\bt}_0 \in T_{n}$ as the base point of fundamental group 
$\pi_1(T_n, \bt_0)$ and fix canonical generators $\gamma_1,  \ldots, \gamma_{n}$ of 
$\Gamma_{n, \bt_0}$.  Again taking  the universal covering map $\tilde{T}_n \lra T_n$, we can 
obtain a trivialization (\ref{eq:trivialization}) and  isomorphisms of the fundamental groups 
(\ref{eq:canonical-isom}).  By using the isomorphisms, for each $\bt \in \tilde{T}_n$, we 
obtain a canonical isomorphism 
$$
\cR(\cP_{n, \bt}) \simeq \cR(\cP_{n,\bt_0}).
$$  
Moreover the group $\pi_1(T_{n}, \bt_0) \simeq PB_{n-1}$  acts on 
the variety  $ \cR(\cP_{n, \bt_0})$ as the group of nonlinear monodromies  and hence defines 
  the action on 
the product $ \cR(\cP_{n, \bt_0}) \times \tilde{T_n}$.  
Define  the subgroup $\Gamma_{n-1}$ of $\pi_1(T_{n}, \bt_0)$ as a kernel of the natural 
homomorphism $ \pi_1(T_{n}, \bt_0)  \lra \Aut (\C[a_1, \ldots, a_n ] )$.    
It is easy to see that $\Gamma_{n-1}$ is a subgroup of $\pi_1(T_{n}, \bt_0)$ of finite index, so defining 
as $T'_n = \tilde{T}_n/\Gamma_{n-1}$ we obtain the finite \'etale covering 
\begin{equation}\label{eq:finite-etale}
T'_n := \tilde{T}_n/\Gamma_{n-1}  \lra T_n.  
\end{equation}
Consider the natural action of  $\Gamma_{n-1}$ 
on the product $\tilde{T}_n \times \cR(\cP_{n, \bt_0})$. 
The natural map $1 \times  p_n: \tilde{T}_n \times \cR(\cP_{n, \bt_0}) \lra  \tilde{T}_n 
\times \cA_n$  is clearly 
 equivariant  with respect to the action of  $\Gamma_{n-1}$, where  $\Gamma_{n-1}$ acts 
on $\cA_n$ as the identity map.  Setting 
\begin{equation}\label{eq:free-quotient}
\cR_n = \tilde{T}_n \times \cR(\cP_{n, \bt_0})/\Gamma_{n-1}, 
\end{equation}
we obtain a morphism 
\begin{equation}\label{eq:family-rep-s4}
\fbox{ \quad $\phi_n: \cR_n \lra  T'_n  \times \cA_n $,  \quad } 
\end{equation}
which is said to be the family of the moduli spaces of $SL_2$-representations of 
the fundamental group. 
The fiber of $\phi_n$ at $(\bt, \ba)$ is given by the affine subscheme of $\cR_n$ 
\begin{equation}\label{eq:mod-rep}
\phi_n^{-1}(\bt, \ba) = \cR(\cP_{n, \bt})_{\ba} := 
\{ [\rho] \in  \cR(\cP_{n, \bt})  \  |  \   \Tr[ \rho(\gamma_i) ] = a_i, 1 \leq i \leq n \}. 
\end{equation}
Since $a_i$ determines the eigenvalues of monodromy matrix $\rho(\gamma_i)$, 
$\ba$ may be considered as the set of spectral of local monodromies. Hence 
the space $\cR(\cP_{n, \bt})_{\ba}$ is said to be the moduli space of {\em isospectral} 
$SL_2$-representations.  Note that though  
the moduli space $M_n^{\balpha}(\bt, \blambda, L)$ is smooth for all $(\bt, \blambda)$ 
if $\ba$ is special in the sense of Definition \ref{def:exponents-intro}
the affine scheme $\cR(\cP_{n, \bt})_{\ba}$ has singularities.

In \S \ref{sec:rep-irred}, 
we will prove the following 
\begin{Proposition} \label{prop: irreducible} For any $\ba \in \cA_n$, 
the scheme $\cR(\cP_{n, \bt})_{\ba}$ in (\ref{eq:mod-rep}) is irreducible.  
\end{Proposition}

\subsection{The case of $n=4$.}

Now we recall the explicit description of the invariant ring for $n = 4$ 
due to Iwasaki (\cite{Iwa02-1}, \cite{Iwa02-2}). 
We denote by $(i, j, k)$ a cyclic permutation  of $(1, 2, 3)$.  Then 
the invariant ring $\left(R_3 \right)^{Ad(SL_2(\C))}$ is generated by 
\begin{equation}
\begin{array}{ccll}
x_i & = & \Tr[M_j M_k] &  \mbox{for $i = 1, 2, 3$} \\
a_i & = & \Tr[M_i]   & \mbox{for $i = 1, 2, 3$} \\
a_4 & = & \Tr[M_1 M_2 M_3]  &  \\
\end{array}
\end{equation}

The following proposition is proved in \cite{Iwa02-2}. 

\begin{Proposition} The invariant 
ring  $ \left( R_3 \right)^{Ad(SL_2(\C))} $ is generated by 
seven elements 
$
x_1, x_2, x_3, a_1, a_2, a_3, a_4
$
and there exists a relation 
\begin{equation}
f(\x, \ba) = x_1 x_2 x_3 + x_1^2 + x_2^2 + x_3^2 
- \theta_1(\ba) x_1 - \theta_2(\ba) x_2- \theta_3(\ba) x_3 +  \theta_4(\ba), 
\end{equation}
where we set
\begin{eqnarray}
\theta_i(\ba) & = & a_i a_4 + a_j a_k , \quad (i, j, k) = \mbox{a cyclic permutation of  }  (1, 2, 3), \\
\theta_4(\ba) & = & a_1a_2 a_3 a_4 + a_1^2 + a_2^2 + a_3^2 + a_4^2 - 4. 
\end{eqnarray}
Therefore we have an isomorphism
\begin{equation}
\left( R_3 \right) ^{Ad(SL_2(\C))} \simeq \C[x_1, x_2, x_3, a_1, a_2, a_3, a_4]/(f(\z, \ba)).
\end{equation}
\end{Proposition}

Recall that fixing canonical generators of the fundamental group, for any $ \bt \in \tilde{T}_4  $, 
the categorical quotient $ \cR_{4,\bt} $ is given by 
$
\cR(\cP_{4, \bt}) :=\Spec[\left( R_3\right)^{Ad(SL_2(\C))}] \simeq 
\Spec[\C[\x, \ba]/(f(\x, \ba))].
$
Setting $\cA_4 = \C^4 = \Spec[\C[a_1, \ldots, a_4]]$, as in (\ref{eq:isospect}) 
we have a surjective morphism 
$$
p_4:\cR(\cP_{4, \bt})  = \Spec[\C[\x, \ba]/(f(\x, \ba))] \lra \cA
$$
whose fiber at $\ba \in \cA$ is an affine cubic hypersurface in 
$ \C^3 $
$$
\cR(\cP_{4, \bt})_{\ba}  \simeq \{  (x_1, x_2, x_3) \in \C^3 \  | \ f(\x, \ba) = 0 \} \subset \C^3.  
$$
Therefore,  the family in (\ref{eq:family-rep-s4}) 
$\phi_4: \cR_4 \lra T'_4 \times \cA_4$
is  a family of affine cubic hypersurfaces in $\C^3$.

The subgroup $\Gamma_3$ of $\pi_1(T_4, \bt_0)$ acts both 
on the space $\cR(\cP_{4, \bt})$ and the space $\cR(\cP_{4, \bt})_{\ba} $
as nonlinear monodromies.   Iwasaki \cite{Iwa02-1} showed  the following 

\begin{Proposition} \label{prop:nonlinear-fix} 
There exists a one-to-one correspondence between the set of 
fixed points 
of the action of $\Gamma_{3}$ on  $\cR(\cP_{4, \bt})_{\ba} $ and the set of 
singular points on the affine cubic hypersurface $\cR(\cP_{4, \bt})_{\ba} $. 
\end{Proposition}

\vspace{1cm} 
\section{Construction of the moduli space $\overline{M_n^{\balpha'\bbeta}}(\bt, \blambda, L)$ 
and Proof of Theorem  \ref{thm:fund}, (1)}
\label{sec:proof}

%%%%%%%%%%%%%%%%%%%%%%%%%%%%% Subsection{section:tr-definition} %%%%%%%%%%%%%%
%%%%%%%%%%%%%%%%%%%%%%%%%%%%%%%%%%%%%%%%%%%%%%%%%%%%%%%%%%%%%%%%%%%%%%%%%%%%%%%
\subsection{Translation of the definition of parabolic $\phi$-connection}
%%%%%%%%%%%%%%%%%%%%%%%%%%%%%%%%%%%%%%%%%%%%%%%%%%%%%%%%%%%%%%%%%%%%%%%%%%%%%%%

In this section, we will translate the definition of
parabolic $\phi$-connection, since it is rather convenient
to generalize the definition for the construction of
the moduli space.

Let $X$ be a smooth projective curve over $\C$
and $D$ be an effective divisor on $X$.

We define an $\cO_X$-bimodule structure on
 $\Lambda^1_D= \cO_X\oplus
 (\Omega^1_X(D))^{\vee}$ by
\begin{align}\label{bimodule}
 (a,v)f&:= (fa+ \langle v,df \rangle ,fv) \\
 f(a,v)&:= (fa,fv) \notag
\end{align}
 for $a,f\in\cO_X$ and $v\in(\Omega^1_X(D))^{\vee}$,
 where $\langle\,,\,\rangle:
 (\Omega^1_X(D))^{\vee}
 \times\Omega^1_X(D)\ra\cO_X$
 is the canonical pairing.

%%%%%%%%% Definition of generalized parabolic $\phi$-connection 
%%%%%%%%%
%%%%%%%%%%%%%%%%%%%%%%%%%%%%%%%%%%%%%%%%%%%%%%%%%%%%%%%%%%%%%%%%%%%%%%%%%
%%%%%
\begin{Definition}\rm
 A parabolic $\Lambda^1_D$-triple $(E_1,E_2,\Phi,F_*(E_1))$ on $X$
 consists of two vector bundles $E_1,E_2$ on $X$, a left 
$\cO_X$-homomorphism
 $\Phi:\Lambda^1_D\otimes_{\cO_X}E_1 \ra E_2$
 and a filtration of coherent subsheaves:
 $E_1= F_1(E_1)\supset F_2(E_1)\supset\cdots\supset F_l(E_1)\supset
 F_{l+1}(E_1)= E_1(-D)$.
\end{Definition}

%%%%%%%%%%%%%%%%%%%%%%%%% Remark 
%%%%%%%%%%%%%%%%%%%%%%%%%%%%%%%%%%%%%%%%%%
\begin{Remark}
Assume that two vector bundles $E_1,E_2$ on $X$ are given.
Then giving morphisms $\phi:E_1\ra E_2$,
$\nabla:E_1\ra E_2\otimes\Omega^1_X(D)$ satisfying
$\phi(fa)= f\phi(a)$, $\nabla(fa)= \phi(a)\otimes df+f\nabla(a)$
for $f\in\cO_X$, $a\in E_1$
is equivalent to giving a left $\cO_X$-homomorphism
$\Phi:\Lambda^1_D\otimes_{\cO_X}E_1 \ra E_2$.
\end{Remark}

%%%%%%%%%%%%%%%%%% Definition of subconnection %%%%%%%%%%%%%%%%%%%
\begin{Definition}\rm
 A parabolic $\Lambda^1_D$-triple $(E'_1,E'_2,\Phi',F_*(E'_1))$
 is said to be a parabolic $\Lambda^1_D$-subtriple of
 $(E_1,E_2,\Phi,F_*(E_1))$ if
 $E'_1\subset E_1$, $E'_2\subset E_2$,
 $\Phi|_{\Lambda^1_D\otimes E'_1}= \Phi'$ and
 $F_i(E'_1)\subset F_i(E_1)$ for any $i$.
\end{Definition}

Fix rational numbers 
$0\leq \alpha'_1<\alpha'_2<\cdots<\alpha'_l<\alpha'_{l+1}= 1$
and positive integers $\beta_1,\beta_2$.  
We write $\balpha'= (\alpha'_1,\ldots,\alpha'_{l})$ and
$\bbeta= (\beta_1,\beta_2)$. 
We also fix an ample line bundle $\cO_X(1)$ and a rational number  $\gamma$
with $\gamma\gg 0$.

%%%%%%%%%%%%%%%%%%%%% Definition of $\mu$ %%%%%%%%%%%%%%%%%%%%%%%%
\begin{Definition}\rm
 For a parabolic $\Lambda^1_D$-triple $(E_1,E_2,\Phi,F_*(E_1))$,
 we put
$$
\begin{array}{ll}
 \mu(E_1,E_2,\Phi,F_*(E_1)):=   
 & \frac{\beta_1\deg E_1(-D) + \beta_2\deg E_2
 -\beta_2\gamma\deg\cO_X(1)\rank E_2 +
 \sum_{i= 1}^{l} \beta_1\alpha'_i \length(F_i(E_1)/F_{i+1}(E_1))}
 {\beta_1\rank E_1 +\beta_2\rank E_2}. \\
 \end{array}
$$
\end{Definition}

%%%%%%%%%%%%%%%%%%%%%% Definition of stability 
%%%%%%%%%%%%%%%%%%%%%%%%%%%%%
\begin{Definition}\rm Assume that $\gamma$ is sufficiently large. 
A parabolic $\Lambda^1_D$-triple $(E_1,E_2,\Phi,F_*(E_1))$ is
 $(\balpha', \bbeta)$-stable (resp.\ $(\balpha', \bbeta)$-semistable) 
 if for any non-zero proper parabolic $\Lambda^1_D$-subtriple
 $(E'_1,E'_2,\Phi',F_*(E'_1))$ of $(E_1,E_2,\Phi,F_*(E_1))$,
 the inequality
 \begin{eqnarray*}
  \mu(E'_1,E'_2,\Phi',F_*(E'_1))
  &<& \mu(E_1,E_2,\Phi,F_*(E_1)) \\
  &{\text{(resp.\ $\leq$)}}&
 \end{eqnarray*}
holds. (If we fix a weight $(\balpha', \bbeta)$, ``$(\balpha', \bbeta)$-stable (resp.\ $(\balpha', \bbeta)$-semistable) " 
may  be abbreviated to ``stable (resp.\ semistable)" for simplicity.) 
\end{Definition}

Let $S$ be a connected noetherian scheme and
$\pi_S:\cX\ra S$ be a smooth projective morphism
whose geometric fibers are curves of genus $g$. 
Let $\cD\subset \cX$ be an effective Cartier divisor which is flat over 
$S$.
A similar formula to (\ref{bimodule}) enables us to consider
the $\cO_{\cX}$-bimodule structure on
$\Lambda^1_{\cD/S}:= \cO_{\cX}\oplus(\Omega^1_{\cX/S}(\cD))^{\vee}$.

Fix rational numbers 
$0\leq\alpha'_1<\alpha'_2<\cdots<\alpha'_l<\alpha'_{l+1}= 1$,
positive integers $r,d,\{d_i\}_{1\leq i\leq l}$,
$\beta_1$, $\beta_2$, $\gamma$ with $\gamma\gg 0$.

%%%%%%%%%%%%%%%%%%%%%%%%% Definition of Moduli functor 
%%%%%%%%%%%%%%%%%%%%%%%%
\begin{Definition}\rm
We define the moduli functor
$\overline{\cM_{\cX/S}^{\cD,\balpha',\bbeta,\gamma}}(r,d,\{d_i\})$
of the category of locally noetherian schemes over $S$ to
the category of sets by
\begin{equation}\label{eq:functor}
 \overline{\cM_{\cX/S}^{\cD,\balpha',\bbeta,\gamma}}(r,d,\{d_i\})(T)
 := \{(E_1,E_2,\Phi,F_*(E_1))\}/\sim ,
\end{equation}
where $T$ is a locally noetherian scheme over $S$ and
\begin{enumerate}
\item
 $E_1,E_2$ are vector bundles on $\cX\times_ST$ such that
 for any geometric point $s$ of $T$,
 $\rank(E_1)_s= \rank(E_2)_s= r$, $\deg (E_1)_s= \deg (E_2)_s= d$,
\item
 $\Phi:\Lambda^1_{\cD/S}\otimes_{\cO_{\cX}}E_1 \to E_2$ is
 a homomorphism of left $\cO_{\cX\times_ST}$-modules,
\item
 $E_1= F_1(E_1)\supset F_2(E_1)\supset\cdots\supset
 F_l(E_1)\supset F_{l+1}(E_1)= E_1(-\cD_T)$
 is a filtration of $E_1$ by coherent subsheaves
 such that each $E_1/F_{i+1}(E_1)$ is flat over $T$
 and for any geometric point $s$ of $T$,
 $\length ((E_1/F_{i+1}(E_1))_s)= d_i$,
\item
 for any geometric point $s$ of $S$,
the parabolic $\Lambda^1_{\cD_s}$-triple
 $((E_1)_s,(E_2)_s,\Phi_s,F_*(E_1)_s)$ is stable (that is, $(\balpha', \bbeta)$-stable) .
\end{enumerate}

$(E_1,E_2,\Phi,F_*(E_1))\sim (E'_1,E'_2,\Phi',F_*(E'_1))$
if there exist a line bundle $\cL$ on $T$ and isomorphisms
$\sigma_j:E_j\stackrel{\sim}\to E'_j\otimes\cL$
for $j= 1,2$ such that
$\sigma_1(F_{i+1}(E_1))= F_{i+1}(E'_1)\otimes\cL$ for any $i$
and the diagram
\[
 \begin{CD}
  \Lambda^1_{\cD/S}\otimes_{\cO_X}E_1 @>\Phi>> E_2 \\
  @V\mathrm{id}\otimes\sigma_1V\cong V 
  @V\sigma_2 V\cong V \\
  \Lambda^1_{\cD/S}\otimes_{\cO_X}E'_1\otimes_S\cL 
  @>\Phi'\otimes\mathrm{id}>>  E'_2\otimes_S\cL
 \end{CD}
\]
commutes.
\end{Definition}

We call $(E_1,E_2,\Phi,F_*(E_1))$ a flat family of parabolic
$\Lambda^1_{\cD_T/T}$-triples on $\cX_T\times T$ over $T$
if it satisfies the above conditions $(1),(2)$ and $(3)$.

%%%%%%%%%%%%%%%%%%%%%%% Subsection %%%%%%%%%%%%%%%%%%%%%%%%%%%%%%%%
%%%%%%%%%%%%%%%%%%%%%%%%%%%%%%%%%%%%%%%%%%%%%%%%%%%%%%%%%%%%%%%%%%%
\subsection{Boundedness and Openness of stability}

%%%%%%%%%%%%%%%%%%%%%%%%%%%%%%%%%%%%%%%%%%%%%%%%%%%%%%%%%%%%%%%%%%%%%

%%%%%%%%%%%%%%%%%%%%%%%%% Proposition %%%%%%%%%%%%%%%%%%%%%%%%%%%%%
\begin{Proposition}\label{bounded}
 The family of geometric points of
 $\overline{\cM_{\cX/S}^{\cD,\balpha',\bbeta,\gamma}}(r,d,\{d_i\})$
 is bounded.
\end{Proposition}

\begin{proof}
Take any geometric point
$(E_1,E_2,\Phi,F_*(E_1))\in
\overline{\cM_{\cX/S}^{\cD,\balpha',\bbeta,\gamma}}(r,d,\{d_i\})(K)$.
By Serre duality, we have
\[
 H^1(\cX_K,E_1(m-1))= \Hom(E_1,\omega_{\cX_K}(1-m))^{\vee}.
\]
Take any nonzero homomorphism
$f:E_1\rightarrow\omega_{\cX_K}(1-m)$.
Then $(\ker f,E_2,\Phi|_{\ker f},F_*(E_1)\cap\ker f)$
becomes a parabolic $\Lambda^1_{\cD_K}$-subtriple of
$(E_1,E_2,\Phi,F_*(E_1))$.
Thus we must have the inequality
\begin{gather*}
 \mu(\ker f,E_2,\Phi|_{\ker f},F_*(E_1)\cap\ker f)
 < \mu(E_1,E_2,\Phi,F_*(E_1)).
\end{gather*}
Since $\deg(\ker f)\geq \deg E_1+m-2g+1$, we can find an integer $m$
which depends only on $r,d,d_i,\bbeta,\balpha',\gamma$, $\cX$ and $\cD$
such that $\Hom(E_1,\omega_{\cX_K}(1-m))= 0$.
Then all $E_1$ become $m$-regular.

Similarly we can find an integer $m'$ such that
$E_2$ are all $m'$-regular.
Then the family of $(E_1,E_2)$ is bounded and
the boundedness of the family of $(E_1,E_2,\Phi,F_*(E_1))$
can be deduced from it.
\end{proof}

We put $\epsilon_i:= \alpha'_{i+1}-\alpha'_i$
for $i= 1,\ldots,l$.
Take an $S$-ample line bundle
$\cO_{\cX}(1)$ on $\cX$.

%%%%%%%%%%%%%%%%%%%%%%%%% Proposition %%%%%%%%%%%%%%%%%%%%%%%%%%%%%
\begin{Proposition}\label{fundamental-lemma}
 There exists an integer $m_0$ such that
 for any geometric point $(E_1,E_2,\Phi,F_*(E_1))\in
 \overline{\cM_{\cX/S}^{\cD,\balpha',\bbeta,\gamma}}(r,d,\{d_i\})(K)$,
 the inequality
\begin{align*}
 &\frac{\beta_1\alpha'_1 h^0(E'_1(m))+\beta_2 h^0(E'_2(m-\gamma))+
 \sum_{i= 1}^l \beta_1\epsilon_i h^0(F_{i+1}(E'_1)(m))}
 {\beta_1\rank(E'_1)+\beta_2\rank(E'_2)} \\
 &<
 \frac{\beta_1\alpha'_1 h^0(E_1(m))+\beta_2 h^0(E_2(m-\gamma))+
 \sum_{i= 1}^l \beta_1\epsilon_i h^0(F_{i+1}(E_1)(m))}
 {\beta_1\rank(E_1)+\beta_2\rank(E_2)}
\end{align*}
 holds for any proper non-zero parabolic
 $\Lambda^1_{\cD_K}$-subtriple $(E'_1,E'_2,\Phi',F_*(E'_1))$
 of $(E_1,E_2,\Phi,F_*(E_1))$ and any integer $m\geq m_0$.
\end{Proposition}

\begin{proof}
By Proposition \ref{bounded}, there exists an integer $N_1$
such that for any geometric point
$(E_1,E_2,\Phi,F_*(E_1))$ of
$\overline{\cM_{\cX/S}^{\cD,\balpha',\bbeta,\gamma}}(r,d,\{d_i\})$,
$h^i(F_j(E_1)(m))= h^i(E_2(m-\gamma))= 0$ for $i>0$,
$1\leq j\leq l+1$ and $m\geq N_1$.
There also exists an integer $e$ such that for any
geometric point $(E_1,E_2,\Phi,F_*(E_1))$ of
$\overline{\cM_{\cX/S}^{\cD,\balpha',\bbeta,\gamma}}(r,d,\{d_i\})$
and for any coherent subsheaf $E'$ of
$E_1^{\oplus\beta_1}\oplus E_2^{\oplus\beta_2}(-\gamma)$,
the inequality
\[
 \deg E' \leq \rank E'
 (\mu(E_1^{\oplus\beta_1}\oplus E_2(-\gamma)^{\oplus\beta_2})+e)
\]
holds.
Note that we write $\mu(E):= \rank(E)^{-1}\deg(E)$
for a vector bundle $E$.
Applying \cite{MY}, Lemma 2.6 to the case
\begin{gather*}
 P(m)= \frac{\beta_1\alpha'_1\chi(E_1(m))+\beta_2\chi(E_2(m-\gamma))
 +\sum_{i= 1}^l \beta_1\epsilon_i\chi(F_{i+1}(E_1)(m))}
 {\beta_1\rank E_1 + \beta_2\rank E_2}-1, \\
 r= \rank(E_1^{\oplus\beta_1}\oplus E_2^{\oplus\beta_2}),
 \quad
 a= \mu(E_1^{\oplus\beta_1}\oplus E_2(-\gamma)^{\oplus\beta_2})+e,
\end{gather*}
we can take integers $L,M$ such that $M\leq a$ and for any
integer $m\geq L$, the inequality
\[
 h^0(E'(m))\leq \rank(E')P(m)
\]
holds for any vector bundle $E'$ on a fiber of $\cX$ over $S$
satisfying $0<\rank(E')<\beta_1\rank(E_1)+\beta_2\rank(E_2)$,
$\mu(E')\leq M$ and
$\deg\tilde{E}' \leq a\rank(\tilde{E}')$
for any proper nonzero coherent subsheaf $\tilde{E}'$ of $E'$.

Now we put
\[
 \cG:= \left\{ E' \left|
 \begin{array}{l}
  \text{there exists a geometric point $(E_1,E_2,\Phi,F_*(E_1))$
  of $\overline{\cM_{\cX/S}^{\cD,\balpha',\bbeta,\gamma}}(r,d,\{d_i\})$} 
\\
  \text{such that $E'$ is a subbundle of
  $E_1^{\oplus\beta_1}\oplus E_2(-\gamma)^{\oplus\beta_2}$
  and $\mu(E')\geq M$}
 \end{array}
 \right\}\right..
\]
Then $\cG$ is bounded.
Thus there exists an integer $L'\geq L$ such that
for any $E'\in\cG$ and any $m\geq L'$,
$E'(m-\gamma)$ is generated by its global sections,
$h^i(E'(m-\gamma))= h^i((F_j(E_1)\cap E')(m))= 0$
for $i>0$ and $1\leq j\leq l+1$.
If we put
\[
 \tilde{\cG}:= \left\{ (E'_1,E'_2) \left|
 \begin{array}{l}
  \text{$E'_1\subset E_1$ (resp.\ $E'_2\subset E_2$)
  is a subbundle such that} \\
  \text{$\Phi(\Lambda^1_{\cD/S}\otimes E'_1)\subset E'_2$ and
  $\mu(E'_1{} ^{\oplus\beta_1}\oplus E'_2(-\gamma)^{\oplus\beta_2})\geq 
M$}
 \end{array}
 \right\}\right.,
\]
then the set of polynomials
\[
 \left\{ \beta_1\alpha'_1\chi(E'_1(m))+\beta_2\chi(E'_2(m-\gamma))
 +\sum_{i= 1}^l\beta_1\epsilon_i\chi((F_{i+1}(E_1)\cap E'_1)(m))
 \right\}_{(E'_1,E'_2)\in\tilde{\cG}}
\]
is finite, because
${E'_1}^{\oplus\beta_1}\oplus E'_2(-\tau)^{\oplus\beta_2}\in\cG$
for any $(E'_1,E'_2)\in\tilde{\cG}$.
Thus there exists an integer $m_0\geq L'$ such that
for any $m\geq m_0$ and for any $(E'_1,E'_2)\in\tilde{\cG}$,
the inequality
\[
 \frac{\beta_1\alpha'_1\chi(E_1'(m))+\beta_2\chi(E'_2(m-\gamma))
 +\sum_{i= 1}^l\beta_1\epsilon_i\chi((F_{i+1}(E_1)\cap E'_1)(m))}
 {\beta_1\rank(E'_1)+\beta_2\rank(E'_2)}
 <P(m)+1
\]
holds.
We can easily see that this $m_0$ satisfies the desired condition.
\end{proof}

%%%%%%%%%%%%%%%%%%%%%%%%%% Proposition 
%%%%%%%%%%%%%%%%%%%%%%%%%%%%%%%%%%%%%
\begin{Proposition}\label{prop:opennes-of-stability}
 Let $T$ be a noetherian scheme over $S$ and
 $(E_1,E_2,\Phi,F_*(E_1))$ be a flat family of parabolic
 $\Lambda^1_{\cD_T/T}$-triples on $\cX\times_S T$ over $T$.
 Then there is an open subscheme $T^s$ of $T$ such that
\[
 T^s(k)= \left\{ t\in T(k) \left|
 \text{$(E_1,E_2,\Phi,F_*(E_1))\otimes k(t)$ is stable}
 \right\}\right.
\]
for any algebraically closed field $k$.
\end{Proposition}

\begin{proof}
We may assume that $T$ is connected.
Put $P_1(m):= \chi((E_1\otimes k(s))(m))$,
$P_2(m):= \chi((E_2\otimes k(s))(m-\gamma))$ and
$P_1^{(i)}(m):= \chi((F_i(E_1)\otimes k(s))(m))$
for a geometric point $s$ of $T$.
Since the family
\[
 \cG= \left\{ E' \left|
 \begin{array}{l}
  \text{$E'$ is a subbundle of
  $(E_1^{\oplus\beta_1}\oplus E_2(-\gamma)^{\oplus\beta_2})\otimes k(s)$
  for some geometric} \\
  \text{point $s$ of $T$ and
  $\mu(E')\geq \mu((E_1,E_2,\Phi,F_*(E_1))\otimes k(s))$}
 \end{array}
 \right\}\right.
\]
is bounded, the family
\[
 \tilde{\cG}= \left\{ (E'_1,E'_2,\Phi',F_*(E'_1)) \left|
 \begin{array}{l}
  \text{$(E'_1,E'_2,\Phi',F_*(E'_1))$ is a parabolic
  $\Lambda^1_{\cD_s}$-subtriple of} \\
  \text{$(E_1,E_2,\Phi,F_*(E_1))\otimes k(s)$
  for some geometric point $s$ of $T$ } \\
  \text{such that $E'_1\subset E_1\otimes k(s)$
  (resp.\ $E'_2\subset E_2\otimes k(s)$) is a subbundle} \\
  \text{and $\mu(E'_1,E'_2,\Phi',F_*(E'_1))\geq
  \mu((E_1,E_2,\Phi,F_*(E_1))\otimes k(s))$}
 \end{array}
 \right\}\right.
\]
is also bounded.
So the set of sequences of polynomials
\[
 \cP:= \left\{
 \left(\chi(E'_1(m)),\chi(E'_2(m-\gamma)),
 (\chi(F_{i+1}(E'_1)(m)))_{1\leq i\leq l}\right) \left|
 (E'_1,E'_2,\Phi',F_*(E'_1))\in\tilde{\cG}
 \right\}\right.
\]
is finite.
For each $\BP':= (P'_1,P'_2,((P'_1)^{(i+1)}))\in\cP$, put
\[
Q:= \Quot^{P_1-P'_1}_{E_1/\cX_T/T}\times_T\Quot^{P_2-P'_2}_{E_2/\cX_T/T}
\]
Let $(E_1)_Q \stackrel{\pi_1}\ra G_1$ and
$(E_2)_Q \stackrel{\pi_2}\ra G_2$
be the universal quotient sheaves.
We put
\[
 Q':= \Quot^{P'_1-(P'_1)^{(2)}}_{\ker\pi_1/\cX_Q/Q}\times_Q\cdots
 \times_Q\Quot^{P'_1-(P'_1)^{(l+1)}}_{\ker\pi_1/\cX_Q/Q}.
\]
Let $(\ker\pi_1)_{Q'}\stackrel{\pi_1^{(i)}}\ra G^{(i)}_1$
($1\leq i\leq l$) be the universal quotient sheaves.
We consider the composite homomorphisms
\begin{align*}
 &\Psi': \Lambda^1_{\cD/S}\otimes(\ker\pi_1)_{Q'}\hookrightarrow
 \Lambda^1_{\cD/S}\otimes(E_1)_{Q'}\stackrel{\Phi_{Q'}}\lra(E_2)_{Q'}
 \stackrel{(\pi_2)_{Q'}}\lra (G_2)_{Q'} \\
 &\psi_i:\ker\pi_1^{(i+1)}\hookrightarrow(\ker\pi_1)_{Q'}
 \stackrel{\pi_1^{(i)}}\lra G^{(i)}_1 \quad (2\leq i\leq l) \\
 &\psi_{l+1}:(\ker\pi_1)_{Q'}\otimes\cO_{\cX}(-\cD) \lra
 (\ker\pi_1)_{Q'} \stackrel{\pi_1^{(l+1)}}\lra G_1^{(l+1)}.
\end{align*}
Let $\tilde{Q}'_{\BP'}$ be the maximal closed subscheme of $Q'$
satisfying $\Psi'_{\tilde{Q}'_{\BP'}}= 0$ and
$(\psi_i)_{\tilde{Q}'_{\BP'}}= 0$ for $2\leq i\leq l+1$.
Since $f_{\BP'}:\tilde{Q}'_{\BP'}\ra T$ is a proper morphism,
\[
 T^s= T\setminus\bigcup_{\BP'\in\cP}f_{\BP'}(\tilde{Q}'_{\BP'})
\]
is an open subscheme which satisfies the desired condition.
\end{proof}

%%%%%%%%%%%%%%%%%%%%%%%% Subsection %%%%%%%%%%%%%%%%%%%%%%%%%%%%%%%
%%%%%%%%%%%%%%%%%%%%%%%%%%%%%%%%%%%%%%%%%%%%%%%%%%%%%%%%%%%%%%%%%%%
\subsection{Construction of the moduli 
space.}\label{moduli-construction}

Now we construct the moduli scheme of
$\overline{\cM_{\cX/S}^{\cD,\balpha',\bbeta,\gamma}}(r,d,\{d_i\})$.
We define a polynomial $P(m)$ in $m$ by
$P(m):= rd_{\cX}m+d+r(1-g)$
where $d_{\cX}= \deg\cO_{\cX_s}(1)$ for $s\in S$
and $g$ is the genus of $\cX_s$.
We take an integer $m_0$ in Proposition \ref{fundamental-lemma}.
By Proposition \ref{bounded}, we may assume, by replacing $m_0$,
that for any $m\geq m_0$,
$h^j(F_i(E_1)(m))= h^j(E_2(m-\gamma))= 0$ for $j>0$, 
$i= 1,\ldots,l+1$
and $E_2(m-\gamma)$, $F_i(E_1)(m)$ ($i= 1,\ldots,l+1$) are generated
by their global sections for any geometric point
$(E_1,E_2,\Phi,F_*(E_1))$ of
$\overline{\cM_{\cX/S}^{\cD,\balpha',\bbeta,\gamma}}(r,d,\{d_i\})$.
Put
$n_1= P(m_0)$ and $n_2= P(m_0-\gamma)$.
Take two free $\cO_S$-modules $V_1,V_2$ such that
$\rank V_1= n_1$, $\rank V_2= n_2$.
Let $Q_1$ be the Quot-scheme
$\Quot_{V_1\otimes\cO_{\cX}(-m_0)/\cX/S}^{P(m)}$
and $V_1\otimes\cO_{\cX_{Q_1}}(-m_0)\to\cE_1$
be the universal quotient sheaf.
Similarly let $Q_2$ be the Quot-scheme
$\Quot_{V_2\otimes\cO_{\cX}(-m_0+\gamma)/\cX/S}^{P(m)}$
and $V_2\otimes\cO_{\cX_{Q_2}}(-m_0+\gamma)\to\cE_2$
be the universal quotient sheaf.
We put $Q_1^{(i)}:= \Quot_{\cE_1/\cX_{Q_1}/Q_1}^{d_i}$.
Let $F_{i+1}(\cE_1)\subset(\cE_1)_{Q_1^{(i)}}$
be the universal subsheaf.
We define $Q$ as the maximal closed subscheme of
$Q_1^{(1)}\times_{Q_1}\cdots\times_{Q_1}Q_1^{(l)}\times Q_2$
such that there are factorizations
\begin{equation}\label{factorization}
 (\cE_1)_Q\otimes\cO_{\cX_Q}(-\cD_Q)\lra F_{i+1}(\cE_1)_Q\hookrightarrow
 F_i(\cE_1)_Q\subset (\cE_1)_Q
\end{equation}
for $i= 1,\ldots,l$, where $F_1(\cE_1)= \cE_1$.
Since $(\cE_2)_Q$ is flat over $Q$, there is a coherent sheaf $\cH$ on 
$Q$
such that there is a functorial isomorphism
\begin{equation}\label{hom-rep}
 \Hom_{\cX_T}(\Lambda^1_{\cD/S}\otimes_{\cO_{\cX}}(\cE_1)_T,
 (\cE_2)_T\otimes\cL)
 \cong \Hom_T(\cH\otimes\cO_T,\cL)
\end{equation}
for any noetherian scheme $T$ over $Q$ and
any quasi-coherent sheaf $\cL$ on $T$.

We denote $\Spec S(\cH)$ by ${\bf V}^*(\cH)$, where $S(\cH)$
is the symmetric algebra of $\cH$ over $\cO_Q$.
Let 
\[
 \tilde{\Phi}:\Lambda^1_{\cD/S}\otimes_{\cO_{\cX}}
 (\cE_1)_{{\bf V}^*(\cH)}\lra(\cE_2)_{{\bf V}^*(\cH)}
\]
be the universal homomorphism.
We define the open subscheme $R^{s}$ of ${\bf V}^*(\cH)$ by
\[
 R^{s}:= \left\{ s\in {\bf V}^*(\cH) \left|
 \begin{array}{l}
  \text{$(V_1)_s\ra H^0((\cE_1)_s(m_0))$,
  $(V_2)_s\ra H^0((\cE_2)_s(m_0-\gamma))$
  are bijective,} \\
  \text{$F_i(\cE_1)_s(m_0),(\cE_2)_s(m_0-\gamma)$
   are generated by their global sections,} \\
  \text{$h^j(F_i(\cE_1)_s(m_0))= h^j((\cE_2)_s(m_0-\gamma))= 0$
  for $j>0$, $1\leq i\leq l+1$} \\
  \text{and $((\cE_1)_s,(\cE_2)_s,\tilde{\Phi}_s,F_*(\cE_1)_s)$
  is stable}
 \end{array}
 \right.\right\}.
\]

For $y\in R^{s}$ and
vector subspaces $V'_1\subset(V_1)_y$, $V'_2\subset(V_2)_y$,
let $E(V'_1,V'_2,y)_1$ be the image of 
$V'_1\otimes\cO_{\cX}(-m_0)\ra(\cE_1)_y$
and $E(V'_1,V'_2,y)_2$ be that of
$V'_1\otimes\Lambda^1_{\cD_y}(-m_0)\oplus 
V'_2\otimes\cO_{\cX}(-m_0+\gamma)
\ra(\cE_2)_y$.
Since the family
\[
 \cF=  \left\{ (E(V'_1,V'_2,y)_1,E(V'_1,V'_2,y)_2)|
 y\in R^{s},V'_1\subset(V_1)_y,V'_2\subset(V_2)_y \right\}
\]
is bounded, there exists an integer $m_1(\geq m_0)$
such that for all $m\geq m_1$,
\begin{gather*}
 V'_1\otimes H^0(\cO_{\cX_y}(m))
 \ra H^0(E(V'_1,V'_2,y)_1(m)), \\
 V'_1\otimes H^0(\cO_{\cX_y}(m_0+m-\gamma)
 \otimes\Lambda^1_{\cD_y}\otimes\cO_{\cX_y}(-m_0))
 \oplus V'_2\otimes H^0(\cO_{\cX_y}(m))
 \ra H^0(E(V'_1,V'_2,y)_2(m-\gamma))
\end{gather*}
are surjective and
$H^i(\cO_{\cX_y}(m_0+m-\gamma)
\otimes\Lambda^1_{\cD_y}\otimes\cO_{\cX_y}(-m_0))= 0$,
$H^i(\cO_{\cX_y}(m))= 0$ for $i>0$
for all members $(E(V'_1,V'_2,y)_1,E(V'_1,V'_2,y)_2)\in\cF$
and the inequality
\begin{align*}
 & (\beta_1\rank E'_1+\beta_2\rank E'_2)d_{\cX}
 \left( \beta_1h^0((\cE_1)_y(m_0))+\beta_2h^0((\cE_2)_y(m_0-\gamma))
 -\sum_{i= 1}^l \beta_1\epsilon_i d_i\right) \\
 & \quad-(\beta_1+\beta_2)rd_{\cX}
 \Big(\beta_1\alpha'_1 h^0(E'_1(m_0))+\beta_2h^0(E'_2(m_0-\gamma))
  + \sum_{i= 1}^l \beta_1\epsilon_i h^0(F_{i+1}(E'_1)(m_0)) \Big) \\
 &>m^{-1}\Big( \beta_1\dim V'_1+\beta_2\dim V'_2
 -\beta_1\chi(E'_1(m_0))-\beta_2\chi(E'_2(m_0-\gamma)) \Big)
 (\beta_1\dim V_1+\beta_2\dim V_2-\sum_{i= 1}^l\beta_1\epsilon_id_i)
\end{align*}
holds for $(0,0)\subsetneq(V'_1,V'_2)\subsetneq((V_1)_y,(V_2)_y)$,
where $E'_1:= E(V'_1,V'_2,y)_1,E'_2:= E(V'_1,V'_2,y)_2$
and $F_{i+1}(E'_1):= E'_1\cap F_{i+1}(\cE_1)_y$ for
$i= 1,\ldots,l$.
{F}rom now on, we fix such a large integer $m_1$.

The composite
\[
 V_1\otimes\Lambda^1_{\cD/S}\otimes
 {\mathcal O}_{\cX}(-m_0)_{R^{s}}
 \ra \Lambda^1_{\cD/S}\otimes(\cE_1)_{R^{s}}
 \stackrel{\tilde{\Phi}}\lra (\cE_2)_{R^{s}}
\]
induces a homomorphism
\[
 V_1\otimes W_1\otimes\cO_{R^{s}}\ra (\pi_{R^{s}})_*
 (\cE_2(m_0+m_1-\gamma)_{R^{s}}),
\]
where
$W_1:= (\pi_S)_*(\cO_{\cX}(m_0+m_1-\gamma)
\otimes\Lambda^1_{\cD/S}\otimes\cO_{\cX}(-m_0))$
and the quotient 
$V_2\otimes\cO_{\cX}(-m_0+\gamma)\ra \cE_2$
induces a homomorphism
\[
V_2\otimes W_2\otimes\cO_{R^{s}}\ra
(\pi_{R^{s}})_*(\cE_2(m_0+m_1-\gamma)_{R^{s}}),
\]
where $W_2:= (\pi_S)_*(\cO_{\cX}(m_1))$.
These homomorphisms induce a quotient bundle
\begin{equation*}
 \left(V_1\otimes W_1\oplus V_2\otimes W_2\right)
 \otimes\cO_{R^{s}}\lra
 (\pi_{R^{s}})_*(\cE_2(m_0+m_1-\gamma)_{R^{s}}).
\end{equation*}
This quotient and the canonical quotient bundles
\begin{align*}
 &V_1\otimes W_2\otimes\cO_{R^s}= 
 V_1\otimes(\pi_S)_*(\cO_{\cX}(m_1))\otimes\cO_{R^{s}}
 \ra (\pi_{R^{s}})_*(\cE_1(m_0+m_1)_{R^{s}}), \\
 &V_1\otimes\cO_{R^{s}}\ra
 (\pi_{R^{s}})_*(\cE_1/F_{i+1}(\cE_1)(m_0)_{R^{s}})
 \quad (i= 1,\ldots,l)
\end{align*}
determine a morphism
\[
 \iota:R^{s}\ra
 \Grass_{r_2}(V_1\otimes W_1 \oplus V_2\otimes W_2)\times
 \Grass_{r_1}(V_1\otimes W_2)\times
 \prod_{i= 1}^l \Grass_{d_i}(V_1),
\]
where
$r_1= h^0(\cE_1(m_0+m_1)_s)$,
$r_2:= h^0(\cE_2(m_0+m_1-\gamma)_s)$
for any point $s\in R^{s}$ and
$\Grass_r(V)$ is the Grassmannian parametrizing
$r$-dimensional quotient vector spaces of $V$.
We can check that $\iota$ is an immersion.

We set $G:= (GL(V_1)\times GL(V_2))/({\bf G}_m\times S)$,
where ${\bf G}_m\times S$ is contained in $GL(V_1)\times GL(V_2)$ as
scalar matrices.
Then $G$ acts canonically on $R^{s}$ and on
$\Grass_{r_2}(V_1\otimes W_1 \oplus V_2\otimes W_2)\times
 \Grass_{r_1}(V_1\otimes W_2)\times
 \prod_{i= 1}^l \Grass_{d_i}(V_1)$.
We can see that $\iota$ is a $G$-equivariant immersion.
There is an $S$-ample line bundle
$\cO_{\Grass_{r_2}(V_1\otimes W_1 \oplus V_2\otimes W_2)}(1)$ on
$\Grass_{r_2}(V_1\otimes W_1 \oplus V_2\otimes W_2)$ induced by
Pl\"{u}cker embedding.
Similarly there are canonical $S$-ample line bundles
$\cO_{\Grass_{r_1}(V_1\otimes W_2)}(1)$, $\cO_{\Grass_{d_i}(V_1)}(1)$,
on $\Grass_{r_1}(V_1\otimes W_2)$,
$\Grass_{d_i}(V_1)$, respectively.
We define positive rational numbers $\nu_1$, $\nu_2$,
$\nu^{(i)}_1$ ($1\leq i\leq l$) by
\begin{align*}
 \nu_1&= \beta_1(\beta_1P(m_0)+\beta_2P(m_0-\gamma)
 - \sum_{i= 1}^l \beta_1\epsilon_i d_i), \\
 \nu_2&= \beta_2(\beta_1P(m_0)+\beta_2P(m_0-\gamma)
 - \sum_{i= 1}^l \beta_1\epsilon_i d_i), \\
 \nu^{(i)}_1&= (\beta_1+\beta_2)\beta_1rd_{\cX}m_1\epsilon_i.
\end{align*}
Let us consider the $\Q$-line bundle
\[
 L:= \iota^*\left(
 \cO_{\Grass_{r_2}(V_1\otimes W_1 \oplus V_2\otimes W_2)}(\nu_1)
 \otimes \cO_{\Grass_{r_1}(V_1\otimes W_2)}(\nu_2)\otimes
 \bigotimes_{i= 1}^l \cO_{\Grass_{d_i}(V_1)}(\nu^{(i)}_1)\right)
\]
on $R^{s}$.
Then for some positive integer $N$,
$L^{\otimes N}$ becomes a $G$-linearized
$S$-ample line bundle on $R^s$.

%%%%%%%%%%%%%%%%%%%%%%%%%%%% GIT-Stability 
%%%%%%%%%%%%%%%%%%%%%%%%%%%%%%%%%%
\begin{Proposition}\label{git-stability}
 All points of $R^s$ are properly stable with respect to the action of 
$G$ and the $G$-linearized $S$-ample line bundle $L^{\otimes N}$.
\end{Proposition}

\begin{proof}
Take any geometric point $x$ of $R^{s}$.
Let $y$ be the induced geometric point of $S$.
We must show that $x$ is a properly stable point
of the fiber $R^s_y$ with respect to the action of
$G_y$ and the polarization $L^{\otimes N}_y$.
So we may assume that $S= \Spec K$ with $K$
an algebraically closed field.
We put
\[
(E_1,E_2,\Phi,F_*(E_1)):= 
((\cE_1)_x,(\cE_2)_x,\tilde{\Phi}_x,F_*(\cE_1)_x).
\]
Let
\[
 \pi_2:V_1\otimes W_1 \oplus V_2\otimes W_2 \ra N_2,
 \quad \pi_1:V_1\otimes W_2 \ra N_1,
 \quad \pi^{(i)}_1:V_1 \ra N^{(i)}_1 \quad (i= 1,\ldots,l)
\]
be the quotient vector spaces corresponding to $\iota(x)$.
We will show that $\iota(x)$ is a properly stable
point with respect to the action of $G$
and the linearization of $L^{\otimes N}$.
Consider the character
\[
 \chi: GL(V_1)\times GL(V_2) \lra \mathbf{G}_m;
 \quad (g_1,g_2)\mapsto \det(g_1)^{\beta_1}\det(g_2)^{\beta_2}.
\]
Then there is an isogeny 
$\ker\chi\longrightarrow G$ and
we may prove the stability with respect to the action of 
$\ker\chi$ instead of $G$.
Take any one parameter subgroup $\lambda$ of 
$\ker\chi$.
For a suitable basis $e^{(1)}_1,\ldots,e^{(1)}_{n_1}$
(resp.\ $e^{(2)}_1,\ldots,e^{(2)}_{n_2}$) of $V_1$ (resp.\ $V_2$),
the action of $\lambda$ on $V_1$ (resp.\ $V_2$) is represented by
\[
e^{(1)}_i \mapsto t^{u^{(1)}_i}e^{(1)}_i \;
(\text{resp.}\: e^{(2)}_i \mapsto t^{u^{(2)}_i}e^{(2)}_i)
\quad (t\in\mathbf{G}_m),
\]
where $u^{(1)}_1\leq \cdots \leq u^{(1)}_{n_1}$
(resp.\ $u^{(2)}_1\leq \cdots \leq u^{(2)}_{n_2}$)
and
$\sum_{i= 1}^{n_1}\beta_1u^{(1)}_i
+\sum_{i= 1}^{n_2}\beta_2u^{(2)}_i= 0$.
Take a basis $f^{(k)}_1,\ldots,f^{(k)}_{b_k}$
of $W_k$ for $k= 1,2$.

We define functions $a_1(p)$ and $a_2(p)$ in 
$p\in\{0,1,\ldots,\beta_1n_1+\beta_2n_2\}$ as follows.
First we put $(a_1(0),a_2(0)):= (0,0)$.
We put
\[
 (a_1(1),a_2(1)):= 
  \begin{cases}
   (1,0) & \text{if $\beta_1u^{(1)}_1\leq\beta_2u^{(2)}_1$} \\
   (0,1) & \text{if $\beta_1u^{(1)}_1>\beta_2u^{(2)}_2$}.
  \end{cases}
\]
Inductively we define
\[
\left\{
\begin{array}{ll}
 (a_1(p+1),a_2(p+1)):= (a_1(p),a_2(p)) & 
  \text{if $p<\beta_1a_1(p)+\beta_2a_2(p)$} \\
 (a_1(p+1),a_2(p+1)):= (a_1(p)+1,a_2(p)) &
 \text{if $p= \beta_1a_1(p)+\beta_2a_2(p)$,
 $\beta_1u^{(1)}_{a_1(p)+1}\leq \beta_2u^{(2)}_{a_2(p)+1}$} \\
 & \text{and $a_1(p)<n_1$} \\
 (a_1(p+1),a_2(p+1)):= (a_1(p),a_2(p)+1)  &
 \text{if $p= \beta_1a_1(p)+\beta_2a_2(p)$,
 $\beta_1u^{(1)}_{a_1(p)+1}> \beta_2u^{(2)}_{a_2(p)+1}$} \\
 & \text{and $a_2(p)<n_2$.} \\
 (a_1(p+1),a_2(p+1)):= (a_1(p)+1,a_2(p)) &
 \text{if $p= \beta_1a_1(p)+\beta_2a_2(p)$ and $a_2(p)= n_2$} \\
 (a_1(p+1),a_2(p+1)):= (a_1(p),a_2(p)+1) &
 \text{if $p= \beta_1a_1(p)+\beta_2a_2(p)$ and $a_1(p)= n_1$}
\end{array}
\right.
\]
Then $a_1(p)$ and $a_2(p)$ are integers with 
$0\leq a_1(p)\leq n_1$, $0\leq a_2(p)\leq n_2$, 
$a_1(p)\leq a_1(p+1)$
and $a_2(p)\leq a_2(p+1)$.
We define $v_1,\ldots,v_{\beta_1n_1+\beta_2n_2}$ and 
$e'_1,\ldots,e'_{\beta_1n_1+\beta_2n_2}$ by
\[
 \left\{
 \begin{array}{ll}
  v_p:= \beta_1u^{(1)}_{a_1(p)}, \;  e'_p:= e^{(1)}_{a_1(p)}  \quad &
  \text{if $a_1(p-1)<a_1(p)$}  \\
  v_p:= \beta_2u^{(2)}_{a_2(p)}, \;  e'_p:= e^{(2)}_{a_2(p)}  \quad &
  \text{if $a_2(p-1)<a_2(p)$} \\
  v_p:= v_{p-1}, \; e'_p:= e'_{p-1} \quad &
  \text{if $a_1(p-1)= a_1(p)$ and $a_2(p-1)= a_2(p)$.}
 \end{array}
 \right.
\]
We put $\delta_p:= (v_{p+1}-v_p)(\beta_1n_1+\beta_2n_2)^{-1}$ 
for $p= 1,\ldots,\beta_1n_1+\beta_2n_2-1$.
Then $\delta_p$ are non-negative rational numbers
and for each $1\leq i\leq n_1$
\[
 \beta_1u^{(1)}_i= \sum_{1\leq p\leq \beta_1n_1+\beta_2n_2-1\atop
a_1(p)< i}
 p\delta_p 
 +\sum_{1\leq p\leq \beta_1n_1+\beta_2n_2-1\atop a_1(p)\geq i}
 (p-\beta_1n_1-\beta_2n_2)\delta_p
\]
and for each $1\leq i\leq n_2$
\[ 
 \beta_2u^{(2)}_i= \sum_{1\leq p\leq \beta_1n_1+\beta_2n_2-1\atop
a_2(p)<i}
 p\delta_p 
 +\sum_{1\leq p\leq \beta_1n_1+\beta_2n_2-1\atop a_2(p)\geq i}
 (p-\beta_1n_1-\beta_2n_2)\delta_p.
\]

For $\mu= 1,\ldots,\beta_1n_1b_1+\beta_2n_2b_2$, we can find
unique integers
$p_0,p_1\in\{0,1,\ldots,\beta_1n_1+\beta_2n_2\}$ such that
\begin{gather*}
 (a_1(p_1),a_2(p_1))= (a_1(p_0+1),a_2(p_0+1))=(a_1(p_0)+1,a_2(p_0)),
 \;{\rm or} \\
 (a_1(p_1),a_2(p_1))= (a_1(p_0+1),a_2(p_0+1)=(a_1(p_0),a_2(p_0)+1)
\end{gather*}
and
\[
 \mu= \begin{cases}
     (a_1(p_0)\beta_1+(p_1-p_0-1))b_1+a_2(p_0)\beta_2b_2+j
      & \begin{array}{l}
         \text{for some $1\leq j\leq b_1$} \\
         \text{if $(a_1(p_1),a_2(p_1))= (a_1(p_0)+1,a_2(p_0))$}
         \end{array} \\
      a_1(p_0)\beta_1b_1+(a_2(p_0)\beta_2+(p_1-p_0-1))b_2+j
      & \begin{array}{l}
         \text{for some $1\leq j\leq b_2$} \\
         \text{if $(a_1(p_1),a_2(p_1))= (a_1(p_0),a_2(p_0)+1)$.}
        \end{array} 
      \end{cases}
\]
For such $\mu$, we put
$s^{(2)}_{\mu}:= v_{p_1}$ and
\[
 h'_{\mu}:= 
 \begin{cases}
  e'_{p_1} \otimes f^{(1)}_j &
  \text{if $(a_1(p_1),a_2(p_1))= (a_1(p_0)+1,a_2(p_0))$} \\
  e'_{p_1} \otimes f^{(2)}_j &
  \text{if $(a_1(p_1),a_2(p_1))= (a_1(p_0),a_2(p_0)+1)$.}
 \end{cases}
\]
Let $U^{(2)}_{\mu}$ be the vector subspace of
$V_1\otimes W_1 \oplus V_2\otimes W_2$ generated by
$h'_1,\ldots,h'_{\mu}$.
We put $U^{(2)}_0=0$.
For $q= 1,\ldots,r_2$, we can find an integer
$\mu^{(2)}_q\in\{1,\ldots,\beta_1n_1b_1+\beta_2n_2b_2\}$ such that
$\dim\pi_2(U^{(2)}_{\mu^{(2)}_q})= q$
and $\dim\pi_2(U^{(2)}_{\mu^{(2)}_q-1})= q-1$.
Then
\begin{align*}
 \sum_{q= 1}^{r_2} s^{(2)}_{\mu^{(2)}_q} &= 
 \sum_{q= 1}^{r_2} s^{(2)}_{\mu^{(2)}_q}
 \left( \dim\pi_2(U^{(2)}_{\mu^{(2)}_q})-
 \dim\pi_2(U^{(2)}_{\mu^{(2)}_q-1}) \right) \\
 &=  \sum_{\mu= 1}^{\beta_1n_1b_1+\beta_2n_2b_2} s^{(2)}_{\mu}
 \left( \dim\pi_2(U^{(2)}_{\mu})-\dim\pi_2(U^{(2)}_{\mu-1}) \right) \\
 &=  r_2 s^{(2)}_{\beta_1n_1b_1+\beta_2n_2b_2} -
 \sum_{\mu= 1}^{\beta_1n_1b_1+\beta_2n_2b_2-1} 
(s^{(2)}_{\mu+1}-s^{(2)}_{\mu})
 \dim\pi_2(U^{(2)}_{\mu}) \\
 &=  r_2 v_{\beta_1n_1+\beta_2n_2} - 
\sum_{p= 1}^{\beta_1n_1+\beta_2n_2-1}
 (v_{p+1}-v_p) \dim \pi_2(U^{(2)}_{\beta_1a_1(p)b_1+\beta_2a_2(p)b_2}) 
\\
 &=  \sum_{p= 1}^{\beta_1n_1+\beta_2n_2-1} \left(
 r_2 p - (\beta_1n_1+\beta_2n_2)
 \dim\pi_2(U^{(2)}_{\beta_1a_1(p)b_1+\beta_2a_2(p)b_2})
 \right) \delta_p.
\end{align*}

For $\mu= (i-1)b_2+j$,
put $h^{(1)}_{\mu}:= e^{(1)}_i\otimes f^{(2)}_j$
for $i= 1,\ldots,n_1$, $j= 1,\ldots,b_2$.
We define integers $s^{(1)}_1,\ldots,s^{(1)}_{b_2n_1}$ by putting
$s^{(1)}_{\mu}:= \beta_1u^{(1)}_i$ for $\mu= (i-1)b_2+j$ with
$1\leq j\leq b_2$.
Let $U^{(1)}_{\mu}$ be the vector subspace of $V_1\otimes W_2$
generated by $h^{(1)}_1,\ldots,h^{(1)}_{\mu}$ for
$\mu= 1,\ldots,b_2n_1$.
We put $U^{(1)}_0=0$.
For $q= 1,\ldots,r_1$, 
let $\mu^{(1)}_q$ be the integer such that
$\dim \pi_1(U^{(1)}_{\mu^{(1)}_q})= q$ and 
$\dim \pi_1(U^{(1)}_{\mu^{(1)}_q-1})= q-1$.
Then
\begin{eqnarray*}
 \sum_{q= 1}^{r_1} s^{(1)}_{\mu^{(1)}_q}
 &= & \sum_{\mu= 1}^{b_2n_1} s^{(1)}_{\mu}
 \left(\dim \pi_1(U^{(1)}_{\mu})-
 \dim \pi_1(U^{(1)}_{\mu-1})\right)  \\
 &= & r_1 s^{(1)}_{b_2n_1} - \sum_{\mu= 1}^{b_2n_1-1}
 (s^{(1)}_{\mu+1}-s^{(1)}_{\mu})\dim\pi_1(U^{(1)}_{\mu}) \\
 &= & r_1 \beta_1 u^{(1)}_{n_1} - \sum_{i= 1}^{n_1-1}
 (u^{(1)}_{i+1}-u^{(1)}_i)\beta_1\dim\pi_1(U^{(1)}_{ib_2}) \\
 &= & r_1\beta_1 u^{(1)}_{n_1}
 +\sum_{a_1(p)<n_1}(v_{p+1}-v_p)\dim\pi_1(U_{a_1(p)b_2}^{(1)}) \\
 &= & r_1 \left(
 \sum_{1\leq p\leq \beta_1n_1+\beta_2n_2-1\atop a_1(p)< n_1}p\delta_p 
 +\sum_{1\leq p\leq \beta_1n_1+\beta_2n_2-1\atop a_1(p)\geq n_1}
 (p-\beta_1n_1-\beta_2n_2)\delta_p  \right)  \\
 &&\quad - \sum_{1\leq p\leq \beta_1n_1+\beta_2n_2-1\atop a_1(p)<n_1}
 (\beta_1n_1+\beta_2n_2)\delta_p\dim\pi_1(U^{(1)}_{a_1(p)b_2}) \\
 &= & \sum_{p= 1}^{\beta_1n_1+\beta_2n_2-1}\left(
 r_1p-(\beta_1n_1+\beta_2n_2)\dim\pi_1(U^{(1)}_{a_1(p)b_2})
 \right)\delta_p
\end{eqnarray*}

Let $V^{(1)}_p$ be the vector subspace of $V_1$
generated by $e^{(1)}_1,\ldots,e^{(1)}_p$.
We put $V^{(1)}_0=0$.
For $i= 1,\ldots,l$ and for $q= 1,\ldots,d_i$,
let $\mu^{i}_q$ be the integer such that
$\dim \pi_1^{(i)}(V^{(1)}_{\mu^{i}_q})= q$ and
$\dim \pi_1^{(i)}(V^{(1)}_{\mu^{i}_q-1})= q-1$.
Then
\begin{align*}
 \sum_{q= 1}^{d_i} \beta_1 u^{(1)}_{\mu^{i}_q}
 &=  \sum_{q= 1}^{d_i} \beta_1 u^{(1)}_{\mu^{i}_q}
 \left( \dim\pi_1^{(i)}(V^{(1)}_{\mu^{i}_q})
 -\dim \pi_1^{(i)}(V^{(1)}_{\mu^{i}_q-1}) \right) \\
 &=  \sum_{p= 1}^{n_1} \beta_1 u^{(1)}_p
 \left( \dim\pi^{(i)}_1(V^{(1)}_p)
 -\dim\pi^{(i)}_1(V^{(1)}_{p-1}) \right) \\
 &=  d_i \beta_1 u^{(1)}_{n_1}
 - \sum_{p= 1}^{n_1-1} \beta_1 (u^{(1)}_{p+1}-u^{(1)}_p)
 \dim\pi_1^{(i)}(V^{(1)}_p) \\
 &= d_i\beta_1u_{n_1}^{(1)}
 -\sum_{a_1(p)<n_1}(v_{p+1}-v_p)\dim\pi^{(i)}_1(V_{a_1(p)}^{(1)}) \\
 &=  d_i \left( \sum_{1\leq p\leq \beta_1n_1+\beta_2n_2-1\atop a_1(p)< 
n_1}
 p\delta_p +
 \sum_{1\leq p\leq \beta_1n_1+\beta_2n_2-1\atop a_1(p)\geq n_1}
 (p-\beta_1n_1-\beta_2n_2)\delta_p \right) \\
 & \quad -\sum_{1\leq p\leq \beta_1n_1+\beta_2n_2-1\atop a_1(p)<n_1}
 (\beta_1n_1+\beta_2n_2)\delta_p\dim\pi_1^{(i)}(V^{(1)}_{a_1(p)}) \\
 &=  \sum_{p= 1}^{\beta_1n_1+\beta_2n_2-1} \left( d_i p -
 (\beta_1n_1+\beta_2n_2)\dim\pi^{(i)}_1(V^{(1)}_{a_1(p)}) \right) 
\delta_p.
\end{align*}
Thus we have
\begin{align*}
 \mu^{L^{\otimes N}}(x,\lambda) &= 
 -\left( \sum_{k= 1}^2\nu_k\sum_{q= 1}^{r_k}s^{(k)}_{\mu^{(k)}_q}+
\sum_{i= 1}^l\nu_1^{(i)}\sum_{q= 1}^{d_i}\beta_1u^{(1)}_{\mu^{i}_q}\right)N \\
 &= \sum_{p= 1}^{\beta_1n_1+\beta_2n_2-1} N\delta_p \Bigg\{\left(
 -p\sum_{i= 1}^l\nu_1^{(i)} d_i +
 (\beta_1n_1+\beta_2n_2)\sum_{i= 1}^l \nu_1^{(i)}
 \dim\pi^{(i)}_1(V^{(1)}_{a_1(p)})\right) \\
 & -(\nu_1r_1+\nu_2r_2)p+
 (\beta_1n_1+\beta_2n_2)\left(\nu_1\dim\pi_1(U^{(1)}_{a_1(p)b_2})+
 \nu_2\dim\pi_2(U^{(2)}_{\beta_1a_1(p)b_1+\beta_2a_2(p)b_2})\right)
 \Bigg\}.
\end{align*}
See \cite{Mum:GIT}, Definition 2.2 for the definition of
$\mu^{L^{\otimes N}}(x,\lambda)$.
By \cite{Mum:GIT}, Theorem 2.1, $x$ is a properly stable point if
\begin{align*}
 & -p(\nu_1r_1+\nu_2r_2)+
 (\beta_1n_1+\beta_2n_2)(\nu_1\dim\pi_1(U^{(1)}_{a_1(p)b_2})+
 \nu_2\dim\pi_2(U^{(2)}_{\beta_1a_1(p)b_1+\beta_2a_2(p)b_2})) \\
 &  \quad  -p\sum_{i= 1}^l \nu_1^{(i)} d_i +
 (\beta_1n_1+\beta_2n_2)\sum_{i= 1}^l \nu_1^{(i)}
 \dim\pi^{(i)}_1(V^{(1)}_{a_1(p)}) > 0
\end{align*}
for all $p= 1,\ldots,\beta_1n_1+\beta_2n_2-1$.

For each $p$ ($1\leq p\leq \beta_1n_1+\beta_2n_2-1$), let $V'_k$ be the 
vector subspace
of $V_k$ generated by $e^{(k)}_1,\ldots,e^{(k)}_{a_k(p)}$ for $k= 1,2$.
Then $U^{(1)}_{a_1(p)b_2}= V'_1\otimes W_2$ and
$U^{(2)}_{\beta_1a_1(p)b_1+\beta_2a_2(p)b_2}= 
V'_1\otimes W_1 \oplus V'_2\otimes W_2$.
Put
\begin{gather*}
 E'_1:= \im (V'_1\otimes\cO_{\cX}(-m_0)\ra E_1), \quad
 F_{i+1}(E'_1):= F_{i+1}(E_1)\cap E'_1, \quad (i= 1,\ldots,l), \\
 E'_2:= \im (V'_1\otimes\Lambda^1_{\cD/S}(-m_0)
 \oplus V'_2\otimes\cO_{\cX}(-m_0+\gamma)\ra E_2),
 \quad \Phi':= \Phi|_{\Lambda^1_{\cD/S}\otimes E'_1}.
\end{gather*}
Then $(E'_1,E'_2,\Phi',F_*(E'_1))$ is a parabolic
$\Lambda^1_{\cD_K}$-subtriple of $(E_1,E_2,\Phi,F_*(E_1))$.
By the choice of $m_1$, we have
$\pi_2(U^{(2)}_{\beta_1a_1(p)b_1+\beta_2a_2(p)b_2})= H^0(E'_2(m_0+m_1-\gamma))$
and $\pi_1(U^{(1)}_{a_1(p)b_2})= H^0(E'_1(m_0+m_1))$.
Put $r'_1:= \rank E'_1$, $r'_2:= \rank E'_2$.
Let ${V'}_1^{(i)}$ be the kernel of the composite
$V'_1\hookrightarrow V_1 \stackrel{\pi_1^{(i)}}\lra N_1^{(i)}$.
Then we have
\begin{align*} 
 & -p(\nu_1r_1+\nu_2r_2)+
 (\beta_1n_1+\beta_2n_2)(\nu_1\dim\pi_1(U^{(1)}_{a_1(p)b_2})+
 \nu_2\dim\pi_2(U^{(2)}_{\beta_1a_1(p)b_1+\beta_2a_2(p)b_2})) \\
 &  \quad  -p\sum_{i= 1}^l \nu_1^{(i)} d_i +
 (\beta_1n_1+\beta_2n_2)\sum_{i= 1}^l \nu_1^{(i)}
 \dim\pi^{(i)}_1(V^{(1)}_{a_1(p)}) \notag \\
 &\geq (\beta_1\dim V_1+\beta_2\dim V_2
 -\sum_{i= 1}^l\beta_1\epsilon_id_i)\times \\
 & \quad \Big\{-(\beta_1\dim V'_1+\beta_2\dim V'_2)
 (\beta_1h^0(E_1(m_0+m_1))+\beta_2h^0(E_2(m_0+m_1-\gamma)))\\
 & \hspace{50pt}+(\beta_1\dim V_1+\beta_2\dim V_2)
 (\beta_1h^0(E'_1(m_0+m_1))+\beta_2h^0(E'_2(m_0+m_1-\gamma)))\Big\} \\
 & \quad -(\beta_1\dim V'_1+\beta_2\dim V'_2) 
 \sum_{i= 1}^l \nu_1^{(i)} d_i + (\beta_1\dim V_1+\beta_2\dim V_2)
 \sum_{i= 1}^l\nu_1^{(i)}(\dim V'_1-\dim {V'}_1^{(i)}) \notag\\
 &= (\beta_1\dim V_1+\beta_2\dim V_2
 -\sum_{i= 1}^l\beta_1\epsilon_id_i)\times \\
 &\quad \Big\{-(\beta_1\dim V'_1+\beta_2\dim V'_2)
 \left(rd_{\cX}(\beta_1+\beta_2)m_1+\beta_1\dim V_1+\beta_2\dim 
V_2\right) \\
 &\quad +(\beta_1\dim V_1+\beta_2\dim V_2)
 \left((\beta_1r'_1+\beta_2r'_2)d_{\cX}m_1+
 \beta_1\chi(E'_1(m_0))+\beta_2\chi(E'_2(m_0-\gamma))\right) \Big\} 
\notag\\
 & \quad -(\beta_1\dim V'_1+\beta_2\dim V'_2) 
 \sum_{i= 1}^l \nu_1^{(i)} d_i + (\beta_1\dim V_1+\beta_2\dim V_2)
 \sum_{i= 1}^l\nu_1^{(i)}(\dim V'_1-\dim {V'}_1^{(i)}) \\
 &=  (\beta_1\dim V_1+\beta_2\dim V_2
  -\sum_{i= 1}^l \beta_1\epsilon_i d_i )\times \\
 &\quad \Big\{-rd_{\cX}(\beta_1+\beta_2)m_1(\beta_1\dim V'_1+\beta_2\dim 
V'_2)
 +(\beta_1r'_1+\beta_2r'_2)d_{\cX}m_1(\beta_1\dim V_1+\beta_2\dim 
V_2)\Big\} \\
 &\quad +(\beta_1\dim V_1+\beta_2\dim V_2
  -\sum_{i= 1}^l \beta_1\epsilon_i d_i)
 (\beta_1\dim V_1+\beta_2\dim V_2)\times \\
 &\quad \Big\{-(\beta_1\dim V'_1+\beta_2\dim V'_2)
 +(\beta_1\chi(E'_1(m_0))+\beta_2\chi(E'_2(m_0-\gamma)))\Big\} \\
 & \quad -(\beta_1\dim V'_1+\beta_2\dim V'_2) 
 \sum_{i= 1}^l (\beta_1+\beta_2)\beta_1rd_{\cX}m_1\epsilon_i d_i \\
 & \quad + (\beta_1\dim V_1+\beta_2\dim V_2)
 \sum_{i= 1}^l (\beta_1+\beta_2)\beta_1rd_{\cX}m_1\epsilon_i
 (\dim V'_1-\dim {V'}_1^{(i)}) \\
&=  -(\beta_1\dim V_1+\beta_2\dim V_2)(\beta_1+\beta_2)
  rd_{\cX}m_1 \Big( \beta_1\dim V'_1+\beta_2\dim V'_2
 -\sum_{i= 1}^l\beta_1\epsilon_i(\dim V'_1-\dim {V'}_1^{(i)})\Big) \\
 & \quad +(\beta_1\dim V_1+\beta_2\dim 
V_2)(\beta_1r'_1+\beta_2r'_2)d_{\cX}m_1
 \Big( \beta_1\dim V_1+\beta_2\dim V_2
 -\sum_{i= 1}^l \beta_1\epsilon_i d_i \Big) \\
 &\quad +\Big(\beta_1\dim V_1+\beta_2\dim V_2
  -\sum_{i= 1}^l \beta_1\epsilon_i d_i \Big)(\beta_1\dim 
V_1+\beta_2\dim V_2)
  \times \\
 & \quad \Big(-(\beta_1\dim V'_1+\beta_2\dim V'_2)
 +(\beta_1\chi(E'_1(m_0))+\beta_2\chi(E'_2(m_0-\gamma)))\Big) \\
 \end{align*}
 \begin{align*}
  &\geq (\beta_1\dim V_1+\beta_2\dim V_2) 
 \Big\{(\beta_1r'_1+\beta_2r'_2)d_{\cX}m_1
 \Big( \beta_1h^0(E_1(m_0))+\beta_2h^0(E_2(m_0-\gamma))
 -\sum_{i= 1}^l \beta_1\epsilon_i d_i \Big) \\
 & \quad -(\beta_1+\beta_2)rd_{\cX}m_1
 \Big( \beta_1h^0(E'_1(m_0))+\beta_2h^0(E'_2(m_0-\gamma))
 -\sum_{i= 1}^l\beta_1\epsilon_i
 (h^0(E'_1(m_0))-h^0(F_{i+1}(E'_1)(m_0))
 \Big) \Big\} \notag \\
 &\quad +\Big(\beta_1\dim V_1+\beta_2\dim V_2
  -\sum_{i= 1}^l \beta_1\epsilon_i d_i \Big)
 (\beta_1\dim V_1+\beta_2\dim V_2)\times \\
 &\quad \Big(-(\beta_1\dim V'_1+\beta_2\dim V'_2)
 +(\beta_1\chi(E'_1(m_0))+\beta_2\chi(E'_2(m_0-\gamma)))\Big) \\
 &>0.
\end{align*}

Note that the last inequality holds by the choice of $m_1$.
Hence $x$ is a properly stable point.
\end{proof}

By Proposition \ref{git-stability}, there exists a geometric quotient
$R^s/G$.
The following proposition follows from a standard argument.

\begin{Theorem}\label{thm-construction}
$\overline{M_{\cX/S}^{\cD,\balpha',\bbeta,\gamma}}(r,d,\{d_i\}):= R^s/G$
 is a coarse moduli scheme of
 $\overline{\cM_{\cX/S}^{\cD,\balpha',\bbeta,\gamma}}(r,d,\{d_i\})$.
\end{Theorem}

\begin{Remark}\rm
 The quotient map
 $R^s\ra \overline{M_{\cX/S}^{\cD,\balpha',\bbeta,\gamma}}(r,d,\{d_i\})$
 is a principal $G$-bundle, which we can see by the following lemma
 and the same argument as \cite{M}, Proposition 6.4.
\end{Remark}

\begin{Lemma}
 Take any geometric point
 $(E_1,E_2,\Phi,F_*(E_1))\in
 \overline{M_{\cX/S}^{\cD,\balpha',\bbeta,\gamma}(r,d,\{d_i\})}(K)$.
 Then for any endomorphisms
 $f_1:E_1\ra E_1$, $f_2:E_2\ra E_2$ satisfying
 $\Phi\circ(1\otimes f_1)= f_2\circ\Phi$ and
 $f_1(F_{i+1}(E_1))\subset F_{i+1}(E_1)$ for $1\leq i\leq l$,
 there exists $c\in K$ such that
 $(f_1,f_2)= (c\cdot\mathrm{id}_{E_1},c\cdot\mathrm{id}_{E_2})$.
\end{Lemma}

\begin{proof}
Take such $(f_1,f_2)$.
Let $c\in K$ be an eigenvalue of $f_1\otimes k(x)$ for some
$x\in\cX_K(K)$.
Then $f_1-c\cdot\mathrm{id}_{E_1}$ becomes an endomorphism
of $E_1$ which is not an isomorphism.
Put $E'_1:= \im(f_1-c\cdot\mathrm{id}_{E_1})$,
$E'_2:= \im(f_2-c\cdot\mathrm{id}_{E_2})$,
$\Phi':= \Phi|_{\Lambda^1_{\cD_K}\otimes E'_1}$ and
$F_{i+1}(E'_1):= (f_1-c\cdot\mathrm{id}_{E_1})(F_{i+1}(E_1))$
for $i= 1,\ldots,l$.
Then $(E'_1,E'_2,\Phi',F_*(E'_1))$ becomes a parabolic
$\Lambda^1_{\cD_K}$-subtriple of $(E_1,E_2,\Phi,F_*(E_1))$.
If we put $G_1:= \ker(E_1\ra E'_1)$, $G_2:= \ker(E_2\ra E'_2)$,
$\Phi_G:= \Phi|_{\Lambda^1_{\cD_K}\otimes G_1}$ and
$F_{i+1}(G_1):= F_{i+1}(E_1)\cap G_1$ for $i= 1,\ldots,l$,
then $(G_1,G_2,\Phi_G,F_*(G_1))$ becomes a parabolic
$\Lambda^1_{\cD_K}$-subtriple of $(E_1,E_2,\Phi,F_*(E_1))$.
If $(E'_1,E'_2)\neq(0,0)$, then, by the stability of
$(E_1,E_2,\Phi,F_*(E_1))$, we must have the inequalities
\begin{align*}
 & \frac{\beta_1\alpha'_1\chi(E_1(m))+\beta_2\chi(E_2(m-\gamma))
 +\sum_{i= 1}^l\beta_1\epsilon_i\chi(F_{i+1}(E_1)(m))}
 {\beta_1\rank(E_1)+\beta_2\rank(E_2)} \\
 &> \frac{\beta_1\alpha'_1\chi(E'_1(m))+\beta_2\chi(E'_2(m-\gamma))
 +\sum_{i= 1}^l\beta_1\epsilon_i\chi(F_{i+1}(E'_1)(m))}
 {\beta_1\rank(E'_1)+\beta_2\rank(E'_2)} \\
 &>\frac{\beta_1\alpha'_1\chi(E_1(m))+\beta_2\chi(E_2(m-\gamma))
 +\sum_{i= 1}^l\beta_1\epsilon_i\chi(F_{i+1}(E_1)(m))}
 {\beta_1\rank(E_1)+\beta_2\rank(E_2)} 
\end{align*}
for $m\gg 0$, which is a contradiction.
Therefore we have $(E'_1,E'_2)= (0,0)$,
which means that
$(f_1,f_2)= (c\cdot\mathrm{id}_{E_1},c\cdot\mathrm{id}_{E_2})$.
\end{proof}

%%%%%%%%%%%%%%%%%%%%%%%% Subsection %%%%%%%%%%%%%%%%%%%%%%%%%%%%
%%%%%%%%%%%%%%%%%%%%%%%%%%%%%%%%%%%%%%%%%%%%%%%%%%%%%%%%%%%%%%%%%
\subsection{Projectivity of the moduli space}
%%%%%%%%%%%%%%%%%%%%%%%%%%%%%%%%%%%%%%%%%%%%%%%%%%%%%%%%%%%%%%%%%%

\begin{Proposition}\label{valuative-criterion}
 Let $R$ be a discrete valuation ring over $S$
with residue field $k= R/m$ and quotient field $K$.
 Let $(E_1,E_2,\Phi,F_*(E_1))$ be a semistable
 parabolic $\Lambda^1_{\cD_K}$-triple on $\cX_K$.
 Then there exists a flat family
 $(\tilde{E}_1,\tilde{E}_2,\tilde{\Phi},F_*(\tilde{E}_1))$
 of parabolic $\Lambda^1_{\cD_R}$-triples on $\cX_R$ over $R$ such that
 $(E_1,E_2,\Phi,F_*(E_1))\cong
 (\tilde{E}_1,\tilde{E}_2,\tilde{\Phi},F_*(\tilde{E}_1))\otimes_RK$
 and that
 $(\tilde{E}_1,\tilde{E}_2,\tilde{\Phi},F_*(\tilde{E}_1))\otimes_R k$
 is semistable.
\end{Proposition}

\begin{proof}
Two surjections
\begin{align*}
 V_1\otimes\cO_{\cX_K}(-m_0)\cong
 H^0(E_1(m_0))\otimes\cO_{\cX_K}(-m_0)\ra E_1, \\
 V_2\otimes\cO_{\cX_K}(-m_0+\gamma)\cong
 H^0(E_2(m_0-\gamma))\otimes\cO_{\cX_K}(-m_0+\gamma)\ra E_2
\end{align*}
and the quotients
$E_1\ra E_1/F_{i+1}(E_1)$ ($i= 1,\ldots,l$)
give a morphism
$f:\Spec K \ra Q$,
where $Q$ is defined by the property (\ref{factorization})
in subsection \ref{moduli-construction}.
Since $Q$ is proper over $S$,
$f$ extends to a morphism
$\tilde{f}:\Spec R \ra Q$.
Thus there are coherent sheaves
$E^{(0)}_1$, $E^{(0)}_2$ on ${\cX}_R$
flat over $R$ and a flat family of filtrations
$F_*(E^{(0)}_1)$ of $E^{(0)}_1$ such that
$E^{(0)}_1\otimes K\cong E_1$,
$E^{(0)}_2\otimes K\cong E_2$
and $F_*(E^{(0)}_1)\otimes_RK= F_*(E_1)$.
The pullback of $\cH$ by 
the morphism $\tilde{f}:\Spec R \ra Q$
is denoted by $\cH_R$.
Recall that $\cH$ is defined by (\ref{hom-rep})
in subsection \ref{moduli-construction}.
The homomorphism $\Phi:\Lambda^1_{\cD/S}\otimes E_1\ra E_2$
corresponds to a homomorphism
$\psi:\cH_R\otimes_RK \ra K$.
There is a non-zero element $t\in K\setminus\{0\}$ and
a homomorphism $\tilde{\psi}:\cH_R \ra R$ such that
$t\psi= \tilde{\psi}\otimes_RK$.
Let $\Phi^{(0)}:\Lambda^1_{\cD/S}\otimes E^{(0)}_1\ra E^{(0)}_2$
be the homomorphism corresponding to $\tilde{\psi}$.
Then we have
$(E_1,E_2,\Phi,F_*(E_1))\cong
(E^{(0)}_1,E^{(0)}_2,\Phi^{(0)},F_*(E^{(0)}_1))\otimes_R K$,
since $(E_1,E_2,\Phi,F_*(E_1))\cong(E_1,E_2,t\Phi,F_*(E_1))$.
Our proposition follows from the following claim:

{\bf Claim} There is a flat family
$(\tilde{E}_1,\tilde{E}_2,\tilde{\Phi},F_*(\tilde{E}_1))$
of parabolic $\Lambda^1_{\cD_R}$-triples on $\cX_R$ over $R$ such that
$\tilde{E}_j\subset E^{(0)}_j$ for $j= 1,2$,
$F_{i+1}(\tilde{E}_1)\subset F_{i+1}(E^{(0)}_1)$
for $i= 1,\ldots,l$,
$\tilde{\Phi}= \Phi^{(0)}|_{\tilde{E}_1\otimes\Lambda^1_{\cD/S}}$,
$(\tilde{E}_1,\tilde{E}_2,\tilde{\Phi},F_*(\tilde{E}_1))\otimes_RK
\cong(E_1,E_2,\Phi,F_*(E_1))$ and
$(\tilde{E}_1,\tilde{E}_2,\tilde{\Phi},F_*(\tilde{E}_1))\otimes_Rk$
is semistable.

Assume that $E^{(0)}_1\otimes k$ or $E^{(0)}_2\otimes k$ have torsions.
In this case let $B^{(0)}_1$ and $B^{(0)}_2$ be the torsion parts of
$E^{(0)}_1\otimes k$ and $E^{(0)}_2\otimes k$, respectively.
Then there are exact sequences
\begin{align*}
 &0 \ra B^{(0)}_1 \lra E^{(0)}_1\otimes k \lra G^{(0)}_1 \ra 0 \\
 &0 \ra B^{(0)}_2 \lra E^{(0)}_2\otimes k \lra G^{(0)}_2 \ra 0,
\end{align*}
where $G^{(0)}_1$ and $G^{(0)}_2$ are vector bundles on $\cX_k$.
Put $E^{(1)}_1:= \ker(E^{(0)}_1\ra ((E^{(0)}_1\otimes k)/B^{(0)}_1))$,
$E^{(1)}_2:= \ker(E^{(0)}_2\ra ((E^{(0)}_2\otimes k)/B^{(0)}_2))$,
$\Phi^{(1)}:= \Phi^{(0)}|_{\Lambda^1_{\cD_R}\otimes E^{(1)}_1}$ and
$F_{i+1}(E^{(1)}_1):= F_{i+1}(E^{(0)}_1)\cap E^{(1)}_1$
for $i= 1,\ldots,l$.
Then there are exact sequences
\begin{align*}
 &0 \ra G^{(0)}_1 \lra E^{(1)}_1\otimes k \lra B^{(0)}_1 \ra 0 \\
 &0 \ra G^{(0)}_2 \lra E^{(1)}_2\otimes k \lra B^{(0)}_2 \ra 0.
\end{align*}
Again let $B^{(1)}_1$ and $B^{(1)}_2$ be the torsion parts of
$E^{(1)}_1\otimes k$ and $E^{(1)}_2\otimes k$, respectively.
Repeating these operations, we obtain sequences
$(E^{(n)}_1,E^{(n)}_2,\Phi^{(n)},F_*(E^{(n)}_1))_{n\geq 0}$,
$(B^{(n)}_1,B^{(n)}_2)_{n\geq 0}$ and
$(G^{(n)}_1,G^{(n)}_2)_{n\geq 0}$.
Then the injections
$B^{(n+1)}_1\hookrightarrow B^{(n)}_1$,
$B^{(n+1)}_2\hookrightarrow B^{(n)}_2$
are induced by the homomorphisms
$E^{(n+1)}_1\otimes k \ra E^{(n)}_1\otimes k$,
$E^{(n+1)}_2\otimes k \ra E^{(n)}_2\otimes k$.
Since $(\length B^{(n)}_1,\length B^{(n)}_2)_{n\geq 0}$ is stationary,
we may assume that it is constant.
Then we have isomorphisms
$B^{(n+1)}_1\stackrel{\sim}\ra B^{(n)}_1$,
$B^{(n+1)}_2\stackrel{\sim}\ra B^{(n)}_2$,
$G^{(n)}_1\stackrel{\sim}\ra G^{(n+1)}_1$,
$G^{(n)}_2\stackrel{\sim}\ra G^{(n+1)}_2$
for all $n$.
Assume that $(B^{(n)}_1,B^{(n)}_2)\neq(0,0)$.
There is an exact sequence
\[
 E^{(n)}_j/m^nE^{(0)}_j \stackrel{u}\lra
 E^{(0)}_j/m^nE^{(0)}_j \lra E^{(0)}_j/E^{(n)}_j \ra 0
\]
for $n\geq 1$ and $j= 1,2$.
We can see that $(E^{(n)}_j/m^nE^{(0)}_j)\otimes k\cong B_j^{(n-1)}$
and that
\[
 u\otimes k:(E^{(n)}_j/m^nE^{(0)}_j)\otimes k\cong B_j^{(n-1)}
 \ra E^{(0)}_j\otimes k
\]
is injective.
Thus $E^{(0)}_j/E^{(n)}_j$ is flat over $R/m^n$ and the
quotient $E^{(0)}_j/m^nE^{(0)}_j \ra E^{(0)}_j/E^{(n)}_j$
determines a morphism
$f_n:\Spec R/m^n \ra \Quot_{E^{(0)}_j/\cX_R/R}$ for $n\geq 1$.
So we obtain a morphism
$f:\Spec\hat{R} \ra \Quot_{E^{(0)}_j/\cX_R/R}$,
where $\hat{R}$ is the completion of $R$.
$f$ corresponds to a quotient sheaf
$E^{(0)}_j\otimes\hat{R}\stackrel{\pi}\ra G$.
Since $(\ker\pi)\otimes R/m\cong B^{(0)}_j$,
$\ker\pi\otimes \hat{K}$ is a torsion submodule of $E^{(0)}_j$,
which is nonzero either for $j= 1$ or $j= 2$,
where $\hat{K}$ is the quotient field of $\hat{R}$.
However, it is a contradiction, because
$E^{(0)}_1\otimes\hat{K}$, $E^{(0)}_2\otimes\hat{K}$ are vector bundles.
Hence we must have $(B^{(n)}_1,B^{(n)}_2)= (0,0)$ for some $n$.
So we may assume without loss of generality that
$E^{(0)}_1\otimes k$ and $E^{(0)}_2\otimes k$ are locally free.

Now assume that the claim does not hold.
Then we can define a descending sequence
of flat families of parabolic $\Lambda^1_{\cD_R}$-triples
\[
 (E^{(0)}_1,E^{(0)}_2,\Phi^{(0)},F_*(E^{(0)}_1))\supset
 (E^{(1)}_1,E^{(1)}_2,\Phi^{(1)},F_*(E^{(1)}_1))\supset
 (E^{(2)}_1,E^{(2)}_2,\Phi^{(2)},F_*(E^{(2)}_1))\supset\cdots
\]
as follows:
Suppose $(E^{(n)}_1,E^{(n)}_2,\Phi^{(n)},F_*(E^{(n)}_1))$
has already been defined.
There exists a maximal destabilizer
$(B^{(n)}_1,B^{(n)}_2,\Phi_{B^{(n)}},F_*(B^{(n)}_1))$ of
$(E^{(n)}_1,E^{(n)}_2,\Phi^{(n)},F_*(E^{(n)}_1))\otimes k$
as in the usual case of semistability of coherent sheaves.
We can see that $B^{(n)}_j$ is a subbundle of $E^{(n)}_j\otimes k$
for $j= 1,2$ and 
$F_{i+1}(B^{(n)}_1)= B^{(n)}_1\cap(F_{i+1}(E^{(n)}_1)\otimes k)$
for $i= 1,\ldots,l$.
We put $G^{(n)}_j:= (E^{(n)}_j\otimes k)/B^{(n)}_j$ for $j= 1,2$.
Then $G^{(n)}_1$ has an induced quotient parabolic structure
$F_*(G^{(n)}_1)$.
A homomorphism
$\Phi_{G^{(n)}}:\Lambda^1_{\cD/S}\otimes G^{(n)}_1\ra G^{(n)}_2$
is induced by $\Phi^{(n)}$ and
$(G^{(n)}_1,G^{(n)}_2,\Phi_{G^{(n)}},F_*(G^{(n)}_1))$
becomes a parabolic $\Lambda^1_{\cD_k}$-triple.
Put
\begin{gather*}
 E^{(n+1)}_j= \ker(E^{(n)}_j\ra G^{(n)}_j), \quad
 \Phi^{(n+1)}:= \Phi^{(n)}|_{\Lambda^1_{\cD/S}\otimes E^{(n+1)}_1}, \\
 F_{i+1}(E^{(n+1)}_1)= \ker(F_{i+1}(E^{(n)}_1)\ra F_{i+1}(G^{(n)}_1))
\end{gather*}
Then $(E^{(n+1)}_1,E^{(n+1)}_2,\Phi^{(n+1)},F_*(E^{(n+1)}_1))$
becomes a flat family of parabolic $\Lambda^1_{\cD_R}$-triples
on $\cX_R$ over $R$.
There are exact sequences
\begin{equation}\label{exact-sequence}
 0\ra B^{(n)}_j \ra E^{(n)}_j\otimes k \ra
 G^{(n)}_j \ra 0 \quad \text{and}\quad
 0\ra G^{(n)}_j \ra E^{(n+1)}_j\otimes k
 \ra B^{(n)}_j \ra 0
\end{equation}
for $j= 1,2$.
Then we can see that
$(G^{(n)}_1,G^{(n)}_2,\Phi_{G^{(n)}},F_*(G^{(n)}_1))$
becomes a parabolic $\Lambda^1_{\cD_k}$-subtriple of
$(E^{(n+1)}_1,E^{(n+1)}_2,\Phi^{(n+1)},F_*(E^{(n+1)}_1))\otimes k$.
We can check that
$F_{i+1}(G^{(n)}_1)= G^{(n)}_1\cap(F_{i+1}(E^{(n+1)}_1)\otimes k)$
for $i= 1,\ldots,l$.
Put 
\begin{gather*}
 C^{(n)}_j:= G^{(n)}_j \cap B^{(n+1)}_j, \quad
 \Phi_{C^{(n)}}:= (\Phi^{(n+1)}\otimes k)|
 _{\Lambda^1_{\cD_k}\otimes C^{(n)}_1}, \\
 F_{i+1}(C^{(n)}_1):= F_{i+1}(G^{(n)}_1)\cap F_{i+1}(B^{(n+1)}_1)
 \quad (i= 1,\ldots,l).
\end{gather*}
Then $(C^{(n)}_1,C^{(n)}_2,\Phi_{C^{(n)}},F_*(C^{(n)}_1))$
becomes a parabolic $\Lambda^1_{\cD_k}$-triple and 
$F_{i+1}(C^{(n)}_1)= 
C^{(n)}_1\cap(F_{i+1}(E^{(n+1)}_1)\otimes k)$
for $i= 1,\ldots,l$.
A quotient parabolic structure $F_*(B^{(n+1)}_1/C^{(n)}_1)$
is induced on $B^{(n+1)}_1/C^{(n)}_1$ and a homomorphism
$\Phi_{B^{(n+1)}/C^{(n)}}:\Lambda^1_{\cD_k}\otimes B^{(n+1)}_1/C^{(n)}_1
\ra B^{(n+1)}_2/C^{(n)}_2$
is induced by $\Phi^{(n+1)}$.
Then
\[
 (B^{(n+1)}_1/C^{(n)}_1,B^{(n+1)}_2/C^{(n)}_2,
 \Phi_{B^{(n+1)}/C^{(n)}},F_*(B^{(n+1)}_1/C^{(n)}_1))
\]
becomes a parabolic $\Lambda^1_{\cD_k}$-triple.
If $(C^{(n)}_1,C^{(n)}_2)\neq(0,0)$, then we have
\begin{gather*}
 \mu((C^{(n)}_1,C^{(n)}_2,\Phi_{C^{(n)}},F_*(C^{(n)}_1)))\leq
 \mu_{\max}((G^{(n)}_1,G^{(n)}_2,\Phi_{G^{(n)}},F_*(G^{(n)}_*))) \\
 <\mu_{\max}((E^{(n)}_1,E^{(n)}_2,\Phi^{(n)},F_*(E^{(n)}_1))\otimes k)
 = \mu((B^{(n)}_1,B^{(n)}_2,\Phi_{B^{(n)}},F_*(B^{(n)}_1))),
\end{gather*}
where $\mu_{\max}$ means the value of $\mu$ at the
maximal destabilizer.
Thus, in any case, we have the inequality
\begin{align*}
 &\mu((B^{(n+1)}_1,B^{(n+1)}_2,\Phi_{B^{(n+1)}},F_*(B^{(n+1)}_1))) \\
 \leq & \mu((B^{(n+1)}_1/C^{(n)}_1,B^{(n+1)}_2/C^{(n)}_2,
 \Phi_{B^{(n+1)}/C^{(n)}},F_*(B^{(n+1)}_1/C^{(n)}_1))) \\
 \leq &\mu((B^{(n)}_1,B^{(n)}_2,\Phi_{B^{(n)}},F_*(B^{(n)}_1)))
\end{align*}
with equality if and only if $(C^{(n)}_1,C^{(n)}_2)= (0,0)$.

The descending sequence
\[
 \{ \mu((B^{(n)}_1,B^{(n)}_2,\Phi_{B^{(n)}},
 F_*(B^{(n)}_1))) \}_{n\in{\bf N}}
\]
must become stationary since it is bounded below.
We may assume without loss of generality that
\[
 \mu((B^{(n)}_1,B^{(n)}_2,\Phi_{B^{(n)}},F_*(B^{(n)}_1)))
\]
is constant for all $n$.
In this case we must have
$(C^{(n)}_1,C^{(n)}_2)= (0,0)$ and
\[
 (B^{(n+1)}_1,B^{(n+1)}_2,\Phi_{B^{(n+1)}},F_*(B^{(n+1)}_1))
\]
becomes a parabolic $\Lambda^1_{\cD_k}$-subtriple of
\[
 (B^{(n)}_1,B^{(n)}_2,\Phi_{B^{(n)}},F_*(B^{(n)}_1))
\]
for all $n$.
Since the descending sequence
$\{ \rank B^{(n)}_1 + \rank B^{(n)}_2 \}_{n\in{\bf N}}$
must be stationary,
we may assume without loss of generality that
$\rank B^{(n)}_1 + \rank B^{(n)}_2$ is constant for all $n$.
Then we must have
\[
 (B^{(n)}_1,B^{(n)}_2,\Phi_{B^{(n)}},F_*(B^{(n)}_1))= 
 (B^{(n+1)}_1,B^{(n+1)}_2,\Phi_{B^{(n+1)}},F_*(B^{(n+1)}_1))
\]
for all $n$.
Thus the sequences (\ref{exact-sequence}) split
and 
\[
 (E^{(n)}_1,E^{(n)}_2,\Phi^{(n)},F_*(E^{(n)}_1))\otimes k
 \cong (B^{(n)}_1,B^{(n)}_2,\Phi_{B^{(n)}},F_*(B^{(n)}_1))\oplus
 (G^{(n)}_1,G^{(n)}_2,\Phi_{G^{(n)}},F_*(G^{(n)}_1)).
\]
Then all the maps $G^{(n)}_j\ra G^{(n+1)}_j$
are isomorphisms.
Since $B^{(n+1)}_j\ra B^{(n)}_j$ are all isomorphic,
every image of
$E^{(n)}_j\otimes k \ra E^{(0)}_j\otimes k$
is $B^{(0)}_j$ for $j= 1,2$.
So we have an isomorphism
$(E^{(0)}_j/E^{(n)}_j)\otimes k\cong G^{(0)}_j$
for any $n$.
On the other hand, every image of
$m^n/m^{n+1}\otimes E^{(0)}_j \ra E^{(n)}_j\otimes k$
is $G^{(n-1)}_j$.
So we have an isomorphism
$(E^{(n)}_j/m^n E^{(0)}_j)\otimes k \cong B^{(n-1)}_j$.
Consider the exact sequence
\[
 0\lra E^{(n)}_j/m^n E^{(0)}_j \stackrel{u}\lra
 E^{(0)}_j/m^n E^{(0)}_j \lra E^{(0)}_j/E^{(n)}_j \lra 0.
\]
Then
$u\otimes k:(E^{(n)}_j/m^n E^{(0)}_j)\otimes k \cong B^{(n-1)}_j
\ra E^{(0)}_j\otimes k$ is injective.
Thus $u$ is injective and $E^{(0)}_j/E^{(n)}_j$
is flat over $R/m^nR$.
Then quotients $E^{(0)}_j\otimes R/m^n \ra E^{(0)}_j/E^{(n)}_j$
define a system of morphisms
$\Spec R/m^n \ra Q'_j:= 
\Quot_{E^{(0)}_j/\cX_R/R}^{\chi(G^{(0)}_j(n))}$
which induces a morphism $f_j:\Spec\hat{R}\ra Q'_j$,
where $\hat{R}$ is the completion of $R$.
If $\tilde{G}_j$ is the quotient sheaf
of $E^{(0)}_j\otimes\hat{R}$
corresponding to $f_j$, then we have
$\tilde{G}_j\otimes R/m^nR \cong E^{(0)}_j/E^{(n)}_j$.
Similarly we can lift the parabolic structure
$F_*(G^{(0)}_1)$ to a flat family $F_*(\tilde{G}_1)$
of parabolic structure on $\tilde{G}_1$ over $\hat{R}$.
We can also lift $\Phi_{G^{(0)}}$ to
$\Phi_{\tilde{G}}:\Lambda^1_{\cD/S}\otimes\tilde{G}_1\ra \tilde{G}_2$ 
and
$(\tilde{G}_1,\tilde{G}_2,\Phi_{\tilde{G}},F_*(\tilde{G}_1))$
becomes a flat family of parabolic $\Lambda^1_{\cD_R}$-triples
which is a quotient of
$(E^{(0)}_1,E^{(0)}_2,\Phi^{(0)},F_*(E^{(0)}_*))\otimes\hat{R}$.
If $\hat{K}$ is the quotient field of $\hat{R}$, then
$(\tilde{G}_1,\tilde{G}_2,\Phi_{\tilde{G}},F_*(\tilde{G}_1))
\otimes\hat{K}$
becomes a destabilizing quotient parabolic
$\Lambda^1_{\cD_{\hat{K}}}$-triple of
$(E_1,E_2,\Phi,F_*(E_1))\otimes\hat{K}$, which contradicts
the semistability of $(E_1,E_2,\Phi,F_*(E_1))$.
\end{proof}

As a corollary of Proposition \ref{valuative-criterion},
we obtain the following proposition:

\begin{Proposition}\label{prop-projective}
 Assume that $\alpha'_1,\ldots,\alpha'_l$ are sufficiently general
 so that all the semistable parabolic
 $\Lambda^1_{\cD/S}$-triples are stable.
 Then the moduli scheme
 $\overline{M_{\cX/S}^{\cD,\balpha',\bbeta,\gamma}}(r,d,\{d_i\})$
 is projective over $S$.
\end{Proposition}

There is another corollary of Proposition \ref{valuative-criterion}
which is used in the proof of the surjectivity of the
Riemann-Hilbert morphism in Lemma \ref{lem:surjective}.
For a parabolic connection $(E,\nabla,\varphi,\{l_i\})$, let
$(0, 0)=(F_0,\nabla_0)\subset(F_1,\nabla_1)\subset\cdots\subset
(F_l,\nabla_l)=(E,\nabla)$
be a Jordan-H\"{o}lder filtration of $(E,\nabla)$,
that is, each $(F_i/F_{i+1},\overline{\nabla_i})$ is irreducible,
where $\overline{\nabla_i}:F_i/F_{i+1}\ra F_i/F_{i+1}\otimes\Omega^1_X(D(\bt))$
is the connection induced by $\nabla_i$.
Then we put
\[
 gr(E,\nabla):=\bigoplus_{i=1}^l (F_i/F_{i+1},\overline{\nabla_i}).
\]

\begin{Corollary}\label{v-c-connection}
Let $R$ be a discrete valuation ring with quotient field $K$ and
residue field $k$.
Let $(E,\nabla,\varphi,\{l_i\})$ be a flat family of connections
with parabolic structures on $X\times\Spec R$ over $R$
such that the generic fiber
$(E,\nabla,\varphi,\{l_i\})\otimes_RK$ is $\balpha$-semistable.
Then there exists a flat family
$(\tilde{E},\tilde{\nabla},\tilde{\varphi},\{\tilde{l}_i\})$
of $\balpha$-semistable parabolic connections such that
$(\tilde{E},\tilde{\nabla},\tilde{\varphi},\{\tilde{l}_i\})\otimes K
\cong (E,\nabla,\varphi,\{l_i\})\otimes K$ and
$gr((\tilde{E},\tilde{\nabla})\otimes k)
\cong gr((E,\nabla)\otimes k)$.
\end{Corollary}

%%%%%%%%%%%%%%%%%%%%%%%%%%%%%%%%%%%%%%%%%%%%%%%%%%%%%%%%%%%%%%%%%%%%%%%%%
%%%%%%%%%%%%%%%%%%%%%%%%%%%%%%%%%%%%%%%%%%%%%%%%%%%%%%%%%%%%%%%%%%%%%%%%%
%%%%%%%%%%%%

\subsection{Proof of Theorem \ref{thm:fund} (1)}

Now we prove the assertion (1) of Theorem \ref{thm:fund}.  

We take $S$ for $T_n\times\Lambda_n$ and
$\cX$ for $\BP^1\times T_n\times\Lambda_n$.

Let $\cD_i\subset\BP^1\times T_n\times\Lambda_n$
be the effective divisor determined by the section
\[
 T_n\times\Lambda_n\hookrightarrow
 \BP^1\times T_n\times\Lambda_n; \quad
 ((t_j)_{1\leq j\leq n},(\lambda_k)_{1\leq k\leq n})
 \mapsto (t_i,(t_j)_{1\leq j\leq n},(\lambda_k)_{1\leq k\leq n})
\]
for $i= 1,\ldots,n$ and put $\cD:= \sum_{i= 1}^n \cD_i$.
Then $\cD$ becomes an effective Cartier divisor on
$\BP^1\times T_n\times\Lambda_n$
which is flat over $T_n\times\Lambda_n$.

We fix a line bundle $L$ on $\BP^1\times T_n\times\Lambda_n$
with a relative connection
\[
 \nabla_L: L\ra L\otimes
 \Omega^1_{\BP^1\times T_n\times\Lambda_n/T_n\times\Lambda_n}(\cD)
\]
over $T_n\times\Lambda_n$.  Let $\balpha'=(\alpha'_1, \ldots, \alpha'_{2n}) $, 
$\bbeta=(\beta_1, \beta_2)$, and $\gamma \gg 0$  be as in 
Theorem \ref{thm:fund}.

We define a moduli functor
$\overline{\cM^{\balpha'\bbeta}_n}(L)$
of the category of locally noetherian schemes over $T_n\times\Lambda_n$ to
the category of sets by
\[
 \overline{\cM^{\balpha'\bbeta}_n}(L)(S)
 := \{(E_1,E_2,\phi,\nabla,\varphi,\{l_i\}_{i= 1}^n)\}/\sim ,
\]
where $S$ is a locally noetherian scheme over $T_n\times\Lambda_n$
corresponding to
$(\bt,\blambda)= (t_1,\ldots,t_n,\lambda_1,\ldots,\lambda_n)\in
T_n(S)\times\Lambda_n(S)$ and
\begin{enumerate}
\item
 $E_1,E_2$ are rank $2$ vector bundles on $\BP^1\times S$,
\item
 $\phi:E_1\ra E_2$ is an $\cO_{\BP^1\times S}$-homomorphism,
 $\nabla:E_1\ra E_2\otimes\Omega^1_{\BP^1}(D(\bt))$
 is a morphism such that
 $\nabla(fa)= \phi(a)\otimes df+f\nabla(a)$
 for $f\in\cO_{\BP^1\times S}$, $a\in E_1$,
\item
 $l_i\subset E_1|_{t_i}$ are rank $1$ subbundles
 such that 
 $(\res_{t_i}(\nabla)-
 \lambda_i\phi|_{t_i})|_{l_i}= 0$
 for $i= 1,\ldots,n$,
\item
 $\varphi:\bigwedge^2 E_2\stackrel{\sim}\lra L\otimes\cL_{\varphi}$
 is an isomorphism such that
 $(\varphi\otimes 1)
 (\nabla(s_1)\wedge\phi(s_2)+\phi(s_1)\wedge\nabla(s_2))
 = (\nabla_L\otimes\mathrm{id}_{\cL_{\varphi}})
 (\varphi(\phi(s_1)\wedge\phi(s_2)))$
 for $s_1,s_2\in E_1$ where $\cL_{\varphi}$ is a line bundle on $S$,
\item
 for any geometric point $s$ of $S$, the fiber
 $((E_1)_s,(E_2)_s,\phi_s,\nabla_s,\varphi_s,
 \{l_i|_{t_i\otimes k(s)}\}_{i= 1}^n)$ is $(\balpha', \bbeta)$-stable
 and $\deg(E_1)_s= \deg L_s$.
\end{enumerate}

Here $(E_1,E_2,\phi,\nabla,\varphi,\{l_i\})\sim
(E'_1,E'_2,\phi',\nabla',\varphi',\{l'_i\})$
if there exist a line bundle $\cL$ on $S$ and isomorphisms
$\sigma_j:E_j\stackrel{\sim}\to E'_j\otimes\cL$
for $j= 1,2$ such that
$\sigma_1|_{t_i\times S}(l_i)= l'_i$ for any $i$, the diagrams
\[
 \begin{CD}
  E_1 @>\phi >> E_2 \\
  @V\sigma_1 V\cong V @V\cong V\sigma_2 V \\
  E'_1\otimes \cL @>\phi' >> E'_2\otimes \cL
 \end{CD}
 \quad \text{and} \quad
 \begin{CD}
  E_1 @>\nabla >> E_2\otimes\Omega^1_{\BP^1}(D(\bt)) \\
  @V\sigma_1 V\cong V @V\cong V\sigma_2\otimes\mathrm{id} V \\
  E'_1\otimes \cL 
  @>\nabla' >>
  E'_2\otimes\Omega^1_{\BP^1}(D(\bt))\otimes \cL 
 \end{CD}
\]
commute and there is an isomorphism
$\sigma:\cL_{\varphi}\stackrel{\sim}\ra
\cL_{\varphi'}\otimes\cL^{\otimes 2}$ such that the diagram
\[
 \begin{CD}
  \bigwedge^2 E_2 @>\varphi>\sim> L\otimes\cL_{\varphi} \\
  @V\wedge\sigma_2 V \cong V @V \cong V\mathrm{id}\otimes\sigma V \\
  \bigwedge^2 E'_2\otimes\cL @>\varphi'\otimes\mathrm{id}_{\cL}>\sim>
  L\otimes\cL_{\varphi'}\otimes\cL^{\otimes 2}
 \end{CD}
\]
commutes.

We can define another weight $\balpha  = (\alpha_1, \ldots, \alpha_{2n})$ with
$ 0 \leq \alpha_1 < \cdots < \alpha_{2n} < \frac{\beta_1}{\beta_1+ 
 \beta_2} < 1 $ by 
 $$
 \balpha = \balpha' \frac{\beta_1}{\beta_1 + \beta_2}. 
 $$

Theorem \ref{thm:fund}, (1) follows from the following theorem:

\begin{Theorem}\label{thm:exist-proj}
 There exists a coarse moduli scheme
 $\overline{M^{\balpha'\bbeta}_n}(L)$ of
 $\overline{\cM^{\balpha'\bbeta}_n}(L)$,
 which is projective over $T_n\times\Lambda_n$ if $\balpha'$ is generic.
 If we put
\[
 M^{\balpha}_n(L):=\left.\left\{
 (E_1,E_2,\phi,\nabla,\varphi,\{l_i\}) \in
 \overline{M^{\balpha'\bbeta}_n}(L)
 \right| \text{\rm $\phi:E_1\ra E_2$ is an isomorphism} \right\},
\]
then $M^{\balpha}_n(L)$ is a Zariski open subset of
$\overline{M^{\balpha'\bbeta}_n}(L)$,
which is a fine moduli scheme of $\balpha$-stable parabolic connections.
\end{Theorem}

\begin{proof}
We put $r= 2$, $d= \deg L_s$ for $s\in T_n\times\Lambda_n$,
$l= 2n$ and $d_i= i$ for $i= 1,\ldots,2n$
and consider the moduli scheme 
$\overline{M_{\BP^1\times T_n\times\Lambda_n/T_n\times\Lambda_n}
^{\cD,\balpha',\bbeta,\gamma}}(2,d,\{d_i\})$.
For each
$(E_1,E_2,\phi,\nabla,\varphi,\{l_i\})
\in\overline{\cM^{\balpha'\bbeta}_n}(L)(S)$,
let $\Phi:\Lambda^1_{\cD_S}\otimes E_1 \ra E_2$
be the left $\cO_{\BP^1\times S}$-homomorphism
corresponding to $(\phi,\nabla)$ and
put $F_{2i+1}(E_1):= E_1(-\sum_{j= 1}^i t_j)$ for
$i= 0,\ldots,l$ and
$F_{2i}(E_1):= \ker(F_{2i-1}(E_1)\ra (E_1|_{t_i}/l_i))$
for $i= 1,\ldots,l$,
where $(t_1,\ldots,t_n,\lambda_1,\ldots,\lambda_n)
\in T_n(S)\times\Lambda_n(S)$
corresponds to the structure morphism
$S\ra T_n\times\Lambda_n$.
Then the correspondence
$(E_1,E_2,\phi,\nabla,\varphi,\{l_i\})\mapsto
(E_1,E_2,\Phi,F_*(E_1))$ determines a morphism of functors
\[
 \iota:\overline{\cM^{\balpha'\bbeta}_n}(L)\ra
 \overline{\cM_{\BP^1\times T_n\times\Lambda_n/T_n\times\Lambda_n}
 ^{\cD,\balpha',\bbeta,\gamma}}(2,d,\{d_i\}).
\]
We can easily see that $\iota$ is represented by a closed immersion.
Recall that
$R^s\ra \overline{M_{\BP^1\times T_n\times\Lambda_n/T_n\times\Lambda_n}
^{\cD,\balpha',\bbeta,\gamma}}(2,d,\{d_i\})$
is a principal $G$-bundle.
Then there exists a closed subscheme $Z\subset R^s$
such that
\[
 h_Z= h_{R^s}\times_{\overline{\cM_{\BP^1\times T_n\times\Lambda_n
 /T_n\times\Lambda_n}^{\cD,\balpha',\bbeta,\gamma}}(2,d,\{d_i\})}
 \overline{\cM^{\balpha'\bbeta}_n}(L).
\]
$Z$ descends to a closed subscheme of
$\overline{M_{\BP^1\times T_n\times\Lambda_n/T_n\times\Lambda_n}
^{\cD,\balpha',\bbeta,\gamma}}(2,d,\{d_i\})$
which is just the coarse moduli scheme of
$\overline{\cM^{\balpha'\bbeta}_n}(L)$. 

If we take $\gamma$ sufficiently large,
we can check that a parabolic connection
$(E,\nabla_E,\varphi,\{l_i\})$ is $\balpha$-stable
if and only if the associated parabolic $\phi$-connection
$(E,E,\mathrm{id}_E,\nabla_E,\varphi,\{l_i\})$ is $(\balpha',\bbeta)$-stable.
Thus the open subscheme
\[
 M^{\balpha}_n(L):=\left\{ \left.
 (E_1,E_2,(\phi,\nabla),F_*(E_1))\in \overline{M^{\balpha'\bbeta}_n}(L)
 \right| \text{$\phi:E_1\ra E_2$ is an isomorphism} \right\}
\]
of $\overline{M^{\balpha'\bbeta}_n}(L)$ is just the moduli space
of $\balpha$-stable parabolic connections with the determinant $L$.

If $\deg L$ is odd, we can see by the same argument as
[\cite{M}, Theorem 6.11] or [\cite{HL}, Theorem 4.6.5] that
$M^{\balpha}_n(L)$ is in fact a fine moduli scheme.
If $\deg L$ is even, then
we can obtain, by an elementary transform, an isomorphism
\[
 \sigma: \cM_n(L) \ra \cM_n(L')
\]
of moduli stacks of parabolic connections without stability condition,
where $\deg L'$ is odd.
Then we can see by the same argument that $\sigma(M^{\balpha}_n(L))$
becomes a fine moduli scheme, and so $M^{\balpha}_n(L)$ is also fine.
\end{proof}

\vspace{1cm}
\section{Tangent spaces of the moduli spaces and Canonical symplectic structure. }
\label{sec:symp}

In this section, we will work over the  finite \'etale covering $T'_n \lra T_n$ defined in 
\eqref{eq:finite-etale}.
Fix $(\bt, \blambda) \in T'_n \times \Lambda_n$ and set 
$a_i = 2 \cos 2 \pi \lambda_i $ and $\ba = (a_1, \ldots, a_n) $.  
Moreover fix a determinant line bundle $ L = (L, \nabla_L) $ on $ \BP^1 $ such that 
$\res_{t_i} (\nabla_L) \in \Z$. We have defined two moduli spaces 
$
M_n^{\balpha}(\bt, \blambda, L),  \cR(\cP_{n, \bt})_{\ba}
$
where $M_n^{\balpha}(\bt, \blambda, L)$ is the moduli space of stable 
\index{something@}{$(\bt, \blambda)$-parabolic connections}  with the determinant $L$ and 
$ \cR(\cP_{n, \bt})_{\ba} $ is the moduli of Jordan equivalence classes of the 
$SL_2(\C)$-representations of $\pi_1(\BP^1 \setminus D(\bt), *)$ with fixed  
local exponents  $\ba =(a_1, \ldots, a_n)$. 
As we show in Theorem \ref{thm:fund},  for a suitable  (or generic) 
weight  $\balpha$,  the moduli space $M_n^{\balpha}(\bt, \blambda, L)$ 
is a non-singular complex scheme.   
In this section, we will describe the tangent space to 
$M_n^{\balpha}(\bt, \blambda, L)$ and a non-degenerate holomorphic $2$-form on the
moduli space $M_n^{\balpha}(\bt, \blambda, L)$.  

Although the moduli space $\cR(\cP_{n, \bt})_{\ba}$ may be singular, 
we can define a 
Zariski dense open set 
$\cR(\cP_{n, \bt})^{\sharp}_{\ba}$ of $\cR(\cP_{n, \bt})_{\ba}$ 
such that $\cR(\cP_{n, \bt})^{\sharp}_{\ba}$ is a non-singular variety.  
(Note that for generic $\ba \in \cA_n$,   
$\cR(\cP_{n, \bt})^{\sharp}_{\ba} =  \cR(\cP_{n, \bt})_{\ba}$).  
Moreover 
on  $\cR(\cP_{n, \bt})^{\sharp}_{\ba}$ we can also define a canonical 
symplectic structure $\Omega_1$.  
In \S \ref{sec:R-H}  we  define  the Riemann-Hilbert correspondence 
${\bf RH}_{\bt, \blambda}: M_{n}^{\balpha}(\bt, \blambda, L)  \lra  \cR(\cP_{n, \bt})_{\ba}$.
We  show that   ${\bf RH}_{\bt, \blambda}$
 is bimeromorphic proper surjective morphism and 
gives an analytic isomorphism between 
$M_{n}^{\balpha}(\bt, \blambda, L)^{\sharp}:= 
{\bf RH}_{\bt, \blambda}^{-1}(\cR(\cP_{n, \bt})^{\sharp}_{\ba}) $  
and $ \cR(\cP_{n, \bt})^{\sharp}_{\ba}$. 
(Again, for a generic $\blambda$, 
$M_{n}^{\balpha}(\bt, \blambda, L)^{\sharp} = M_{n}^{\balpha}(\bt, \blambda, L)$).  
Note that ${\bf RH}_{\bt, \blambda}$ 
is not an algebraic morphism, and hence the algebraic structures of 
$M_{n}^{\balpha}(\bt, \blambda, L)^{\sharp}$ 
and $\cR(\cP_{n, \bt})^{\sharp}_{\ba}$ are completely different).  
The canonical symplectic structures on both moduli spaces can be identified via 
${\bf RH}_{\bt, \blambda|M_{n}^{\balpha}(\bt, \blambda, L)^{\sharp}}$, 
that is,  
$\left({\bf RH}_{\bt, \blambda|M_{n}^{\balpha}(\bt, \blambda, L)^{\sharp}}\right)^{*}(\Omega_1) 
= \Omega$.

\subsection{Tangent space to $M_n^{\balpha}(\bt, \blambda, L)$.}

Consider the base extension of the family of moduli spaces in (\ref{eq:fam-moduli}) by the 
\'etale covering $T'_n \lra T_n$:
\begin{equation}
\pi_n:M^{\balpha}_{n}(L) \lra T'_n \times \Lambda_n, 
\end{equation}
such that for every $(\bt, \blambda) \in T'_n \times \Lambda_n $, we have    $
\pi_n^{-1}((\bt, \lambda)) \simeq  M^{\balpha}_n(\bt, \blambda, L)$. 
For simplicity,  we will omit 
$L$ from now on, so we write as  $M^{\balpha}_n = M^{\balpha}_n(L)$ , $M^{\balpha}_n(\bt, \blambda) =  M^{\balpha}_n(\bt, \blambda, L)$.    
We assume that  $\balpha$ is  generic so that $\pi_n$ is a smooth morphism (cf. Theorem 
\ref{thm:fund}). 

Let us consider  natural projection maps 
$$
\begin{array}{rcccl}
&&  M_n^{\balpha} &&  \\
&& \quad  \downarrow \pi_n &&  \\
&&  T'_n \times \Lambda_n &&  \\
&  \stackrel{p_1}{\swarrow}  &    &  
\stackrel{p_2}{\searrow} &  \\
T'_n &   &     & & \Lambda_n\\
\end{array}
$$
and set $\varphi_i = p_i \circ \pi_n $.  
Since $\varphi_1:M^{\balpha}_n \lra T'_n $ 
is smooth, 
we have the following exact sequence of 
tangent sheaves on $M_{n}^{\balpha}$
%%%%%%%%%%%%%%
\begin{equation}
 0 \lra \Theta_{M_{n}^{\balpha}/T'_n \times \Lambda_n} \lra \Theta_{M_{n}^{\balpha}/T'_n} 
 \lra \pi_n^{*} (\Theta_{T'_n\times \Lambda_n / T'_n}) \lra 0.   
\end{equation}
%%%%%%%%%%%%
We will describe this exact sequence in terms of the infinitesimal deformation of 
the stable parabolic connections.  
Let us consider the natural projection map
$ q_2 : \BP^1 \times T'_n \lra T'_n$ and defines a divisor 
$\cD \subset \BP^1 \times T'_n $ such that 
$q_2^{-1} (\bt) \cap \cD =D(\bt) =t_1 + \cdots + t_n \subset \BP^1 $.

Let $(\tilde{E}, \tilde{\nabla}, \tilde{\varphi}, \{ \tilde{l_i} \})$ be 
a universal family on 
$\BP^1 \times M_{n}^{\balpha}$.  
Consider  the following commutative diagram:
\begin{equation}
\begin{array}{ccccc}
  & & \BP^1 \times M_n^{\balpha} &  & \\
  &\stackrel{\tilde{\varphi}_1}{\swarrow}  &   \downarrow \tilde{\pi}_n  & &   \\
   &  &  \BP^1 \times T'_n \times \Lambda_n &&  \\
  & \swarrow &  &   &   \\
   \BP^1 \times T'_n & & & & \\
  q_2 \downarrow \quad & & & & \\
   T'_n & & &  \\
 \end{array}
\end{equation}

For a coherent sheaf $\cG$ on 
$\BP^1 \times M_{n}^{\balpha}$ and a closed point 
$\x \in M_{n}^{\balpha}$, we set 
$\cG_{\x} := \cG_{| \BP^1 \times \x }$. 

We define coherent sheaves  on  $\BP^1\times M^{\balpha}_{n}$ as follows.
\begin{equation}
\label{eq:end-p0}
 \cF^0 :=\left\{ s \in {\mathcal End}(\tilde{E})  |
 \Tr(s)=0, ;  \  (s|_{t_i\times M^{\balpha}_n})(\tilde{l}_i)
 \subset \tilde{l}_i \right\}
 \end{equation}
\begin{equation}\label{eq:end-p1}
\cF^1 :=\left\{ s \in {\mathcal End}(\tilde{E})\otimes
\tilde{\varphi}_1^{*}( \Omega^1_{\BP^1 \times T'_n/T'_n}(\cD))
 \ | \ \Tr(s)=0,\; \ (s|_{t_i\times M^{\balpha}_n})(\tilde{l}_i)=0
 \right\}
\end{equation}
\begin{equation}
\label{eq:end-p1s}
 \cF^{1,+} :=\left\{ s \in {\mathcal End}(\tilde{E}) \otimes\tilde{\varphi}_1^{*}( \Omega^1_{\BP^1 \times T'_n/T'_n}(\cD)) | 
 \Tr(s)=0, \;  \  \res_{(t_i\times M^{\balpha}_n)}(s)(\tilde{l}_i)
 \subset\tilde{l}_i \right\}
\end{equation}
For a local section $s$ of $\cF^0$, define
$\nabla_1(s) := \tilde{\nabla} s - s \tilde{\nabla}$.  
Then it is easy to see that  $\nabla_1(s)$ is a local section of $\cF^1$. Since we have a natural inclusion of sheaves 
$\iota: \cF^1 \hookrightarrow \cF^{1,+}$, we can define two complexes of sheaves on $\BP^1 \times M_{n}^{\balpha}$:
\begin{equation}\label{eq:complex-org}
\cF^{\bullet}:= \left[\nabla_1: \cF^0 \lra \cF^1 \right], 
\end{equation} 
\begin{equation}\label{eq:complex+}
\cF^{\bullet, +}: =\left[\nabla_{1}^{+}: \cF^0 \lra \cF^{1, +}\right].  
\end{equation} 

Let $\x \in M^{\balpha}_{n}$ be a closed point and set $\pi(\x) = (\bt, \blambda)$. 
Setting $\cT_1 = \cF^{1,+}/\cF^1$, we have the following exact sequences of the 
complexes on $\BP^1 \times M_{n}^{\balpha}$
and $\BP^1 \times \{ \x \}$.  
\begin{equation} \label{eq:complex}
\begin{array}{cccc}
    &                                & 0           &   \\  
    &                                & \downarrow  &    \\
\cF^0       &  \stackrel{\nabla_1}{\longrightarrow} & \cF^1       &  \\
|| &                                & \downarrow  &    \\
\cF^0 & \stackrel{\nabla_{1}^{+}}{\lra} & \cF^{1,+} &   \\
 \downarrow&                                & \downarrow  &    \\
0 & \lra                               & \cT_1  &     \\
 &                                & \downarrow  &    \\
 &                                & 0  &   \\
\end{array}
 \mbox{on $\BP^1 \times M_{n}^{\balpha}$},  
\hspace{0.5cm}
\begin{array}{cccc}
    &                                & 0           &   \\  
    &                                & \downarrow  &    \\
\cF^0_{\x}       &  \stackrel{\nabla_{1, \x}}{\longrightarrow} & \cF^1_{\x}       &  \\
|| &                                & \downarrow  &    \\
\cF^0_{\x} & \stackrel{\nabla_{1, \x}^{+}}{\lra} & \cF^{1,+}_{\x} &   \\
\downarrow &                                & \downarrow  &    \\
0 & \lra                              &\cT_{1, \x}  &     \\
 &                                & \downarrow  &    \\
 &                                & 0 &.    \\
\end{array}  \mbox{on $\BP^1 \times \{\x\}$} 
\end{equation}

Note that at each point  $t_i$, $1 \leq i \leq n$, 
the stalk $(\cT_{1, \x})_{t_i}$ is isomorphic to $ 
\C((t_i, \x))$, hence 
$H^0(\cT_{1,\x}) \simeq \oplus_{i=1}^n \C((t_i, \x)) \simeq \C^n$.

\begin{Lemma}\label{lem:tangent}   
At each closed point $ \x \in M^{\balpha}_n(\bt, \blambda) \subset M^{\balpha}_{n}$ the 
tangent spaces can be  given as follows.  
\begin{equation}\label{eq:tangentall}
(\Theta_{M^{\balpha}_{n}/T'_n})_{\x} \simeq 
{\bf H}^1(\BP^1, [\cF^0_{\x} \stackrel{\nabla_{1,\x}^+}{\lra} \cF^{1,^+}_\x ]), 
\end{equation}
\begin{equation} \label{eq:tangentrel}
( \Theta_{M^{\balpha}_n/T'_n \times \Lambda_n})_{\x} 
\simeq {\bf H}^1(\BP^1, 
[\cF^0_{\x} \stackrel{\nabla_{1,\x}}{\lra} \cF^{1}_{\x} ] ),  
\end{equation}
\begin{equation}
(\Theta_{T'_n \times \Lambda_n/T'_n})_{\pi(\x)} \simeq  H^0(\cT_{1, \x}) \simeq \C^n.  
\end{equation}
Under these isomorphisms, we have the following identification of the natural 
exact sequences of the tangent spaces with the exact sequences of the hypercohomologies:
\begin{equation}
\begin{array}{ccccccc}
 0 \lra & 
( \Theta_{M^{\balpha}_n/T'_n \times \Lambda_n})_{\x} &  
\lra &
(\Theta_{M^{\balpha}_{n}/T'_n})_{\x} &
\lra & 
 (\Theta_{T'_n \times \Lambda_n/T'_n})_{\pi(\x)}&
 \lra 0 \\
 & & & & & & \\
&
|| 
&   
& 
||  
& 
& 
||
&  \\
& & & & & &  \\ 
 0 \lra & 
  {\bf H}^{1}([\cF^0_{\x}  \stackrel{\nabla_{1,\x}}{\lra} \cF^1_{\x}
  ]) & 
   \lra & 
   {\bf H}^{1}([\cF^0_{\x}  \stackrel{\nabla^+_{1, \x}}{\lra} \cF^{1,+}_{\x}] ) &
\lra & 
H^0(\cT_{1, \x})&  \lra 0.  \\
\end{array} 
\end{equation}
\end{Lemma}

\begin{proof} The smoothness of 
the natural map $\pi_n:M^{\balpha}_{n} \lra T'_n \times  \Lambda_n$ 
  follows from Theorem \ref{thm:fund}.  (Actually, one can show that  
  ${\bf H}^{2}([\cF^0_{\x}  \stackrel{\nabla_{1,\x}}{\lra} \cF^1_{\x}
  ]) = \{ 0 \}$ (cf. Lemma \ref{lem:h2})).   
The space of the infinitesimal deformations of logarithmic parabolic connection 
with fixing the  eigenvalues of the residue matrix  of $\tilde{\nabla}_{\x, t_i}$ at $t_i$ 
is given by 
the hypercohomology
$$
{\bf H}^1(\BP^1, [\cF^0_{\x}  \stackrel{\nabla_{1,\x}}{\lra} \cF^1_{\x}] ).
$$
(Cf.\ Arinkin \cite{A}).  Moreover it is easy to see that 
${\bf H}^1(\BP^1, [\cF^0_{\x}  \stackrel{\nabla^+_{1,\x}}{\lra} \cF^{1,+}_{\x}] )$ 
is the set of infinitesimal deformations of logarithmic parabolic connections without 
fixing the eigenvalues of the residues of $\tilde{\nabla}_{\x}$.   

Since $\cT_{1, \x}$ is a skyscraper sheaf supported on 
$D(\bt) \subset \BP^1 \times \{\x\}$, we see that 
${\bf H}^0 (0 \ra \cT_{1,\x}) = \{0 \}$, ${\bf H}^1 (0 \ra \cT_{1,\x}) = H^0(\cT_{1, \x}) \simeq \C^n$.  
 Local calculations of the maps $\nabla_1, \nabla_{1,+}$ in 
  the commutative diagram in (\ref{eq:complex}) show that 
  the natural map 
  $$
d\pi_{n, \x}: {\bf H}^1([\nabla^+_{1,\x}:\cF^0_{\x} \lra \cF^{1,+}]) \lra H^0(T_1)  
  $$
  gives the differential of the map $\pi_n:M^{\balpha}_n \lra \Lambda_n $ at 
  $\x$. Since  ${\bf H}^{2}([\cF^0_{\x}  \stackrel{\nabla_{1,\x}}{\lra} \cF^1_{\x}
  ]) = \{ 0 \}$ or equivalently $\pi_n$ is smooth at $\x$, the map $d \pi_{\x}$ is surjective.  
\end{proof}

\subsection{The relative symplectic 
form $\Omega$ for $\pi_n: M^{\balpha}_n \lra T'_n \times 
\Lambda_n $.}

Let us consider the smooth family of moduli spaces of stable parabolic connections: 
\begin{equation}
\pi_n:M^{\balpha}_n \lra T'_n \times \Lambda_n. 
\end{equation}
Now we will  show that each closed 
fiber $\pi_n^{-1}(\bt, \blambda) = M^{\balpha}_n(\bt, \blambda, L)$ 
admits a canonical {\em symplectic structure} $\Omega$, which induces 
a non-degenerate skew symmetric bilinear form on its tangent sheaf:
\begin{equation}\label{eq:form}
\Omega_{|M^{\balpha}_n(\bt, \blambda, L)}: \Theta_{M^{\balpha}_n(\bt, \blambda, L)}
 \otimes \Theta_{M^{\balpha}_n(\bt, \blambda, L)} \lra \cO_{M^{\balpha}_n(\bt, \blambda, L)}. 
\end{equation}

First, a local calculation shows the following 

\begin{Lemma}\label{lem:duality}
For each point  $ \x \in M^{\balpha}_n(\bt, \blambda, L) = \pi_n^{-1}(\bt, \blambda) 
\subset M^{\balpha}_n $,  set 
 $ \cF^i_{\x} = \cF^i_{|\BP^1 \times \x}$ for 
$ i = 0, 1$.  Then we have  isomorphisms 
\begin{equation}\label{eq:duality}
\cF^{1}_{\x} \simeq \cF^{0 \vee}_{\x}
\otimes \Omega^1_{\BP^1}, \quad 
\cF^{0}_{\x} \simeq \cF^{1 \vee}_{\x} \otimes \Omega^1_{\BP^1}. 
\end{equation}
where $\cF^{i \vee}_{\x} = {\mathcal Hom}(\cF^i_{\x}, \cO_{\BP^1})$. 
\end{Lemma}

The following lemma is a key of proof of  
the smoothness of the moduli space $M^{\balpha}_n(\bt, \blambda, L)$. 
The stability assumption on the objects in $M^{\balpha}_n(\bt, \blambda, L)$ 
is essential in this lemma.

\begin{Lemma}\label{lem:h2} Under the notation as above, we have 
\begin{equation} \label{eq:h2}
{\bf H}^2 (\BP^1, \cF_{\x}^{\bullet} ) =  \{ 0 \}. 
\end{equation}
\end{Lemma}

\begin{proof} 
Consider the dual complex $(\nabla^1)^{\vee}: (\cF^1)^{\vee} \otimes 
\Omega^1_{\BP^1} 
\lra (\cF^0)^{\vee} \otimes \Omega^1_{\BP^1}$ which  can be identified 
with the original 
complex $\nabla^1$ by Killing form  (cf.  Lemma \ref{lem:duality}).  
Therefore 
$$
\begin{array}{ccc}
{\bf H}^2(\cF^{\bullet})  &  \simeq  &  
\coker \left[ H^1( \cF_{\x}^0) \stackrel{\nabla^1}{\lra} H^1( 
\cF_{\x}^1) \right]  \\
&  &   \\
& \simeq &  \ker \left[ H^1( \cF_{\x}^1)^{\vee} 
\stackrel{(\nabla^1)^{\vee}}{\lra} H^1(\cF_{\x}^0)^{\vee} \right]^{\vee} 
\\ 
&  &  \\
&\simeq& \ker \left[ H^0(\cF_{\x}^0) \stackrel{\nabla^1}{\lra} 
H^0(\cF_{\x}^1) \right]^{\vee}. 
\end{array}
$$ 
Since $\cF_{\x}^0$ is in the trace free part of the endomorphisms, 
it suffices to show that any $s \in H^0(\cF_{\x}^0) $ such that
 $s \nabla =  \nabla s $ 
is  a scalar.  
For any $\lambda \in \C$, let us set 
$E^0_{\lambda} =  \ker (s - \lambda)$  and 
$E^1_{\lambda} =  \im (s - \lambda) $.  
Then both $E^0_{\lambda}$ and $E^1_{\lambda}$ are  
 subsheaves of $E$ stable under $\nabla$.  
 If $E^0_{\lambda} $ is locally free of rank 1, 
 one  can see that either  $E^0_{\lambda} $ or $E^1_{\lambda}$ 
violates the stability of $E$.  Hence $E^0_{\lambda}$ is zero or 
coincides with $E$.  Therefore $s$ is scalar.  
\end{proof}

\begin{Proposition}\label{prop:symplectic}
 There exists a  global relative 2-form
\begin{equation}\label{eq:symplectic}
\Omega \in H^0(M^{\balpha}_n, \Omega^2_{M^{\balpha}_n/T'_n\times\Lambda_n} ). 
\end{equation}
which induces a symplectic structure on each fiber of $\pi_n$.
\end{Proposition}

\begin{proof} Let us consider the following commutative diagram:
\begin{equation}
\begin{array}{ccc}
\BP^1 \times M^{\balpha}_n & \stackrel{p_2}{\lra} & M^{\balpha}_n \\
\downarrow &    &  \downarrow \\
\BP^1 \times T'_n \times \Lambda_n  & \stackrel{q_2}{\lra} & T'_n \times \Lambda_n \\ 
\end{array}
\end{equation}
Let $\cF^{\bullet}:=\left[\nabla_{1}:\cF^0 \lra \cF^1 \right] $ be 
the complex of sheaves defined in (\ref{eq:complex-org}). From Lemma \ref{lem:tangent},  
we have a natural isomorphism of sheaves:
\begin{equation}
{\bf R}^1 p_{2,*} ( \cF^{\bullet}) \stackrel{\simeq}{\lra} \Theta_{M^{\balpha}_n/T'_n\times\Lambda_n}. 
\end{equation}
By this isomorphism, it suffices  to define a non-degenerate skew-symmetric form 
\begin{equation}\label{eq:rel-symp}
\Omega:{\bf R}^1 p_{2,*} ( \cF^{\bullet}) \otimes {\bf R}^1 p_{2,*}
 ( \cF^{\bullet}) \lra {\bf R}^2 p_{2,*} (
 \Omega^{\bullet}_{\BP^1 \times M^{\balpha}_n/M^{\balpha}_n})
 \cong \cO_{M^{\balpha}_n}. 
\end{equation}
Let us fix a point $ \x \in M^{\balpha}_n(\bt, \blambda, L) = \pi_n^{-1}(\bt, \blambda) 
\subset M^{\balpha}_n $ and define the restriction as 
$\cF^{\bullet}_{\x} = \cF^{\bullet}_{|\BP^1 \times \{\x\}} $. 
{F}rom the following definition of $\Omega(\x)$ 
at the stalk level of (\ref{eq:rel-symp}), it is obvious the definition of the 
global relative 2-form $\Omega$ in (\ref{eq:rel-symp}), and the non-degeneracy of 
$\Omega$ will be  checked at  the stalk of  each closed point $\x$.  

Take an affine open covering $ \{ U_{\alpha} \} $ of $ \BP^1 $ and
consider the following pairing
\begin{equation} \label{eq:mod-sym}
\begin{array}{cccc}
 \Omega(\x): &  \bH^1(\BP^1, \cF^{\bullet}_{\x} ) 
 \otimes \bH^1(\BP^1, \cF^{\bullet}_{\x}) &  
 \lra  &  \bH^2( \Omega_{\BP^1}^{\bullet}) \simeq \C(\x) \\
 & & & \\
 &  ([ \{ v_{\alpha\beta}\},\{u_{\alpha}\}],
 [\{v'_{\alpha\beta}\},\{u'_{\alpha}\}]) &  \mapsto  & 
 [\{\Tr(v_{\alpha\beta}\circ u'_{\beta})
 -\Tr(u_{\alpha}\circ v'_{\alpha\beta})\}]
 -[\{\Tr(v_{\alpha\beta}\circ v'_{\beta\gamma})\}]  \\
\end{array}
\end{equation}
where we consider in \v{C}ech cohomology and
$\{v_{\alpha\beta}\}\in C^1(\cF^0_{\x})$, $\{u_{\alpha}\}\in C^0(\cF^1_{\x})$,
$\{\tilde{\nabla}_{\x}v_{\alpha\beta}
-v_{\alpha\beta}\tilde{\nabla}_{\x}\}=\{u_{\beta}-u_{\alpha}\}$
and so on.
We can check that $\Omega(\x)$ is a skew symmetric pairing. 
Let us show   that $\Omega(\x)$ is non-degenerate for any point
$\x \in M^{\balpha}_n(\bt,\blambda, L)$.
 From Lemma \ref{lem:h2},  one can show that   $\mathbf{H}^2( \cF^{\bullet}_{\x})=0$ for
any $\x\in M^{\balpha}_n(\bt,\blambda, L)$.
$\Omega(\x)$ induces a homomorphism
\[
 \mathbf{H}^1(\cF^{\bullet}_{\x})
 \stackrel{\xi}\lra
 \mathbf{H}^1(\cF^{\bullet}_{\x})^{\vee}.
\]
 From the spectral sequence
$H^p\left(H^q(\cF^0_{\x})\ra H^q(\cF^1_{\x})\right)
\Rightarrow\mathbf{H}^{p+q}(\cF^{\bullet}_{\x})$,
we obtain the following exact sequence
\begin{equation}\label{eq:spectral}
 0 \lra H^0(\cF^0_{\x}) \lra  H^0(\cF^1_{\x})
 \lra \mathbf{H}^1(\cF^{\bullet}_{\x})\lra
 H^1(\cF^0_{\x}) \lra H^1(\cF^1_{\x}) \lra 0.
\end{equation}
Then we obtain the exact commutative diagram
\[
 \begin{CD}
  H^0(\cF^0_{\x}) @>>> H^0(\cF^1_{\x}) @>>>
  \mathbf{H}^1(\cF^{\bullet}_{\x})
  @>>> H^1(\cF^0_{\x}) @>>> H^1(\cF^1_{\x}) \\
  @V b_1 VV @V b_2 VV @V \xi VV @V b_3 VV @V b_4 VV \\
  H^1(\cF^1_{\x})^{\vee} @>>> H^1(\cF^0_{\x})^{\vee}
  @>>> \mathbf{H}^1(\cF^{\bullet}_{\x})^{\vee} @>>>
  H^0(\cF^1_{\x})^{\vee} @>>> H^0(\cF^0_{\x})^{\vee},
 \end{CD}
\]
where $b_1,\ldots,b_4$ are isomorphisms induced by
the isomorphisms $\cF^0_{\x} \cong (\cF^1)_{\x}^{\vee}\otimes\Omega^1_{\BP^1}$,
$\cF^1_{\x}\cong(\cF^0)_{\x}^{\vee}\otimes\Omega^1_{\BP^1}$
and Serre-duality.
Thus $\xi$ becomes an isomorphism by five lemma.
\end{proof}

\subsection{Smoothness of $M^{\balpha}_n(\bt, \blambda, L)$ and its dimension}

In this subsection, we prove that  the morphism 
$\pi_n: M^{\balpha}_n(L) \lra T_n \times \Lambda_n$ is smooth of equidimension $2n-6$, 
which is stated in Theorem \ref{thm:fund}, (2).

\begin{Proposition}\label{prop:smoothness}
\begin{enumerate}
\item The morphism $\pi_n:  M_{n}^{\balpha}(L)  \lra T_n \times \Lambda_n $ is smooth. 
\item 
For any closed point $\x \in M^{\balpha}_n(\bt, \blambda, L)$, we have 
\begin{equation}\label{eq:dimension}
\dim_{\C} {\mathbf H}^1( \cF^{\bullet}_{\x}) = 2n -6. 
\end{equation}
In particular, the moduli space   $M^{\balpha}_n(\bt, \blambda, L)$ is 
smooth of equidimension $2n -6$.  
\end{enumerate}
\end{Proposition}

\begin{proof}
(1): By a standard argument as in [Lemma 4,  \cite{A}], 
the smoothness of $\pi_n$ at $\x$ follows from Lemma 
\ref{lem:h2}. 
(2): First, by (\ref{eq:duality}), we have 
$\cF^1 \simeq (\cF^0)^{\vee} \otimes \Omega_{\BP^1}^{1}$, 
and 
hence Serre duality implies that 
$\chi(\BP^1, \cF^1_{\x}) = 
\chi( (\cF^0_{\x})^{\vee} \otimes \Omega_{\BP^1}^1) 
= - \chi (\BP^1, \cF^0)$. 
Together with  the exact sequence (\ref{eq:spectral}), we obtain   
\begin{equation}\label{eq:chi}
\dim \mathbf{H}^1(\cF^{\bullet}_{\x}) 
=-  \chi(\BP^1, \cF^0_{\x}) +  \chi(\BP^1, \cF^1_{\x})
 = -2 \chi(\BP^1,  \cF^0_{\x}) .
\end{equation}
Setting 
${\mathcal End}^{0} (\E_{\x})
 = \{ s \in {\mathcal End}(\tilde{E}_{\x}) \ | \ \Tr(s) = 0 \}$, 
by definition of $\cF^0_{\x}$ (\ref{eq:end-p0}), 
we obtain the following exact sequence
%%%%%%%%%%%%%%%%%%
\begin{equation}\label{eq:end-exact}
0 \lra \cF^0_{\x} \lra {\mathcal End}^{0} (\tilde{E}_{\x})  \lra \otimes_{i=1}^n \C((t_i, \x)) \lra 0.
\end{equation}
%%%%%%%%%%%%%%%%%%
Since ${\mathcal End}^{0} (\E_{\x}) $ is a self-dual  locally free sheaf of rank $3$ on 
$\BP^1$, Riemann-Roch theorem implies that  $\chi(\BP^1,  {\mathcal End}^{0} (\E_{\x})) =   
3 + \deg {\mathcal End}^{0} (\E_{\x})  = 3 $.  Then the exact sequence (\ref{eq:end-exact}) 
together with  (\ref{eq:chi})
shows that 
$$
\chi (\BP^1, \cF^0_{\x}) = \chi({\mathcal End}^0(\tilde{E}_{\x}))  - n  = 3-n, 
$$
which implies the assertion (\ref{eq:dimension}) . 
\end{proof}

\begin{Remark} One can also show that 
\begin{equation} \label{eq:h2+}
 {\bf H}^2 (\BP^1, \cF_{\x}^{\bullet   +}) =  \{ 0 \},  
 \end{equation}
 which implies that the morphism $M_n^{\balpha} \lra T_n $ is smooth.  
\end{Remark}

\subsection{Tangent space to $\cR(\cP_{n, \bt})_{\ba}$}

Let $ (E, \nabla, \varphi, \{ l_i \})$  be a stable parabolic connection on $\BP^1$  
corresponding to a point $\x$ in  
$M^{\balpha}_n( \bt, \blambda, L)$.  Let us consider the inclusion 
$j: \BP^1 \setminus D(\bt) \hookrightarrow \BP^1$ and define 
\begin{equation}\label{eq:loc-costant}
{\bf E} = \ker \left[\nabla^{an} : E
\lra (E \otimes \Omega^1_{\BP^1})\right]_{\BP^1 \setminus D(\bt)}.
\end{equation}
Then ${\bf E}$ becomes a locally constant sheaf on $\BP^1 \setminus D(\bt)$. 
The correspondence $(E, \nabla, \varphi, \{ l_i \}) \mapsto {\bf E}$ induces an analytic 
morphism 
\begin{equation}\label{eq:r-h}
{\bf RH}_{\bt, \blambda} :M^{\balpha}_n( \bt, \blambda, L) \lra \cR(\cP_{n, \bt})_{\ba}
\end{equation} 
which is  called the {\em Riemann-Hilbert correspondence}.  
(Here we set $\ba=(a_i), a_i = 2 \cos 2\pi \lambda_i$ ).  For the precise definition, see    
Definition \ref{def:def-rh} in \S \ref{sec:R-H}. 

The morphism  ${\bf RH}_{\bt, \blambda}$   
 will be studied in detail in the next section.   

Define another locally constant sheaf  on $\BP^1 \setminus D(\bt)$ by 
\begin{equation}\label{eq:trace-free}
{\bf V} := \{ s \in Hom({\bf E}, {\bf E}) |  \Tr(s) = 0 \}.
\end{equation}
Note that for each point $ u \in \BP^1 \setminus D(\bt)$ the fiber of ${\bf V}_u$ is 
isomorphic to the Lie algebra $sl_2(\C)$.  Therefore 
 ${\bf V}$ admits the natural non-degenerate pairing 
$
q: {\bf V} \otimes {\bf V}  \lra  \C_{\BP^1\setminus D(\bt)} 
$
induced by the Killing form on each fiber ${\bf V}_u, u \in \BP^1 \setminus D(\bt)$.  
Now consider the constructible sheaf 
$
j_*({\bf V}) 
$
and the following exact sequence induced by the 
Leray spectral sequence for the inclusion 
$ j:\BP^1 \setminus D(\bt) \hookrightarrow \BP^1$:
\begin{equation}\label{eq:Leray}
0 \ra H^1(\BP^1, j_* \V) \ra H^1(\BP^1 \setminus D(\bt), \V) \ra 
H^0(\BP^1, R^1j_*({\bf V}) )\ra H^2(\BP^1, j_*(\V)) \ra 
H^2 (\BP^1 \setminus D(\bt), \V).  
\end{equation}
Recall that in \S \ref{sec:monodromy} we have obtained the morphism 
\begin{equation}\label{eq:family-rep}
\phi_n:\cR_n\lra T'_n \times \cA_n 
\end{equation}
such that  $\phi_n^{-1}((\bt, \ba)) = \cR(\cP_{n, \bt})_{\ba}$.  Fixing $\bt \in T'_n$, 
we can also define 
\begin{equation}
\phi_{n, \bt} : \cR(\cP_{n, \bt} )  \lra {\bt} \times \cA_n. 
\end{equation}

\begin{Lemma} \label{lem:coh-ex}
 Let $ (E, \nabla, \varphi, \{ l_i \}) \in M^{\balpha}_n( \bt, \blambda, L)$
 be a stable parabolic connection, and ${\bf E}:= 
\ker \nabla_{|\BP^1-D(\bt)}$ 
the corresponding local system.   Moreover let $\V $ be the trace 
free part of $\End ({\bf E})$. 
 Let us fix  a monodromy representation 
$\rho_{{\bf E}}: \pi_1(\BP^1 \setminus D(\bt), \ast) \lra SL_2(\C)$ 
associated to the local system 
${\bf E}$. Fix  canonical generators $\gamma_i, 1 \leq i \leq n$ of 
$ \pi_1(\BP^1 \setminus D(\bt), \ast) $ and 
 set $M_i = \rho_{{\bf E}}(\gamma_i) \in SL_2(\C) $ for $ 1 \leq i \leq n$.  
Consider the following conditions.  
\begin{equation}\label{eq:smooth-cond1}
\mbox{The representation $\rho_{{\bf E}}$ is irreducible.} 
\end{equation}
\begin{equation}\label{eq:smooth-cond2}
\mbox{ For each $i, 1 \leq i \leq n$, the local monodromy matrix $M_i$ around 
$t_i$ is not equal  to $\pm I_2$.}  
\end{equation}

\begin{enumerate}
\item Under the condition $(\ref{eq:smooth-cond1})$, we have 
\begin{equation}\label{eq:mon-h2}
H^2(\BP^1 \setminus D(\bt), \V) = \{  0 \}. 
\end{equation}

\item Under the conditions $(\ref{eq:smooth-cond1})$  and
$(\ref{eq:smooth-cond2})$,  we have a sheaf isomorphism 
\begin{equation}\label{eq:local-r1}
R^1j_*({\bf V}) \simeq  \oplus_{i=1}^n 
\C(t_i), 
\end{equation}
and the exact sequence of cohomology groups.
\begin{equation}\label{eq:fundrep-ex}
0 \lra
  H^1(\BP^1, j_* \V) \lra 
   H^1(\BP^1 \setminus D(\bt), \V) \lra  
H^0(\BP^1, R^1j_*(\V))   \lra 0.  
\end{equation}
\end{enumerate}

\end{Lemma}

\begin{proof}
Since we have a canonical non-degenerate pairing 
$$
j_*( \V ) \otimes j_*(\V)  \stackrel{\mbox{\footnotesize Killing}}{\lra} \C_{\BP^1},  
$$
we have  a self-duality $(j_*\V)^{\vee} \simeq j_* \V$ and hence a duality isomorphism  
$$
H^2 (\BP^1, j_*(\V))  \simeq H^0(\BP^1, j_*\V)^{\vee} \simeq H^0(\BP^1 \setminus 
D(\bt), \V)^{\vee}.  
$$
Since by \eqref{eq:smooth-cond1} the monodromy representation $\rho_{{\bf E}}$ 
is irreducible, $H^0(\BP^1 \setminus D(\bt), 
\End({\bf E})) \simeq \C \cdot Id_{\bf E}$ by Schur's lemma  and hence its trace free part 
$H^0(\BP^1 \setminus D(\bt), \V) $ is $\{ 0 \} $, thus
\begin{equation}\label{eq:vanishing-h2}
H^2(\BP^1, j_*(\V)) = \{ 0 \}. 
\end{equation}
Moreover $H^1(\BP^1, R^1j_{*} \V) = \{ 0 \}$, 
for the sheaf $R^1j_* \V $
 is supported only on $D(\bt) = t_1 + \cdots + t_n$.  Then the assertion  (\ref{eq:mon-h2})  
now easily follows from the Leray spectral sequence for $j:\BP^1 \setminus D(\bt) \hookrightarrow 
\BP^1$.

{F}rom    (\ref{eq:Leray}),   we obtain the exact sequence 
(\ref{eq:fundrep-ex})  because of (\ref{eq:vanishing-h2}).  

For the assertion (\ref{eq:local-r1}), we first remark that the 
sheaf  $R^1 j_* \V$ is   supported on $D(\bt) = t_1 + \cdots +t_n$.  We will determine  the 
 stalk $R^1 j_* \V_{t_i} $ at each $t_i$.  Let us take a small neighborhood  $U_i $ of $t_i$ and 
 $ u_i \in U_i - \{ t_i \}$.  Then one can identify the fiber $\V_{u_i}$ with  the symmetric tensor
$ Sym^2({\bf E}_{u_i})$.  Consider the $\V_{u_i} \simeq Sym^2 ({\bf E}_{u_i}) $ as 
the vector space with the action of $M_i$.  Then define the invariant part as 
$$
\V_{u_i} ^{< M_i>} = Sym^2({\bf E}_{u_i})^{<M_i>} : = \ker(Sym^2( M_i) - Id : Sym^2({\bf E}_{u_i})\ra Sym^2({\bf E}_{u_i})). 
$$
Then it is easy to see that 
$$
R^1 j_* \V_{t_i} \simeq \left(\V_{u_i}^{<M_i>} \right)^{\vee}.  
$$ 
Choose a suitable basis of ${\bf E}_{u_i}$ and  write $M_i$ 
 as  $M_i = \left(\begin{array}{cc} a & b \\ c & d \end{array} \right) $ with $ad -bc = 1$.  
 Then the action of $M_i $ on $Sym^2({\bf E}_{u_i})$ has the following matrix representation. 
 \begin{equation}
 Sym^2 (M_i) := \left(\begin{array}{ccc}
 a^2 & ab &  b^2 \\
 2ac & ad + bc & 2bd \\
 c^2 & cd & d^2 \\
 \end{array}
 \right) 
 \end{equation}
 Then it is easy to check that the eigenvalues of $Sym^2(M_i) $ are given by the roots of 
 $$
 (x - 1) (x^2 -((a+d)^2 - 2) x +1) = 0. 
 $$
If neither of the roots of $x^2 - ((a+d)^2 -2) x + 1 = 0$ is $ 1 $,
then $\dim \ker (Sym^2(M_i) - Id) = 1$.  
If one of the roots of    $x^2 - ((a+d)^2 -2) x + 1 = 0$ is one, 
then we have $(a+d)^2 =4$, which implies that $a+d = \pm 2$.  
For those cases,  the eigenvalues of $M_i$ are $ 1 $ or $-1$ respectively. 
We may assume that $M_i  \not= \pm I_2$. Then we can assume that
$M_i = \left( \begin{array}{cc} 1 &  b \\ 0 & 1 \end{array} \right) $ or
$M_i =  \left( \begin{array}{cc} 
-1 &  b \\ 0 & -1 \end{array} \right) $ with $b \not=0$. 
For these cases, we can write 
$$
Sym^2 (M_i) = \left(\begin{array}{ccc}
 1 & b &  b^2 \\
 0 & 1 & 2b \\
 0 & 0 & 1 \\
 \end{array}
 \right )   \quad \mbox{or}  \quad  \left(\begin{array}{ccc}
 1 & -b &  b^2 \\
 0 & 1  & -2b \\
 0 & 0  & 1  \\
 \end{array}
 \right)  .  
$$
Now it is easy to check that $\dim \ker ( Sym^2 (M_i) - Id ) = 1$.  
\end{proof}

\begin{Lemma}\label{lem:resol}
Let us fix $\bt \in T'_n$. 
The notation being as in Lemma \ref{lem:coh-ex}, let us take a point  
$\y:=[{\bf E}] \in \cR(\cP_{n, \bt})_{\ba} \subset \cR(\cP_{n, \bt})$.  
\begin{enumerate}

\item Assume that the condition (\ref{eq:smooth-cond1}) holds for ${\bf E}$. 
Then the total space  $\cR(\cP_{n, \bt})$ is smooth at $\y = [{\bf E}]$ and we have the 
isomorphism 
$$
\Theta_{\cR(\cP_{n, \bt}), \y} \simeq H^1(\BP^1 \setminus D(\bt),  \V). 
$$

\item Assume that the conditions  (\ref{eq:smooth-cond1}) and (\ref{eq:smooth-cond2}) hold for 
${\bf E}$.
Then, the map $\phi_{n, \bt} : \cR(\cP_{n, \bt})  \lra {\bt} \times \cA_n$ is also 
smooth at $\y = [{\bf E}]$. Hence the fiber $\phi_{n, \bt}^{-1}(\ba) =   
\cR(\cP_{n, \bt})_{\ba} $ is smooth at $\y$ where 
$\ba = \phi_{n, \bt}(\y) $.  
Moreover we have the following linear isomorphisms: 
$$ 
(\Theta_{\cR(\cP_{n, \bt})_{\ba}})_{\y} \simeq H^1(\BP^1, j_* \V)
$$
$$
(\Theta_{\cR(\cP_{n, \bt})})_{\y} \simeq H^1(\BP^1 \setminus D(\bt), \V)
$$
$$
(\Theta_{{\bt} \times \cA_n})_{\phi_{n, \bt}(\y)} \simeq H^0(\BP^1, R^1j_*(\V))
$$
   Under the isomorphisms above, we have the following identification of the natural 
exact sequences of the tangent spaces with the sequence (\ref{eq:fundrep-ex}) 
\begin{equation}
\begin{array}{ccccccc}
 0 \lra & 
(\Theta_{\cR(\cP(\cP_{n, \bt}))_{\ba}})_{\y} &  
\lra &
(\Theta_{\cR(\cP_{n, \bt})})_{\y} &
\stackrel{d\phi_{n, \bt *}}{\lra} & 
(\Theta_{{\bt} \times \cA_n})_{\phi_{n, \bt}(\y)}
 &
 \lra 0. \\
 & & & & & & \\
&
|| 
&   
& 
||  
& 
& 
||
&  \\
& & & & & &  \\ 
 0 \lra & 
  H^1(\BP^1, j_* \V) & 
   \lra & 
   H^1(\BP^1 \setminus D(\bt), \V) &
\lra & 
H^0(\BP^1, R^1j_*(\V)) &  \lra 0.  
\end{array}
\end{equation}
\end{enumerate}
\end{Lemma}

\begin{proof}
1. Since $\E$ is irreducible, 
 it is easy to see that the Zariski tangent space $\Theta_{\cR(\cP_{n, \bt}), \y}$  
of $\cR(\cP_{n, \bt}) $ at $\y = [{\bf E}]$ 
is given by $ H^1(\BP^1\setminus D(\bt), \V) $ and the obstructions to deformations lie in 
$H^2(\BP \setminus D(\bt), \V)$.  Since we assume that $\rho_{\bf E}$ is irreducible, we have 
$ H^2(\BP \setminus D(\bt), \V) = \{ 0 \} $ (cf. (\ref{eq:mon-h2})) , 
from which the assertion follows.

2. From Lemma \ref{lem:coh-ex},  under the assumptions, we can see that 
the differential  $(\Theta_{\cR(\cP_{n, \bt})})_{\y} 
\stackrel{d\phi_{n, \bt *}}{\lra} 
(\Theta_{{\bt} \times \cA_n})_{\phi_{n, \bt}(\y)}$ can be identified with the 
linear map 
$$
H^1(\BP^1 \setminus D(\bt), \V) 
\lra 
H^0(\BP^1, R^1j_*(\V)) \simeq \C^n, 
$$
which is surjective because of $H^2(\BP^1, j_* \V) = \{ 0 \}$. 
Therefore the map $\phi_{n, \bt} $ is smooth at $\y = [{\bf E}]$ and the fiber 
$\phi_{n, \bt}^{-1}(\ba) =   \cR(\cP_{n, \bt})_{\ba} $ is smooth at $\y$. Other assertions 
now easily follow from the exact sequence  (\ref{eq:fundrep-ex}).  
\end{proof}

\begin{Lemma}\label{lem:inf-rh}
Under the conditions  (\ref{eq:smooth-cond1}) and (\ref{eq:smooth-cond2}) for ${\bf E}$, 
we have an isomorphism of locally constant sheaves
\begin{equation}\label{eq:quasi-isom}
j_* \V \simeq \ker \nabla_1 \simeq \left[\nabla_1:\cF^0 \lra \cF^1\right],  
\end{equation}
which induces  a canonical isomorphism 
\begin{equation} \label{eq:rh-tangent}
H^1(\BP^1, j_* \V) \stackrel{\simeq}{\lra}  {\bf H}^1( \left[\nabla_1:\cF^0 \lra \cF^1 \right]). 
\end{equation}
Moreover we have the canonical non-degenerate pairing 
\begin{equation} \label{eq: rh-sym}
H^1(\BP^1, j_* \V)  \otimes H^1(\BP^1, j_* \V)  \lra H^2 (\BP^1, \C_{\BP^1}) \simeq \C, 
\end{equation}
which induces the non-degenerate pairing  $\Omega_1(\y)$ at $\y = {\bf E}$ 
\begin{equation}\label{eq:tangent-mon}
\Omega_1(\y): 
(\Theta_{\cR(\cP_{n, \bt})_{\ba}})_{\y} \otimes (\Theta_{\cR(\cP_{n, \bt})_{\ba}})_{\y} \lra (\cO_{\cR(\cP_{n, \bt})_{\ba}})_{ \y}. 
\end{equation}
This pairing can be identified with (\ref{eq:mod-sym}) via the isomorphism (\ref{eq:rh-tangent}). 
\end{Lemma}

\begin{proof}
The assertion (\ref{eq:quasi-isom}) is trivial at  the point $u \in \BP^1 \setminus D(\bt)$.
At each point $t_i$ $ i = 1, \ldots, n$,  we will describe the connection $\nabla$ and $\nabla_1$ 
locally around $t_i$.  Let us set $n=n_i = \res_{t_i}(\nabla_L) \in \Z $. 
We separate the proof into two cases.

i) Let $\lambda, n - \lambda$ be the eigenvalues of $\res_{t_i} (\nabla)$. 
First assume that  $ 2 \lambda \not\in \Z $. Then  $\lambda \not= n- \lambda$.   
By a standard reduction theory of 
connection near a regular singularity, we can 
choose a suitable local coordinate $ z $ around $ t= t_i$ and 
write down the connection matrix of  $\nabla$ by 
$$
\nabla = \frac{dz}{z - t} \left( \begin{array}{cc} 
 \lambda  & 0 \\ 0 & n-\lambda 
\end{array} \right). 
$$
Then for a local section $ s = 
 \begin{pmatrix}
  s_1 & s_2 \\
  s_3 & -s_1
 \end{pmatrix} \in \End(E) $, 
the connection $\nabla_1 s = \nabla_E s - s \nabla_E $ is given by 
$$
 \begin{pmatrix}
  s_1 & s_2 \\
  s_3 & -s_1
 \end{pmatrix}
 \mapsto
 \begin{pmatrix}
  ds_1 & ds_2+(2\lambda-n) s_2(z-t)^{-1}dz \\
  ds_3+(n-2\lambda) s_3(z-t)^{-1}dz & -ds_1
 \end{pmatrix}.
$$
Solving  $\nabla_1 = 0$  locally near $z = t$, 
we obtain the local solutions for $z \not= t$ as follows. 
\begin{gather*}
  c_1
 \begin{pmatrix}
  1 & 0 \\
  0 & -1
 \end{pmatrix}
  + c_2
 \begin{pmatrix}
  0 & (z-t)^{n-2\lambda} \\
  0 & 0
 \end{pmatrix} 
  + c_3
 \begin{pmatrix}
  0 & 0 \\
  (z-t)^{2\lambda-n} & 0
 \end{pmatrix}. 
\end{gather*}
Here,  $c_1,c_2,c_3 \in \C$ are constants.  These solutions have to be 
 single-valued well-defined  section  around $z = t$, hence $\ker \nabla_1$ 
 is generated by $\begin{pmatrix}
  1 & 0 \\
  0 & -1
 \end{pmatrix}.
 $ 
  (Note that this local section lies in $ \cF^0 $. ) 
 On the other hand, the stalk $(j_*\V)_{t}$ is the space of monodromy  invariant 
 trace-free endomorphisms of ${\bf E}_{u}$,  
 which is also generated by  
 $\begin{pmatrix}
  1 & 0 \\
  0 & -1
 \end{pmatrix}.  
 $
 Hence we have an isomorphism $(j_*\V)_t \simeq \ker(\nabla_1)_t$.  

ii) Again,  let $\lambda, n-\lambda$ be the eigenvalues of 
$ \res_{t_i} (\nabla) $ and assume that $ 2 \lambda \in \Z $.  
Since we assume that the local monodromy $M_i$ is not $\pm I_2$, 
by a reduction theory of a connection near 
a regular singularity, we can  choose 
a suitable local coordinate $z$ around $t=t_i$ and 
 write down the connection matrix of  $\nabla$ by 
$$
\nabla = \frac{dz}{z-t}
\begin{pmatrix}
  m_1 & (z-t)^{m_2-m_1} \\
  0                 & m_2 
 \end{pmatrix}, 
$$
where $ 2m_1, 2m_2 \in \mathbf{Z}$,  $ m_2-m_1 \in \mathbf {Z}  $ and $ m_1\leq m_2 $. 
For local section $ \begin{pmatrix}
  s_1 & s_2 \\
  s_3 & -s_1
 \end{pmatrix} \in \End(E) $ , the connection $\nabla_1 s $ can be given by 
$$
\begin{pmatrix}
 s_1 & s_2 \\
 s_3 & -s_1
\end{pmatrix}
\mapsto
\begin{pmatrix}
 ds_1+s_3(z-t)^{m_2-m_1-1}dz &
 ds_2-2s_1(z-t)^{m_2-m_1-1}dz+s_2(m_1-m_2)(z-t)^{-1}dz \\
 ds_3+s_3(m_2-m_1)(z-t)^{-1}dz &
 -ds_1-s_3(z-t)^{m_2-m_1-1}dz
 \end{pmatrix}.
$$ 

Solving $\nabla_1 s = 0$ locally for $z \not= t$, we obtain the solutions
\begin{equation*}
\begin{array}{ll}
s =  c_1
\begin{pmatrix}
 0 & (z-t)^{m_2-m_1} \\
 0 & 0
\end{pmatrix}
 & + c_2
\begin{pmatrix}
 1 & 2(z-t)^{m_2-m_1}\log(z-t) \\
 0 & -1
\end{pmatrix} \\
\quad  & + c_3
\begin{pmatrix}
 \log(z-t) & (z-t)^{m_2-m_1}(\log(z-t))2 \\
 -(z-t)^{m_1-m_2}  & -\log(z-t)
\end{pmatrix}
\end{array}
\end{equation*} 
where $c_1, c_2, c_3 \in \C$ are the constants.  Then we can see that all 
 single valued solutions for $\ker \nabla_1$ are  
$$
c_1
 \begin{pmatrix}
  0 & (z-t)^{m_2-m_1} \\
  0 & 0
 \end{pmatrix}, 
$$ 
which are also sections of  $(j_* \V)_t$. Therefore we have an isomorphism 
$(\ker \nabla_1)_t  \simeq  (j_* \V)_t$. Hence we have proved the assertion (\ref{eq:quasi-isom}) 
which shows also (\ref{eq:rh-tangent}).  

It is easy to see that the pairing of sheaves 
$ j_* \V \otimes j_* \V  \lra j_* \C_{\BP^1 \setminus D(\bt)}  \simeq \C_{\BP^1}$ is non-degenerate 
at each point of $\BP^1$.  Therefore, the pairing  (\ref{eq: rh-sym})  is also non-degenerate.

\end{proof}

Summarizing all results in this section, we have the following 

\begin{Proposition}\label{prop:symplectic-rep-1}
  Let $\phi_n: \cR_n \lra T'_n \times \cA_n$ be a family of 
moduli spaces of representations of the fundamental group
$\pi_1(\BP^1 \setminus D(\bt), *)$ as in 
(\ref{eq:family-rep}).  Let  $\cR_n^{\sharp} $ be the subset of $\cR_n$ whose closed points 
satisfy the conditions (\ref{eq:smooth-cond1}) and (\ref{eq:smooth-cond2}).  

Then $\cR_n^{\sharp}$ is a non-singular variety and the restricted morphism 
\begin{equation}\label{eq:rest-morphism}
\phi_n: \cR_n^{ \sharp} \lra T'_n \times \cA_n
\end{equation} 
 is smooth, so that all fibers  $\cR_{n, (\bt, \ba)}^{\sharp} = \cR(\cP_{n, \bt})^{\sharp}_{\ba}$ are non-singular 
 varieties.  
On $\cR_n^{\sharp}$, there exists a relative symplectic form 
\begin{equation}\label{eq:sympl-rep}
\Omega_1 \in \Gamma(\cR_n^{\sharp}, \Omega^2_{\cR_n^{\sharp}/T'_n \times \cA_n })
\end{equation}
induced by (\ref{eq:tangent-mon}). 
\end{Proposition}

\begin{Remark}
{\rm 
\begin{enumerate}
\item 
Since  $p_2 \circ \phi_n:\cR_n^{\sharp} \lra T'_n  \times \cA_n \lra T_n'$ is locally trivial, one can lift 
$\Omega_1 \in \Gamma(\cR_n^{\sharp}, \Omega^2_{\cR_n^{\sharp}/T'_n \times \cA_n }) $ 
 to a relative regular $2$-form  
\begin{equation}\label{eq:relative-two-form}
\tilde{\Omega}_1\in \Gamma(\cR_n^{\sharp}, \Omega^2_{\cR_n^{' \sharp}/\cA_n}).  
\end{equation}
In \S \ref{sec:R-H}, we can define the Riemann-Hilbert 
correspondence ${\bf RH}_n:M_n^{\balpha} \lra \cR_n$ 
which is a  surjective holomorphic map.  Set $(M_n^{\balpha})^{\sharp} = {\bf RH}_n^{-1}(\cR_n^{\sharp})$.  
 From Lemma \ref{lem:inf-rh}, one can see that ${\bf RH}^{*}_{n|\cR_n^{\sharp}} (\Omega_{1}) $ coincides 
with the two form $ \Omega_{|(M^{\balpha}_{n})^{\sharp}} \in 
\Gamma((M^{\balpha}_{n})^{\sharp}, \Omega^2_{(M^{\balpha}_n)^{\sharp}/T'_n \times \Lambda_n})$ defined in 
(\ref{eq:mod-sym}).  Pulling back $\tilde{\Omega}_1 $ via ${\bf RH}_{n|\cR_n^{\sharp}}$, we obtain 
\begin{equation}\label{eq:global-twoform}
\tilde{ \Omega} \in \Gamma((M^{\balpha}_{n})^{\sharp}, \Omega^2_{(M^{\balpha}_n)^{\sharp}/ \Lambda_n})
\end{equation}
which is a lift of $\Omega_{|(M^{\balpha}_{n})^{\sharp}}$
via  the canonical morphism
 $\Gamma((M^{\balpha}_{n})^{\sharp}, \Omega^2_{(M^{\balpha}_n)^{\sharp}/ \Lambda_n})
\lra 
\Gamma((M^{\balpha}_{n})^{\sharp}, \Omega^2_{(M^{\balpha}_n)^{\sharp}/T'_n \times  \Lambda_n}).
$
(Note that a lift $\tilde{ \Omega}$ can be induced  from the splitting homomorphism 
(\ref{eq:descend-split}) and the splitting homomorphism can be defined algebraically). 
Since the codimension of $M_n^{\balpha} \setminus (M_n^{\balpha})^{\sharp}$ 
in $M_n^{\balpha}$ is at least two, the two form $\tilde{ \Omega}$ can be extended to 
a regular relative two form on $M_n^{\balpha}$ which will be denoted also by $\tilde{\Omega}$. 
This extended two form 
$\tilde{\Omega}  \in \Gamma(M^{\balpha}_{n}, \Omega^2_{M^{\balpha}_n/ \Lambda_n}) $ 
is a lift of $\Omega \in 
\Gamma(M^{\balpha}_{n}, \Omega^2_{M^{\balpha}_n/T'_n \times  \Lambda_n})$ 
in (\ref{eq:mod-sym}) on the 
whole total space $M_n^{\balpha}$.

\item 
The closedness of $\tilde{\Omega}$, $(d_{M^{\balpha}_n/\Lambda_n}(\tilde{\Omega}) = 0)$,  can be 
proved as follows. It is easy to see that the two form  $\tilde{\Omega}$ here 
coincides with the symplectic two form introduced in \cite{Iwa91} and \cite{Iwa92}
on a Zariski dense open subset   $(M_n^{\balpha})'$ of $M_n^{\alpha} $.  
As proved in  \cite{Iwa91}, \cite{Iwa92},  there exists a 
suitable affine open covering  $\{ U_{i} \}_{i} $ of $(M_n^{\balpha})'$ with local 
coordinates (for $U_i$ ) 
$$
(q^i_1, \ldots, q^i_r, p^i_1, \ldots, p^i_r, t_1, \ldots, t_n, \lambda_1, \ldots, \lambda_n)
$$
such that $\tilde{\Omega}_{|U_i}$ can be written as 
\begin{equation}\label{eq:canonical-form}
\tilde{\Omega}_{|U_i} =\sum_{k=1}^r  dq^i_k \wedge dp^i_k - \sum_{l=1}^n  dt_l 
\wedge dH_l^{i} ({\bf p},{\bf q}, \bt, \blambda). 
\end{equation} 
where $r = n-3 $ (= the half of the relative dimension of $\pi_n$)  
and $H^{i}_l({\bf p},{\bf q}, \bt, \blambda)$ are 
regular algebraic functions on $U_i$.  
The closedness $d_{M^{\balpha}_n/\Lambda_n}(\tilde{\Omega}) = 0$ on  $U_i$ easily follows 
from the expression (\ref{eq:canonical-form}),  hence by analytic continuation we see that 
$d_{M^{\balpha}_n/\Lambda_n}(\tilde{\Omega}) = 0$ on the total space $M_n^{\balpha}$.   

\item The regular functions 
$H^i_l({\bf p},{\bf q}, \bt, \blambda)$ on $U_i$ in (\ref{eq:canonical-form})
are called {\em Hamiltonians} for Painlev\'e or Garnier systems  
with respect to the time variable $t_l$.  
 Actually on an affine open set $U_i$ 
 one  can obtain the Hamiltonian systems $($Cf. \cite{Iwa91}, \cite{Iwa92}$)$. 
\begin{equation}\displaystyle{
\frac{\partial q^i_k}{\partial t_l}   =   \frac{\partial H^i_l }{ \partial p^i_k},  \quad
 \frac{\partial p^i_k}{\partial t_l}   =  - \frac{\partial H^i_l }{ \partial q^i_k} }   
 \quad ( 1 \leq k \leq n-3, 1\leq l \leq n). 
\end{equation}
Although these Hamiltonian systems are defined on a Zariski open subset 
$(M_n^{\balpha})'$ of $M_n^{\balpha}$,  these Hamitonian systems can be 
extended to Hamitonian systems on  the total space  $M_n^{\balpha}$. 
This is because 
 global vector fields on $M_n^{\balpha}$ induced from the isomonodromic flows 
coincide with these Hamitonian systems on the Zariski open set $(M_n^{\balpha})'$ 
and the  global vector fields on $M_n^{\balpha}$ also  
preserves the symplectic form $\tilde{\Omega}$.  
\end{enumerate}
}
\end{Remark}

\vspace{1cm}

\section{The Riemann-Hilbert correspondence}
\label{sec:R-H}

In this section, we also work over $T'_n$ (cf. \eqref{eq:finite-etale}).  
Fix $(\bt, \blambda) \in T'_n \times \Lambda_n$ and set 
$D(\bt)= t_1+ \cdots + t_n  \subset \BP^1$, $a_i = 2 \cos 2 \pi \lambda_i$ and 
$\ba = (a_1, \ldots, a_n) \in \cA_n$.  
 Moreover fix  a determinant line bundle 
$L =(L, \nabla_L)$ on $\BP^1$ such that $\res_{t_i}(\nabla_L) \in \Z$ for every 
$ 1\leq i \leq n$. 
We have defined two  moduli spaces 
$M^{\balpha}_n ( \bt, \blambda, L) $ in \eqref{eq:modulispace} and  $\cR(\cP_{n, \bt})_{\ba}$ 
in \eqref{eq:mod-rep}.  In this section, we define 
 the {\em Riemann-Hilbert correspondence}
$
\RH_{\bt, \blambda}: M^{\balpha}( \bt, \blambda, L) \lra \cR(\cP_{n, \bt})_{\ba},  
$ and  show our main results for the Riemann-Hilbert correspondence (Theorem \ref{thm:RH}).

\subsection{Definition of ${\bf RH}_{\bt, \blambda}$}
\label{subsec:def-rh}
As in  \eqref{eq:r-h}, take  
$E = (E, \nabla, \varphi, l) \in   M^{\balpha}_n ( \bt, \blambda, L)$ and 
define the local system on $\BP^1 \setminus D(\bt)$ as 
$\E = \ker\left(\nabla_{|\BP^1 \setminus  D(\bt)}\right)^{an}$. 
(Here we denote by $\left(\nabla_{|\BP^1 \setminus  D(\bt)}\right)^{an}$ the 
analytic connection associated to $\left(\nabla_{|\BP^1 \setminus  D(\bt)}\right)$.)  
Choosing a suitable flat basis for the fiber $\E_{*}$ at the base point $* \in \BP^1 \setminus D(\bt)$, 
 one can define a  monodromy representation
$
\rho_{\E}: \pi^1(\BP^1 \setminus D(\bt), *) \lra SL_2 (\C)$. 
The difference of choices of  flat basis can be given 
by the adjoint action of 
$SL_2(\C)$, and hence one has a correspondence
\begin{equation}\label{eq:object-rh}
E = (E, \nabla, \varphi, l) \mapsto [\rho_{\E}].
\end{equation}
Here$ [\rho_{\E}]$ denotes the Jordan equivalence class of $\rho_{\E}$. 

Fix canonical generators $\gamma_i, 1 \leq i \leq n$ of $ \pi^1(\BP^1 \setminus D(\bt), *)$. 
For a monodromy representation $\rho_{\E}$ of  $(E, \nabla, \varphi, l)$, set 
$M_i = \rho_{\E}(\gamma_i)$ as in  \S \ref{sec:monodromy}. 
Since eigenvalues of $\res_{t_i}(\nabla)$ can be given by 
$\lambda_i, \res_{t_i} (\nabla_L)  - \lambda_i$ and 
$\res_{t_i}(\nabla_L) \in \Z$, we see that the eigenvalues of 
$M_i$ are given by 
$\exp( \mp  2 \pi \sqrt{-1} \lambda_i)$. 
Therefore,  we have   local exponents for $\rho_{\E}$ 
\begin{equation}\label{eq:localexp2}
a_i := \Tr[M_i] =  \exp( -2 \pi \sqrt{-1} \lambda_i) + \exp(  2 \pi \sqrt{-1} \lambda_i) = 2 \cos ( 2
 \pi \lambda_i), \end{equation}
which are  invariant under the adjoint action.  
\begin{Definition}{\rm \label{def:def-rh}
Under the relation \eqref{eq:localexp2}, 
the correspondence  \eqref{eq:object-rh} gives 
an analytic morphism
\begin{equation}\label{eq:RHF}
       \RH_{\bt, \blambda}:  M_n^{\balpha}(\bt, \blambda, L) \lra \cR(\cP_{n, \bt})_{\ba},  
\end{equation}
which is called the {\em Riemann-Hilbert correspondence}. 
}
\end{Definition}

\subsection{Fundamental properties  of Riemann-Hilbert correspondence}

Let us assume that $n \geq 4$.  
In \S \ref{sec:monodromy}, \eqref{eq:family-rep-s4},  we have  defined 
the family of moduli spaces of representations of 
fundamental group $\phi_n: \cR_n \lra T'_n \times \cA_n$ and   
we also have a smooth family 
$\pi_n: M^{\balpha}_n(L) \lra T'_n \times \Lambda_n$ whose geometric fibers are 
$M^{\balpha}(\bt, \blambda, L)$ (cf. Theorem \ref{thm:fund}).
{F}rom Definition \ref{def:def-rh} we obtain the following commutative diagram:
\begin{equation}\label{eq:RH}
\begin{CD}
M^{\balpha}_n(L)  @> \RH_n >> \cR_n \\
@V \pi_n VV  @V \phi_n VV \\
T'_n \times \Lambda_n @> id \times \mu_n >> T'_n \times \cA_n. 
\end{CD}
\end{equation} 
Here   $\mu_n:\Lambda_n \lra \cA_n$ is given by 
\begin{equation}\label{eq:cor-local}
\mu_n (\lambda_1, \ldots, \lambda_n) = (a_1, \ldots, a_n)  =( 2 \cos ( 2 \pi \lambda_1), \ldots, 2 \cos ( 2 \pi \lambda_n)).
\end{equation}
Of course, for each $(\bt, \blambda) \in T'_n \times  \Lambda_n$, 
the morphism $\RH_{n| M_n^{\balpha}(\bt, \blambda, L)}$ is equal to 
$\RH_{\bt, \blambda}$ in (\ref{eq:RHF}).

\begin{Theorem} \label{thm:RH}  Under the notation and the assumption  
as above, we have 
the following assertions. 
\begin{enumerate}
\item For all $(\bt, \blambda) \in T'_n \times \Lambda_n$, 
the Riemann-Hilbert correspondence $
\RH_{\bt,\blambda}: M_n^{\balpha}(\bt, \blambda, L) \lra \cR(\cP_{n, \bt})_{\ba}$ in \eqref{eq:RHF} 
is a \underline{\em bimeromorphic  proper surjective}  morphism.   
\item For any $(\bt, \blambda)$, let $\cR(\cP_{n, \bt})^{\sharp}_{\ba}$ be the Zariski open
subset of $\cR(\cP_{n, \bt})_{\ba}$ whose closed points satisfy the conditions 
(\ref{eq:smooth-cond1}) and (\ref{eq:smooth-cond2}) in \S \ref{sec:symp}, 
and set $M_n^{\balpha}(\bt, \blambda, L)^{\sharp} = \RH_{\bt, \blambda}^{-1}( \cR(\cP_{n, \bt})^{\sharp}_{\ba})$. Then the Riemann-Hilbert correspondence gives an analytic 
isomorphism 
\begin{equation}\label{eq:rh-isom}
\RH_{\bt,\blambda, | M_n^{\balpha}(\bt, \blambda, L)^{\sharp}} : M_n^{\balpha}(\bt, \blambda, L)^{\sharp}  \stackrel{\simeq}{\lra}   
\cR(\cP_{n, \bt})^{\sharp}_{\ba}.  
\end{equation}
$($Note that if $\blambda$ is generic (cf. Definition \ref{def:exponents}, (\ref{eq:special1}),  
(\ref{eq:special2})),    $\cR(\cP_{n, \bt})^{\sharp}_{\ba} = \cR(\cP_{n, \bt})_{\ba}$, hence $\RH_{\bt, \blambda}$ gives an analytic isomorphism between 
$M_n^{\balpha}(\bt, \blambda, L)$ and $\cR(\cP_{n, \bt})_{\ba}$.) 
\item Let us set $\cR(\cP_{n, \bt})^{sing}_{\ba} = \cR(\cP_{n, \bt})_{\ba} \setminus 
\cR(\cP_{n, \bt})^{\sharp}_{\ba}$.  Then the codimension of  $\cR(\cP_{n, \bt})^{sing}_{\ba}$ 
in $\cR(\cP_{n, \bt})_{\ba}$ is at least $2$.

\item The symplectic structures  $\Omega$ on 
$M_n^{\balpha}(\bt, \blambda, L)$ and 
$ \Omega_1$ on $\cR(\cP_{n, \bt})^{\sharp}_{\ba}$ can be identified with each other 
via $\RH_{\bt, \blambda}$, that is, 
\begin{equation} \label{eq:pull-back}
\Omega = \left( \RH_{\bt, \blambda, | \cR(\cP_{n, \bt})^{\sharp}_{\ba}}\right)^*(\Omega_1) \quad 
\mbox{on}  \ \  M_n^{\balpha}(\bt, \blambda, L)^{\sharp}. 
\end{equation} 
 \end{enumerate}
\end{Theorem}

\begin{Remark}{\rm 
\begin{enumerate}
\item The moduli spaces 
$M_n^{\balpha}(\bt, \blambda, L)$  and 
$\cR(\cP_{n, \bt})_{\ba}$ are irreducible.  (See \S \ref{sec:rep-irred} and \S  
\ref{sec:irr-stable}). 
\item The  statement (\ref{eq:pull-back}) is originally shown by Iwasaki 
 in  \cite{Iwa91}, \cite{Iwa92}.  
\end{enumerate}
}
\end{Remark}

Let us denote  ${\bf RH}_{\bt, \blambda}$ in (\ref{eq:RHF}) simply by 
$\RH$. We first show the following

\begin{Lemma}\label{lem:surjective}
 Assume that $n\geq 4$ and that
 $\alpha_i$ ($i=1,\ldots,n$) are so general
 that all the semistable parabolic connections are stable.
 Then the morphism
 $\RH: M^{\balpha}_n(\bt, \blambda, L) \lra \cR(\cP_{n, \bt})_{\ba}$
 is a bimeromorphic surjective morphism.
\end{Lemma}

\begin{proof}
Let $\cR^{\rm irr}(\cP_{n,\bt})_{\ba}$ be the open subscheme of
$\cR(\cP_{n,\bt})_{\ba}$ whose points correspond to the irreducible
representations.
First we will show that
$\cR^{\rm irr}(\cP_{n,\bt})_{\ba}$ is contained in the image of $\RH$.

Let $M^{\rm irr}_n(\bt,\blambda,L)$
be the open subscheme of $M^{\balpha}_n(\bt, \blambda,L)$
consisting of the points corresponding to the irreducible connections.
Note that if $(E,\nabla_E,\varphi,\{l_i\})$ is a parabolic connection
such that $(E,\nabla_E)$ is an irreducible connection,
we have $(E,\nabla_E,\varphi,\{l_i\})\in M^{\balpha}_n(\bt,\blambda,L)$.
We consider the isomorphism of the moduli spaces
\[
 Elm^{-}_{t_i}: M^{\rm irr}_n(\bt,\blambda,L)
 \stackrel{\sim}\longrightarrow
 M^{\rm irr}_n(\bt,\blambda',L(-t_i)); \quad
 (E,\nabla_E,\varphi,\{l_i\})\mapsto
 (E',\nabla_{E'},\varphi',\{l'_i\}),
\]
where $E'=\ker(E\ra E_{t_i}/l_i)$,
$\nabla_{E'}$ is a connection on $E'$ induced by $\nabla_E$,
$l'_i=\ker(E'_{t_i}\ra E_{t_i})$, $l'_j=l_j$ for $j\neq i$,
$\lambda'_i=1+\res_{t_i}(\nabla_L)-\lambda_i$,
$\lambda'_j=\lambda_j$ for $j\neq i$ and
$\varphi':\bigwedge^2 E'\stackrel{\sim}\ra L(-t_i)$
is the horizontal isomorphism induced by $\varphi$.
We also consider the isomorphisms of the moduli spaces
\begin{gather*}
 \otimes\cO(t_i):M^{\rm irr}_n(\bt,\blambda,L)\longrightarrow
 M^{\rm irr}_n(\bt,\blambda',L\otimes\cO(2t_i)); \\
 (E,\nabla_E,\varphi,\{l_i\})\mapsto
 (E\otimes\cO(t_i),\nabla_{E\otimes\cO(t_i)},
 \varphi\otimes1,\{l'_i\otimes\cO(t_i)|_{t_i}\}),
\end{gather*}
where $\lambda'_i=\lambda_i-1$ and $\lambda'_j=\lambda_j$ for $j\neq i$ and
we consider for $\lambda_i\neq \res_{t_i}(\nabla_L)-\lambda_i$, an isomorphism
\[
 s_i:M^{\rm irr}_n(\bt,\blambda,L)
 \stackrel{\sim}\longrightarrow
 M^{\rm irr}_n(\bt,\blambda',L); \quad
 (E,\nabla_E,\varphi,\{l_i\})\mapsto
 (E,\nabla_E,\varphi,\{l'_i\}),
\]
where $\lambda'_i=\res_{t_i}(\nabla_L)-\lambda_i$,
$\lambda'_j=\lambda_j$ for $j\neq i$,
$l'_i=\ker(\res_{t_i}(\nabla_E)-\lambda'_i)$
and $l'_j=l_j$ for $j\neq i$.
Note that these isomorphisms all commute with
the Riemann-Hilbert morphism $\RH$.

Now we fix $(\lambda_1,\ldots,\lambda_n)\in\mathbf{C}^n$ and put
$\lambda_i^+:=\lambda_i$, $\lambda_i^-:=\res_{t_i}(\nabla_L)-\lambda_i$
for $i=1,\ldots,n$.
Applying a certain composition of $Elm^-_{t_i}$,
$\otimes\cO_{\BP^1}(t_i)$ and $s_i$
for $i=1,\ldots,n$, we obtain an isomorphism
\[
 \tau:M^{\rm irr}_n(\bt,\blambda,L)\stackrel{\sim}\lra
 M^{\rm irr}_n(\bt,\blambda',L'),
\]
where $\lambda'_i:=\lambda_i+m_i^+$ and
$\res_{t_i}(\nabla_{L'})=\res_{t_i}(\nabla_L)+m_i^++m_i^-$
for some integers $m^+_i,m^-_i\in\mathbf{Z}$ such that
$0\leq \mathrm{Re}(\lambda^{+}_i+m^{+}_i)<1$,
$0\leq \mathrm{Re}(\lambda^{-}_i+m^{-}_i)<1$
for $1\leq i\leq n$.

Let $N^{\rm irr}_n(\bt,\blambda',L')$ be the moduli space of
rank $2$ irreducible connections $(E,\nabla_E)$ with a horizontal
isomorphism $\bigwedge^2 E\stackrel{\sim}\ra L'$ such that
$\det(\res_{t_i}(\nabla_E)-\lambda'_i)=0$ for $i=1,\ldots,n$.
By [\cite{Deligne:70}, Proposition 5.4], we obtain an isomorphism
\begin{equation}\label{eq:rh-irr}
 \mathrm{rh}:N^{\rm irr}_n(\bt,\blambda',L')
 \stackrel{\sim}\lra \cR^{\rm irr}(\cP_{n,\bt})_{\ba}.
\end{equation}
There is a canonical surjective morphism
\begin{equation}\label{forget}
 M^{\rm irr}_n(\bt,\blambda',L')
 \rightarrow N^{\rm irr}_n(\bt,\blambda',L')
\end{equation}
which is obtained by forgetting parabolic structures.
Composing $\tau$, (\ref{forget}) and $\mathrm{rh}$,
we obtain a surjective morphism
\begin{equation}\label{RH}
 \RH: M^{\rm irr}_n(\bt,\blambda,L)
 \longrightarrow \cR^{\rm irr}(\cP_{n,\bt})_{\ba}.
\end{equation}

Note that the morphism (\ref{forget}) is isomorphic except
on the locus where the parabolic structures are not uniquely
determined by $(E,\nabla_E)$, namely,
\[
 M^{\rm app}_n(\bt,\blambda',L')=
 \left\{
 (E,\nabla_E,\varphi,\{l_j\})\in M^{irr}_n(\bt,\blambda',L') \left|
 \text{$\mathrm{Res}_{t_i}(\nabla_E)= O $ or
 $\frac{1}{2}\mathrm{id}_{E_{t_i}}$ for some $i$}
 \right\}\right.
\]
whose image in
$\cR(\cP_{n,\bt})_{\ba}$ is
\[
 \cR^{\rm app}(\cP_{n,\bt})_{\ba}
 =\left.\left\{ \rho\in \cR^{irr}(\cP_{n,\bt})_{\ba}
 \right|
 \text{$\rho(\gamma_i)=\pm \mathrm{id}$ for some $i$}
 \right\}.
\]
Thus the restriction of $\RH$
\[
 M_n^{\balpha}(\bt, \blambda, L)^{\sharp}=
 M^{\rm irr}_n(\bt,\blambda,L)
 \setminus \tau^{-1}(M^{\rm app}_n(\bt,\blambda',L'))
 \stackrel{\mathrm{RH}}\longrightarrow
 \cR^{\rm irr}(\cP_{n,\bt})_{\ba}\setminus
 \cR^{\rm app}(\cP_{n,\bt})_{\ba}=
 \cR(\cP_{n, \bt})^{\sharp}_{\ba}
\]
is an isomorphism.
Since $\dim \cR^{\rm app}(\cP_{n, \bt})_{\ba}
<\dim \cR^{\rm irr}(\cP_{n,\bt})_{\ba}$
for $n\geq 4$,
$\RH$ is a bimeromorphic morphism.

Next we will show that
$M^{\rm red}_n(\bt,\blambda,L)\ra\cR^{\rm red}(\cP_{n,\bt})_{\ba}$
is surjective.
Take any point $[\rho]\in\cR^{\rm red}(\cP_{n,\bt})_{\ba}$.
Then the representation $\rho$ is Jordan equivalent to the
representation $\rho_1\oplus\rho_2$ for some one dimensional
representations $\rho_1,\rho_2$ of
$\pi_1(\BP^1\setminus D(\bt), *)$.
Put
\[
 U_{n,\ba}:=\left\{ (M_1,\ldots,M_{n-1})\in SL_2(\C)^{n-1}  \left|
 \Tr(M_i)=a_i\;(1\leq i\leq n-1), \quad
 \Tr((M_1M_2\cdots M_{n-1})^{-1})=a_{n} \right\}\right..
\]
Then $U_{n,\ba}$ is irreducible by Proposition \ref{prop:irr-2}.
Let $\Phi_n: U_{n,\ba}\ra\cR(\cP_{n,\bt})_{\ba}$
be the quotient map.
Then there exists a point $p_0\in U_{n,\ba}$ such that
$\Phi_n(p_0)=[\rho]$.
Since $U_{n,\ba}$ is irreducible,
there exists a smooth irreducible curve $C$,
a point $p$ of $C$ and a morphism
$f:C\ra U_{n,\ba}$ such that
$f(p)=p_0$ and that
$\Phi_n(f(C))\cap \cR^{\rm irr}(\cP_{n,\bt})_{\ba}\neq\emptyset$.
{F}rom [\cite{Deligne:70}, Proposition 5.4], there exists an analytic
flat family $(\tilde{E},\nabla_{\tilde{E}},\tilde{\varphi})$
of connections such that
$\ker\nabla_{\tilde{E}}|_{(\BP^1\setminus D(\bt))\times C}$
is equivalent to the flat family of local systems on
$(\BP^1\setminus D(\bt) )\times C$ over $C$
induced by the morphism $f$.
Applying certain elementary transformations and tensoring
line bundles to $(\tilde{E},\nabla_{\tilde{E}},\tilde{\varphi})$,
we may assume that the eigenvalues of
$\res_{t_i}(\tilde{E})$ are $\lambda_i$ and
$\res_{t_i}(\nabla_L)-\lambda_i$ for $i=1,\ldots,n$.
We can construct a flat family of parabolic structures
$\{\tilde{l}_i\}$ and
$(\tilde{E},\nabla_{\tilde{E}},\tilde{\varphi},\{\tilde{l}_i\})$
becomes a flat family of parabolic connections.
Taking the completion at $p$,
we obtain a flat family of parabolic connections
$(E,\nabla_E,\varphi,\{l_i\})$ on $\BP^1_{\C[[x]]}$ over $\C[[x]]$.
By Corollary \ref{v-c-connection}, there exists a flat family
$(E',\nabla_{E'},\varphi',\{l'\})$
of $\balpha$-semistable parabolic connections such that
$(E,\nabla_E,\varphi,\{l_i\})\otimes \C((x)) \cong
(E',\nabla_{E'},\varphi',\{l'\})\otimes \C((x))$ and
$gr((E',\nabla_{E'})\otimes \C[[x]]/(x)) \cong
gr((E,\nabla_E)\otimes \C[[x]]/(x))$.
Then $(E',\nabla_{E'},\varphi',\{l'\})$ determines a morphism
$\Spec\C[[x]] \ra M^{\balpha}_n(\bt,\blambda,L)$.
If $q$ is the image of the closed point by this morphism,
then we have $\RH(q)=[\rho]$.
\end{proof}

\subsubsection{Proof of Theorem \ref{thm:RH} except for the properness of 
$\RH_{\bt,\blambda}$}

\begin{proof} Here we prove the assertions in Theorem  \ref{thm:RH}
except for the properness of $\RH_{\bt,\blambda}$ which will be proved in 
Proposition \ref{prop:properness}. 
The first assertion except for the properness  follows from 
Lemma \ref{lem:surjective}  and 
the second assertion is proved in  the proof of Lemma \ref{lem:surjective}.  
The last assertion follows from these assertions and Lemma \ref{lem:inf-rh}. 
For the third assertion recall the definition of 
$\cR^{irr}(\cP_{n, \bt})_{\ba}$ and $\cR^{app}(\cP_{n, \bt})_{\ba}$ in the proof of Lemma \ref{lem:surjective}. 
Let us set $\cR^{red}(\cP_{n, \bt})_{\ba} = \cR(\cP_{n, \bt})_{\ba} \setminus 
\cR^{irr}(\cP_{n, \bt})_{\ba}$.  Then we see that 
$$
\cR(\cP_{n, \bt})^{sing}_{\ba}  = \cR^{red}(\cP_{n, \bt})_{\ba} \cup \cR^{app}(\cP_{n, \bt})_{\ba}. 
$$
If $[\rho] \in \cR^{red}(\cP_{n, \bt})_{\ba}$, then $\rho$ is a reducible representation. 
Then the semisimplification of $\rho$ is a direct sum of one dimensional representation 
$\rho_1, \rho_2$. Since $\wedge^2 \rho $ is trivial, $\rho_2 \simeq \rho_1^{-1}$.  
Moreover since  $\Tr[\rho(\gamma_i) ] = a_i$ are fixed for all $1 \leq i \leq n$, we see that 
Jordan equivalence class of $\rho$, which is equal to the Jordan equivalence class of 
 $ \rho_1 \oplus \rho_1^{-1}$, has finitely many possiblity. Hence 
 $ \cR^{red}(\cP_{n, \bt})_{\ba}$ is a zero dimensional subscheme.  
 Moreover for a closed point 
 $[\rho] \in \cR^{app}(\cP_{n, \bt})_{\ba}$, $\rho$ is irreducible and 
 $\rho(\gamma_i) = \pm id$ for some $i$ by definition. 
 This means that $\rho$ is determined by $\rho(\gamma_j)$ for $j \not= i$ and 
 so  $\dim \cR^{app}(\cP_{n, \bt})_{\ba} = \dim \cR(\cP_{n-1, \bt})_{\ba'}$. 
Noting  that $\dim \cR(\cP_{n, \bt})_{\ba} = 2n -6$ for $n \geq 3$, we have 
$ \dim \cR^{app}(\cP_{n, \bt})_{\ba} = 2n- 8$ for $n \geq 4$. In both cases, 
if $n \geq 4$, the codimensions of the subschemes are at least $2$. 
\end{proof}

\subsection{The case of $n=4$}
In  the case of  $n = 4$, let us recall the 
isomorphism 
$$
T'_4/PGL_2 \simeq B = \BP^1 - \{ 0, 1, \infty\}, 
$$
where $B$ is one-dimensional  space of time variables as usual.  
Here the group $PGL_2$ acts on the base space $\BP^1$ by 
linear fractional transformations. Therefore 
the  family and the morphism can be descended and one 
obtains the following commutative diagram:
\begin{equation}
 \begin{CD}
\cS & \stackrel{u}{\simeq} &  M^{\balpha}_4 @> \RH_4 >> \cR_4 \\
@V \pi VV @V \pi_4 VV  @V \phi_4 VV \\
B \times \Lambda_4 @=B \times \Lambda_4 @> id \times \mu >> B \times \cA_4. 
\end{CD}
\end{equation} 
Here the family $\pi:\cS \lra B \times \Lambda_4$ is the family of Okamoto space 
of initial conditions.  The isomorphism $u$ will be constructed in 
\cite{IIS-2}.

\vspace{1cm} 
\section{Irreducibility of  $\cR(\cP_{n, \bt})_{\ba}$}
\label{sec:rep-irred}

As in Lemma \ref{lem:cat-quot}, we have   the natural quotient morphism  
$$
\begin{array}{ccl}
\Phi_n: SL_2(\C)^{n-1} &  \lra &  \cR(\cP_{n, \bt}) \simeq \Spec [
\left(R_{n-1}\right)^{Ad(SL_2(\C))}] \\
   &   &   \\
   (M_1, M_2, \ldots, M_{n-1}) & \mapsto &  [M_1, M_2, \ldots, M_{n-1}] \\
\end{array} 
$$
where $R_{n-1}$  denotes the affine coordinate ring of $SL_2(\C)^{n-1}$.  
Under this quotient morphism, for   $\ba = (a_1, \ldots, a_n) \in \cA_n = \C^n $,  
the  subscheme $\cR(\cP_{n, \bt})_{\ba}$  in \eqref{eq:mod-rep} 
is isomorphic to  
$$
\cR(\cP_{n, \bt})_{\ba} = \{ [M_1, \ldots, M_{n-1}] \in 
\cR(\cP_{n, \bt}) \ | \ \Tr(M_i) = a_i,  (1 \leq i \leq n-1),  
 \Tr(M_1M_2 \cdots M_{n-1})^{-1}  = a_n \}. 
$$
\begin{Proposition} \label{prop:irr-p2} Assume that $n \geq 4$. 
The affine scheme $\cR(\cP_{n, \bt})_{\ba}$ is irreducible.  
\end{Proposition}

Set $U_{n, \ba} := \Phi_{n}^{-1} ( \cR(\cP_{n, \bt})_{\ba} )$ so that we have  
a surjective morphism $U_{n, \ba} \lra \cR(\cP_{n, \bt})_{\ba}$ of schemes. 
Because $\Tr (M_1M_2 \cdots M_{n-1})^{-1}= \Tr (M_1 M_2 \cdots M_{n-1})$ for 
$M_i \in SL_2(\C)$,   we have  
\begin{equation} \label{eq:isospectral} 
  U_{n, \ba}  = \{ (M_1, \ldots, M_{n-1})  
   \in SL_2(\C)^{n-1} \ | \ \Tr (M_i) = a_i,  1 \leq i \leq n-1,   
 \Tr (M_1M_2 \cdots M_{n-1})  = a_n  \} . 
\end{equation}
Then it suffices to show the following

\begin{Proposition}\label{prop:irr-2} 
The scheme $U_{n, \ba}$ is irreducible.   
\end{Proposition}

Let us prove some easy lemma which we will use later. The proof of the following lemma is easy and we omit it.  

\begin{Lemma}\label{lem:quadric} 
Fix    $a \in \C$ and define  
$$
V_{ a} =\{  A= \left( \begin{array}{cc} s & t \\  u & v \\ \end{array} \right) 
\in SL_2(\C) \ | \ \Tr (A) = a \}.  
$$
\begin{enumerate}
\item Then $V_{a}$ is an irreducible affine subscheme of $\C^3$.  
\item Let us define a quadratic hypersurface in $\BP^3_{\C}$ as:
\begin{equation}\label{eq:quadric}
\overline{V}_a := \{ [x: y: z: w] \in \BP^3_{\C} \ | \   x^2 - a x w +w^2 + yz = 0 \} .
\end{equation}  
Then we have an isomorphism 
$V_a \simeq \overline{V}_a \setminus \{ w = 0 \}$, that is, $\overline{V}_a$ 
is a compactification of $V_a$.  If $a \not= \pm 2$, 
$\overline{V}_a$ is a smooth quadric hypersurface, and if $a = \pm 2$, $\overline{V}_a$ is 
a cone over a conic and have a unique singular point at $p_a=[x:y:z:w]=[a/2:0:0:1]$.  
\end{enumerate}
\end{Lemma}

Fix  $\ba=(a_1, \ldots, a_{n}) \in \cA_n$ and set 
$\ba' = (a_1, \ldots, a_{n-1}) $.   Using the notation in Lemma \ref{lem:quadric}, we 
set  
\begin{equation}\label{eq:product}
V_{\ba'} :=V_{a_1} \times V_{a_2} \times \cdots \times V_{a_{n-1}}
\subset 
\overline{V}_{\ba'}:=\overline{V}_{a_1} \times \overline{V}_{a_2}  \cdots \times \overline{V}_{a_{n-1}}
\end{equation}
It is obvious that  $U_{n, \ba}$ is a Cartier divisor  of the scheme 
$$
 V_{\ba'} = \{ (M_1, \ldots, M_{n-1} ) \in SL_2(\C) \ | \ \Tr (M_i) = a_i, 1 \leq i \leq n-1 \}
$$
defined by the equation
\begin{equation}\label{eq:traceeq}
\Tr (M_1M_2 \cdots M_{n-1}) = a_n . 
\end{equation}
Again from Lemma \ref{lem:quadric}, 
we can introduce a homogeneous coordinates $[x_i: y_i: z_i: w_i] \in \BP_{\C}^3$ such that
$$
\overline{V}_{a_i} = \{ [x_i: y_i: z_i: w_i] \in \BP_{\C}^3 \ | \  F_{a_i} = x_i^2 - a_i x_i w_i +w_i^2 + y_iz_i = 0 \} 
$$
Let us denote by 
$\overline{U}_{n, \ba} $ the closure of $U_{n, \ba} \subset V_{\ba'}$ in $\overline{V}_{\ba'} 
\subset (\BP_{\C}^3)^{n-1}$.  It is easy to see that $\overline{U}_{n, \ba}$ is also a Cartier divisor 
in $\overline{V}_{\ba'}$.

For $ 1 \leq i \leq n-1$, 
 set $T_{n-2, i} = V_{a_1} \times \cdots \times \widehat{{V}_{a_i}} \times \cdots \times V_{a_{n-1}} $ and 
$ \overline{T}_{n-2, i} = \overline{V}_{a_2} \times\cdots \times \widehat{\overline{V}_{a_i}} \times \cdots  \times \overline{V}_{a_{n-1}}$ (omitting 
$i$-th factors) and consider the $i$-th projections
\begin{equation}\label{eq:family}
\begin{array}{ccc}
U_{n, \ba} &  \hookrightarrow  &V_{a_1} \times V_{a_2} \times\cdots \times V_{a_{n-1}}  \\
\pi_i \downarrow \quad  &        &   \pi'_i \downarrow  \quad \\
T_{n-2, i} & =  & 
T_{n-2, i}  \\
\end{array}
\quad 
\begin{array}{ccc}
\overline{U}_{n, \ba} &  \hookrightarrow  & \overline{V}_{a_1} \times \overline{V}_{a_2} \times\cdots \times \overline{V}_{a_{n-1}}  \\
\pi_i \downarrow \quad  &        &   \pi'_i \downarrow  \quad \\
\overline{T}_{n-2, i} & =  & 
\overline{T}_{n-2, i}  \\
\end{array}
\end{equation}

\begin{Lemma}\label{lem:muti-1}
  For each $ 1 \leq i \leq n-1$, 
the family $\pi_i  : \overline{U}_{n, \ba} \lra \overline{T}_{n-2, i}$ can be considered as a 
family of hyperplane sections of 
$\overline{V}_{a_i} \subset \BP_{\C}^3$ parametrized by $\overline{T}_{n-2, i}$. Therefore 
$ \overline{U}_{n, \ba} \subset \overline{V}_{\ba'} $ 
is a hypersurface defined by a multi-homogeneous polynomial
\begin{equation}\label{eq:deftrace}
H_{\ba}= H_{\ba}(x_1, y_1, z_1, w_1, \ldots, x_{n-1}, y_{n-1}, z_{n-1}, w_{n-1}) 
\end{equation}
 in the homogeneous coordinate ring of $(\BP_{\C}^3)^{n-1}$ of multi-degree $(1, \ldots, 1)$. 
 \end{Lemma}
 \begin{proof} 
 
 First we prove the assertion for $i=1$. For simplicity, we set $T = T_{n-2, 1}$ and $\overline{T} = \overline{T}_{n-2, 1}$ and we write as  $\pi_1:U_{n, \ba} \lra T$, $\overline{\pi}_1:\overline{U}_{n, \ba} 
 \lra \overline{T}$.   Take an element $(M_{2},  M_3,  \ldots,  M_{n-1}) \in T$ and set 
 $$
 M_{2} M_3 \cdots M_{n-1} = \left( \begin{array}{cc} f_1 & f_2 \\ f_3 & f_4 \end{array} \right)
 $$ 
 and $M_{1} = \left( \begin{array}{cc} s & t \\ u & a_{1}-s  \end{array} \right)
 \in V_{a_1} $ with $ s(a_1-s) -tu =1$.   
Then we can write as  
 $$
 M_{1} (M_2 \cdots M_{n-1}) =\left( \begin{array}{cc} s & t \\ u & a_{1}-s  \end{array} \right) \left( \begin{array}{cc} f_1 & f_2 \\ f_3 & f_4 \end{array} \right) = 
 \left( \begin{array}{cc} s f_1 + t f_3  & s f_2 + t f_4 \\ f_3(a_1 -s) + u f_1 & 
(a_1-s) f_4 + u f_2   \end{array} \right).
 $$
 Hence the  Cartier divisor $U_{n, \ba} \subset V_{\ba'}$ is defined by  the polynomial 
\begin{equation}\label{eq:trace-equation}
\begin{array}{lcl}
 \Tr (M_{1} M_2 \cdots M_{n-1}) -a_n & = &  s f_1 + t f_3  + (a_1-s) f_4 + u f_2 -  a_n  \\
&  = &  
 (f_1- f_4)s + f_3 t + f_2 u +(a_1 f_4 - a_n) . \\
\end{array}
 \end{equation}
First,  we will show that any irreducible component of 
$U_{n, \ba}$ is not a pullback divisor  via $\pi_1$. 
Consider the subscheme $Z $ of $T$ defined by the ideal generated by the following 
 elements:
 $$
 f_1 -  f_4,  \  f_3, \  f_2, \ a_1 f_4 - a_n.
 $$ 
 Then, 
it suffices to show  that the codimension of $Z$ in $T$ is at least $2$.  Recall that 
$T$ is a product of $V_{a_i}$'s for $i =2, \ldots, n-1$.   

If $ n \geq 4$, let us consider the natural projection  $\varphi: Z \lra 
V_{a_2} \times \cdots \times V_{a_{n-2}}$. We will show  that every closed 
fiber of $\varphi$ consists of  a finite number of points or becomes empty, which means that codimension of $Z$ in $T$ is 
at least $2$.  
For this purpose, let us set
$$
M_2 \cdots M_{n-2} = \left( \begin{array}{cc} g_1 & g_2 \\ g_3 & g_4 \end{array} \right)
$$
and 
$$
M_{n-1} = \left( \begin{array}{cc} s & t \\ u & b - s  \end{array} \right).  
$$ 
(Note that $g_1 g_4 -  g_2 g_3 = 1, s (b - s) - t u = 1$,  $s = s_{n-1}, t= t_{n-1}, u = t_{n-1},  b = a_{n-1} $).       
Since 
 $$
\left( \begin{array}{cc} f_1 & f_2 \\ f_3 & f_4 \end{array} \right) = 
 M_2 \cdots M_{n-2} M_{n-1} = \left( \begin{array}{cc} g_1 & g_2 \\ g_3 & g_4 \end{array} \right) 
  \left( \begin{array}{cc} s & t \\ u & b - s  \end{array} \right) =  \left( \begin{array}{cc} g_1 s + g_2 u  & g_1 t + g_2 (b-s)  \\ g_3 s + g_4 u  & g_3 t + g_4(b-s)   \end{array} \right), 
 $$
the ideal of $Z$ contains the following element 
 $$
 \begin{array}{lcl}
f_1 - f_4 & = &  (g_1 + g_4) s + g_2 u - g_3 t - g_4 b  \\
f_2 & = & g_1 t + g_2 (b - s)  \\
f_3 & = & g_3 s + g_4 u.     
\end{array}
 $$
 Using the relations  $ g_1 g_4 - g_2 g_3  = 1, s(b-s) - t u = 1$,  from these elements 
 we can obtain the following  elements of the ideal  of $Z$ 
$$
 s^2 - g_4^2, \quad t^2 - g_2^2,  \quad  u^2 -  g_3^2. 
$$
This means that every closed fiber of 
the projection $\varphi:Z \lra V_{a_2} \times \cdots \times V_{a_{n-2}}$  
consists of  finitely many points or becomes the empty set as desired.  
Let us recall the natural projection $\pi_1:U_{n, \ba} \lra  T$. 
The assertion implies that the Cartier divisor $U_{n, \ba}$ defined by the 
polynomial (\ref{eq:trace-equation}) has no irreducible component 
which is a pullback Cartier divisor  by $\pi_1$.  Then,  
from the expression in  (\ref{eq:trace-equation}), 
we conclude  that the polynomial (\ref{eq:trace-equation}) is of  degree 1 with respect to $s, t, u$
and hence the fibers of  the compactifications  
$\overline{\pi}_1:\overline{U}_{n, \ba} \lra \overline{T} $ of morphism $\pi_1$ 
 are hyperplane sections of the quadric 
hypersurface $\overline{V}_{a_1} \subset \BP^3_{\C}$.   This proves the assertion for $i=1$.  
Since 
\begin{equation}\label{eq:cyclic}
\Tr( M_i M_{i+1} \cdots M_{n-1} M_1 \cdots M_{i-1}) = 
\Tr (M_1 M_2 \cdots M_{n-1}), 
\end{equation}
the same is true for the $i$-th factor.  Now we can  conclude 
that $\overline{U}_{n, \ba}$ is defined by the multi-homogeneous 
polynomial of $H_{\ba}$ of multi-degree $(1, \ldots, 1)$.  
\end{proof} 

Now we prove
\begin{Lemma}\label{lem:fiber-irr}
For any $\ba = (a_1, \ldots, a_n) \in \cA_{n}$ and $i,  1 \leq i \leq n-1$,  the general  
fiber of $\pi_i : \overline{U}_{n, \ba} \lra \overline{T}_{n-2, i}$ is irreducible and reduced.
\end{Lemma}

\begin{proof}  By  (\ref{eq:cyclic}), we only have to prove the assertion for $i =1$.  {F}rom now on, we set 
$T_{n-2} := T_{n-2, 1}, \overline{T}_{n-2} = \overline{T}_{n-2, 1}$ 
and $\pi= \pi_1$, etc.

For  $(M_1, M_2, M_3, \ldots, M_{n-1}) \in U_{n, \ba}$, 
write 
$$
M_1  = \left( \begin{array}{cc} s_1 & t_1 \\ u_1 & a_{1}-s_1  \end{array} \right), \quad 
M_2 M_3 \cdots M_{n-1} = \left( \begin{array}{cc} f_1 & f_2 \\ f_3 & f_4 \end{array} \right). 
$$
Then for a fixed $(M_2, M_3, \ldots, M_{n-1}) \in T_{n-2}$, 
the fiber of $\pi_{|U_{n, \ba}}:U_{n, \ba} \lra T_{n-2}$ is defined by the equations 
$$
s_1 f_1 + t_1 f_3  + (a_1-s_1) f_4 + u_1 f_2 -  a_n = 0, \quad s_1^2 -a_1 s_1 + 1 + 
 t_1 u_1 = 0. 
$$
of $(s, t, u) \in \C^3$.  Therefore again for a fixed $(M_2, M_3, \ldots, M_{n-1}) \in T_{n-2}$
the equations of the fiber $\overline{U}_{\bf f} $
of $\pi:\overline{U}_{n, \ba} \lra \overline{T}_n$   is given by 
\begin{equation}\label{eq:equations}
 (f_1- f_4) x_1 + f_3 y_1 + f_2 z_1 + (a_1 f_4 - a_n) w_1 =0, \quad 
 x_1^2 - a_1 x_1 w_1 + w_1^2 + y_1 z_1 = 0. 
\end{equation} 
Recall that $[x_1:y_1:z_1:w_1] $ 
is the homogeneous coordinate for $\BP_{\C}^3$ and 
$\overline{V}_{a_1} =\{   x_1^2 - a_1 x_1 w_1 + w_1^2 - y_1 z_1 = 0 \} 
\subset \BP^3_{\C}$ is a quadric hypersurface.  Therefore the fiber $\overline{U}_{\bf f}$ is
a complete intersection of the quadric and the hyperplane defined by 
the equations above.  In order to see that the general 
fiber $\overline{U}_{\bf f}$  is irreducible, it suffices to show that for a special 
choice of $(M_2, \ldots, M_n)$ the fiber is irreducible.  
First we consider the Zariski open subset defined by  $z_1 \not= 0$. 
Setting $X = \frac{x_1}{z_1}, Y = \frac{y_1}{z_1}, W = \frac{w_1}{z_1}$, 
the equation above can be reduced to 
\begin{equation}
(f_1 - f_4) X + f_3 Y + f_2 + (a_1 f_4 - a_n) W = 0,  \quad Y = -X^2 + a_1 X W - W^2. 
\end{equation}
Therefore the above open subscheme is isomorphic to 
an affine quadratic curve in $\C^2$ defined by 
\begin{equation}\label{eq:affine-curve} 
- f_3 X^2 + a_1 f_3 X W - f_3 W^2 + (f_1 - f_4) X + (a_1 f_4 - a_n) W + f_2 = 0.  
\end{equation}
It is easy to see that this affine quadric curve  is reducible if and only if 
\begin{equation}\label{eq:determinant}
\det \left(
\begin{array}{ccc}
-f_3 & a_1f_3/2 &  (f_1 - f_4)/2 \\
a_1f_3/2 &  -f_3 &  (a_1 f_4 - a_n)/2 \\
(f_1 - f_4)/2 & (a_1 f_4 - a_n)/2 & f_2 \\
\end{array}
\right) = 0.  
\end{equation}
Then in order to prove the Lemma, we only have to find 
$(M_2, M_3, \ldots, M_{n-1}) \in T_{n-2} $ such that 
the equation (\ref{eq:determinant}) has no solution.  
For this purpose we choose the following matrices. 
$$
M_2 = \left(  \begin{array}{cc} 0 &  1 \\ -1 & a_2 
\end{array} \right),  \quad  M_3 M_4 \cdots M_{n-1} = 
\left(  \begin{array}{cc} c &  p c^{-1} \\  0  & c^{-1} 
\end{array} \right). 
$$
Here $c$ is a non-zero constant and $p$ is a free parameter.  
Then we have  
$$
M_2 M_3 \cdots M_{n-1} = \left(  \begin{array}{cc} f_1 & f_2 \\  f_3  & f_4 \end{array} \right) = 
\left(  \begin{array}{cc} 0 &   c^{-1} 
\\  -c  & -p c^{-1} + a_2 c^{-1}  
\end{array} \right). 
$$
Then setting $f_1 = 0, f_2 = c^{-1}, f_3=  -c , f_4 = (-p + a_2)c^{-1} $, the equation (\ref{eq:determinant})  becomes 
\begin{equation}\label{eq:determinant2}
\det \left(
\begin{array}{ccc}
c & -a_1c/2 &  (p- a_2)c^{-1} /2 \\
-a_1 c/2 &  c &  (a_1(a_2-p) c^{-1}- a_n)/2 \\
(p- a_2)c^{-1}/2 & (a_1(a_2-p) c^{-1}- a_n)/2 & c^{-1} \\
\end{array}
\right) = 0.  
\end{equation}
Expanding the determinant and multiplying a non-zero constant, 
 we obtain the explicit equation
\begin{equation}\label{eq:final-eq}
p^2 +(a_1 a_{n-1} c-2 a_2) p + 
(a_2^2 -a_1 a_2 a_n c - 4c^2 + a_1^2 c^2 + a_n^2 c^2)=0. 
\end{equation}
Since $a_1, a_2, a_{n-1}, c$ are fixed constant and  
$p$ is  a free parameter, we can choose  a value of $p$ so that 
the equation (\ref{eq:final-eq}) has no solution.  Therefore 
for such a choice, the affine quadratic curve (\ref{eq:affine-curve}) 
is irreducible and reduced, and then it is easy to conclude that the fiber 
$\overline{U}_{\bf f}$ defined by (\ref{eq:equations}) is irreducible and reduced.  

\end{proof}

\noindent
{\bf Proof of Proposition \ref{prop:irr-2}. } \quad 

\begin{proof}
In the proof of Lemma \ref{lem:muti-1}, we have shown that any irreducible component of 
$U_{n, \ba}$ is not a pullback divisor  via $\pi_i $ for $1 \leq i \leq n-1$. 
Since   the  Cartier divisor  $\overline{U}_{n, \ba}$ in  
$\overline{V}_{\ba'}= \overline{V}_{a_1} \times  \cdots \times \overline{V}_{a_{n-1}}$  has the multi-degree   
$(1, \ldots, 1)$ with respect to the embeddings $\overline{V}_{a_i} \hookrightarrow \BP^3_{\C}$, 
each closed fiber of 
 $\overline{\pi}_i: \overline{U}_{n, \ba} \lra \overline{T}_{n-2, i} $ is isomorphic to 
 a hyperplane section of 
 the quadric hypersurface $\overline{V}_{a_i}$.  If $\overline{U}_{n, \ba}$ is a reducible Cartier divisor in 
 $\overline{V}_{\ba'}$, there exists an integer $i$, $1 \leq i \leq n-1$ such that 
 all closed fiber $\overline{\pi}_i$ is a reducible conic or a double line, and this contradicts to Lemma \ref{lem:fiber-irr}.  
\end{proof}

\vspace{1cm}
\section{Irreducibility of  $M^{\balpha}_n(\bt,\blambda,L)$}
\label{sec:irr-stable}

\begin{Proposition}\label{prop:stable-mod-irr} 
 $M^{\balpha}_n(\bt,\blambda,L)$ is irreducible.
\end{Proposition}

\begin{proof}
 From the proof of Lemma \ref{lem:surjective},
\[
 M^{\balpha}_n(\bt,\blambda,L)^{\sharp}
 \stackrel{\RH}\longrightarrow
 \cR(\cP_{n, \bt})_{\ba}^{\sharp}
\]
is an analytic isomorphism.
Since $\cR(\cP_{n, \bt})_{\ba}$ is irreducible by
Proposition \ref{prop:irr-p2},
$M^{\balpha}_n(\bt,\blambda,L)^{\sharp}$
is also irreducible.
If $M^{\balpha}_n(\bt,\blambda,L)=
M^{\balpha}_n(\bt,\blambda,L)^{\sharp}$,
there is nothing to prove.
So assume that $2\lambda_i\in\Z$ for some $i$ or
$\sum_{i=1}^n \epsilon_i \lambda_i \in \Z$ for
some $(\epsilon_i)\in \{\pm 1\}^n$.
Note that $M^{\balpha}_n(\bt,\blambda,L)$ is smooth over $\C$
of equidimension $2n-6$ (cf. Proposition \ref{prop:smoothness}). 
So it suffices to show that the dimension of
$M^{\balpha}_n(\bt,\blambda,L)\setminus
M^{\balpha}_n(\bt,\blambda,L)^{\sharp}$
is less than $2n-6$.

First we consider the case $2\lambda_i\in\mathbf{Z}$ for some $i$.
By composing elementary transforms at $t_i$, we can obtain an isomorphism
\[
 \tau:\cM_n(\bt,\blambda,L)\stackrel{\sim}\lra
 \cM_n(\bt,\blambda',L')
\]
of moduli stacks of parabolic connections without stability condition,
where $\lambda'_i=\res_{t_i}(\nabla_{L'})/2$
and $\lambda'_j=\lambda_j$ for $j\neq i$.
Put
\[
 A_i:=\left\{ \tau^{-1}(E,\nabla_E,\varphi,\{l_j\})\in
 M^{\balpha}_n(\bt,\blambda,L) \left|
 \text{$(E,\nabla_E)$ is an irreducible connection and
 $\res_{t_i}(\nabla_E)=\lambda'_i\mathrm{id}$} \right.\right\}.
\]
Tensoring a certain line bundle with a connection
having the residue $m/2$ at $t_i$ for some integer $m$,
we may assume that $\lambda'_i=0$.
Put
\begin{align*}
 \tilde{\bt}&:=(t_1,\ldots,t_{i-1},t_{i+1},\ldots,t_n), \\
 \tilde{\blambda}&:=
(\lambda_1,\ldots,\lambda_{i-1},\lambda_{i+1},\ldots,\lambda_n).
\end{align*}
Then $A_i$ consists of the irreducible
$(\tilde{\bt},\tilde{\blambda})$-parabolic connections
$(E,\nabla_E,\varphi,\{l_j\})$ with a one dimensional subspace
$l_i\subset \tau(E)|_{t_i}$.
So we have $\dim A_i=2(n-1)-6+1<2n-6$.

Next consider the case
$\sum_{i=1}^n \epsilon_i\lambda_i\in\mathbf{Z}$
with $\epsilon_i\in\{\pm 1\}$.
As in the proof of Lemma \ref{lem:surjective}, we can obtain,
by composing elementary transforms, an isomorphism
\[
 \tau:\cM_n(\bt,\blambda,L)\stackrel{\sim}\lra
 \cM_n(\bt,\blambda',L')
\]
of moduli stacks of parabolic connections without stability condition,
where $\deg L'=0$, $\epsilon_i\lambda_i-\lambda'_i\in\mathbf{Z}$ for any $i$,
$\sum_{i=1}^n\lambda'_i=0$ and $\lambda'_i=\res_{t_i}(\nabla_{L'})/2$
if $2\lambda_i\in\mathbf{Z}$.
Put
\[
 B:=\left\{ \tau^{-1}(E,\nabla_E,\varphi,\{l_i\})\in
 M^{\balpha}_n(\bt,\blambda,L) \left|
 \begin{array}{l}
 \text{there exists a subconnection
 $(F,\nabla_F)\subset(E,\nabla_E)$} \\
 \text{satisfying $\res_{t_i}(\nabla_F)=\lambda'_i$ for any $i$.}
 \end{array}
 \right.\right\}
\]
Take the sections $\omega_1,\omega_3\in H^0(\Omega^1_{\BP^1}(D(\bt)))$
satisfying $\res_{t_i}(\omega_1)=\lambda'_i$ and
$\res_{t_i}(\omega_3)=\res_{t_i}(\nabla_{L'})-\lambda'_i$
for $i=1,\ldots,n$.
Take any member $\tau^{-1}(E,\nabla_E,\varphi,\{l_j\})\in B$.
Then there exists a subbundle $F\subset E$ satisfying
$\nabla_E(F)\subset F\otimes\Omega^1_{\BP^1}(D(\bt))$
and $\res_{t_i}(\nabla_F)=\lambda'_i$ for any $i$.
Since $\sum_{i=1}^n\lambda'_i=0$,
$F\cong\cO_{\BP^1}$ and $E/F\cong\cO_{\BP^1}$.
Thus we have $E\cong\cO_{\BP^1}\oplus\cO_{\BP^1}$
and $\nabla_E$ can be given by
\begin{gather*}
 \nabla_E:\cO_{\BP^1}\oplus\cO_{\BP^1}\lra
 (\cO_{\BP^1}\oplus\cO_{\BP^1})\otimes\Omega^1_{\BP^1}(D(\bt)) \\
 \begin{pmatrix}
  f_1 \\ f_2
 \end{pmatrix}
 \mapsto
 \begin{pmatrix}
  df_1 \\ df_2
 \end{pmatrix}
 +
 \begin{pmatrix}
  \omega_1 & \omega_2 \\
  0 & \omega_3
 \end{pmatrix}
 \begin{pmatrix}
  f_1 \\ f_2
 \end{pmatrix}
\end{gather*}
So put
\[
 Z:=\left\{ (\omega_2,\{l_i\}_{i=1}^n) \left|
 \begin{array}{l}
 \text{$\omega_2\in H^0(\Omega^1_{\BP^1}(D(\bt)))$ and
 $l_i\subset(\cO_{\BP^1}\oplus\cO_{\BP^1})|_{t_i}$
 is a line such that} \\
 \text{$l_i=(\cO_{\BP^1}\oplus 0)|_{t_i}$ for $i$ satisfying
 $\res_{t_i}(\omega_2)\neq 0$ or $2\lambda'_i\notin\mathbf{Z}$,} \\ 
 \text{the parabolic connection
 $\tau^{-1}(\cO_{\BP^1}\oplus\cO_{\BP^1},\nabla,\varphi,\{l_i\})$
 is $\balpha$-stable,} \\
 \text{where $\nabla$ is given by
 $\nabla(f_1,f_2)=(df_1,df_2)+(\omega_1f_1+\omega_2f_2,\omega_3f_2)$} \\
 \text{for $f_1,f_2\in\cO_{\BP^1}$ and
 $\varphi:\bigwedge^2(\cO_{\BP^1}\oplus\cO_{\BP^1})\stackrel{\sim}\ra L'$}
 \end{array}
 \right.\right\}
\]
Then we obtain a morphism
\[
 f: Z\lra B
\]
which is surjective.

Assume that $2\lambda_i\notin\mathbf{Z}$ for some $i$.
In this case, $\omega_1-\omega_3\neq 0$ and the group
\[
 G=\left\{\left.
 \begin{pmatrix}
  c & a \\
  0 & c^{-1}
 \end{pmatrix}
 \right|
 c\in\C^{\times}, \; a\in\C \right\}\left/
 \left\{ \pm
 \begin{pmatrix}
  1 & 0 \\
  0 & 1
 \end{pmatrix}
 \right\}\right.
\]
acts freely on
$Z\setminus\{\text{the locus $\omega_2\in\C(\omega_1-\omega_3)$}\}$
satisfying $f(gx)=f(x)$ for $g\in G$ and $x\in Z$.
Thus we can see that
$\dim 
f(Z\setminus\{\text{the locus $\omega_2\in\C(\omega_1-\omega_3)$}\})
\leq\dim H^0(\Omega^1_{\BP^1}(D(\bt)))-2=n-3$.
Take a point $t_j$ satisfying $\res_{t_j}(\omega_1-\omega_3)\neq 0$
and put
\begin{align*}
 Y&:=\{(\omega_2,\{l_i\})\in Z|\omega_2\in\C(\omega_1-\omega_3)\}, \\
 X&:=\left\{(\omega_2,\{l_i\})\in Y\left|
 l_i=\C
 \begin{pmatrix}
  \res_{t_j}(\omega_2) \\ \res_{t_j}(\omega_3-\omega_1)
 \end{pmatrix}
 \; \text{if $2\lambda_i\in\mathbf{Z}$ and $\res_{t_i}(\omega_2)\neq 0$}
 \right\}\right.. 
\end{align*}
Then $\dim X\leq 1$ and $G$ acts freely on $Y\setminus X$.
So we have $\dim f(Y)\leq n-3$.
Thus we have $\dim B<2n-6$.

Assume that $2\lambda_i\in\mathbf{Z}$ for all $i$.
Take any member $\tau^{-1}(E,\nabla_E,\varphi,\{l_j\})$ of $B$.
If $\res_{t_i}(\nabla_E)\neq\lambda'_i\mathrm{id}$ for any $i$,
the object $\tau^{-1}(E,\nabla_E,\varphi,\{l_j\})$
can not be stable for any choice of $\balpha$
because it has a nontrivial endomorphism.
Thus $\res_{t_i}(\nabla_E)=\lambda'_i\mathrm{id}$ for some $i$.
Then we can check that $G$ acts freely on $Z$ and
$\dim f(Z)\leq \dim Z-\dim G=n-3$.
Since $M^{\balpha}_n(\bt,\blambda,L)\setminus
M^{\balpha}_n(\bt,\blambda,L)^{\sharp}$
is a union of the subsets like $A_i,B$, we have
$\dim (M^{\balpha}_n(\bt,\blambda,L)\setminus
M^{\balpha}_n(\bt,\blambda,L)^{\sharp})<2n-6$.
\end{proof}

\section{Properness of the Riemann-Hilbert correspondence}
\label{sec:properness}

We set
\[
 \cR(\cP_{n, \bt})^{irr}_{\ba}:=
 \{ [\rho]\in\cR(\cP_{n, \bt})_{\ba} |
 \text{$\rho$ is irreducible} \}.
\]
For $\blambda\in\Lambda_n$ such that
$2\cos2\pi\lambda_i=a_i$ for $i=1,\ldots,n$,
we put
\[
 M^{irr}_n(\bt, \blambda, L):=
 \{ p\in M^{\balpha}_n(\bt, \blambda, L) |
 \RH(p)\in \cR(\cP_{n, \bt})^{irr}_{\ba} \}.
\]

\begin{Lemma}\label{lem:proper-on-irred}
 The restriction
 $M^{irr}_n(\bt, \blambda, L)\stackrel{\RH}\longrightarrow
 \cR(\cP_{n, \bt})^{irr}_{\ba}$
 is a proper morphism.
\end{Lemma}

\begin{proof}
As in the proof of Lemma \ref{lem:surjective},
we can obtain an isomorphism
\[
 \tau: M^{irr}_n(\bt,\blambda,L)\stackrel{\sim}\lra
 M^{irr}_n(\bt,\blambda',L')
\]
by composing elementary transforms and other transforms,
where $0\leq \mathrm{Re}(\lambda'_i)<1$
and $0\leq \mathrm{Re}(\res_{t_i}(\nabla_{L'})-\lambda'_i)<1$
for $i=1,\ldots,n$.
Let $N^{irr}_n(\bt,\blambda',L')$ be the moduli space of
irreducible connections $(E,\nabla_E)$ with a horizontal isomorphism
$\det(E)\cong L'$ satisfying $\det(\res_{t_i}(\nabla_E)-\lambda'_i)=0$
for $i=1,\ldots,n$.
As in the proof of Lemma \ref{lem:surjective},
we obtain a factorization
\[
 \RH:M^{irr}_n(\bt,\blambda,L)\stackrel{\sim}\lra
 M^{irr}_n(\bt,\blambda',L')\stackrel{\pi}\lra
 N^{irr}_n(\bt,\blambda',L')\stackrel{\sim}\lra
 \cR(\cP_{n, \bt})^{irr}_{\ba}.
\]
So it is sufficient to show that
$\pi:M^{irr}_n(\bt,\blambda',L')\ra N^{irr}_n(\bt,\blambda',L')$
is proper.
Let $R$ be any discrete valuation ring and $K$
its quotient field.
Assume that a commutative diagram
\[
 \begin{CD}
  \Spec K @>\eta >> M^{irr}_n(\bt,\blambda',L') \\
  @VVV @VV\pi V \\
  \Spec R @>\xi >> N^{irr}_n(\bt,\blambda',L')
 \end{CD}
\]
is given.
Then $\xi$ corresponds to a flat family of
connections $(E,\nabla_E,\varphi)$ on $\BP^1_R$ over $R$
and $\eta$ corresponds to parabolic structures
$\{l_j\subset (E\otimes K)|_{t_j\otimes K} \}$ satisfying
$(\res_{t_j\otimes K}(\nabla_E\otimes K)-\lambda'_j)|_{l_j}=0$.
We can take subbundles $\tilde{l}_j\subset E|_{t_j\otimes R}$
such that $\tilde{l}_j\otimes K=l_j$ as subspaces of
$(E\otimes K)|_{t_j\otimes K}$.
Then we have
$(\res_{t_j\otimes R}(\nabla_E)-\lambda'_j)|_{\tilde{l}_j}=0$
and $(E,\nabla_E,\varphi,\{l_j\})$ becomes an $R$-valued point
of $M^{irr}_n(\bt,\blambda',L')$ whose image by $\pi$
is $\xi$ and the restriction to $\Spec K$ is $\eta$.
Thus $\pi:M^{irr}_n(\bt,\blambda',L')\ra N^{irr}_n(\bt,\blambda',L')$
becomes a proper morphism by the valuative criterion of properness.
\end{proof}

Put
\begin{align*}
 \cR(\cP_{n, \bt})^{red}_{\ba}&:=
 \cR(\cP_{n, \bt})_{\ba}\setminus
 \cR(\cP_{n, \bt})^{irr}_{\ba} \\
 M^{red}_n(\bt, \blambda, L) &:=
 M^{\balpha}_n(\bt, \blambda, L) \setminus
 M^{irr}_n(\bt,\blambda,L).
\end{align*}

\begin{Lemma}\label{lem:red-compact}
$M^{red}_n(\bt, \blambda, L)$ is a closed subset of 
$M^{\balpha}_n(\bt, \blambda, L)$ which is proper over $\C$.
\end{Lemma}

\begin{proof}
First we will show that $M^{red}_n(\bt,\blambda,L)$
is a constructible subset of $M^{\balpha}_n(\bt, \blambda, L)$.
Take any member $(E,\nabla_E,\varphi,\{l_i\})$ of
$M^{\balpha}_n(\bt, \blambda, L)$.
Then
\begin{align*}
 & (E,\nabla_E,\varphi,\{l_i\}) \in M^{red}_n(\bt,\blambda,L) \\
 \Leftrightarrow &
 \text{$\exists F\subset E$: subbundle of rank $1$ such that
 $\nabla_E(F)\subset F\otimes\Omega^1_{\BP^1}(D(\bt))$}.
\end{align*}
Note that such subbundles $F$ must satisfy
\begin{equation}\label{choice-subconnection}
 \res_{t_i}(\nabla_E|_F)=\lambda_i \; \text{or} \;
 \res_{t_i}(\nabla_E|_F)=\res_{t_i}(\nabla_L)-\lambda_i
 \quad \text{for $i=1,\ldots,n$},
\end{equation}
from which the choice of isomorphism classes of $(F,\nabla_E|_F)$
must be finite. 
Let $L_1,L_2,\ldots,L_m$ be all the line bundles with connections
satisfying the condition (\ref{choice-subconnection}).
Then we can see that the set
\[
 \left\{ (p,L_i\stackrel{\iota}\ra E) \left|
 \begin{array}{l}
  \text{$p=(E,\nabla_E,\varphi,\{l_i\})\in
  M^{\balpha}_n(\bt, \blambda, L)$ and $\iota$ is an injective} \\
  \text{homomorphism such that
  $\nabla_E(\iota(L_i))\subset \iota(L_i)\otimes\Omega^1_{\BP^1}(D(\bt))$}
 \end{array}
 \right\}\right.
\]
can be parameterized by a scheme of finite type over
$M^{\balpha}_n(\bt, \blambda, L)$.
Thus $M^{red}_n(\bt,\blambda,L)$ is a constructible subset of
$M^{\balpha}_n(\bt, \blambda, L)$.

So it suffices to show by \cite{Hart}, Chapter II, Lemma 4.5 that
$M^{red}_n(\bt,\blambda,L)$ is stable under specialization
in the compactification
$\overline{M^{\balpha'\bbeta}_n}(\bt, \blambda, L)$.
Take any scheme point $x_1\in M^{red}_n(\bt,\blambda,L)$ and
$T$ be the closure of $\{x_1\}$ in
$\overline{M^{\balpha'\bbeta}_n}(\bt, \blambda, L)$.
Take any point $x_0\in T$.
Put $K:=k(x_1)$.
Then there exists a discrete valuation ring $R$
with quotient field $K$ which dominates $\cO_{x_0}$.
A morphism
$\iota:\Spec R \ra \overline{M^{\balpha'\bbeta}_n}(\bt, \blambda, L)$
satisfying $\iota(\eta)=x_1$ and $\iota(\xi)=x_0$ is induced,
where $\eta$ is the generic point of $\Spec R$
and $\xi$ the closed point of $\Spec R$.
$\iota$ corresponds to a flat family of
$(\bt,\blambda)$-parabolic $\phi$-connections
$(\tilde{E}_1,\tilde{E}_2,\tilde{\phi},\tilde{\nabla},
\tilde{\varphi},\{\tilde{l}_i\})$ on $\BP^1_R$ over $R$.
Put
\begin{align*}
 (E_1,E_2,\phi,\nabla,\varphi,\{l_i\}) & :=
 (\tilde{E}_1,\tilde{E}_2,\tilde{\phi},\tilde{\nabla},
 \tilde{\varphi},\{\tilde{l}_i\})\otimes k(\eta) \\
 (E'_1,E'_2,\phi',\nabla',\varphi',\{l'_i\}) & :=
 (\tilde{E}_1,\tilde{E}_2,\tilde{\phi},\tilde{\nabla},
 \tilde{\varphi},\{\tilde{l}_i\})\otimes k(\xi),
\end{align*}
Since
$(E_1,E_2,\phi,\nabla,\varphi,\{l_i\})\in
M^{red}_n(\bt,\blambda,L)(K)$,
$\phi$ is isomorphic and there exist commutative diagrams
\begin{equation}\label{comm-phi-connection-1}
 \begin{array}{ccc}
  0 & & 0 \\
  \downarrow & & \downarrow \\
  F_1 & \xrightarrow[\sim]{\phi_1} & F_2 \\
  \downarrow & & \downarrow \\
  E_1 & \xrightarrow[\sim]{\phi} & E_2 \\
  \downarrow & & \downarrow \\
  G_1 & \xrightarrow[\sim]{\phi_2} & G_2 \\
  \downarrow & & \downarrow \\
  0 & & 0
 \end{array}
 \hspace{30pt}
 \begin{array}{ccc}
  0 & & 0 \\
  \downarrow & & \downarrow \\
  F_1 & \xrightarrow{\nabla_1}
  & F_2\otimes\Omega^1_{\BP^1}(D(\bt)) \\
  \downarrow & & \downarrow \\
  E_1 & \xrightarrow{\nabla}
  & E_2\otimes\Omega^1_{\BP^1}(D(\bt)) \\
  \downarrow & & \downarrow \\
  G_1 & \xrightarrow{\nabla_2}
  & G_2\otimes\Omega^1_{\BP^1}(D(\bt)) \\
  \downarrow & & \downarrow \\
  0 & & \; 0,
 \end{array}
\end{equation}
where $F_1,F_2,G_1,G_2$ are line bundles.
There exist quotient coherent sheaves
$\tilde{E}_1\ra \tilde{G}_1$,
$\tilde{E}_2\ra \tilde{G}_2$
which are flat over $R$ and
whose fibers over $\eta$ are isomorphic to
$E_1\ra G_1$ and $E_2\ra G_2$, respectively.
Put $\tilde{F}_1:=\ker(\tilde{E}_1\ra\tilde{G}_1)$
and $\tilde{F}_2:=\ker(\tilde{E}_2\ra\tilde{G}_2)$.
Then we obtain commutative diagrams
\[
 \begin{array}{ccc}
  0 & & 0 \\
  \downarrow & & \downarrow \\
  \tilde{F}_1 & \xrightarrow{\tilde{\phi}_1} & \tilde{F}_2 \\
  \downarrow & & \downarrow \\
  \tilde{E}_1 & \xrightarrow{\tilde{\phi}} & \tilde{E}_2 \\
  \downarrow & & \downarrow \\
  \tilde{G}_1 & \xrightarrow{\tilde{\phi}_2} & \tilde{G}_2 \\
  \downarrow & & \downarrow \\
  0 & & 0
 \end{array}
 \hspace{30pt}
 \begin{array}{ccc}
  0 & & 0 \\
  \downarrow & & \downarrow \\
  \tilde{F}_1 & \xrightarrow{\tilde{\nabla}_1}
  & \tilde{F}_2\otimes\Omega^1_{\BP^1}(D(\bt)) \\
  \downarrow & & \downarrow \\
  \tilde{E}_1 & \xrightarrow{\tilde{\nabla}}
  & \tilde{E}_2\otimes\Omega^1_{\BP^1}(D(\bt)) \\
  \downarrow & & \downarrow \\
  \tilde{G}_1 & \xrightarrow{\tilde{\nabla}_2}
  & \tilde{G}_2\otimes\Omega^1_{\BP^1}(D(\bt)) \\
  \downarrow & & \downarrow \\
  0 & & \; 0,
 \end{array}
\]
whose fibers over $\eta$ are the commutative diagrams
(\ref{comm-phi-connection-1}).
Put $G'_1:=\tilde{G}_1(\xi)/(\tilde{G}_1(\xi))_{tor}$,
$G'_2:=\tilde{G}_2(\xi)/(\tilde{G}_2(\xi))_{tor}$,
$F'_1:=\ker(\tilde{E}_1(\xi) \ra G'_1)$ and
$F'_2:=\ker(\tilde{E}_2(\xi) \ra G'_2)$,
where $(\tilde{G}_1(\xi))_{tor}$ and $(\tilde{G}_1(\xi))_{tor}$
are the torsion parts of $\tilde{G}_1(\xi)$ and $\tilde{G}_2(\xi)$,
respectively.
Then we obtain commutative diagrams
\begin{equation}\label{comm-phi-connection-2}
 \begin{array}{ccc}
  0 & & 0 \\
  \downarrow & & \downarrow \\
  F'_1 & \xrightarrow{\phi'_1} & F'_2 \\
  \downarrow & & \downarrow \\
  \tilde{E}_1(\xi) & \xrightarrow{\phi'} & \tilde{E}_2(\xi) \\
  \downarrow & & \downarrow \\
  G'_1 & \xrightarrow{\phi'_2} & G'_2 \\
  \downarrow & & \downarrow \\
  0 & & 0
 \end{array}
 \hspace{30pt}
 \begin{array}{ccc}
  0 & & 0 \\
  \downarrow & & \downarrow \\
  F'_1 & \xrightarrow{\nabla'_1}
  & F'_2\otimes\Omega^1_{\BP^1}(D(\bt)) \\
  \downarrow & & \downarrow \\
  \tilde{E}_1(\xi) & \xrightarrow{\nabla'}
  & \tilde{E}_2(\xi)\otimes\Omega^1_{\BP^1}(D(\bt)) \\
  \downarrow & & \downarrow \\
  G'_1 & \xrightarrow{\nabla'_2}
  & G'_2\otimes\Omega^1_{\BP^1}(D(\bt)) \\
  \downarrow & & \downarrow \\
  0 & & \; 0.
 \end{array}
\end{equation}
Assume that $x_0\notin M^{\balpha}_n(\bt, \blambda, L)$.
Then we have either
\[
 \left\{
 \begin{array}{l}
 \text{$\phi'_1=0$ and $\deg F'_1\geq\deg F'_2$ or} \\
 \text{$\phi'_2=0$ and $\deg G'_1\geq\deg G'_2$.}
 \end{array}\right.
\]
Assume $\phi'_1=0$ and $\deg F'_1\geq \deg F'_2$.
Then $\res_{t_i}(\nabla'_1)=\lambda'_i\phi'_1(t_i)=0$
for $i=1,\ldots,n$, where $\lambda'_i=\lambda_i$ or
$\lambda'_i=\res_{t_i}(\nabla_L)-\lambda_i$.
We must have $\nabla'_1=0$,
since $H^0(\BP^1,(F'_1)^{\vee}\otimes F'_2\otimes\Omega^1_{\BP^1})=0$.
Then
$(F'_1,0)$ collapses the stability of
$(\tilde{E}_1(\xi),\tilde{E}_2(\xi),\tilde{\phi}(\xi),
\tilde{\nabla}(\xi),\tilde{\varphi}(\xi),\{\tilde{l}_i(\xi)\})$.
Assume $\phi'_2=0$ and $\deg G'_1\geq\deg G'_2$.
Then $\res_{t_i}(\nabla'_2)=\lambda'_i\phi'_2(t_i)=0$
for $i=1,\ldots,n$, where $\lambda'_i=\lambda_i$ or
$\lambda'_i=\res_{t_i}(\nabla_L)-\lambda_i$.
Again we have $\nabla'_2=0$, because $\deg G'_1\geq\deg G'_2$.
Then $(\tilde{E}_1(\xi),F'_2)$ collapses the stability of
$(\tilde{E}_1(\xi),\tilde{E}_2(\xi),\tilde{\phi}(\xi),
\tilde{\nabla}(\xi),\tilde{\varphi}(\xi),\{\tilde{l}_i(\xi)\})$.
Hence we must have $x_0\in M^{\balpha}_n(\bt, \blambda, L)$
and so $x_0\in M^{red}_n(\bt,\blambda,L)$.
\end{proof}

Now we are ready to prove the following 

\begin{Proposition}\label{prop:properness}
 The analytic morphism 
 $$
 \RH_{\bt, \blambda}: M^{\balpha}_n(\bt, \blambda, L) \lra
 \cR_n (\cP_{n, \bt})_{\ba}
$$
defined in \eqref{eq:RHF} is  proper. 
\end{Proposition}

\begin{proof}  Recall that we have proved the assertions in 
Theorem \ref{thm:RH} except for the properness of  $\RH_{\bt, \blambda}$. 
Therefore, we see  
that $\RH_{\bt, \blambda}$ induces the an analytic isomorphism 
$$
M^{\balpha}_n(\bt, \blambda, L)^{\sharp} \stackrel{\simeq}{\lra}
 \cR (\cP_{n, \bt})^{\sharp}_{\ba}.
$$
(See \eqref{eq:rh-isom}). 
 From the third assertion  of Theorem \ref{thm:RH}, we see that 
 $\cR(\cP_{n, \bt})^{sing}_{\ba} = \cR(\cP_{n, \bt})_{\ba} \setminus 
\cR(\cP_{n, \bt})^{\sharp}_{\ba}$ has codimension $\geq 2$ in $\cR(\cP_{n, \bt})_{\ba}$.  Moreover Lemma \ref{lem:proper-on-irred} and Lemma 
\ref{lem:red-compact} shows that every  fiber of  $\RH_{\bt, \blambda}$ 
at each closed point of  $\cR_n (\cP_{n, \bt})_{\ba}$ is compact.  
Therefore the assertion  follows from the following lemma due to A. Fujiki.
\end{proof}

\begin{Lemma}\label{lem:fujiki}
 Let $f:X\ra Y$ be a surjective holomorphic mapping of
 irreducible analytic varieties.
 Assume that an analytic closed subset $S$ of $Y$ exists such that
 $\codim_Y S\geq 2$, 
 $X^{\sharp}:=f^{-1}(Y^{\sharp})$ is dense in $X$,
 where $Y^{\sharp}=Y\setminus S$ and that the restriction
 $f|_{X^{\sharp}}:X^{\sharp}\ra Y^{\sharp}$ is an analytic isomorphism.
 Moreover assume that the fibers $f^{-1}(y)$ are compact for all $y\in Y$.
 Then $f$ is a proper mapping.
\end{Lemma}

\begin{proof}
Since normalization morphism is proper, we may assume that
both $X$ and $Y$ are normal by replacing them by their normalizations.
Take any point $y\in Y$.
Since $f^{-1}(y)$ is compact, there is an open neighborhood
$U$ of $f^{-1}(y)$ and $V$ of $y$ such that
$f(U)\subset V$ and the restriction
$f|_U:U\ra V$ is a proper mapping.
Since $Y^{\sharp}$ is normal, we may assume that
$Y^{\sharp}\cap V$ is connected.
$f(X^{\sharp}\cap U)$ is open in $Y^{\sharp}\cap V$,
because $X^{\sharp}\stackrel{f}\lra Y^{\sharp}$
is an isomorphism.
$f(X^{\sharp}\cap U)$ is also closed in $Y^{\sharp}\cap V$,
because $X^{\sharp}\cap U\stackrel{f}\lra Y^{\sharp}\cap V$
is proper.
Thus we have $f(X^{\sharp}\cap U)=Y^{\sharp}\cap V$,
because $X^{\sharp}$ is dense in $X$.

Assume that $f^{-1}(V)\setminus U\neq\emptyset$ and take
a point $b\in f^{-1}(V)\setminus U$.
Take an open neighborhood $W$ of $b$ such that
$f(W)\subset V$.
$X^{\sharp}\cap W$ is nonempty and dense in $W$.
{F}rom the injectivity of $f|_{X^{\sharp}}$ and the
fact $f(X^{\sharp}\cap U)=Y^{\sharp}\cap V$,
we have $X^{\sharp}\cap W \subset X^{\sharp}\cap U$,
since $f(X^{\sharp}\cap W)\subset Y^{\sharp}\cap V$.
Taking closures in $W$, we have
$b\in W\subset U$, which is a contradiction.
Thus we have $f^{-1}(V)=U$.
\end{proof}

\section{List of Notation}
\begin{tabular}{lcc}
Notation  &   number of equation &  page \\
    &  &  \\
   $T_n$  & \eqref{eq:config-space}  & \pageref{eq:config-space}  \\
    $\Lambda_{n}$ & \eqref{eq:exponents-space} & 
    \pageref{eq:exponents-space}  \\
    $D(\bt) $ & \eqref{eq:divisor}  & \pageref{eq:divisor} \\
   $ \pardeg_{\balpha} E$ &  
   \eqref{eq:pardeg} & \pageref{eq:pardeg} \\
$M_n^{\balpha}(\bt, \blambda, L)$ &  \eqref{eq:modulispace} &  \pageref{eq:modulispace}\\
$\overline{M_n^{\balpha'\bbeta}}(\bt, \blambda, L) $ & \eqref{eq:coarse-moduli-compact} & \pageref{eq:coarse-moduli-compact}\\
$\overline{\pi}_n:
\overline{\cM_n^{\balpha' \bbeta}}(L) \lra T_n \times \Lambda_n$ & 
\eqref{eq:fam-moduli} & \pageref{eq:fam-moduli}  \\
 $\pi_n: \cM_n^{\balpha}(L) \lra T_n \times \Lambda_n$ & \eqref{eq:fam-moduli}
 & \pageref{eq:fam-moduli} \\
 $Elm_{t_i}^+: M_n^{\balpha}(\bt, \blambda, L)   \lra  
 M_n^{\balpha}(\bt, \blambda', L(t_i))$  & 
 \eqref{eq:upper-elm}  & 
\pageref{eq:upper-elm} \\
$Elm_{t_i}^-: M_n^{\balpha}(\bt, \blambda, L)   \lra  
 M_n^{\balpha}(\bt, \blambda', L(-t_i))$  & 
 \eqref{eq:lower-elm}  & 
\pageref{eq:lower-elm} \\
$ R_i(E)$   & \eqref{eq:eigen-change} & \pageref{eq:eigen-change} \\
$BL_n$ &  \eqref{eq:backlund-group} & \pageref{eq:backlund-group} \\
$\cR(\cP_{n, \bt}) = \Spec[ \left( R_{n-1} \right)^{Ad(SL_2(\C))}]$ & 
\eqref{eq:categorical-quotient-1} & \pageref{eq:categorical-quotient-1}  \\
$T'_n := \tilde{T}_n/\Gamma_{n-1}  \lra T_n$ & \eqref{eq:finite-etale} & \pageref{eq:finite-etale} \\
$\phi_n: \cR_n \lra  T'_n  \times \cA_n $ &\eqref{eq:family-rep-s4}   & \pageref{eq:family-rep-s4} \\
 $\cR(\cP_{n, \bt})_{\ba}$ &  \eqref{eq:mod-rep} & \pageref{eq:mod-rep} \\
$\cF^0$ & \eqref{eq:end-p0} & \pageref{eq:end-p0} \\
$\cF^1$ & \eqref{eq:end-p1} & \pageref{eq:end-p1} \\
$ \cF^{1,+}$ & \eqref{eq:end-p1s} & \pageref{eq:end-p1s} \\
$\cF^{\bullet}$ & \eqref{eq:complex-org} & \pageref{eq:complex-org} \\
$\cF^{\bullet, +} $&  \eqref{eq:complex+} & \pageref{eq:complex+} \\ 
$\Omega \in H^0(M^{\balpha}_n, \Omega^2_{M^{\balpha}_n/T'_n\times\Lambda_n} )$   & \eqref{eq:symplectic} & \pageref{eq:symplectic} \\ 
$ \Omega_1 \in \Gamma(\cR_n^{\sharp}, \Omega^2_{\cR_n^{\sharp}/T'_n \times \cA_n })$ & \eqref{eq:sympl-rep} & \pageref{eq:sympl-rep} \\
${\bf RH}_{\bt, \blambda} :M^{\balpha}_n( \bt, \blambda, L) \lra 
\cR(\cP_{n, \bt})_{\ba}$ & \eqref{eq:RHF} & \pageref{eq:RHF} \\
$\RH_n: M_n^{\balpha}(L) \lra \cR_n$ & \eqref{eq:RH} & \pageref{eq:RH} \\
$  \mu_n $ &\eqref{eq:cor-local}& \pageref{eq:cor-local} \\ 
$\RH_{\bt,\blambda, | M_n^{\balpha}(\bt, \blambda, L)^{\sharp}} : M_n^{\balpha}(\bt, \blambda, L)^{\sharp}  \stackrel{\simeq}{\lra}   
\cR(\cP_{n, \bt})^{\sharp}_{\ba}$   & \eqref{eq:rh-isom} &  \pageref{eq:rh-isom} \\
 $\overline{\cM_{\cX/S}^{\cD,\balpha,\tau}}(r,d,\{d_i\})$ & \eqref{eq:functor} & 
 \pageref{eq:functor} \\
\end{tabular}

\end{document}